\documentclass{amsart}
\usepackage{graphicx,amssymb,amsfonts,amsmath,amsthm,newlfont}

\usepackage{epsfig}
\newtheorem{thm}{Theorem}[section]
\newtheorem{cor}[thm]{Corollary}
\newtheorem{lem}[thm]{Lemma}
\newtheorem{prop}[thm]{Proposition}

\theoremstyle{definition}
\newtheorem{defn}[thm]{Definition}
\newtheorem{propdef}[thm]{Definition-Proposition}

\newtheorem{example}[thm]{Example}
\theoremstyle{remark}
\newtheorem{rem}[thm]{Remark}
\numberwithin{equation}{section}

\font \rus= wncyr10

\newcommand{\xs}{\mathsf x}
\newcommand{\Z}{\mathbb Z}
\newcommand{\C}{\mathbb C}
\newcommand{\Pro}{\mathbb P}

\newcommand{\R}{\mathbb R}
\newcommand{\N}{\mathbb N}
\newcommand{\Q}{\mathbb Q}

\newcommand{\Hom}{\hbox{Hom}}
\newcommand{\Reg}{\mathrm{Reg}}
\newcommand{\Hy}{\mathfrak{H}}

\newcommand{\phil}{{\phi}_*}
\newcommand{\Upi}{\mathcal{U}_M}

\newcommand{\Up}{\mathfrak{up}}
\newcommand{\Ut}{\mathfrak{ut}}

\newcommand{\Urel}{U_{\widehat{R}/R}}
\newcommand{\PSL}{\mathsf{PSL}}

\newcommand{\A}{\mathbb {A}}

\newcommand{\Or}{\mathcal{O}}

\newcommand{\uij}{\underline{u}}
\newcommand{\ddim}{\ell}
\newcommand{\mdim}{\ell}

\newcommand{\MT}{\mathrm{MT}}

\newcommand{\gr}{\mathrm{gr}}

\newcommand{\Spec}{\mathrm{Spec\,}}

\newcommand{\ord}{\mathrm{ord}}
\newcommand{\Modf}{\mathfrak{M}^{\,\delta}}
\newcommand{\an}{\mathrm{an}}

\newcommand{\MZV}{\mathcal{Z}}
\newcommand{\I}{\mathbb{I}}
\newcommand{\ke}{k\{\epsilon\}}

\newcommand{\kes}{k\{\epsilon_1,\ldots, \epsilon_\ell\}}
\newcommand{\Ues}{U\{\epsilon_1,\ldots, \epsilon_\ell\}}

\newcommand{\Fo}{\mathcal{F}}

\newcommand{\DB}{\mathfrak{B}}

\newcommand{\Li}{\mathrm{Li}}

\newcommand{\Mo}{\mathcal{M}}
\newcommand{\Mod}{\mathfrak{M}}

\newcommand{\Fa}{\mathcal F}

\newcommand{\To}{\longrightarrow}
\newcommand{\sha}{\, \hbox{\rus x} \,}

\newcommand{\ld}{ \langle \langle}
\newcommand{\rd}{ \rangle \rangle}

\newcommand{\Gal}{\mathrm{Gal}}
\newcommand{\Sym}{\mathfrak{S}}

\newcommand{\Image}{\mathrm{Im\,}}
\newcommand{\Real}{\mathrm{Re\,}}
\newcommand{\ad}{\mathrm{ad\,}}

\newcommand{\x}{ \sim_{\mathsf{x}}}

\newcommand{\DR}{\mathrm{DR}}

\begin{document}

\title[Multiple zeta values and periods of moduli spaces.]{}%
\author{Francis C.S. Brown}%
\address{45 Rue d'Ulm}
\email{brown@clipper.ens.fr}

\maketitle
\begin{center}
\Large{\textbf{Multiple zeta values and periods of moduli spaces
$\overline{\Mod}_{0,n}(\R)$ .}}
\end{center}

\vspace{1in}

%\tableofcontents

\begin{center}
{\sc Contents.} \end{center}

\vspace{0.2in}

\noindent
1. Introduction \hfill{2}\\
2. Dihedral coordinates on $\overline{\Mod}_{0,n}(\R)$ \hfill{11}\\
3. The reduced Bar construction and Picard-Vessiot theory \hfill{33}\\
4. Manifolds with corners and Fuchsian differential equations \hfill{53}\\
5. Hyperlogarithms \hfill{61}\\
6.  The universal algebra of polylogarithms on $\Mod_{0,n}$ \hfill{67}\\
7. Period integrals on $\overline{\Mod}_{0,n}(\R)$ and generalised shuffle products \hfill{84}\\
8. Calculation of the periods of $\Mod_{0,n}$ \hfill{97}\\
9.  Appendix \hfill{104}\\
10. Bibliography \hfill{105}\\
11. Index of Notations \hfill{109}\\
\newpage

\section{Introduction}
Let $n=\ell+3\geq 4$, and let   $\Mod_{0,n}$ denote the moduli space of curves of genus $0$ with $n$ marked points.
There is a smooth compactification $\overline{\Mod}_{0,n}$, defined by Deligne, Knudsen and Mumford, such that the complement
$$\overline{\Mod}_{0,n} \backslash \Mod_{0,n}$$
is a  normal crossing divisor. Let $A, B \subset \overline{\Mod}_{0,n} \backslash \Mod_{0,n}$ denote two sets of boundary divisors
which share no irreducible components. In [G-M], Goncharov and Manin show that the relative   cohomology group
\begin{equation} \label{modularmotive} H^{\ell} (\overline{\Mod}_{0,n}\backslash A, B\backslash B\cap A)
\end{equation}
%with an additional structure called the framing, can be used to
defines a mixed Tate motive which is unramified over $\Z$.

On the other hand, let $n_1,\ldots, n_r\in \N$, and suppose that $n_r\geq 2$. The multiple zeta value $\zeta(n_1,\ldots, n_r)$
 is
the real number  defined by the convergent sum
\begin{equation} \label{mzv}
\zeta(n_1,\ldots, n_r) = \sum_{0<k_1<\ldots<k_r} {1 \over k_1^{n_1} \ldots k_r^{n_r}}\ .\end{equation}
Its weight is the quantity $n_1+\ldots +n_r$, and its depth is the number of indices $r$. We will say that
the period $2i\pi$ has weight $1$.  A very general conjecture [Go1] claims that the
periods of any mixed Tate motive unramified  over $\Z$ are  multiple zeta values. In the case of the motives $(\ref{modularmotive})$ arising from moduli spaces, this says the following. Consider a
real smooth compact submanifold  $X_B\subset \overline{\Mod}_{0,n}$ of dimension $\ell$, whose boundary is contained in $B$ and
%$$[X_B] \in \gr^W_0 H_\ell(\overline{\Mod}_{0,n}, B)\ ,$$
which does not meet $A$. It represents a class in $H_\ell(\overline{\Mod}_{0,n},B)$.  Let $\omega_A\in \Omega^\ell(\overline{\Mod}_{0,n}\backslash A)$ denote an algebraic form with singularities
 contained in $A$. In [G-M], Goncharov and Manin
conjecture  that the   integral
\begin{equation} \label{introIdefn}
I = \int_{X_B} \omega_A\end{equation} is a linear combination of
multiple zeta values, and proved that every multiple zeta value
can occur as such a period integral. In this paper, we develop
some general methods for computing periods and  prove this
conjecture as
 an application. %in this paper, we prove this conjecture.
\begin{thm}
The integral $I$ is a $\Q[2\pi i]$-linear combination of multiple zeta values  of weight
  at most $\ell$.
\end{thm}

This theorem thus
lends significant
 weight to the conjecture on the periods of all mixed Tate motives which are unramified over $\Z$.

The rough idea of our method is as follows. The set of real points $\Mod_{0,n}(\R)$ is
tesselated by a number of open cells $X_n$ which can naturally be identified with a Stasheff polytope, or associahedron.
We can reduce to the  case where the domain of integration in $(\ref{introIdefn})$ is a single  cell $X_n$.
%Thus
% $B$ is the set of divisors which bound $X_n$, and $\omega_A$ is an arbitrary algebraic form which has no poles along $B$.
The key  is then to apply  a version of  Stokes' theorem  to the closed polytope $\overline{X}_n\subset \overline{\Mod}_{0,n}(\R)$. Since
each face of $\overline{X}_n$ is itself  a product of associahedra $\overline{X}_a\times \overline{X}_b$,  we repeatedly
 take primitives  to obtain
 a cascade of integrals over associahedra
 of smaller and smaller dimension. In order to do this, we need to  construct a graded algebra $L(\Mod_{0,n})$ of multiple polylogarithm
 functions
   on $\Mod_{0,n}$ in which  primitives exist. At each stage of the induction, the dimension of the domain
  of integration decreases by one, and the weight of the integrand increases by one. At the final stage, we evaluate a multiple
  polylogarithm at the point $1$, and this gives a linear combination of multiple zeta values.
This gives an effective algorithm for computing such integrals. Our approach also  works in  greater generality, and our results should extend without
difficulty, for example, to the case of configuration spaces related to other Coxeter groups.

\subsection{General overview}
This paper is essentially  a  study of  the de Rham theory of the
motivic fundamental group of $\Mod_{0,n}$. Previously, the focus
has mainly been on the projective line minus roots of unity and
its products: in particular, $\Mod_{0,4} \cong
\Pro^1\backslash\{0,1,\infty\}$ and $\Mod_{0,4}^\ell$ ([De1],
[D-G], [Go1-2], [Ra]). The advantage of considering the moduli
spaces $\Mod_{0,n}$ is that we can bring to bear the full richness
of their geometry.  We show, for example, that the double shuffle
relations for multiple zeta values are just two special cases of
generalised product relations arising naturally from  functorial
maps between moduli spaces.
 An essential part of this work is devoted to     multiple polylogarithms, which are functions
first defined by Goncharov for
all $n_1,\ldots, n_\ell\in \N$ by the power series:
\begin{equation} \label{intromultpolysdef}\Li_{n_1,\ldots, n_\ell} (x_1,\ldots, x_\ell) =
\sum_{0<k_1<\ldots<k_\ell} {x_1^{k_1}\ldots x_\ell^{k_\ell} \over k_1^{n_1} \ldots k_\ell^{n_\ell}}\ ,\qquad \hbox{where  } \quad |x_i|<1\ .\end{equation}
By analytic continuation, they define multi-valued functions on $\Mod_{0,n}$, where $n=\ell+3$.
One of our main objects of study in this paper  is the larger set  $L(\Mod_{0,n})$ of all homotopy-invariant iterated integrals on $\Mod_{0,n}$.
It forms a differential algebra of
multi-valued functions on $\Mod_{0,n}$, in which the set of functions $(\ref{intromultpolysdef})$ is strictly contained.
From the point
of view of differential Galois theory, $L(\Mod_{0,n})$  defines   a  maximal unipotent Picard-Vessiot theory on $\Mod_{0,n}$.
We then define the universal algebra of multiple polylogarithms $B(\Mod_{0,n})$ to be a modified version of Chen's reduced bar construction.
It is a differential graded Hopf algebra which is an
abstract algebraic version of $L(\Mod_{0,n})$.
%which is given
% by a variant of Chen's reduced bar construction $B(\Mod_{0,n})$.
%It is a differential graded Hopf algebra.
 One of our key results  states that  the de Rham cohomology
of  $B(\Mod_{0,n})$ is trivial. From this we deduce the existence of primitives in $L(\Mod_{0,n})$. We also need to understand the
regularised  restriction
of polylogarithms to the faces of $\overline{X}_n$. This requires a canonical regularisation theorem, and amounts to studying what happens
when singularities of an iterated integral collide. We are thus led to work on  certain   blow-ups of $\Mod_{0,n}$, described below.  It follows that
the structure
 of $L(\Mod_{0,n})$,  and hence the function theory of multiple polylogarithms,  is  intimately related to the combinatorics
of the associahedron.%, and the geometry of $\Mod_{0,n}$.
%The appearance of the multiple zeta values comes from a monodromy computation in which they are coefficients
%of an associator.

\subsection{Detailed summary of results}
In section 2 we review some aspects  of the geometry of the moduli
spaces $\Mod_{0,n}$, and study certain blow-ups obtained from
them. Let $S$ denote a set with $n$ elements, each labelling a
marked point on the projective line $\Pro^1$, and write
$\Mod_{0,S}=\Mod_{0,n}$. A \emph{dihedral structure} on $S$ is an
identification of $S$ with the set of edges (or vertices) of an
unoriented $n$-gon. For each such  dihedral structure $\delta$, we
embed $\Mod_{0,S}$ in the affine space $\A^\ell$, where
$\ell=n-3$, and blow up
 parts of the boundary in $\A^\ell\backslash \Mod_{0,S}$ to obtain an intermediary space
 $$\Mod_{0,S} \subset \Modf_{0,S} \subset \overline{\Mod}_{0,S}\ ,$$
 where $\Modf_{0,S}$ is an affine scheme defined over $\Z$.
We prove that the set of $\Modf_{0,S}$, for varying $\delta$,  form a set of smooth affine charts on $\overline{\Mod}_{0,S}$. In order to define them, we introduce
 \emph{dihedral coordinates}, which are one of the key %technical
 tools used throughout this paper. These are functions
 $$u_{ij} : \Mod_{0,S} \rightarrow \Pro^1\backslash\{0,1,\infty\}\ , \quad \hbox{ where } \quad  \{i,j\} \in \chi_{S,\delta}\ ,$$
  indexed by the set of chords $\chi_{S,\delta}$ in the $n$-gon defined by $\delta$. Together, they define an embedding
$(u_{ij})_{\chi_{S,\delta}}:\Mod_{0,S} \rightarrow \A^{n(n-3)/2}$, and the scheme $\Modf_{0,S}$ is
the Zariski closure of the image of this map.
   %See figure 1 for the case $n=5$.
  For example, in the case $n=5$, we
 can identify $\Mod_{0,S} = \{(t_1,t_2)\in \Pro^1 \times \Pro^1 : t_1t_2(1-t_1)(1-t_2)(t_1-t_2)\neq 0, \, t_1,t_2\neq \infty\}.$ The pentagon $(S,\delta)$
 has five chords, labelled $\{13,24,35,41,52\}$ (fig. 1), and  we have
$$
u_{13}= 1-t_1,\quad u_{24}={t_1\over t_2},\quad u_{35}={t_2-t_1\over t_2(1-t_1)}, \quad
u_{41}={1-t_2\over 1-t_1}, \quad u_{52}=t_2\ .$$

\begin{figure}[h!]
  \begin{center}
%    \leavevmode
    \epsfxsize=13.0cm \epsfbox{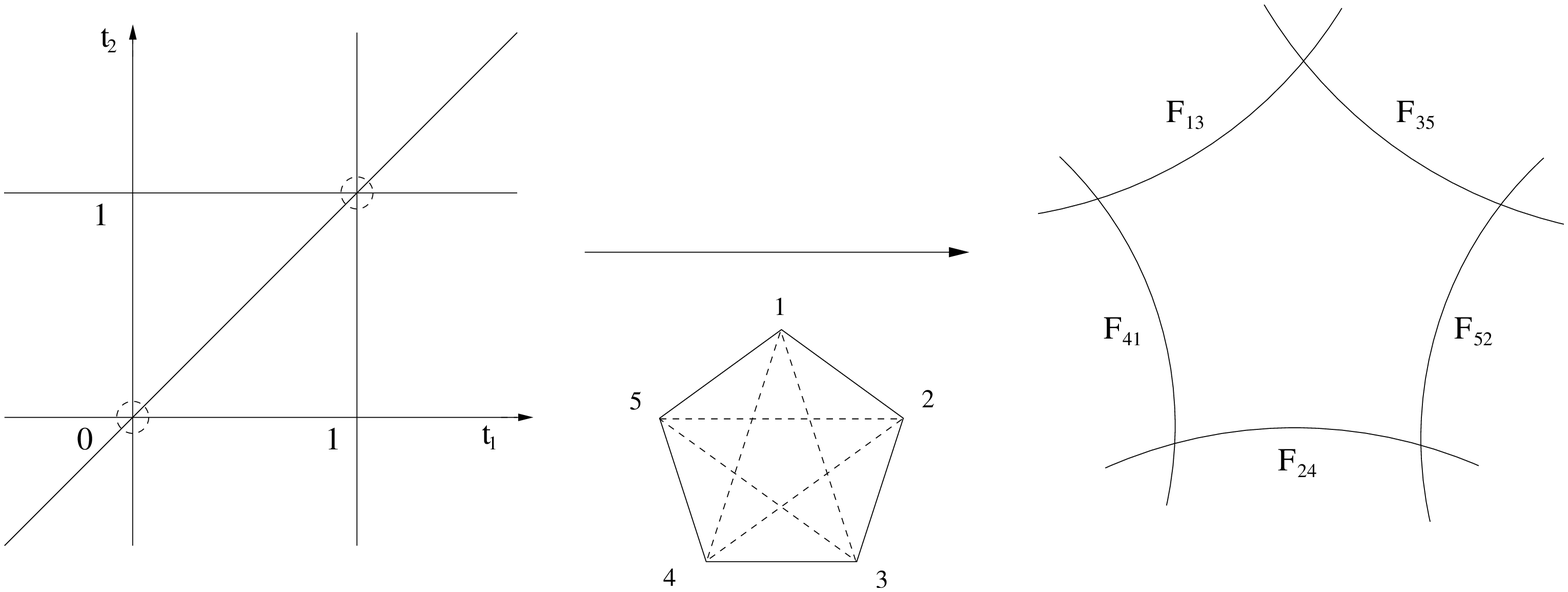}
\put(-252,115){\small{$(t_1,t_2) \mapsto (u_{13}, u_{24}, u_{35}, u_{41}, u_{52})$}}
  \label{pentagons}
  \caption{Dihedral coordinates on $\Mod_{0,5}$. The scheme $\Modf_{0,5}$ (right) is defined to be the Zariski closure of the image of the embedding
 $\{u_{ij}\}: \Mod_{0,5} \hookrightarrow \A^5$ defined by the set of dihedral coordinates, which are indexed by chords in a pentagon (middle). This map has the effect of blowing up
the points $(0,0)$ and $(1,1)$. A cell $X_{S,\delta}$ is given by the region
$0< t_1< t_2<1$ (left). After blowing-up it becomes a pentagon with sides $F_{ij}=\{u_{ij}=0\}$.}
  \end{center}
\end{figure}

\noindent
Now consider the set of real points $\Mod_{0,S}(\R)$. There is a   bounded cell
$X_{S,\delta} \subset \Mod_{0,S}(\R)$
defined by the region $\{0<u_{ij}<1\}$. One shows that $\Mod_{0,S}(\R)$ is the disjoint
 union of the open  cells $X_{S,\delta}$ of dimension $\ell=n-3$, as $\delta$ runs over the set of dihedral structures on $S$, so
a dihedral structure corresponds to choosing a connected component of $\Mod_{0,S}(\R)$.
The closure of the  cell $\overline{X}_{S,\delta}$ satisfies
\begin{equation}\label{introclosedStasheff}
\overline{X}_{S,\delta}=\{0\leq u_{ij}\leq 1\} \subset \Modf_{0,S}(\R)\ ,\end{equation}
and  $\Modf_{0,S}\backslash \Mod_{0,S}$ is the union of all divisors meeting the boundary of $X_{S,\delta}$.
%obtained by blowing up divisors
%which bound it.
 Therefore $\overline{X}_{S,\delta}$ is a convex polytope, and its boundary divisors give an explicit algebraic model of
the  associahedron.  It is well-known that the combinatorics of the associahedron is given by triangulations of polygons. But because dihedral coordinates are already defined in terms of polygons,
the main combinatorial properties
of the associahedron, and its dihedral symmetry, follow immediately from properties of the coordinates $u_{ij}$. In particular, the face $F_{ij}=\{u_{ij}=0\}$
of $\overline{X}_{S,\delta}$ is a product
\begin{equation} \label{introfacedecomp}
 F_{ij} \cong \overline{X}_{T_1,\delta_1} \times \overline{X}_{T_2,\delta_2}\ ,
\end{equation}
where $(T_1,\delta_1), (T_2,\delta_2)$ are two smaller polygons obtained by cutting the $n$-gon $S$ along the chord $\{i,j\}$
(fig. 3, $\S\ref{sect22}$).
%
%Now the polytope $X_{S,\delta}$ is a manifold with corners.
In this way, a vertex $v $ of $\overline{X}_{S,\delta}$
corresponds to a complete triangulation $\alpha$  of the $n$-gon by $\ell$ chords.   We  also introduce explicit
 \emph{vertex coordinates}
  $x_1^\alpha,\ldots, x_\ell^\alpha$  which are
 a certain subset of the set of all dihedral coordinates. These form a system of normal parameters
in the
neighbourhood of  the vertex $v\in \overline{X}_{S,\delta}$ corresponding to $\alpha$, such that  $\Mod_{0,S}(\R)\subset\overline{\Mod}_{0,S}(\R)$ is locally the complement of the normal crossing
divisor $x^\alpha_1\ldots x^\alpha_\ell=0$.
 These systems of coordinates (in the sense of differential geometry) are precisely what is needed for solving differential equations on $\Mod_{0,S}$ and regularising logarithmic singularities of
multiple polylogarithms.

In section 3 we define an abstract algebra of iterated integrals on $\Mod_{0,S}$ using a variant of Chen's reduced bar construction.
Since this construction exists in far greater generality,
we consider the complement of an arbitrary affine hyperplane arrangement defined over a field $k$ of characteristic $0$. So let
$$M= \A^\ell \backslash \bigcup_{i=1}^N H_i\ ,$$
where $H_1,\ldots, H_N$ is any set of hyperplanes in $\A^\ell$. Let $t_1,\ldots, t_\ell$ denote coordinates
on $\A^\ell$, and  let $\Or_M$ denote the ring of regular functions
on $M$. It  is a differential $k$-algebra with respect to the coordinate derivations $\partial/\partial t_i$, for $1\leq i\leq \ell$.
Let
$$\omega_i = {d\alpha_i \over \alpha_i}\ , \quad \hbox{ for } 1\leq i\leq N\ , $$
denote the logarithmic $1$-form corresponding to $H_i$, where $\alpha_i$ is a defining equation for $H_i$.
The version of the bar construction $B(M)$  we consider %is similar to the zero$^{th}$
%cohomology group of Chen's bar complex. To
is defined as follows. Let $V_m(M)$ denote
 the  $k$-vector space  generated by %finite
 linear combinations of symbols
\begin{equation}\label{introgenforbm}
 \sum_{I=(i_1,\ldots, i_m)} c_I [\omega_{i_1}|\ldots |\omega_{i_m}]\ ,\qquad c_I \in k\ ,
\end{equation}
which satisfy the integrability condition:
\begin{equation} \label{introintcond}
\sum_I c_I \,\omega_{i_1}\otimes \ldots \otimes (\omega_{i_j}\wedge \omega_{i_{j+1}}) \otimes \ldots \otimes \omega_{i_m}=0 \qquad \hbox{ for all }
1\leq j<m \ .
\end{equation}
We then set $B(M) = \Or_M \otimes_k \bigoplus_{m\geq 0} V_m(M)$, where $V_0(M)=k$. This is a graded Hopf algebra over $\Or_M$ which is similar to the zero$^{th}$ cohomology group of the  bar complex studied by Chen [Ch1], except that it consists
of 1-forms only (see also [Ha1]).
Using the 1-part of the coproduct on $B(M)$, we define the action of $\ell$ commuting derivations $\partial_i$ on $B(M)$,
and show that $(B(M), \partial_i)$ defines a differential extension of $(\Or_M, \partial/\partial t_i)$.
The possibility of using iterated integrals to construct a Picard-Vessiot theory on manifolds was first suggested
by Chen [Ch2].

%, which we denote $B(\Mod_{0,S})$. This is a differential graded Hopf algebra.
\begin{thm} $B(M)$ is an infinite unipotent Picard-Vessiot extension of $\Or_M$. In other words, it has no non-trivial
differential ideals and its ring of constants is $k$. It is therefore a polynomial algebra. Furthermore, $B(M)$   contains 1-primitives:
$$H^1_\DR(B(M))=0\ .$$
\end{thm}

\noindent
It follows that every unipotent extension of $B(M)$ is trivial, and it is the smallest extension of $\Or_M$ with
this property. Equivalently, $B(M)$ is the limit
$$B(M) = \lim_{\rightarrow}\, U \ , $$
where $U$ ranges over all unipotent extensions $U$ of $\Or_M$. In
this sense it is universal, and  it follows that its differential
Galois group is  a pro-unipotent  group. Now if we identify
$\Mod_{0,S}$ with the affine hyperplane configuration
$$\Mod_{0,S}\cong \{(t_1,\ldots, t_\ell)\in \A^\ell: \quad t_i\neq 0,1\ ,  \quad t_i-t_j \neq 0\}\ ,$$
then we define the \emph{universal algebra of polylogarithms} on $\Mod_{0,S}$ to be $B(\Mod_{0,S})$. In general, it is difficult
to construct words $(\ref{introgenforbm})$ satisfying the integrability condition $(\ref{introintcond})$ since they rapidly
become very complicated as the weight increases. In order
 to overcome
this problem, we consider  two affine hyperplane arrangements, one of which fibers linearly over the other. Therefore, let $M\subset \A^\ell$ and
$M'\subset \A^{\ell-1}$ denote two affine arrangements, and consider a linear projection
$$\pi: M \rightarrow M'$$
with constant fibres $F$, where $F$ is  the affine line $\A^1$
minus a number of marked points. We then prove that  there is a
tensor product decomposition
$$B(M) \cong B(M') \otimes_{\Or_{M'}} B_{M'}(F)\ ,$$
where $B_{M'}(F)$ is a free shuffle algebra which  can be
described explicitly. In the case of moduli spaces $\Mod_{0,S}$,
we
 apply this argument to
the fibration map:
$$ \Mod_{0,n} \To \Mod_{0,n-1}$$
and use induction to deduce that $B(\Mod_{0,S})$ is a tensor
product of free shuffle algebras. As a result, one can  write down
a basis for $B(\Mod_{0,S})$, and one deduces  that the higher
cohomology groups of $B(\Mod_{0,S})$ vanish.

\begin{thm} The de Rham cohomology of $B(\Mod_{0,S})$ is trivial:
$$H^i_\DR (B(\Mod_{0,S})) = 0 \qquad \hbox{ for all }\quad  i \geq 1\ .$$
\end{thm}

\noindent A similar result holds
 for
any hyperplane arrangement of fiber type, {\it i.e.}, one which
can be obtained as a sequence of such fibrations. In an appendix
we also prove that
 $H^i_\DR(B(M))$, for $i\geq1$, vanishes for  all  arrangements $M$ which have quadratic cohomology. The proofs  only
use simple arguments of differential algebra. Theorem $1.3$ holds
 because $\Mod_{0,S}$ is a $K(\pi,1)$-rational space. An equivalent theorem is due to Hain and MacPherson ([H-M], [Kh2]).

Given any point $z_0\in \Mod_{0,S}(\C)$ we define a realisation
\begin{eqnarray} \label{introreal}
\rho_{z_0}: B(\Mod_{0,S}) &\overset{\sim}{\To}& L_{z_0}(\Mod_{0,S}) \\
\sum_I f_I\, [\omega_{i_1}|\ldots |\omega_{i_m}] & \mapsto & \sum_I f_I\, \int_{z_0}^z \omega_{i_m} \ldots \omega_{i_1}\ , \nonumber
\end{eqnarray}
given by iterated integration along any path
$\gamma:[0,1]\rightarrow \Mod_{0,S}(\C)$ which begins at $z_0$ and
ends at a variable point $z\in \Mod_{0,S}(\C)$. The integrability
condition $(\ref{introintcond})$ ensures that the iterated
integral $(\ref{introreal})$ only depends on the homotopy class of
$\gamma$. It therefore defines a multi-valued function of the
parameter $z$, {\it i.e.}, a holomorphic function on the universal
covering space of $\Mod_{0,S}(\C)$. Here,  $L_{z_0}(\Mod_{0,S})$
is a differential graded algebra of multi-valued functions on
$\Mod_{0,S}$.
 We deduce from the previous theorem that $\ell$-forms with coefficients
in $L_{z_0}(\Mod_{0,S})$ have primitives in $L_{z_0}(\Mod_{0,S})$.
%which is the key result for applying Stokes' formula.

The realisation $\rho_{z_0}$ is not quite good enough, however. We
actually need a realisation $\rho_{z_0}: B(\Mod_{0,S}) \rightarrow
L_{z_0}(\Mod_{0,S})$, where the base point $z_0$ does not lie in
$\Mod_{0,S}(\C)$. The point  $z_0$ can be replaced with
 a tangential base point in the sense of
  [De1], but our approach consists  instead of  viewing $z_0$ as the corner of a manifold with corners. This gives rise to divergent integrals,
and to deal with this requires a regularisation procedure. The
best approach is to consider the generating series of all such
iterated integrals, and regularise them  all simultaneously. Such
a generating series satisfies a formal differential equation, and
to solve it requires a generalised Fuchs' theorem in several
variables in the unipotent case. For want of a suitable reference,
we develop the necessary theory from scratch in section 4. We also
study the regularisation of logarithmic singularities along the
boundary of any manifold with corners. Section 5 is devoted to a
detailed study of the case of a one dimensional arrangement
$\Pro^1\backslash \{\sigma_0,\ldots, \sigma_N, \infty\}$. In this
case, the bar construction can be written down explicitly (it is a
free shuffle algebra), and the corresponding iterated integrals
are known as hyperlogarithms, which go back to  Poincar\'e and
Lappo-Danilevsky.

In section 6, we apply all the results developed previously to the
case of the moduli spaces $\Mod_{0,S}$ to obtain the necessary
regularisation results. The generating series of multiple
polylogarithms can be described as follows. To each dihedral
coordinate, or chord, is associated a logarithmic one-form
$$\omega_{ij} = d \log u_{ij}\ ,\qquad \hbox { for } \quad  \{i,j\} \in \chi_{S,\delta}\ .$$
It is symmetric in $i$ and $j$. Let $\delta_{ij}$, for $\{i,j\}\in \chi_{S,\delta}$, denote a set of symbols
 satisfying $\delta_{ij}=\delta_{ji}$, and consider the formal 1-form
\begin{equation} \label{introOmegadefn}
\Omega_{S,\delta} = \sum_{\{i,j\} \in \chi_{S,\delta}}
\delta_{ij}\, \omega_{ij} \ . \end{equation} This is a homogeneous
version of the Knizhnik-Zamolodchikov form [Dr, Ka, K-Z]. The
integrability of $\Omega_{S,\delta}$ is equivalent to certain
quadratic relations in the $\delta_{ij}$, which we call the
dihedral braid relations. In the special case $\Mod_{0,5}$, these
reduce to the relations:
$$
[\delta_{ij}, \delta_{kl} ] = 0\ ,$$
for any pair of  chords  $\{i,j\}, \{k,l\}$ which  do not cross, and the pentagonal relation
$$[\delta_{13},\delta_{24}] + [\delta_{24}, \delta_{35}] +[\delta_{35}, \delta_{41}]+[\delta_{41}, \delta_{52}] +[\delta_{52}, \delta_{13}]=0\ .$$
Let us fix a dihedral structure $\delta$ on $S$, and let
$\widehat{\DB}_{S,\delta}(\C)$ denote the ring of non-commutative
formal power series in the symbols $\delta_{ij}$ with coefficients
in $\C$, modulo the dihedral braid relations. Then
%If the dihedral braid relations hold, then
 we can consider the formal differential equation
\begin{equation} \label{intromainequation} d L = \Omega_{S,\delta} L\ ,\end{equation}
where $L$ takes values in $\widehat{\DB}_{S,\delta}(\C)$.

\begin{thm} Let $v$ denote a vertex  of the associahedron $\overline{X}_{S,\delta}$,
and let $F_v$ denote the set of faces meeting $v$. Then
there is a unique solution $L_{v,\delta}$ of $(\ref{intromainequation})$
such that
$$L_{v,\delta}(z) = f_{v,\delta}(z)\, \exp ( \sum_{\{i,j\}: F_{ij}\in F_v} \delta_{ij} \log u_{ij} ) \ ,$$
where  $f_{v,\delta}(z)$ is holomorphic in  a  neighbourhood of
$v\in \Modf_{0,S}$, and $f_{v,\delta}(v)=1$.
\end{thm}

In other words, the function $L_{v,\delta}(z)$ is holomorphic on
an open set of $\Modf_{0,S}(\C)$ which contains the  real cell
$X_{S,\delta}$, and has explicitly given  monodromy around each
face $u_{ij}=0$ of the associahedron $\overline{X}_{S,\delta}$
which meets the vertex  $v$. The differential equation
$(\ref{intromainequation})$ is closely related to the
Knizhnik-Zamolodchikov equation. Solutions to the latter equation
are usually constructed by induction using fibration maps between
configuration spaces. The previous theorem, however, is proved
directly using the generalised Fuchs' theorem developed in section
$4$. This approach has many advantages: firstly, there are no
coherence conditions to verify; secondly, we obtain a direct
geometric interpretation of Drinfeld's asymptotic zones, which
were studied by Kapranov;
 and thirdly, the functoriality of the solution $L_{v,\delta}(z)$ with respect to maps between moduli spaces follows automatically.
As a result, we obtain a direct definition of an associator on
$\Mod_{0,S}$
 by considering the quotient of two different solutions:
$$Z^{v,v'} = (L_{v,\delta}(z))^{-1} L_{v',\delta}(z)\, \in \widehat{\DB}_{S,\delta}(\C)\ .$$
Here, $z$ is any point in an open neighbourhood of $X_{S,\delta}$ in $\Modf_{0,S}$. The quotient is necessarily constant.
 The main properties of Drinfeld's associator can be derived immediately.
Using the previous theorem, we  deduce an expression for the
monodromy of $L_{v,\delta}(z)$ and its regularisation in terms of
the series $Z^{v,v'}$ ($\S\ref{sect65}$). Then, using explicit
expressions for hyperlogarithms, we deduce the following result,
which was first proved by Le and Murakami, following Kontsevich.

\begin{thm} The coefficients of the series $Z^{v,v'}$ are multiple zeta values.
\end{thm}
\noindent It follows that  the holonomy of the moduli spaces
$\Mod_{0,S}$ can be expressed using  multiple zeta values and the
constant $2\pi i$. Now define $L^{v,\delta}(\Mod_{0,S})$ to be the
differential algebra generated by the coefficients of the series
$L_{v,\delta}(z)$. We can then define  the sought-after
realisation $\rho_{v,\delta}$ which is regularised at the vertex
$v$ of $\overline{X}_{S,\delta}$:
$$\rho_{v,\delta}: B(\Mod_{0,S})\overset{\sim}{\To} L^{v,\delta}(\Mod_{0,S})\ ,$$
and which is defined over the field $k=\Q$.
From this we deduce the main regularisation theorem, which describes   the regularised restriction of a multiple polylogarithm
to the face of the associahedron in terms of multiple zeta values.

\begin{thm} Let $F_{ij}$ denote a face of $\overline{X}_{S,\delta}$ isomorphic to a product $\overline{X}_{T_1,\delta_1}\times \overline{X}_{T_2,\delta_2}$ as in
$(\ref{introfacedecomp})$ above. Then if the vertex $v$ corresponds to the pair $(v_1,v_2)$,
 $$  \Reg \,( L^{v,\delta}(\Mod_{0,S}), F_{ij}) \otimes_\Q\MZV \cong L^{v_1,\delta_1} (\Mod_{0,T_1}) \otimes_{\Q} L^{v_2,\delta_2} (\Mod_{0,T_2})
  \otimes_\Q \MZV\ .$$
\end{thm}
\noindent
In other words, the regularisation of multiple polylogarithms along   divisors at infinity  is completely determined by the combinatorics
of the %corresponding
 associahedron.

In section 7, we study  period integrals on $\Mod_{0,S}(\R)$ in
terms of dihedral coordinates.% and compute their singularities.
 We
first show that, up to multiplication by a rational number, there
is a unique algebraic $\ell$-form $\omega_{S,\delta}$, which has
neither zeros nor poles on $\Modf_{0,S}$. %the closed polytope $\overline{X}_{S,\delta}$.
 This form is invariant under the
natural action of the dihedral group. We deduce that one can write
an arbitrary integral $(\ref{introIdefn})$ as a linear combination
of integrals
\begin{equation} \label{introdihedralintegral}
I_{S,\delta}(\alpha_{ij})=\int_{\overline{X}_{S,\delta}} \prod_{\{i,j\}\in\chi_{S,\delta} } u^{\alpha_{ij}}_{ij}  \, \omega_{S,\delta}\ ,
\end{equation}
for some fixed dihedral structure $\delta$, where the indices $\alpha_{ij} \in \Z$. Such an integral converges if and only if the coefficients
$\alpha_{ij}$ are all non-negative. In explicit coordinates, $(\ref{introdihedralintegral})$  can be written as a generalized Selberg integral
$$ I_{S,\delta}(\alpha_{ij}) =\int_{ [0,1]^\ell} \prod_{i=1}^\ell
x_i^{a_i} (1-x_i)^{b_i} \prod_{i<j} (1-x_ix_{i+1}\ldots x_j
)^{c_{ij}}\, dx_1\ldots dx_\ell \ .$$ Particular subfamilies  of
these kinds of integrals have been considered by various authors
in connection with the diophantine approximation of zeta values
(see, {\it e.g.}, [Fi2, Zl, Zu]). Terasoma has also computed the
Taylor expansions (with respect to the exponents)  of  certain
families of such integrals, and proved they are multiple zeta
values [Ter]. The advantage of the blown-up integral
representation $(\ref{introdihedralintegral})$ is that all poles
of the integrand have been pushed to infinity, which allows an
algebraic interpretation of the integrals as periods, and a
systematic procedure for computing them, which is detailed in
section $8$ and summarised below. As a further application of
dihedral coordinates, we give an explicit formula for the order of
vanishing of any form
$$\prod_{\{i,j\}\in\chi_{S,\delta} } u^{\alpha_{ij}}_{ij}  \, \omega_{S,\delta}\ ,$$
along the divisors at infinity in $\overline{\Mod}_{0,S}$. Using
this formula we retrieve a result, due to Goncharov and Manin,
which gives the singular locus of a certain family of forms
% considered by Kontsevich
which correspond directly to  multiple zeta values.
Our method exploits the action of the symmetric group on
$\Mod_{0,S}$, and completely avoids the delicate calculation of
blow-ups and the cancellation of singularities studied in [G-M].
In $\S\ref{sect75}$, we  show how functorial maps
$$f: \Mod_{0,S} \To \Mod_{0,T_1}\times \Mod_{0,T_2}\ ,$$
where $T_1$ and $T_2$ satisfy certain conditions (\S $\ref{sect27}$),
give rise to generalised product formulae between multiple zeta values. More precisely,
given any such map $f$, there is a set  of dihedral structures $G_f$ on $S$ such that the following formula holds:
\begin{equation}\label{introproductmap}
 \int_{X_{T_1,\delta_1}} \omega_1 \times \int_{X_{T_2,\delta_2}} \omega_2 = \sum_{\gamma \in G_f} \int_{X_{S,\gamma}} f^* (\omega_1 \otimes \omega_2)\ .
\end{equation}
This expresses a product of periods as a $\Q$-linear combination
of  other periods.
 We compute two explicit examples of such maps $f$; one where $G_f$ is as large as possible, and the other when $G_f$ reduces
to a single element. In the first case, $G_f$ is the set of $(p,q)$ shuffles where $p=\dim \Mod_{0,T_1}$ and $q=\dim \Mod_{0,T_2}$,
 and $(\ref{introproductmap})$ gives rise to the shuffle product for multiple zeta values.
In the second case, we show that $(\ref{introproductmap})$, on
applying an identity due to Cartier, gives rise to the stuffle
relations for multiple zeta values. Thus both shuffle and stuffle
relations can be regarded as two extreme cases of generalised
product relations of geometric origin on moduli spaces.

The above results are put together in section 8, where  we give a proof of theorem $1.1$ using  Stokes' formula  as described
above. We summarise the main points of the argument here. The regularisation results of section $6$ provide the existence of a graded algebra
of multi-valued functions $L(\Mod_{0,S})$ with the following properties:
\begin{enumerate}
  \item The graded part of weight $0$ of $L(\Mod_{0,S})$ consists of all regular algebraic functions on $\Mod_{0,S}$ with coefficients in
  $\Q$.
  \item Primitives of $\ell$-forms exist in $L(\Mod_{0,S})$, and increase the weight by one.
  \item The restriction of a function  $f\in L(\Mod_{0,S})$ to a face of  $\overline{X}_{S,\delta}$
 is a product of multiple zeta values with functions in
  $L(\Mod_{0,T_1}) L(\Mod_{0,T_2})$. %in such a way that the total weight is preserved.
\end{enumerate}
The argument for computing the period integrals is then by an inductive application of Stokes' theorem over the associahedron $\overline{X}_{S,\delta}$.
At each stage, we must compute
$$I=\int_{\overline{X}_{S,\delta}} f\, \omega_{S,\delta} \ ,$$
where $f\in L(\Mod_{0,S})$ is a function  which is allowed
logarithmic singularities along the boundary $\partial
\overline{X}_{S,\delta}$, but which has no polar singularities.
Such an integral necessarily converges, and it follows from
property (2)  that there exists a primitive $F$ with coefficients
in $ L(\Mod_{0,S})$ such that $dF=f$. However, such primitives are
not unique, and we may inadvertently have introduced extra poles.
 We show, however, that there exists a primitive $F$ with no poles along
$\overline{X}_{S,\delta}$, and it then follows that $F$ extends
continuously to the boundary $\partial \overline{X}_{S,\delta}$.
The essential remark is that the one-form
$$\log x \, dx \qquad \hbox{ where } x\geq 0\ ,$$
has a logarithmic singularity at the point $x=0$, but that its primitive $x\log x-x$ extends continuously to the point $0$. We can
therefore restrict the primitive $F$ to the faces of the associahedron  by  property $(3)$, and proceed by induction using Stokes' formula
and  $(\ref{introfacedecomp})$ without any further difficulty.
 In $\S8.5$ we show how
the same strategy can be used to compute all relative periods of
moduli spaces $\Mod_{0,S}$, and finish with some  simple examples
in $\S8.6$. The paper is completely self-contained, apart from
some  properties of iterated integrals which are  very clearly
presented in [Ha1], and some remarks on framed motives in
$\S\ref{sectMHS}$.

We expect that the ideas and methods introduced in this paper
should have applications in the following situations. First of
all, one can consider more general hyperplane configurations
associated to other root systems or Coxeter groups, and consider
the corresponding polylogarithm algebras, periods and associators.
Notably, one can introduce $N^{\mathrm{th}}$ roots of unity to
obtain a tower of spaces over $\Pro^1\backslash\{0, e^{2i\pi
k/N},\infty\}$ which are finite covers of $\Mod_{0,S}$ and
construct a similar theory giving a higher dimensional version of
[Ra, D-G]. Furthermore, in perturbative quantum field theory, it
is generally believed that certain period integrals one derives
from a large class of Feynman diagrams should give multiple zeta
values. After blowing up, these are integrals of rational
algebraic forms over  algebraic convex polytopes. It would be very
interesting to try to apply the methods of this paper to such
integrals.
\\

This paper was written during  my doctoral thesis at the
university of Bordeaux. I am very grateful to Richard Hain for his
many detailed comments regarding an earlier version of this
manuscript, and especially to Pierre Cartier, without whose many
suggestions, good humour, and continuous encouragement, this paper
would not have reached its present form.

% Richard Hain for his many detailed
%comments regarding an earlier version of this manuscript. I also
%wish to
%thank Pierre Cartier f

%Here- state stuff about defining partial compactifaction,
%introduce Deligne-mumford-knudsen..

 %mo5 computations/examples here?

%Subtlety of $d^2$ for bar construction.

%Emphasize UNITY: combinatorics of stasheff polytopes -
%differential equations/multiple polylogs - Euler integrals.
%\\

%Ring of periods coincides with the holonomy ring.
%Application to generalised products on mzvs. Methods work for general Hyperplane arrangements - should generalise easily (see below).
%Definitive technique - Previous approaches are ad-hoc.

%K%ontsevich Zagier -  We have at our disposition: Stokes,
%Symmetries (dihedral), extra birational bastards, product
%structures.
%\\

%Idea is simple : apply stokes theorem - spaces of decreasing
%dimension, complexity of the function (weight) goes up by one.
%Needs algebra with three properties: existence of primitives,
%regularisation, contains algebraic forms.
%\\

%Works for general hyperplanes??
%\\

%Taylor expansions.
%\\
%theorem 1.11 which states that whenever an affine hyperplane arrangement fibres
%over another one, the corresponding reduced bar constructions
%decomposes as a tensor product.
%

%Mention Rhin and Viola.
\newpage

\newpage
\section{Dihedral coordinates on $\overline{\Mod}_{0,n}(\R)$.}

\subsection{}\label{sect21}
Let  $n\geq 4$, and  let $S$ denote a set with $n$ elements. Let $\Mod_{0,S}$
 denote the moduli
space of Riemann spheres with $n$ points labelled with elements of $S$. If
$(\mathbb{P}^1)_*^S$ denotes the set of all $n$-tuples of distinct
points $z_s\in \mathbb{P}^1$,  for $s\in S$, then
$$\mathfrak{M}_{0,S}
=\mbox{PSL}_2 \backslash (\mathbb{P}^1)_*^S\ ,$$
 where $\mbox{PSL}_2$ is the algebraic group of automorphisms of $\Pro^1$ and acts  by
M\"{o}bius transformations. The quotient $\mathfrak{M}_{0,S}$ is an
affine variety of dimension $\ell= n-3$.
 A point in $\Mod_{0,S}(\C)$ is therefore an injective map $S\hookrightarrow \Pro^1(\C)$
considered up to the
action of $\PSL_2(\C)$. If $S=\{s_1,\ldots, s_n\}$, then  we  frequently write $i$ instead of  $s_i$, and
%in that case
 denote $\Mod_{0,S}$ by $\Mod_{0,n}$.

 We wish to write down the set
of regular functions on $\mathfrak{M}_{0,S}$, or, equivalently,
the set of $\mbox{PSL}_2$-invariant regular functions
on $(\mathbb{P}^1)_*^n.$
%It is convenient to
%write $\Delta_{ij}= z_i-z_j$
 %for any pair of indices $1\leq i\neq j \leq
%n$. %, where $z_i$ are coordinates on $\mathbb{P}^1$.
 Let $i,j,k,l$ denote any distinct indices in $S$. Recall that
 the cross-ratio is
defined by the formula:
$$\big[ij\,|\,kl\big]={ (z_i-z_k)(z_j-z_{l}) \over (z_i -z_l)(z_j - z_k)}\ .$$
 The cross-ratios do not depend on the choice of coordinates $z_i$ and are %indeed
  $\mbox{PSL}_2$-invariant. We therefore have a set of
  maps
$[i j\,|\,k l]:\Mod_{0,S} \rightarrow \Mod_{0,4}\cong
\Pro^1\backslash\{0,1,\infty\}$.
   The symmetric group on four letters $\mathfrak{S}_4$ acts on
each cross-ratio via the group of anharmonic substitutions $
\langle z\mapsto 1-z, \,\,z\mapsto 1/z\rangle\cong
\mathfrak{S}_3\cong \Sym_4/V$, where $V$ is the Vierergruppe. We have:
\begin{equation}\label{S4}
\big[ij\,|\,kl\big] = 1- \big[ik\,|\,jl\big] \ ,  \quad \hbox{and}
\quad  \big[ij\,|\,lk\big] = \big[ij\,|\,k l
\big]^{-1}=\big[ji\,|\,kl\big]\ ,
\end{equation}
$$\hbox{ and  } \quad \big[ij\,|\,kl\big] = \big[kl\,|\,ij\big] =
\big[ji\,|\,lk\big]=\big[lk\,|\,ji\big]\ .$$
%from which it follows that $[\big[ij\,|\,kl\big]=\big[kl\, |\,ij\big]$.
For any five
distinct indices $i,j,k,l,m\in S$ there is also the multiplicative
relation:
\begin{equation} \label{mult}
\big[ij\,|\,k\,l\big] = \big[ij\,|\,k\, m\big] .
\big[ij\,|\,m\,l\big] \ .
\end{equation}

In order to make explicit computations, it will be convenient to
fix a system  of coordinates on $\Mod_{0,S}$  from the beginning.
This breaks the symmetry, so  we  assume here
that $S=\{1,\ldots,
n\}$. Since the action of $\PSL_2(\C)$ is triply transitive on $\Pro^1(\C)$, we
can place the coordinates $z_1$ at 1, $z_2$ at $\infty$, and $z_3$
at $0$. We define explicit \emph{simplicial coordinates}
$t_1,\ldots, t_\ell$ on $\Mod_{0,S}$ by setting
$$t_1=z_4\ ,\ \ldots \ ,  \ t_\ell=z_n\ .$$
This identifies $\Mod_{0,S}$ with the complement of the affine
hyperplane configuration:
\begin{equation} \label{ExpSimpCoords}
\Mod_{0,S} \cong \{(t_1,\ldots, t_\ell) \in \A^\ell: \, t_i\notin
\{0,1\},\quad  t_i\neq t_j \hbox{ for all } i \neq j\}\
.\end{equation}If we now perform the change of variables
\begin{equation}\label{changeofvars}
t_1=x_1\ldots x_\ell\ , \  t_2=x_2\ldots x_\ell\ , \  \ldots \ ,\
t_\ell =x_\ell\ ,\end{equation} then we can identify $\Mod_{0,S}$
with the open complement of hyperboloids:
\begin{equation} \label{ExpCubeCoords}
\Mod_{0,S} \cong \{(x_1,\ldots, x_\ell) \in \A^\ell: \, x_i\notin
\{0,1\},\quad  x_i\ldots x_j\neq 1 \hbox{ for all } i < j\}\
.\end{equation} The coordinates $x_1,\ldots, x_\ell$ will be
referred to as \emph{cubical coordinates} and are %particularly
well-suited to the study of polylogarithms on $\Mod_{0,S}$ (\S6).
Simplicial and cubical coordinates are two extremal cases of more
general systems of coordinates which we define in an invariant
manner in $\S\ref{sect27}$. We shall pass freely between the two
systems, especially when making comparisons with formulae existing
in the
literature. %Observe that
The change of coordinates $(\ref{changeofvars})$ %from $t_i$ to the $x_i$
 has the effect
of blowing up the origin; the boundary divisors at the origin in
$(\ref{ExpCubeCoords})$ cross normally, but do not in
$(\ref{ExpSimpCoords})$.
\subsection{Dihedral coordinates on $\Mod_{0,S}$}\label{sect22}
% In order to define dihedral coordinates on $\Mod_{0,S}$,
Let $S$
be a finite set with $n\geq 4$ elements.
\begin{defn}
 A \emph{cyclic structure} $\gamma$ on $S$ is  a cyclic ordering of the elements of $S$, or equivalently, an
 identification of the elements of $S$ with the edges %(or vertices)
of an oriented $n$-gon modulo  rotations.
 A \emph{dihedral structure} $\delta$ on $S$ is an identification with the edges %(or vertices)
  of an unoriented $n$-gon
modulo  dihedral symmetries.
\end{defn}

When we write $S=\{s_1,\ldots, s_n\}$,
it will carry the obvious dihedral structure unless stated otherwise.
%and we shall frequently write $k$ for the
%element $s_k$ of $S$.
 In this case,  the group of permutations $\Sym_S$
can be identified with
the symmetric group $\Sym_n$. The set of cyclic (resp. dihedral) structures  on $S$ is then indexed
by the set of cosets $\Sym_n/C_n$ (resp. $\Sym_n/D_{2n}$), where $C_n$ and $D_{2n}$ denote the cyclic and dihedral
groups of orders $n$ and $2n$ respectively.
%Since both points of view are equally
%useful,
We  will often represent a dihedral structure  as a  regular
$n$-gon $(S,\delta)$ with edges labelled $1,2,\ldots, n$ in order.
A number in parentheses $(i)$, where $i\in \Z/n\Z$, will denote
the pair of adjacent edges $\{i,i+1\}$.  We will represent this on
the $n$-gon $(S,\delta)$ by labelling the vertices with the
elements $(1),(2),\ldots, (n)$% in parentheses
; the convention is that the vertex labelled $(i)$ meets the edges
labelled $i$ and $i+1$ modulo $n$ (figures 2 and 3).

 %If we write $S=\{s_1,\ldots, s_n\}$, combinatorial $n$-gon with an element of $S$.

% We  define a
%\emph{cyclic structure} on $S$ to be the choice of a coset $\gamma
%\in \Sym_n/C_n$, where $C_n$ is the cyclic group of order $n$.
%This is equivalent to choosing a   labelling
%$$\gamma: S \overset{\sim}{\To} C_n\ ,$$
%where $\gamma$ is well-defined up to  cyclic permutations.
%Likewise, we define a \emph{dihedral structure} on $S$ to be the
%choice of a coset
%$$\delta\in \Sym_n/D_{2n}\ ,$$
%where $D_{2n}$ is the dihedral group of order $2n$.
% This is
%equivalent to labelling  each vertex, or each edge, of a %regular
%A cyclic structure
%on $S$ gives a planar embedding, or equivalently, an orientation
%of this $n$-gon.

Given a dihedral structure $\delta$ on $S$, we define coordinates
on $\Mod_{0,S}$ using a certain subset of the set of all
cross-ratios as follows.
 Let $\chi_{S,\delta}$ denote the set of all $n(n-3)/2$
unordered pairs $\{i,j\}$, $1\leq i,j\leq n$ such that
$i,j,i\!+\!1,j\!+\!1$ are distinct modulo $n$ ({\it i.e.}, $i, j$ are not
consecutive modulo $n$). Each element $\{i,j\} \in
\chi_{S,\delta}$ will be depicted as a chord joining the vertices
$i$ and $j$ in the regular $n$-gon (fig. 2). We set
\begin{equation}\label{udef}
u_{ij} =\big[i \,\, i\!+\!1 \,|\,  j\!+\!1 \,\, j \big]\quad
\hbox{ for each } \quad \{i,j\} \in \chi_{S,\delta}\ .
\end{equation}
 Using the definition
of the cross ratio, one can check that $u_{ij}$ is symmetric in
$i$ and $j$, and is therefore well-defined. The set of
cross-ratios $\{u_{ij}: \{i,j\} \in \chi_{S,\delta}\}$ only
depends on $\delta$. Consequently, we obtain a regular morphism
\begin{equation} \label{Embed}
\big(u_{ij}\big)_{\{i,j\}\in \chi_{S,\delta}} : \Mod_{0,S} \To
\Mod_{0,4}^{n(n-3)/2} \subset \A^{n(n-3)/2}\ .
\end{equation}
\begin{figure}[h!]
  \begin{center}
%    \leavevmode
    \epsfxsize=12.0cm \epsfbox{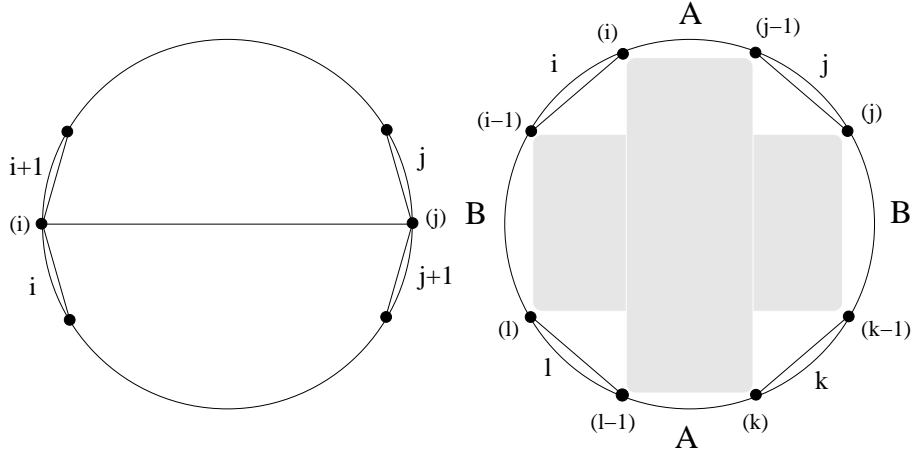}
  \label{Chord}
  \caption{Part of an oriented regular $n$-gon inscribed in a circle. Its edges are labelled with the elements of $S$, and its vertices
  are labelled with elements of $S$ in parentheses. Left - a chord $\{i,j\}\in \chi_{S,\delta}$ meets four edges $i,i+1,j,j+1$
  which define the dihedral coordinate $u_{ij}=[i\,\,i\!+\!1\,|\,j\!+\!1\,\, j]$. Changing the orientation of the $n$-gon  does not alter
  $u_{ij}$ by the last equation in $(\ref{S4})$. Right - a set of four edges $i,j,k,l$ breaks the $n$-gon into four regions as in lemma $\ref{lemmacrossratio}$, and defines
   a pair $A, B \subset \chi_{S,\delta}$ of completely crossing chords, depicted by the shaded rectangles (corollary $\ref{corcrossingchords}$).}
  \end{center}
\end{figure}

A simple calculation in simplicial coordinates shows that
\begin{equation} \label{tcoords}
t_1=u_{24}\ldots u_{2n} \ ,\quad \ldots \ ,\quad  t_{\ell-1} = u_{2\, n-1} u_{2\, n} \ ,\quad  t_\ell = u_{2n}\ ,
\end{equation}
$$
1-t_1=u_{13} \ ,\quad  1-t_2= u_{13}u_{14} \ ,\quad  \ldots \ , \quad 1-t_{\ell} = u_{13}\ldots u_{1\, n-1} \ .
$$
Likewise, the set of cubical coordinates $(x_1,\ldots, x_\ell) = (u_{24},\ldots, u_{2n})$ are completely determined by the functions
$u_{i j}$, and therefore
$(\ref{Embed})$ is  an embedding.
 It follows that every cross-ratio %$[ij|kl]$
 can be written in terms of the functions $u_{ij}$,
$\{i,j\}\in \chi_{S,\delta}$.
\begin{lem}\label{lemmacrossratio} Let $i,j,k,l$ be distinct indices modulo $n$ in dihedral order.
Then
$$\big[ij\,|\,kl\big]= \prod_{a=i}^{j-1}  \prod_{b=k}^{l-1} u^{-1}_{ab}\ . $$
Using $(\ref{S4})$, we can write any cross-ratio as a product of
$u_{ab}$ or their inverses.
\end{lem}

\begin{proof}
Suppose first that $1\leq i<j<k<l\leq n$. Using the definition of
$u_{ab}$,
$$u_{ak}\ldots u_{a\, l-1} = \big[a\, a+1\,|\,k+1\, k\big] \big[a\, a+1\,|\,k+2\,
k+1\big]\ldots \big[a\, a+1\,|\,l\, l-1\big]
 = \big[a\, a+1\,|\,l\, k\big]\ ,
$$
by repeated application of $(\ref{mult})$. Likewise, using  $(\ref{S4})$ and $(\ref{mult})$,
$$\prod_{a=i}^{j-1}  \prod_{b=k}^{l-1} u_{ab} =
\big[i\, i+1\,|\,l\, k\big] \big[i+1\, i+2\,|\,l\, k\big]
  \ldots \big[j-1\, j\,|\,k+1\, k\big]=\big[ij\,|\,lk\big]\ .$$ The
formula is clearly invariant under cyclic rotations.
 Therefore, given
any four  indices $i,j,k,l$ in arbitrary position, we can reduce
to this case by applying the inversion $(\ref{S4})$, which allows
us to interchange $i, j$, or $k,l$ or both pairs $(i,j), (k,l)$.
\end{proof}
It follows from invariant theory that every $\PSL_2$-invariant
regular function on $(\Pro^1)_*^n$ is a polynomial in the
cross-ratios $[ij\,|\,kl]$. We deduce from    lemma    $\ref{lemmacrossratio}$ that the
ring of regular functions on $\Mod_{0,S}$ is
\begin{equation}\label{regring}
\Or(\Mod_{0,S})=\Q\big[u_{ij}, {1\over u_{ij}} : \{i,j\} \in
\chi_{S,\delta} \big]\ .\end{equation}
 We can write down a generating set for  all algebraic relations between the coordinates
$u_{ij}$ in  a dihedrally-invariant manner. %Pl\"{u}cker relations..
Consider any chord $\{i,j\} \in \chi_{S,\delta}$. Then the set $S$ (considered as vertices)
with the elements $i$ and $j$ removed  $S\backslash\{i,j\}$ is
partitioned into two connected pieces $S_1$ and $S_2$. We say that
two chords $\{i,j\}$ and $\{k,l\}\in \chi_{S,\delta}$ \emph{cross}
if and only if $k\in S_1$ and $l \in S_2$. We write this
$$\{i,j\} \x \{k,l\}\ .$$
Given a subset $A \subset \chi_{S,\delta}$, let $A^{\xs}$ denote
the set of chords in $\chi_{S,\delta}$ which cross every chord in
$A$. We say that two sets $A, B\subset \chi_{S,\delta}$
\emph{cross completely} if $A^\xs= B $ and $B^\xs = A$, %consists
%of the set of all chords which cross every chord in $B$,
{\it
i.e.},
$$a \in A \quad \Leftrightarrow \quad a\x b \quad \hbox{ for all }\quad b \in B\ ,$$
 and vice versa (fig. 2). If, for example, $A$ is the single chord $\{i,j\}$, and $B$ is
the set of all chords crossing $\{i,j\}$, then $A$ and $B$ cross
completely.

\begin{cor}\label{corcrossingchords} For every two sets of chords $A,B \subset \chi_{S,\delta}$  which cross
completely,
\begin{equation}\label{id}
 u_A + u_B  =1 \ ,
\end{equation}
where $u_A=\prod_{a\in A} u_a$ and $u_B=\prod_{b\in B}u_b$.
\end{cor}
\begin{proof}
One can verify that $A$ and $B$ cross completely if and only if
there exist four elements $\{i,j,k,l\}\subset S$ in dihedral order
(fig. 2) such that
\begin{eqnarray}
A&= &\{\{p,q\}\in \chi_{S,\delta} :\quad  i\leq p< j \quad \hbox{
and } \quad k\leq q< l \}\ , \nonumber \\
B&= &\{ \{p,q\} \in \chi_{S,\delta} :\quad  j\leq p< k \quad
\hbox{ and } \quad l\leq q< i \}\ . \nonumber
\end{eqnarray}
%$A$ is  the set of all chords with one vertex in the segment
%$[i,j]$ and another in the segment $[k,l]$, and $B$ is the set of
%all chords with one vertex in the segment $[j,k]$ and the other in
%the segment $[l,i]$.
By  lemma $\ref{lemmacrossratio}$ and $(\ref{S4})$,  $u_A= [ij|kl]^{-1}=[ij|lk]$.
Likewise, $u_B = [li|jk]^{-1}= [il|jk]$. It follows that $u_A+u_B= [ij|lk]+[il|jk]=1$ by $(\ref{lemmacrossratio})$.
\end{proof}

%\begin{rem} Observe that $u_{kl}$ occurs as a factor of  $1-u_{ij}$ in the
%right hand side of $(\ref{id})$ if and only if the chords
%$\{k,l\}$ and $\{i, j\}$ cross each other. This relation is
%symmetric in $\{i,j\}$ and $\{k,l\}$. $\mathsf{x}$.
%\end{rem}

\begin{defn}
 Let $I^{\chi}_{S,\delta}\subset \Z[u_{ij}]$ %, \{i,j\} \in \chi_{S,\delta}]$
 denote the ideal
generated by the  identities $(\ref{id})$. Let the
\emph{dihedral extension} $\Modf_{0,S}$ of
$\Mod_{0,S}$ be the affine scheme
\begin{equation} \label{finitecomp}
\Modf_{0,S}= \Spec \,\Z[u_{ij}: \{i,j\}\in \chi_{S,\delta}]/I^{\chi}_{S,\delta}\ .
\end{equation}
%In particular, $\overline{\Mod}_{0,n}$ has a model over $\Z$.
By lemma $\ref{lemmacrossratio}$, we could also define $\Modf_{0,S}$ as follows:
\begin{equation} \label{ringoffunctionsascrossratios}
%\Z[u_{ij}:\,\{i,j\} \in \chi_{S,\delta}]=
\Modf_{0,S} =\Spec \Z\big[ \, [i \, j|l\,k], \hbox{ where } i,j,k,l\in S \hbox{ are in dihedral order}\big]/I_{S,\delta} \ ,
\end{equation}
where $I_{S,\delta}$ is the ideal generated by the identities $[ i \, j|l\,k]+[i\,l|j\, k]=1$ and  $(\ref{mult})$.
\end{defn}

%---------------------------------------------------------------------------------
%---------------------------------------------------------------------------------
%-------------------------INVALID ARGUMENT REMOVED -------------------------------
%------------------NOT CLEAR THAT I IS A PRIME IDEAL AT THIS STAGE----------------
%---------------------------------------------------------------------------------

%The scheme $\Modf_{0,S}$ is of dimension  $\ell=n-3$. To see this,
%suppose we are given the $\ell$ variables $u_{24},\ldots, u_{2n}$.
%Then by $(\ref{id})$,
%\begin{eqnarray}
%1-u_{13}&=& u_{24}\ldots u_{2n}\ , \nonumber \\
%1-u_{13}u_{14} & = & u_{25}\ldots u_{2n} \ , \nonumber \\
%\vdots & & \vdots \nonumber \\
%1-u_{13}\ldots u_{1\,n-1} & = & u_{2n} \ . \nonumber
%\end{eqnarray}
%This determines the variables $u_{13},\ldots, u_{1\, n-1}$
%corresponding to all chords with one vertex equal to $1$.
%Proceeding in this way, we can write all coordinates $u_{ij}$ as
%rational functions of $u_{24}, \ldots, u_{2n}$, which proves that
%the dimension is at most $\ell$.

%---------------------------------------------------------------------------------
%---------------------------------------------------------------------------------
%---------------------------------------------------------------------------------
%---------------------------------------------------------------------------------

Let $i_{\delta}: \Mod_{0,S} \hookrightarrow \Modf_{0,S}$ denote
the embedding obtained from $(\ref{Embed})$.
%We shall prove in the next section that the ideal $I$ is prime.
%The claim follows on
%observing that $\dim \Mod_{0,S}=\ell$.

\begin{lem}\label{lemmadense}
The scheme $\Modf_{0,S}$ is the Zariski closure of the image of $
\Mod_{0,S}$ in the affine space $\A^{n(n-3)/2}$. In particular, it is of dimension $\ell$.\end{lem}

\begin{proof}
  By $(\ref{tcoords})$, the dihedral coordinates $u_{24},\ldots, u_{2n}$ are
  equal to the cubical coordinates $x_1,\ldots, x_\ell$ respectively.  Therefore, $i_\delta(\Mod_{0,S})$ is contained in $\{u_{24},\ldots, u_{2n} \neq 0\}$. It follows using dihedral
symmetry that $i_{\delta}(\Mod_{0,S})
\subset \{u_{ij} \neq 0\,: \{i,j\}\in \chi_{S,\delta} \}$.  Using
$(\ref{id})$, we can write $x_i$, $1-x_i$, and $1-x_i\ldots x_j$ as
a product of dihedral coordinates $u_{a  b}$. This implies that
$$i_{\delta}(\Mod_{0,S}) = \{u_{ij}\neq 0\, : \{i,j\}\in \chi_{S,\delta}\} \subset \Modf_{0,S}\ .$$
By $(\ref{ExpCubeCoords})$, each divisor $\{x_i=0\}=\{u_{2\,
i+3}=0\}$ is in the closure of $\Mod_{0,S}$. Dihedral symmetry
implies that $\{u_{ab}=0\}$ is in the closure of $\Mod_{0,S}$ for
all $\{a,b\}\in \chi_{S,\delta}$.
\end{proof}

The complement $\Modf_{0,S}\backslash i_{\delta} (\Mod_{0,S})$ is
therefore a union of divisors
$$D_{ij} \subset \Modf_{0,S} \ ,\quad \hbox{ for each } \quad  \{i,j\} \in \chi_{S,\delta}\ ,$$
where  $D_{ij}$ is defined by the equation $u_{ij}=0$. In order to
describe the configuration of the divisors $D_{ij}$, consider
cutting the regular $n-$gon along the chord $\{i,j\}$ joining
vertices $i$ and $j$. This partitions the set of
\emph{edges} of $S$ into two sets $S_1$ and $S_2$ and breaks the
$n$-gon into two smaller polygons.
 Their sets of edges are $S_1 \cup \{e\}$ and
$S_2\cup \{e\}$, where $e$ is the new edge given by the chord
$\{i,j\}$ (fig. $3$). Each set inherits a dihedral structure
$\delta_k$ for $k=1,2$, and $\chi_{S,\delta}$ is a disjoint union:
\begin{equation} \label{chidecomp}
\chi_{S,\delta} = \chi_{S_1\cup\{e\}, \delta_1} \sqcup
\chi_{S_2\cup\{e\}, \delta_2}\sqcup \{i,j\} \sqcup
\bigcup_{\{k,l\}\x\{i,j\} } \{k,l\}\ .\end{equation}

\begin{figure}[h!]
  \begin{center}
%    \leavevmode
    \epsfxsize=10.0cm \epsfbox{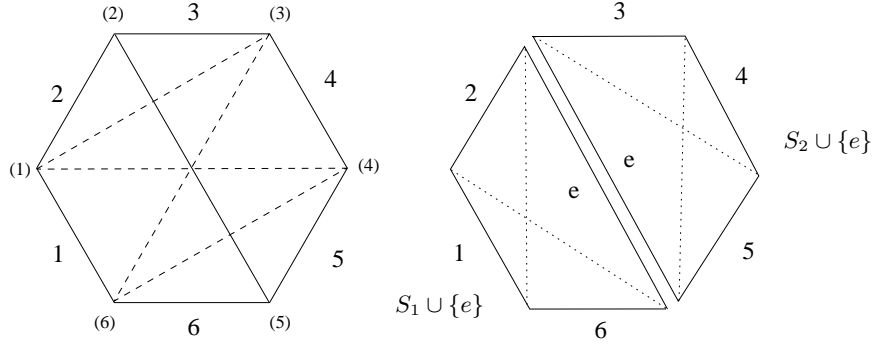}
  \label{Hexdecomp}
  \put(-142,10){\small{$S_1\cup\{e\}$}}
\put(5,72){\small{$S_2\cup\{e\}$}}
%\put(-280,10){\small{$S$}}
  \caption{  Decomposition of the hexagon on setting $u_{25}=0$. The variables corresponding
to  chords which cross $\{2,5\}$, namely $u_{13}$, $u_{46}$,
$u_{14}$, $u_{36}$, are all equal to $1$ (left). The system
$(\ref{id})$ splits into the pair of equations, $u_{15}=1-u_{26}$
and $u_{35}=1-u_{24}$, which identifies $D_{25}$ with
$\Modf_{0,4}\times \Modf_{0,4}$.}
  \end{center}
\end{figure}

\begin{lem}\label{lemmadecomp} The decomposition $(\ref{chidecomp})$ gives
 a canonical isomorphism
 $$D_{ij} \cong \Mod^{\delta_1}_{0,S_1\cup\{e\}}\times \Mod^{\delta_2}_{0,S_2\cup \{e\}}\ .$$
%which is canonical up to the action of  dihedral symmetries
%$D_{\delta_1}$, $D_{\delta_2}$
 % on each component.
\end{lem}
\begin{proof} Equations $(\ref{id})$ imply in particular that
\begin{equation} \label{thing}
u_{ab} + \prod_{\{c,d\}\x \{a,b\}} u_{cd} =1 \quad \hbox{ for all }
\{a,b\} \in \chi_{S,\delta}\ .\end{equation}
 Therefore, setting $u_{ij}=0$ implies that $u_{kl}=1 $ for all chords $\{k, l\}$ which cross
$\{i, j\}$. The system of equations $(\ref{id})$ then decomposes
into two disjoint sets, each one containing all variables $u_{ab}$,
where $\{a,b\}\in \chi_{{S_1\cup\{e\}}, \delta_1}$,  or
$\chi_{{S_2\cup\{e\}}, \delta_2}$ respectively. To see this,
consider the equation
\begin{equation} \label{mainrelforproof}
  u_A+   u_B=1 \ ,
\end{equation}
where $A,B\subset \chi_{S,\delta}$ cross completely, and where we write $u_I=\prod_{i\in I} u_i$ for any subset $I\subset \chi_{S,\delta}$.
 Consider the decomposition
$(\ref{chidecomp})$ above, and set $A_i=A\cap \chi_{{S_i\cup\{e\},\delta_i}}$ for $i=1,2$. It follows from the calculation above that,
since $u_{ij}=0$,
\begin{equation}\label{Acases}
u_A =
  \begin{cases}
    0 & \text{if } \{i,j\} \in A \ , \\
    u_{A_1} u_{A_2} & \text{otherwise}.
  \end{cases}
\end{equation}
A similar formula holds for $u_B$.
The picture below  depicts the three possible cases which can occur, up to exchanging $i,j$ or $A,B$.
 If neither $A_1$ nor $A_2$ is empty,  the set  $A$ contains chords on either side of the chord
$\{i,j\}$ (case I). It follows that $\{i,j\} \in A$, and therefore $u_A=0$.
Since $B=A^\xs$, it
follows that  $B\subset \{i,j\}^\xs$, and so $u_{k l} =1 $ for all
$\{k,l\} \in B$. Thus $(\ref{mainrelforproof})$ reduces to $0+1=1$.
Therefore we can assume without loss of generality that $A_1=\emptyset$ (see cases II and III), and so $u_A=u_{A_2}$
 by $(\ref{Acases})$.
It is clear that $B_1=\emptyset$, and so $u_B=u_{B_2}$.
 %Now if $ B_1$ contained a chord $b$, then
% $b$ could not cross all the chords in $A$ (as is  verified by looking at cases II and III).
%
\begin{figure}[h!]
  \begin{center}
%    \leavevmode
    \epsfxsize=12.0cm \epsfbox{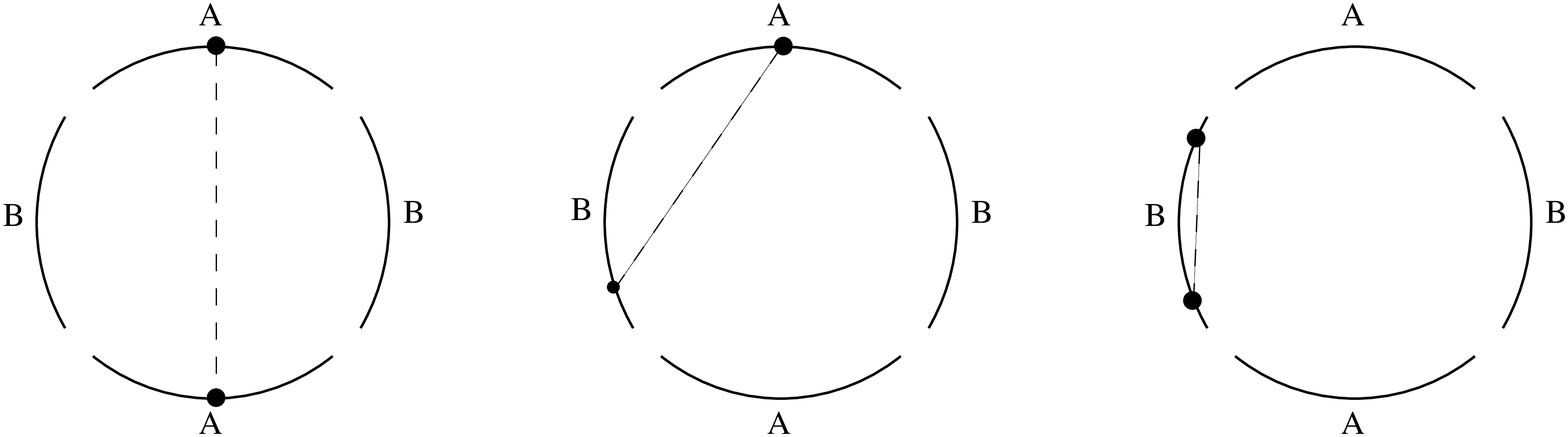}
  \label{threecasepic}
  \put(-260,10){I}\put(-130,10){II}\put(-5,10){III}
  %\caption{ Left - three cases which occur in the proof of lemma $(\ref{lemmadecomp})$ up to exchanging $i,j$ or $A, B$.}
  \end{center}
\end{figure}

\noindent  It follows that
equation $(\ref{mainrelforproof})$  reduces to
$u_{A_2}+ u_{B_2}=1\ ,$
%
%of the sets $A,B$ is contained in a
%$\chi_{{T_i\cup\{e\}}, \delta_i}$ for $i=1,2$. Suppose that $A\subset \chi_{T_1\cup\{e\},\delta_1}$. By the
%decomposition $(\ref{chidecomp})$ above,
%$$B= \big(B\cap \chi_{T_1\cup\{e\},\delta_1}\big) \sqcup \big(B\cap \{i,j\}^\xs\big)\ .$$
%But $u_{kl}=1$ for all $\{k,l\} \in B\cap \{i,j\}^\xs$. Therefore
% if we set $A'= A\cap \chi_{{T_1\cup\{e\}}, \delta_1}= A$ and
%$B'=B\cap \chi_{{T_1\cup\{e\}}, \delta_1}$,
 % becomes
%$$ 1- \prod_{a\in A'} u_a = \prod_{b\in B'} u_b\ ,$$
which  is a defining equation for
$\Mod^{\delta_2}_{0,S_2\cup\{e\}}$.  All pairs of
completely crossing sets in each smaller polygon $S_k\cup\{e\}$, for $k=1,2$, arise in this
way. This proves the result. % to dihedral symmetries.
%
%
% This is
%because two   chords cross in each of the smaller polygons if and
%only if they crossed in the original $n-$gon.
\end{proof}

It follows from the proof of the lemma that $D_{ij}$ and $D_{kl}$
have non-empty intersection if and only if the chords $\{i,j\}$
and $\{k,l\}$ do not cross. By $(\ref{thing})$, $u_{ij}$
and $u_{kl}$ cannot simultaneously be zero if $\{i,j\}\x\{k,l\}$.
  We are therefore led to consider partial
decompositions of the $n-$gon $(S,\delta)$ by $k$ non-crossing chords.

\begin{defn} For each integer $1\leq k\leq \ell$, let  $\chi^k_{S,\delta}$  denote the set of
$k$ distinct chords $\alpha=\{\{i_1,j_1\},\ldots, \{i_k,j_k\}\}$
in the $n$-gon $(S,\delta)$, such that no pair of chords in $\alpha$   cross. For
each such $\alpha \in \chi^k_{S,\delta}$,  let $D_\alpha$ denote
the subvariety defined by the equations $u_{i_1 \,j_1}=\ldots
=u_{i_k \,j_k}=0$, {\it i.e.},
$D_\alpha = \bigcap_{m=1}^k
D_{i_m\,j_m} \ .$
\end{defn}
\noindent  It follows by induction using the previous lemma that
the codimension of $D_{\alpha}$, for  $\alpha \in
\chi^k_{S,\delta}$,
 is exactly $k$, and that every codimension-$k$
intersection of divisors $D_{ij}$ arises in this manner.
%By
%applying this lemma inductively,
% it follows that for
Any set of $k$ chords  $\alpha\in \chi^k_{S,\delta}$ splits the polygon into $k+1$ pieces, and we have:
\begin{equation} \label{Dalphadecomp}
D_{\alpha} \cong \prod_{m=1}^{k+1} \Mod^{\delta_m}_{0,S_m}\ ,\end{equation}
where $(S_m, \delta_m)$ are  given by the set of all edges of each
%fundamental
 small polygon in  the $k$-decomposition  $\alpha$, with the
induced dihedral structures (fig. 4).

\begin{figure}[h!]
  \begin{center}
%    \leavevmode
    \epsfxsize=3.5cm \epsfbox{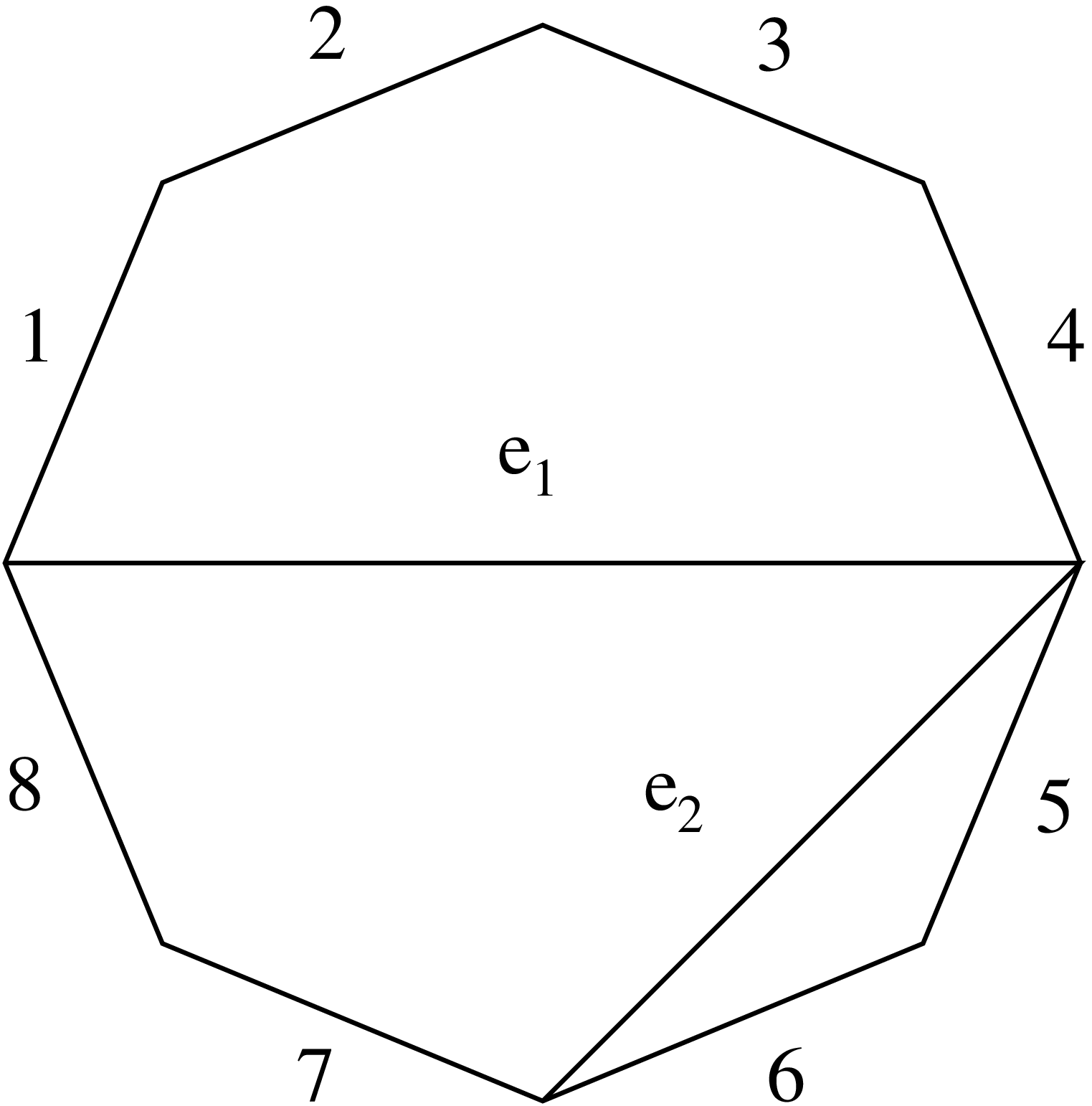}
  \label{Proofpic}
  \caption{  A partial decomposition $\alpha \in \chi^2_{8,\delta}$ of an octagon gives an isomorphism of $D_\alpha$ with
$   \Mod^{\delta_1}_{0,5}  \times \Mod^{\delta_2}_{0,4} \times \Mod^{\delta_3}_{0,3} = \Mod^{\delta_1}_{0,5} \times \A^1\times\{\mathrm{pt}\}.$ }
  \end{center}
\end{figure}

\begin{rem}
The set of all polygons equipped  with the operation of glueing
sides together forms %a cyclic operad
what is  known as the mosaic operad [Dev1]. This says that, given two
polygons with edges labelled $S_1\cup \{e\}$ and $S_2\cup\{e\}$
respectively, there is an  operation of glueing along the common
edge $e$, which gives rise to a map
$$\chi^k_{S_1\cup\{e\},\delta_1}\times \chi^l_{S_2\cup\{e\},\delta_2} \rightarrow
\chi^{k+l+1}_{{S_1\cup S_2},\delta}\ .$$ This corresponds to  the
decomposition of lemma $\ref{lemmadecomp}$.
% By lemma $2.4$, each divisor
%$D_{ij}$ is a product of partial compactifications of moduli
%spaces.
\end{rem}

\subsection{Forgetful maps between moduli spaces and projections}\label{sect25}
%We will need to consider  projection maps  between various moduli
%spaces. %of different dimensions.
Let $T$ denote any subset of $S$ such that $|T|\geq 3$. There is a
natural map
\begin{equation} \label{fTdef}
f_T: \Mod_{0,S} \To \Mod_{0,T} \end{equation}
obtained by forgetting  the marked points of $ S$ which do not lie
in $T$. Now suppose that $S$ has dihedral structure $\delta$. Then $T$ inherits a dihedral structure
which we denote $\delta_T$. If we view $S$ as the set of edges of the $n$-gon $(S,\delta)$, we obtain a map:
$$f_T:\chi_{S,\delta} \rightarrow \chi_{T,\delta_T}$$
which contracts all edges in $S\backslash T$ and combines the corresponding chords (fig. 5).

\begin{lem}\label{lemforgetmap} The map $(\ref{fTdef})$ extends to give a map
$f_{T} : \Modf_{0,S} \To \Mod^{\delta_T}_{0,T}$ such that
\begin{equation} \label{ft}
f_T^*(u_{kl}) = \prod_{\{a,b\}\in f_T^{-1}(\{k,l\})} u_{ab}\ .
\end{equation}
\end{lem}

\begin{proof} By  $(\ref{ringoffunctionsascrossratios})$, $f_T^*$ is induced by the map:
$$f_T^*: \Z\big[ [i\,j|l\,k]: \,\, i,j,k,l\in T^{\delta_T}\big]/I_{T,\delta_T} \hookrightarrow
\Z\big[ [i\,j|l\,k]: \,\, i,j,k,l\in S^{\delta}\big]/I_{S,\delta}\ ,$$
where $i,j,k,l\in T^{\delta_T}$ (resp. $S^\delta$) denotes  four elements in $T$ (resp. $S$) in dihedral order.
Formula $(\ref{ft})$ follows immediately from lemma $\ref{lemmacrossratio}$.
  \end{proof}

\begin{figure}[h!]
  \begin{center}
%    \leavevmode
    \epsfxsize=10cm \epsfbox{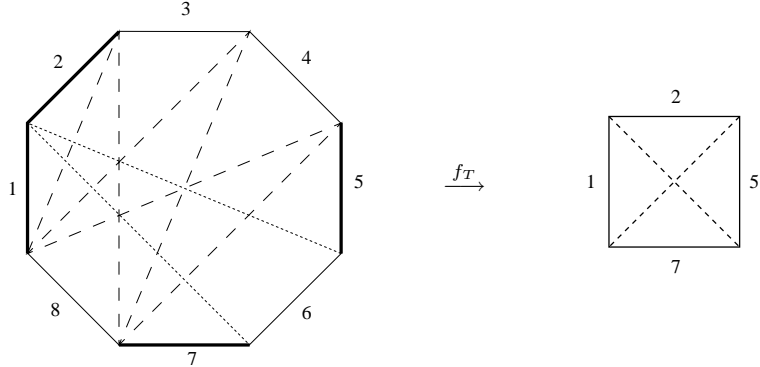}
\put(-120,65){$\overset{f_T}{\To}$}
  \label{Forgetful}
  \caption{The forgetful map $f_T$ contracts edges labelled $3,4,6,8$. The dihedral coordinates corresponding
to the  two chords in the square are pulled back by $f_T^*$ to $u_{15}u_{16}$ and
$u_{27}u_{37}u_{47}u_{28}u_{38}u_{48}$. }
  \end{center}
\end{figure}
\begin{rem} If $\{i,j\} \in \chi_{S,\delta}$, let $T$ denote the
four element set $T=\{i,i+1,j,j+1\}$. Then the dihedral coordinate
$u_{ij}$ is by definition a forgetful map:
$$\Modf_{0,S} \overset{f_T}{\To} \Mod_{0,T}^{\delta_T} \cong
\A^1\ .$$
\end{rem}

If $T_1,T_2$ are two subset of $S$ such that $|T_1\cap T_2|\geq 3$, we obtain  a
map
\begin{equation} \label{Generalfibrationdeltas}
f_{T_1}\times f_{T_2} : \Modf_{0,S}\To \Mod^{\delta_1}_{0,T_1}\times_{\Mod^{\delta'}_{0,T_1\cap T_2}} \Mod^{\delta_2}_{0,T_2}\ ,
\end{equation}
where $\delta_1=\delta_{T_1}$, $\delta_2=\delta_{T_2}$, and $\delta'=\delta_{T_1\cap T_2}$.

We  now consider  what happens when $|T_1\cap T_2|=2$. Suppose that the elements of $T_k$  are consecutive with respect
to $\delta$, for $k=1,2$. Then two cases can occur: either $T_1\cap T_2$
 consists of two consecutive elements  $\{i,i+1\}$ where $i\in S$, or $T_1\cap T_2$ has two components
and  $T_1\cap T_2=\{a,b\}$, where $a,b\in S$ are non-consecutive.  We only consider the first case here. This corresponds to
choosing a directed chord $\{i,j\} \in \chi_{S,\delta}$, and cutting along it.
Let $S_1$ and $S_2$ denote the corresponding
partition of the set $S$ viewed as edges of the $n$-gon (fig. 3), and consider the larger overlapping sets defined by
$T_1=S_1\cup\{i,i+1\}$ and $T_2=S_2\cup\{i,i+1\}$.

%  which arise when we cut along a chord $\{i,j\}\in \chi_{S,\delta}$.
%Let $S_1$ and $S_2$ denote the corresponding
%partition of the set $S$ viewed as \emph{vertices} of the $n$-gon (fig. 3), and consider the larger overlapping sets defined by
%$T_1=S_1\cup\{i,j\}$ and $T_2=S_2\cup\{i,j\}$.

 A product of forgetful maps gives
\begin{equation} \label{compactgenfib}
f_{T_1}\times f_{T_2} : \Modf_{0,S} \To
\Mod^{\delta_1}_{0,T_1}\times
\Mod^{\delta_2}_{0,T_2} \ .
\end{equation}
The dimension of the product on the right hand side is
$(|T_1|-3) + (|T_2|-3)$, which is $\ell-1$, one less than on the left.
Suppose that the chord $\{i,j\}= \{i,i+2\}$ is short. In that case, one of the sets, say $T_2$, has
just three elements, and $\Mod^{\delta_2}_{0,T_2}$ reduces to a point.
The complement $S\backslash T_1$ is a single point.
We write $S=\{s_1,\ldots, s_n\}$, and let $\{s_n\}=S\backslash T_1$. In that case, the restriction of  $f_{T_1}$ to $\Mod_{0,S}$:
$$f_{T_1}:\Mod_{0,\{s_1,\ldots, s_n\}} \rightarrow \Mod_{0,\{s_1,\ldots,s_{n-1}\}}$$
is a
fibration with one-dimensional fibres which are isomorphic to the
punctured projective line $\Pro^1\backslash\{s_1,\ldots,
s_{n-1}\}$.  In general, the restriction of the map $(\ref{compactgenfib})$ to the open set $\Mod_{0,S}$ is not
a fibration, but almost.
Let us compute it in cubical coordinates. By applying a dihedral symmetry, we can assume $i=2$.
 By $(\ref{tcoords})$,  we have $u_{2j}= x_{m}$, where $m=j-3$.
One verifies that
 \begin{eqnarray}
\Mod_{0,T_1} & = & \{(x_1,\ldots, x_{m-1}): x_i \notin\{0,1\}\ ,\quad x_i\ldots x_j\neq 1 \quad \hbox{for } i<j \}\ , \nonumber\\
\Mod_{0,T_2} &=& \{(x_{m+1},\ldots, x_{\ell}): x_i \notin\{0,1\}\ ,
\quad x_i\ldots x_j\neq 1 \quad \hbox{for } i<j \} \ ,\nonumber
\end{eqnarray}
and the map $ f_{T_1}\times f_{T_2} : \Mod_{0,S} \rightarrow
\Mod_{0,T_1}\times \Mod_{0,T_2}$ is just projection onto $x_m=0$ :
$$f_{T_1}\times f_{T_2}: (x_1,\ldots, x_\ell)
\mapsto  \big((x_1,\ldots,x_{m-1}), (x_{m+1},\ldots, x_\ell)\big)\
 .
$$
We can therefore think of $(\ref{compactgenfib})$ as a coordinate
projection in cubical coordinates. Referring to figure $10$, we
see that in $\Mod_{0,6}$ there are two such types of map,
one given by projection onto $\Mod_{0,5}$ %with $N=3$
(set $x_1=0$
or $x_3=0$), and the other given by projection onto
$\Mod_{0,4}\times \Mod_{0,4}$ %with $N=4$
(set $x_2=0$).
Restricting $(\ref{compactgenfib})$ to the divisor $u_{ij}=0$,
%Since this implies that $u_{kl}=1$ for all chords $\{k,l\}$ which
%cross $\{i,j\}$,
 we retrieve the isomorphism $D_{ij} \cong
\Mod^{\delta_1}_{0,T_1}\times \Mod^{\delta_2}_{0,T_2}$ which was defined in  lemma
$\ref{lemmadecomp}$.

\begin{rem}\label{remmakefib}
One can make the map $f_{T_1}\times f_{T_2}: \Mod_{0,S} \rightarrow \Mod_{0,T_1}\times \Mod_{0,T_2}$ into a fibration
by restricting it to an
 open subset $U_S\subset\Mod_{0,S}$. One obtains a map $f_{T_1}\times f_{T_2}:U_S\rightarrow V_S$, where $V_S \subset
\Mod_{0,T_1}\times \Mod_{0,T_2}$, whose fibers are isomorphic to $\A^1$ with $N$ points removed.  The $N$ removed points
correspond to the set of chords $\{k,l\}$ which cross
$\{i,j\}$, plus the chord $\{i,j\}$ itself.
%$N=(\ell-m+1)m+1$
%One can make $(\ref{genfibmap})$ into a fibration as follows.
%Define an open subset $U_S \subset \Mod_{0,T_1}\times \Mod_{0,T_2}$ by
%$$U_{S}=\{(x_1,\ldots, x_{m-1}), (x_{m+1},\ldots, x_\ell)\in \Mod_{0,T_1}\times \Mod_{0,T_2} \hbox{ such that } \qquad \qquad$$
%$$
%x_i\ldots x_j\neq (x_{p}\ldots x_{q})^{\pm1}  \hbox{ for all } 1\leq i<j\leq m-1\ ,  m+1\leq p<q\leq \ell\}\ .$$
%Let $V_S= (f_{T_1}\times f_{T_2})^{-1} U_S$. Then one can verify that the map $f_{T_1}\times f_{T_2}:V_S \rightarrow U_S$
%is a fibration whose  fibres are
%precisely
%$$\{x_m \in \C: x_m\neq 0 \hbox{ and } x_i\ldots x_m\ldots x_j\neq 1, \hbox{ for all } 1\leq i\leq m\leq j\leq \ell\}\ . $$
%One can check that the number of removed points is
%$N=(\ell-m+1)m+1$. This is the  number of chords $\{k,l\}\in \chi_{S,\delta}$ which cross
%$\{i,j\}$ plus the chord $\{i,j\}$ itself.

For example, consider the case $\Mod_{0,6}$, where we write the cubical coordinates $(x_1,x_2,x_3)$ as $(x,y,z)$.
Then $U_S=\{(x,y,z)\in \Mod_{0,6}: x\neq z\}$, and $V_S = \{(x,z)\in \C^2: x,z \neq 0,1: x\neq z^{\pm 1}\}$.
Then the fibration map $(x,y,z)\mapsto (x,z): U_S\rightarrow V_S$ has fibers $\{ y\in \C: y\notin\{0,1,x^{-1},z^{-1},(xz)^{-1}\}\}$.
The removed points in the fiber are given by the five equations $u_{25}=0$, $u_{13}=u_{14}=u_{36}=u_{46}=1$ (see figure 3).
\end{rem}

%The proof shows that
%$$\Q(\Mod_{0,S}) = \Q(F) \otimes f_{T_1}^* \Q(\Mod_{0,T_1})
%\otimes f_{T_2}^* \Q(\Mod_{0,T_1})\ ,$$ which in dihedral
%coordinates is equivalent to the decompostion
%$$\Q(\Mod_{0,S} ) = \Q\Big[ u_{ij},
%{1\over u_{ij}}, \Big({1\over u_{kl}}\Big)_{\{k,l\}\x \{i,j\}}
%\Big] \otimes  $$

\subsection{Normal vertex coordinates on $\Mod_{0,S}(\C)$.}\label{sect23} We wish to show that $\Modf_{0,S}$ is  smooth
and that $\Modf_{0,S}\backslash \Mod_{0,S}$ is normal crossing. In order to do this, we
construct normal coordinates in the neighbourhood of each
subvariety $D_{\alpha}$, for $\alpha \in \chi^k_{S,\delta}$.
We  first consider two further
relations satisfied by the dihedral coordinates $u_{ij}$ on $\Mod_{0,S}$.
% satisfy
%many different relations over $\Z$. It will be useful to consider
%an alternative presentation of the ideal $I$ of definition $2.3$.

We frequently use the following notation: for any two sets $I, J \subset S$, we write
\begin{equation} \label{productnotation} u_{IJ}=\prod_{i\in I,\, j\in J} u_{ij}\ .\end{equation}
Also, given two consecutive indices $i, i+1 $ modulo $n$, we adopt the convention
that $u_{i \,i+1 }=0.$
This is compatible with the decomposition of lemma $\ref{lemmadecomp}$: after cutting along a chord
$\{i,j\} \in \chi_{S,\delta}$, the vertices $i$ and $j$ become adjacent in each small polygon, and indeed, $u_{ij}=0$ is the equation of the corresponding
divisor (fig. 4).

\begin{lem} Let $\{p,q\} \in \chi_{S,\delta}$. Then
any three of the four coordinates $u_{pq}$, $u_{p\, q+1}$,
$u_{p+1\, q}$, $u_{p+1\, q+1}$ determine the fourth, and we have
the butterfly relation on $\Mod_{0,S}$:
\begin{equation} \label{sqid}
{u_{pq} (1-u_{p\, q+1}) \over 1-u_{pq}u_{p\, q+1}} = {1-u_{p+1\,
q+1} \over 1-u_{p+1\, q+1} u_{p+1\, q}}\ ,
\end{equation}
 where  $u_{p+1\,
q}=0$ ($u_{p\, q+1}=0$) if  $p+1$ and $q$,
(respectively $q+1$, $p$) are consecutive.
\end{lem}

\begin{proof}  Let $A$ and $B$ be the subsets
of vertices $S$ pictured in the diagram below (left). Then $(\ref{id})$
implies the following equations:
\begin{eqnarray}
1-u_{p\, q+1} &= & u_{A\, p+1}\, u_{AB} \, u_{A q} \ ,\nonumber \\
1- u_{pq}\,u_{p\,q+1}& = & u_{A\, p+1} \, u_{AB} \ ,\nonumber \\
1-u_{p+1\, q+1} &= & u_{p B} \, u_{AB}\, u_{Aq}\, u_{pq}\ ,
\nonumber\\
1- u_{p+1\,q+1}\,u_{p+1\,q}& = & u_{pB} \, u_{AB} \ .\nonumber
\end{eqnarray}
Identity $(\ref{sqid})$ follows by substitution.
\end{proof}
\begin{figure}[h!]
  \begin{center}
%    \leavevmode
    \epsfxsize=10.0cm \epsfbox{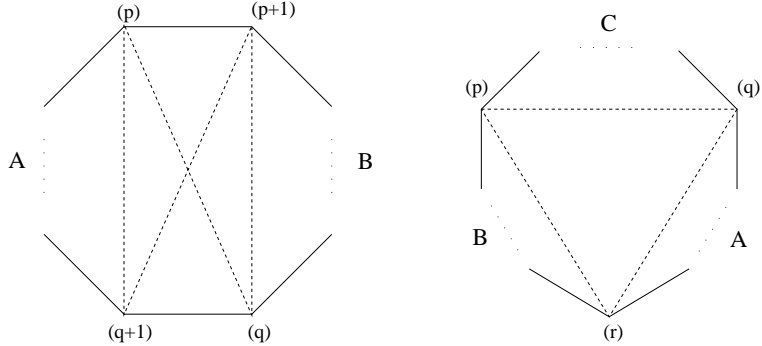}
  \label{Butterfly}
  \caption{Proof of the butterfly (left) and triangle relation (right).}
  \end{center}
\end{figure}

\begin{lem} \label{lemtriangleproof} Let $p,q,r$ denote three non-consecutive elements of $S$, and set
$\pi_r = \prod_{p<i<q} u_{ir}\ .$ Then the
triangle  relation holds on $\Mod_{0,S}$:
\begin{equation} \label{crossbow}
{1-u_{pq} \over (1-u_{pr}) (1-u_{qr}) } = {\pi_r \over (1-\pi_r
u_{pr}) (1-\pi_r u_{qr})}\ .
\end{equation}
If we regard this as a quadratic equation for $\pi_r$ in $\Q (
u_{pq}, u_{pr}, u_{qr})$, then the discriminant is non-zero in  a Zariski-open neighbourhood of $u_{pr}=u_{qr}=u_{pq}=0$. %(or $u_{qr}=0)$.
\end{lem}

\begin{proof}
Let $A$, $B$ and $C$ be the subsets of vertices $S$ pictured in
the diagram above (right). Then $(\ref{id})$ implies the following
equations:
\begin{eqnarray}
1- u_{pq} & = & u_{BC}\, \pi_r \,u_{AC}\ , \nonumber \\
1- u_{pr} & = & u_{BC}\, u_{Bq} \,u_{AB}\ , \nonumber \\
1- u_{qr} & = & u_{AC}\, u_{Ap} \,u_{AB}\ , \nonumber \\
1- u_{pr}\pi_{r} & = &  u_{Bq} \,u_{AB}\ , \nonumber \\
1- u_{qr}\pi_r & = &  u_{Ap} \,u_{AB}\ . \nonumber
\end{eqnarray}
The identity $(\ref{crossbow})$ follows by substitution. One
verifies by straightforward computation that the discriminant of
$(\ref{crossbow})$ is
\begin{equation} \label{Discintdefn}
\Delta_{pq,r}=(1 -u_{pq}u_{qr} +u_{qr}u_{rp}-u_{pr}u_{pq})^2 -4 (1-u_{pq})^2
u_{pr} u_{qr}\ ,\end{equation} from which the last statement follows.
\end{proof}

Let $\alpha \in \chi^\ell_{S,\delta}$ denote a
triangulation of the $n$-gon $(S,\delta)$. An \emph{internal triangle} of
$\alpha$ is a triple $p,q,r\in S$ such that $p,q,r$ are
non-adjacent, and $\{p,q\}$, $\{p,r\}$, and $\{q,r\}$ are in
$\alpha$. A \emph{free vertex} of $\alpha$ is a vertex $i\in
S$ such that $\{i,k\}\notin \alpha$ for all $\{i,k\}\in
\chi_{S,\delta}$. If $t$ denotes the number of internal triangles
in $\alpha$, and $v$ denotes the number of free vertices in
$\alpha$, then it is easy to show that $v=2+t.$ A triangulation of
the $n$-gon has no internal triangles if and only if it has
exactly two free vertices.

% It is easy to show that the number of free
%vertices minus the number of internal triangles is always equal to
%$2$.

\begin{defn}
Let $\alpha \in \chi^\ell_{S,\delta}$, and choose
an ordering on the set of chords  $\{i_1,j_1\}$, \ldots,$\{i_\ell,j_\ell\}$ in $\alpha$. Then the set of
\emph{vertex coordinates}\footnote{The reason for this
 terminology will become apparent in $\S\ref{sect24}$. A triangulation $\alpha$ corresponds to the point $D_{\alpha}$ which is a  vertex (corner)
 of the Stasheff polytope $\overline{X}_{S,\delta}\subset \Modf_{0,S}(\R)$ (see fig. 10).}
 corresponding to the ordered triangulation $\alpha$ is the set of variables:
$$x_1^{\alpha},\ldots, x_\ell^{\alpha}\ ,$$
defined by setting $x^{\alpha}_k=u_{i_kj_k}$ for $1\leq k\leq \ell$.\end{defn}

If $\alpha=\{\{2,4\},\ldots, \{2,n\}\}$ with the natural ordering,  then
$x^{\alpha}_k=x_k=u_{2\,k+3}$ for $1\leq k \leq \ell$, and we retrieve the cubical coordinates defined
in $(\ref{ExpCubeCoords})$ as a special case.

Let $\{i,j\} \in \chi_{S,\delta}$.  Recall from $(\ref{chidecomp})$ that there is a decomposition
$$\chi_{S,\delta} = \chi'\sqcup \{i,j\} \sqcup \bigcup_{\{k,l\}\x\{i,j\}} \{k,l\}\ ,$$
where $\chi'$ consists of all chords $\{a,b\}$ which do not cross $\{i,j\}$. The following lemma
 states  that we can eliminate all dihedral coordinates $u_{kl}$, where $\{k,l\}$ crosses $\{i,j\}$.

\begin{lem} Let $\{i,j\}\in \alpha$. On the open set $\Modf_{0,S}\backslash\{u_{ij}= 1\}$, every variable $u_{kl}$, where
 $\{k,l\} \in \chi_{S,\delta}$ crosses $\{i,j\}$,
can be expressed as a rational function  of $u_{ij}$, and the variables $u_{ab}$ where $\{a,b\}\in \chi'$.
\end{lem}
\begin{proof}
The easiest way to see this is on the example  $\alpha \in \chi^5_{8,\delta}$ depicted in figure 7 (left), where $\{i,j\}=\{1,5\}$. Consider the following  equations given by $(\ref{id})$:
\begin{eqnarray}
u_{28} &=& 1- u_{17}u_{16}u_{15}u_{14}u_{13}\ , \nonumber \\
u_{28}u_{27} &=& 1- u_{16}u_{15}u_{14}u_{13}\ , \nonumber \\
u_{28}u_{27}u_{26} &=& 1- u_{15}u_{14}u_{13}\ , \nonumber \\
u_{38} u_{28} &=& 1- u_{17}u_{16}u_{15}u_{14}\ , \nonumber \\
u_{38} u_{37} u_{28}u_{27} & = & 1- u_{16}u_{15} u_{14} \ , \nonumber \\
u_{38} u_{37} u_{36} u_{28}u_{27}u_{26} & = &  1- u_{15}u_{14} \ , \nonumber \\
 \ldots \nonumber
\end{eqnarray}
The identities $(\ref{id})$ imply that
$$1-u_{ij} = \prod_{\{k,l\}\x \{i,j\}} u_{kl}\ ,$$
and therefore all the variables on the left hand side in the equations above are invertible on $\Modf_{0,S}\backslash\{u_{ij}=1\}$. All the variables on the right hand side
lie in  $\chi'\cup \{i,j\}$. We can
therefore solve for $u_{28}, u_{27}, u_{26}, u_{38}, u_{37}, u_{36}$ and so on, in turn. The general case is similar.
\end{proof}

%
%In order to prove that $\Mod_{0,S}^\delta$ is smooth, we decompose the open sets $U_{\alpha}$ into  products and show
%that each component in the product is smooth. Recall that when we set $u_{ij}=0$, there is a decomposition of the $n$-gon into two pieces by cutting
%along the chord $\{i,j\}$, and the divisor $D_{ij}$ decomposes as a product (lemma $\ref{lemmadecomp}$). A similar decomposition
%holds  when $u_{ij}$ is any fixed value $r\in \A^1\backslash \{1\}$, except that we must now add a new vertex along  the edge $\{i,j\}$, labelled $\star$,  to each of  the smaller polygons  (see
%figure 6 below). In other words, if we view $S$ as a set of vertices, then the decomposition $S=T_1\cup\{i\} \cup \{j\} \cup T_2$ gives rise
%to two smaller polygons $P_1=T_1\cup\{i\}\cup\{j\} \cup \{\star\}$ and $P_2=T_2\cup\{i\}\cup\{j\}\cup\{\star\}$ with corresponding dihedral structures
%$\delta_1, \delta_2$.
%Let us set:
%\begin{eqnarray}
%u_{a \star}& =& \prod_{k\in T_2} u_{a k}\qquad \hbox{ for all } a \in T_1\ , \nonumber \\
%u_{\star b} &=& \prod_{l \in T_1} u_{l b} \qquad \hbox{ for all } b \in T_2 \ .\nonumber
%\end{eqnarray}
%This defines sets of dihedral coordinates $u^{P_1}, u^{P_2}$ on each small polygon $P_1,P_2$. Now it follows from their definition
%that the sets of relations $(\ref{id})$ are compatible: any relation which holds between dihedral coordinates on each smaller polygon $P_1$, $P_2$,
%actually holds on the original polygon $S$ (fig. 6). We obtain a map

Let $\{i,j\} \in \chi_{S,\delta}$ denote a chord. Then $\{i,j\}$
partitions the set of edges of $(S,\delta)$ into two sets, $S_1$
and $S_2$. The chord itself corresponds to the four edges
$E=\{i,i+1,j,j+1\}$. The sets $T_1=S_1\cup E$, and $T_2=S_2\cup E$
overlap in precisely the set $E$, and therefore $|T_1\cap T_2|=
4$. Let $\delta_E$ denote the induced dihedral structure on $E$.
By definition of the dihedral coordinate $u_{ij}$, there is an
isomorphism $u_{ij}: \Mod_{0,E}^{\delta_E} \cong \A^1.$ Therefore
$(\ref{Generalfibrationdeltas})$ defines a map:
\begin{equation} \label{stardecompmap}
\Modf_{0,S} \To \Mod^{\delta_1}_{0,T_1}\times_{\A^1} \Mod^{\delta_2}_{0,T_2}\ .\end{equation}
The chord $\{i,j\}$ is in both $\chi_{T_1,\delta_1}$ and $\chi_{T_2,\delta_2}$ (fig. 7).

\begin{prop}
The map $(\ref{stardecompmap})$ defines an embedding
\begin{equation}\label{Ualphadecomp}
\Modf_{0,S}\backslash \{u_{ij}=1\}\,\,{\hookrightarrow}\,\, \Mod^{\delta_1}_{0,T_1}\backslash \{u_{ij}=1\}\times_{\A^1\backslash \{1\}}
\Mod^{\delta_2}_{0,T_2}\backslash \{u_{ij}=1\}\ .
\end{equation}
\end{prop}

\begin{proof}
The map $(f_{T_1}\times f_{T_2})^*$ is given by:
$$\Z\big[[i\,j|l\,k]: i,j,k,l\in T_1^{\delta_1}\big]/I_{T_1,\delta_1} \otimes_{\Z[u_{ij}]}
 \Z\big[[i\,j|l\,k]: i,j,k,l\in T_2^{\delta_2}\big]/I_{T_1,\delta_1}
$$
$$\To
\Z\big[[i\,j|l\,k]: i,j,k,l\in S^{\delta}\big]/I_{S,\delta}\ ,$$
with the  notation used in the proof of lemma $\ref{lemforgetmap}$.
It  can also be regarded as a map:
$$\Z[ u_{ab}: \{a,b\} \in \chi_{T_1,\delta_1} ]/I^{\chi}_{T_1,\delta_1} \otimes_{\Z[u_{ij}]}
\Z[ u_{ab}: \{a,b\} \in \chi_{T_2,\delta_2} ]/I^{\chi}_{T_2,\delta_2}
 %\Z\big[[i\,j|l\,k]: i,j,k,l\in T_2^{\delta_2}\big]/I_{T_1,\delta_1}
$$
$$\To
\Z[ u_{ab}: \{a,b\} \in \chi_{S,\delta}]/I^{\chi}_{S,\delta}\ ,$$
%Since we know that the map
%$$f_{T_1}\times f_{T_2} : \Mod_{0,S} \To \Mod_{0,T_1}\times_{\A^1} \Mod_{0,T_2}\ ,$$
%is an isomorphism, the coordinate rings are isomorphic after localisation.%, and it follows that
%$(f_{T_1}\times f_{T_2})^*$ is injective.
 The previous lemma implies that the map $(f_{T_1}\times f_{T_2})^*$ is surjective
when we invert the coordinates $u_{kl}$ such that $\{k,l\}$ crosses $\{i,j\}$, {\it i.e.}, on the open set $u_{ij}\neq 1$.
\end{proof}

\begin{figure}[h!]
  \begin{center}
%    \leavevmode
    \epsfxsize=12cm \epsfbox{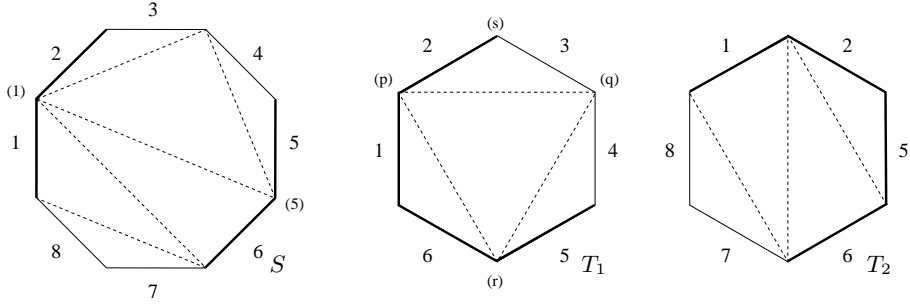}
\put(-124,10){\small{$T_1$}}
\put(-16,10){\small{$T_2$}}
\put(-242,10){\small{$S$}}
  \label{vertexproof}
  \caption{The induction step in the proof of  proposition $\ref{prop211}$.
The chord $\{1,5\}$ in the octagon on the left distinguishes the four  thick edges $E=\{1,2,5,6\}$.
The two sets $T_1=\{1,2,3,4,5,6\}$ and $T_2=\{5,6,7,8,1,2\}$ intersect in $E$ and define the hexagons on the right.
The right hand hexagon can be further decomposed  into a pair of pentagons, but the
  middle hexagon has an internal
triangle, so we have to  invoke the triangle lemma. }
  \end{center}
\end{figure}

Each set of vertex coordinates defines an \'etale map on a certain Zariski-open subset obtained by iterating
the map in the previous proposition.
Let $\alpha \in\chi^k_{S,\delta}$ denote any partial $k$-decomposition of the $n$-gon and define:
\begin{equation} U_\alpha=\bigcap_{\{i,j\}\in \alpha} \{ u_{ij} \neq 1 \}\subset \Modf_{0,S}\ .
\end{equation}
%Using $(\ref{id})$, one can show that
Since $D_{\alpha} = \{u_{ij}=0: \hbox{ for } \{i,j\}\in \alpha\}$,
% is also given by $D_{\alpha}=\{u_{kl}=1: \hbox{ for } \{k,l\} \notin \alpha\}$.
 it follows that
 $D_\alpha \subset U_{\alpha}$, and   $U_{\alpha}$ is an open neighbourhood of the subvariety $D_\alpha$ which contains the open set
$\Mod_{0,S}$. %Furthermore, it follows from $(\ref{id})$ that $u_{ij}\neq 1$ for all $\{i,j\}\in \alpha$ on the set $U_{\alpha}$.

Now let $\{i,j\} \in \alpha$, and consider the map
$(\ref{Ualphadecomp})$. We obtain two decompositions
 $\alpha_1$, $\alpha_2$, and Zariski-open sets $U_{\alpha_1}$, $U_{\alpha_2}$ on $T_1$, $T_2$ respectively.

%\begin{prop}   The map $(\ref{stardecompmap})$ is an embedding on the open set $U_\alpha$:
%\begin{equation} \label{Ualphadecomp}
%U_{\alpha} \hookrightarrow U_{\alpha_1} \times U_{\alpha_{2}}
%\end{equation} whose image is the subvariety defined by
%$u^{P_1}_{ij}= u^{P_2}_{ij}$.
%\end{prop}

%\begin{proof}
%This follows immediately from the previous lemma.
%\end{proof}

\begin{rem}
The embedding $(\ref{Ualphadecomp})$  extends the isomorphism $D_{\alpha}\cong D_{\alpha_1}\times D_{\alpha_2}$.
The product structure on each boundary stratum  of $\Modf_{0,S}$ therefore extends over a Zariski-open subset of the variety.
%This is not obvious from  the usual description
%of $\Mod_{0,S}$
%in terms of stable curves.
\end{rem}
%By taking the product with the projection $u_{ij}: U_{\alpha} \rightarrow \A^1$, we obtain a map
%\begin{equation} \label{parametrisedUdecomp}
%U_{\alpha} \hookrightarrow \A^1 \times U_{\alpha_1}\times U_{\alpha_2}\end{equation}
%whose image is $\{(r, u^{P_1}_{ij}=r, u^{P_2}_{ij}=r),\, r\in   \A^1\backslash\{1\}\}$.

Now for each $\alpha\in \chi^k_{S,\delta}$, we define a Zariski-open set
\begin{equation} U'_\alpha=\bigcap_{\{p,q,r\}\in \alpha} \{\Delta_{pq,r}\neq 0 \}\cap U_\alpha \subset \Modf_{0,S}\ ,
\end{equation}
where the intersection is over all sets of (ordered) internal triangles $\{p,q,r\}\in \alpha$ and $\Delta_{pq,r}$ is defined by $(\ref{Discintdefn})$.
It follows from $(\ref{Discintdefn})$ that $D_{\alpha} \subset U'_{\alpha}$.

\begin{prop}\label{prop211}
 Let $\alpha\in \chi^\ell_{S,\delta}$ denote any ordered triangulation of the $n$-gon $(S,\delta)$.
The set of  vertex coordinates  $\{ x_1^{\alpha},\ldots,
x_\ell^{\alpha}\}$ defines a map
$$(x^\alpha_1,\ldots, x^\alpha_\ell) : U'_{\alpha} \To \A^\ell$$
which is \'etale.
 It therefore defines a  system of coordinates on   $U'_\alpha(\R)$ or $U'_\alpha(\C)$.
\end{prop}

\begin{proof} By iterating the decomposition $(\ref{Ualphadecomp})$, we obtain an embedding
\begin{equation}\label{Utotaldecomp}
U_\alpha \hookrightarrow \prod_{i=1}^N U_{\alpha_i}\ ,\end{equation}
where each $\alpha_i$ is a triangulation of a $k_i$-gon which
 cannot be decomposed any further. We can assume that each decomposition is strict, {\it i.e.}, $k_i\geq 5$ for each $i$.

 Two cases can occur.  If there are no internal triangles,
then necessarily $k_i=5$, and we can assume that $x^{\alpha_i}_1= u_{13}$, and $x^{\alpha_i}_2= u_{14}$. In that case, $(\ref{id})$ gives:
$$ u_{25} = 1- u_{13}u_{14} \ ,\qquad
u_{24}u_{25}  =  1- u_{13} \ , \quad \hbox{ and }\quad
u_{35} u_{25}  =  1 - u_{14}\ .
$$
%\begin{eqnarray}
%u_{25}& = & 1- u_{13}u_{14} \ ,\nonumber \\
%u_{24}u_{25} & = & 1- u_{13} \ , \nonumber \\
%u_{35} u_{25} & = & 1 - u_{14}\ . \nonumber
%\end{eqnarray}
The variables on the left are invertible on $U_{\alpha_i} = \{u_{25}\neq 0\} \cap \{u_{24}\neq 0\} \cap \{u_{35}\neq 0\}$, so it follows
that $x^{\alpha_i}_1, x^{\alpha_i}_2$ is a coordinate system on this open set, {\it i.e.}, the map
$$(x_1^{\alpha_i},x_2^{\alpha_i}): U_{\alpha_i} \To \A^2$$
is injective, and certainly \'etale.
On the other hand,
if there is an internal
triangle $\{p,q,r\}$, and the $k_i$-gon cannot be decomposed any
further, then we are in the situation corresponding to
$\Mod_{0,6}$ pictured above (fig. 6, middle). By  symmetry, we can assume that $x_1^{\alpha_i}=u_{13}$,
$x_2^{\alpha_i}=u_{15}$, and $x_3^{\alpha_i}=u_{35}$.
By equation $(\ref{crossbow})$,  the variable $\pi_r=u_{25}$ depends quadratically on $x_1^{\alpha_i}, x_2^{\alpha_i}, x_3^{\alpha_i}$.
It follows from the triangle lemma $(\ref{lemtriangleproof})$ and the definition of $U'_{\alpha_i}$ that the map
$$(x_1^{\alpha_i},x_2^{\alpha_i},x^{\alpha_i}_3): U'_{\alpha_i} \To \A^3$$
is \'etale. It is in fact two to one. This is because all remaining dihedral coordinates $u_{ij}$ are uniquely
determined by $x_1^{\alpha_i}=u_{15}, \pi_r=u_{25}$, and $x_3^{\alpha_i} = u_{35}$
 by applying the relation $(\ref{id})$ and inverting coordinates which do not vanish on $U'_{\alpha_i}$ in much the same way as above.
From $(\ref{Utotaldecomp})$ we obtain an embedding $U'_{\alpha} \hookrightarrow \prod_{i=1}^N U_{\alpha_i}'$, which in turn  gives rise to a commutative diagram
$$\begin{array}{ccc}
 U'_\alpha  & \xrightarrow{(x_1^\alpha,\ldots, x_\ell^\alpha)}
  & \A^\ell  \\
 \downarrow &  &  \downarrow\\
   \prod_{i=1}^N U'_{\alpha_i} &\xrightarrow{\prod_{i=1}^N (x_1^{\alpha_i},\ldots, x^{\alpha_i}_{k_i})} & \prod_{i=1}^N \A^{k_i}
\end{array}
$$
The vertical maps on the left and on the right are diagonal  embeddings.
We have shown that the horizontal map along the bottom is \'etale. It follows that the horizontal map along the top is \'etale, which
completes the proof.
\end{proof}

If the triangulation $\alpha$ contains no  internal triangles, then $U'_\alpha= U_{\alpha}$, and
the functions $x_1^\alpha,\ldots, x_\ell^\alpha$ give an isomorphism of $\Mod_{0,S}$ with a Zariski open subset of $\A^\ell$.
\begin{cor} \label{lemnotrivertexcoords} Let  $\alpha\in \chi^\ell_{S,\delta}$,
such that   $\alpha$ has  no  internal triangles (and therefore two free vertices).
 Then the set $x_1^\alpha,\ldots, x_\ell^\alpha$ is
a system of coordinates everywhere on $\Mod_{0,S}$, and every
cross-ratio $u_{ij}$ is a $\Q$-rational function of the
$x_i^{\alpha}$.
\end{cor}

%\begin{proof}
%% We prove by induction that every coordinate
%%$u_{ij}$ can be written as a $\Q$-rational function of the
%%$x_i^{\alpha}$.
%We begin with the  triangulation of the $n$-gon corresponding to
%the variables $x_i^{\alpha}$ themselves, and
%  represent functions $u_{ij}$ which have
%already been expressed in terms of the functions $x_i^{\alpha}$ by
%drawing in the chord $\{i,j\}$ in the $n$-gon in turn. Given any
%square $p,p+1,q,q+1\subset S$ with three of its long edges or
%diagonals drawn in, we can use the butterfly lemma $(\ref{sqid})$
%to fill in the remaining edge or diagonal.
% It is a nice exercise
%to show that, starting from any triangle-free $\alpha$,  we
%eventually obtain the fully connected graph on $n$ points by
%applying this procedure repeatedly  (or see the inductive argument  used in the following
%proposition). The idea is to add new chords
%until there is a vertex ({\it e.g.}, 2) with the property  that all
%chords emanating from it ($\{2,4\}$,\ldots, $\{2,n\}$) are drawn
%in. Whenever this situation occurs, we apply the butterfly lemma
%to draw in chords $\{1,n-1\}$, $\{1,n-2\}$, \ldots, $\{1,3\}$ in
%turn, and so on until we have a complete triangulation.
%\end{proof}

We could also define $\Modf_{0,S}$ using the set
of equations $(\ref{sqid})$. One can verify that if $\alpha \in \chi^\ell_{S,\delta}$ has no internal triangles,
then all dihedral coordinates can be expressed in terms of
the vertex
coordinates $\{x_i^\alpha\}$  by  repeatedly applying the butterfly lemma.

% Observe that there is a large number
%of ordered internal triangle-free triangulations of an $n$-gon.
%Such coordinate systems $x_1^{\alpha},\ldots, x_\ell^\alpha$ can
%be used to generate many  birational transformations on period
%integrals.
%

\begin{lem} \label{lemUalphascover}
The sets $U'_{\alpha}$, for $\alpha \in \chi^{\ell}_{S,\delta}$, cover $\Modf_{0,S}$.
\end{lem}
\begin{proof}
For each partial decomposition $\beta \in \chi^{k}_{S,\delta}$,   set
$$N_\beta = \{u_{ij}=0 \hbox{ for all } \{i,j\} \in \beta\} \cap \{ u_{pq} \neq 0  \hbox{ for all } \{p,q\}\notin \beta\}\ ,$$
which is an open subset of $D_{\beta}= \{u_{ij}=0 \hbox{ for all }  \{i,j\} \in \beta\}$. It follows immediately from this definition  that
$\Modf_{0,S}$ decomposes as a disjoint union:
\begin{equation} \label{wierdstrat}
\Modf_{0,S} = \Mod_{0,S} \cup \bigcup_{k=1}^{\ell}  \bigcup_{\beta \in \chi^k_{S,\delta}} N_\beta \ .
\end{equation}
Let $\beta\in \chi^k_{S,\delta}$ denote any partial decomposition
of the $n$-gon $(S,\delta)$. By adding   chords, we can find a complete
triangulation $\alpha\in \chi^\ell_{S,\delta}$ which contains
$\beta$, without creating any new internal triangles in $\alpha$.
It follows from   $(\ref{Discintdefn})$ and lemma $\ref{lemtriangleproof}$ that
$N_{\beta}\subset U'_{\alpha}$. Note that if $\beta$ is the empty triangulation, then
$\Mod_{0,S} \subset U'_{\alpha}$ for any $\alpha$ which has no internal triangles. The decomposition $(\ref{wierdstrat})$ then implies that
$$\Modf_{0,S} \subset \bigcup_{\alpha \in \chi^\ell_{S,\delta}} U'_\alpha\  . $$
\end{proof}

\begin{thm} \label{partialissmooth} The affine variety $\Modf_{0,S}$ is smooth, and the
 divisors $D_{ij}$, for $\{i,j\}\in \chi_{S,\delta}$, are smooth and normal crossing.
\end{thm}
\begin{proof} Let $\alpha\in \chi^\ell_{S,\delta}$.
 Proposition $\ref{prop211}$ states that the vertex coordinates
$x_1^{\alpha},\ldots, x_\ell^{\alpha}$ corresponding to
$\alpha$ define an \'etale map
$U'_\alpha \rightarrow \A^\ell .$
The image of $\Mod_{0,S}\cap U'_{\alpha}$ in $U'_{\alpha}$ is precisely the complement of the
normal crossing divisor
$$x_1^{\alpha}\ldots x_\ell^{\alpha}=0$$
%Now consider any partial  decomposition $\beta\in \chi^k_{S,\delta}$.
%The divisors $D_{ij}\cap U_{\alpha}$ for $\{i,j\} \in \beta$
 %are therefore smooth and normal crossing for every complete triangulation $\alpha\in\chi^\ell_{S,\delta}$  which contains $\beta$
 (see fig. 8).
The theorem follows since the sets $U'_{\alpha}$ cover $\Modf_{0,S}$.
\end{proof}

%Note that since an $n$-gon is triangulated by exactly $n-3$
%chords, exactly $\ell$ divisors $D_{ij}$ can intersect in a point.

% The divisors $D_{\alpha}$ form the following simplicial scheme of affine varieties:
%\begin{equation} \label{simpscheme}
%\Modf_{0,S} \leftarrow \bigsqcup_{\alpha_1\in \chi^1_{S,\delta}}
%D_{\alpha_1} \leftleftarrows \bigsqcup_{\alpha_2\in
%\chi^2_{S,\delta}} D_{\alpha_2} \leftleftarrows \ldots
%\leftleftarrows \bigsqcup_{\alpha_\ell\in \chi^\ell_{S,\delta}}
%D_{\alpha_\ell}\ .
%\end{equation}
%The corresponding spectral sequence computes the cohomology of
%$\Mod_{0,S}$. It follows that  $H^1(\Mod_{0,S}) = H^1(\Modf_{0,S})
%\oplus \bigoplus_{\alpha\in \chi_{S,\delta}} H^0(D_\alpha)$, and
%hence $H^1(\Modf_{0,S})=0$.

\begin{figure}[h!]
  \begin{center}
%    \leavevmode
    \epsfxsize=6.0cm \epsfbox{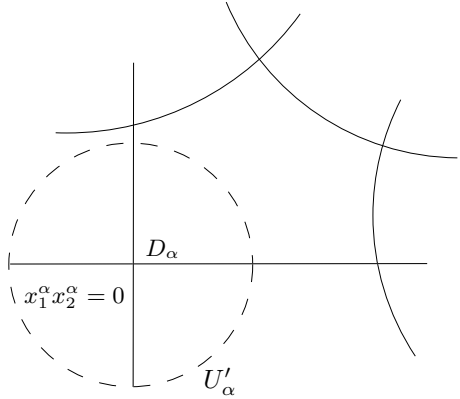}
  \label{Covering}
\put(-122,50){\small{$D_\alpha$}}
\put(-168,32){\small{$x_1^{\alpha}x_2^{\alpha}=0$}}
%\put(-170,50){\small{$x_1^{\alpha}$}}
%\put(-124,20){\small{$x_2^{\alpha}=0$}}
\put(-100,0){$U'_{\alpha}$}
  \caption{The covering of $\Modf_{0,5}$. To each vertex $\alpha \in \chi^\ell_{S,\delta}$, there is an open set $U'_{\alpha}$ on which
  the set of vertex coordinates $\{x_i^{\alpha}\}$  cross normally. The sets $U'_{\alpha}$ form a covering, which implies that the divisors
  $D_{ij}$ are smooth and normal crossing.
  }
  \end{center}
\end{figure}

\subsection{The real moduli space $\Modf_{0,S}(\R)$} \label{sect24}
% Consider the set of real points
%$\Modf_{0,S}(\R)$. % and $\overline{\Mod}_{0,S}(\R)$.
Consider the moduli space of   projective circles with
$n$ ordered marked points:
$$\Mod_{0,S}(\R)= \PSL_2(\R)\backslash\Pro^1(\R)_*^n  \ .$$
The space $\Mod_{0,S}(\R)$ is not connected, but is a disjoint
union of  open cells  which we define as follows. First, we fix
a dihedral structure $\delta$ on $S$, which defines a set of
dihedral coordinates $u_{ij}$ for $\{i,j\} \in \chi_{S,\delta}$.
Let %We define
\begin{equation} \label{Xbardef}
\overline{X}_{S,\delta}= \{ u_{ij} \geq 0\,\, : \,\, \{i,j\}\in
\chi_{S,\delta}\}\ \subset \Modf_{0,S}(\R)\ . \end{equation}
 By $(\ref{id})$,
$\overline{X}_{S,\delta}$ is also defined by  the equations $0\leq
u_{ij} \leq 1$ for all $\{i,j\}\in \chi_{S,\delta},$  and is
therefore compact. We define the  open cell $X_{S,\delta}$
 to be the
interior of $\overline{X}_{S,\delta}$:
\begin{equation} \label{Xdef}
X_{S,\delta}= \overline{X}_{S,\delta} \cap \Mod_{0,S}=\{ u_{ij} >
0\,\, : \,\, \{i,j\}\in \chi_{S,\delta}\}\  .
\end{equation}
The sets $X_{S,\delta}$ and $\overline{X}_{S,\delta}$ are clearly
preserved by  the dihedral symmetries of $\delta$,  so there is
%write the corresponding
an action of the dihedral group
% deduce that there
%is an action of $D_{2n}$ on $\overline{X}_{S,\delta}$.
$D_{2n}\times \overline{X}_{S,\delta}\rightarrow
\overline{X}_{S,\delta} .$
 Using explicit simplicial
coordinates $(\ref{tcoords})$, one  checks that the open set
$X_{S,\delta}$ is homeomorphic to the % unit
 simplex
\begin{equation} \label{Xsimplex}
X_{S,\delta}\cong \{(t_1,\ldots, t_\ell)\quad : \quad
0<t_1<\ldots< t_\ell<1\}\ .\end{equation}
It follows that $X_{S,\delta}$ is
contractible, and moreover, that $X_{S,\delta}$ is a %maximal
connected
component of $\Mod_{0,S}(\R)$. After changing to cubical
coordinates, we see that $X_{S,\delta}$ is the unit hypercube
$\{(x_1,\ldots, x_\ell):
x_i\in(0,1)\}=(0,1)^\ell,$ which explains the nomenclature of each %these
coordinate system (fig. 9).

%---------------------------------------------------------
%-----------PICTURE CAN ALSO GO IN HERE-------------------
%---------------------------------------------------------

Each cell $X_{S,\delta}$ consists of the set of points
$s_1,\ldots, s_n \in \Pro^1(\R)$ such that $s_1,\ldots, s_n$ are in
the dihedral order determined by $\delta$.
Two components $X_{S,\delta}$, $X_{S,\delta'}$ are disjoint if $\delta$ and $\delta'$ are distinct, and the set of dihedral structures
are
permuted transitively by the symmetric group $\Sym_n$.
%  Since the symmetric
%group $\Sym_n$ permutes the set of all dihedral structures
%$\delta$, it follows that the components $X_{S,\delta}$ and
%$X_{S,\delta'}$ for distinct dihedral structures $\delta, \delta'$
%are disjoint. Futhermore, since every set of $n$ points in
%$\mathbb{RP}^1$ has a natural dihedral order, the cells
%$X_{S,\delta}$ cover $\Mod_{0,S}$.
 This implies the  following tiling lemma. Devadoss  has studied
the exact glueing relations between the %various
 cells
$X_{S,\delta}$ in this tiling [Dev1].

 % decomposition efinitions above.
%Symmetric group acts on them, disjoint, and cover

\begin{lem}\label{lemmatess}
The space $\Mod_{0,S}(\R)$ is the disjoint union of the $n!/2n$
open cells $X_{S,\delta}$, as $\delta\in \Sym_n/D_{2n}$ ranges
over the set of all dihedral structures on $S$.
\end{lem}

\noindent It is now clear that the choice of a dihedral structure
$\delta$ on $S$ is equivalent to  the choice of a fundamental cell
$X_{S,\delta} \in \pi_0(\Mod_{0,S}(\R))\cong \Sym_n/D_{2n}$.  The
set of dihedral coordinates $u_{ij}$ corresponding to $\delta$ can
be regarded as a natural set of functions
 which is stable under the action of the  symmetry group
of $X_{S,\delta}$.
%
%By corollary 2.5, $\Modf_{0,S}(\R)$ is obtained by blowing up and
%compactifying along  the divisors which bound the cell
%$X_{S,\delta}$.
%
%Blow up divisors surrounding this cell. It follows that
%$X_{S,\delta}$ is the unique cell in $\Modf_{0,S}(\R)$ which is
%bounded.

\begin{figure}[h!]
  \begin{center}
%    \leavevmode
    \epsfxsize=13.0cm \epsfbox{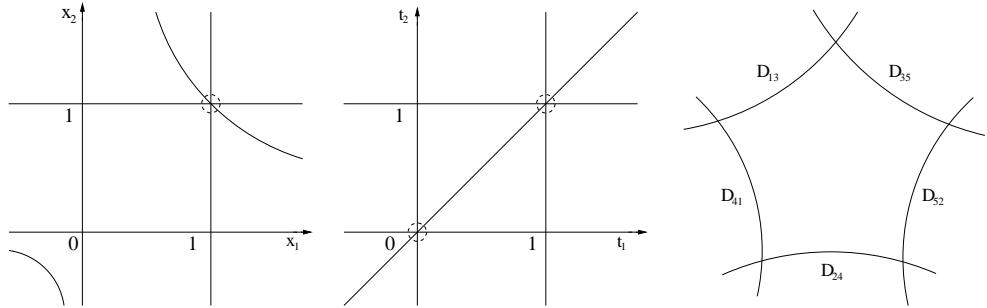}
  \label{blowups}
  \caption{The set of real points $\Mod_{0,5}(\R)$ in cubical
(left),  simplicial (middle), and dihedral coordinates (right). The
dotted circles denote points which are blown up when passing to
dihedral coordinates. %The cell $X_{5,\delta}$ is shaded.
There are
5!/10=12  regions $X_{S,\delta}$ when $|S|=5$. }
  \end{center}
\end{figure}

\begin{defn}
For each chord $\{i,j\}\in \chi_{S,\delta}$, we define the  face
$F_{ij}(\overline{X}_{S,\delta})$ of $\overline{X}_{S,\delta}$  to be  the
closed subset  $F_{ij}(\overline{X}_{S,\delta}) = D_{ij} \cap
\overline{X}_{S,\delta}\subset \Modf_{0,S}(\R).$ Likewise, for
each $\alpha \in \chi^k_{S,\delta}$, we define the codimension-$k$
face of $\overline{X}_{S,\delta}$ to be
$F_\alpha(\overline{X}_{S,\delta})=D_\alpha \cap \overline{X}_{S,\delta}$.
\end{defn}

\noindent  It follows from lemma 2.5  that
\begin{equation}\label{Facestruc}
F_{ij}(\overline{X}_{S,\delta}) \cong
\overline{X}_{T_1\cup\{e\},\delta_1}\times
\overline{X}_{T_2\cup\{e\},\delta_2}\ , \end{equation} where
$T_1\cup T_2$ is the partition of (the set of edges) $S$
corresponding to the chord $e=\{i,j\}$. By equation $(\ref{Dalphadecomp})$, each
codimension-$k$ face $F_{\alpha}(\overline{X}_{S,\delta})$ is  a product
$$F_{\alpha}(\overline{X}_{S,\delta})\cong \prod_{m=1}^{k+1}
\overline{X}_{S_m,\delta_m}\ .$$
By  repeatedly taking boundaries
% or simply by intersecting the
%simplicial scheme $(\ref{simpscheme})$ with the real cell
%$\overline{X}_{S,\delta}$,
we
obtain a stratification:
\begin{equation} \label{Firststrat}
\overline{X}_{S,\delta} \supseteq \partial \overline{X}_{S,\delta}
\supseteq \partial^2 \overline{X}_{S,\delta} \supseteq \ldots
\supseteq \partial^\ell \overline{X}_{S,\delta}\ ,\end{equation}
%
%\begin{defn} For each integer $1\leq k\leq n$, let  $\chi^k_{S,\delta}$  denote the set of
%$k$ distinct chords $\alpha=\{\{i_1,j_1\},\ldots, \{i_k,j_k\}\}$
%in an $n$-gon which do not cross pairwise. For each such $\alpha
%\in \chi^k_{S,\delta}$  we define the codimension-$k$ face of
%$\overline{X}_{S,\delta}$ to be
%$$F_\alpha(\overline{X}_{S,\delta}) = \bigcap_{m=1}^k
%F_{i_m\,j_m} (\overline{X}_{S,\delta})\ ,$$ which is the subset of
%$\overline{X}_{S,\delta}$ defined by $u_{i_1 \,j_1}=\ldots =u_{i_k
%\,j_k}=0$.
%\end{defn}
%
where the  codimension-$k$ boundary of $\overline{X}_{S,\delta}$
is the union of its codimension-$k$ faces:
$$\partial^k \overline{X}_{S,\delta} = \bigcup_{\alpha \in
\chi^k_{S,\delta}} F_\alpha(\overline{X}_{S,\delta}) \qquad \hbox{ for }
1\leq k\leq \ell\ .
$$

For each $n\geq 4$, the associahedron $K_{n-1}$, or Stasheff
polytope [St], is a convex polytope of dimension $n-3$  whose
codimension-$k$ faces are indexed by the partially ordered set of   compatible bracketings
on a set of $n-1$ elements.
% Blah. It was proved by Maclane that
%the lattice of .. can indeed be realised as a convex polytope.

\begin{cor} The lattice of faces $\overline{X}_{S,\delta}$ is
combinatorially equivalent to the associahedron $K_{n-1}$.
\end{cor}

\begin{proof} The set of all codimension-$k$ faces $F_\alpha(\overline{X}_{S,\delta})$
is indexed by $k$-triangulations of a regular $n$-gon, and the
inclusion of one face in another is given by removing a chord. By
taking the dual graph of a partial triangulation of an $n$-gon we
obtain a planar tree. If we fix an edge $s_1$ of $S$, then each
such tree is rooted, and defines, in a standard way, a
bracketing of the ordered set $\{s_2,\ldots, s_n\}$ (fig. 11). We
obtain in this way a bijection between faces of
$\overline{X}_{S,\delta}$ and bracketings on a set of $n-1$
elements (this is beautifully illustrated in  [Dev2]).
\end{proof}

\begin{figure}[h!]
  \begin{center}
%    \leavevmode
    \epsfxsize=9.0cm \epsfbox{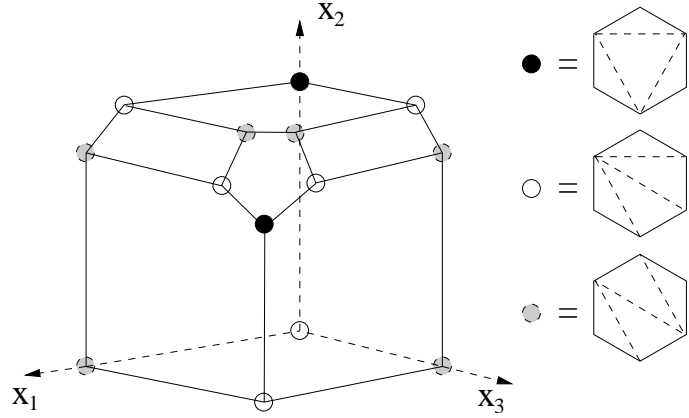}
  \label{Stasheff}
  \caption{The associahedron or Stasheff polytope $\overline{X}_{6,\delta}=K_5$ in  $\Modf_{0,6}(\R)$ obtained by
truncating the unit cube in $\R^3$, or blowing up along
$x_1=x_2=x_3=1$, and then $x_1=x_2=1$ and $x_2=x_3=1$. It has six
faces $F_{13},F_{24},F_{35}, F_{46}, F_{51}, F_{62}$ which are pentagons $\overline{X}_{5,\delta}$, and three
faces $F_{14}, F_{25}, F_{36}$ which are quadrilaterals
$\overline{X}_{4,\delta_1}\times\overline{X}_{4,\delta_1}$. These
are permuted by the group $D_{12}$. There are three types of
vertices corresponding to three kinds of triangulation of a
hexagon. The vertex coordinates defined in \S2.3 for each triangulation provide local
affine charts in the neighbourhood of each vertex.
}  % in cubical coordinates. }
  \end{center}
\end{figure}

\noindent Since each face $F_\alpha$ is  contractible, we can view
$\overline{X}_{S,\delta}$ in the coordinates $u_{ij}$ as a
dihedrally-symmetric algebraic model of the associahedron
$K_{n-1}$. The fact that the divisors $D_{ij}$ cross normally
implies  that the associahedron is a simple polytope, {\it i.e.},
each vertex is the intersection of exactly $\ell$ distinct faces.

\begin{rem}
As remarked earlier, $\Modf_{0,S}$ can be obtained by blowing up a
set of
divisors bounding $X_{S,\delta}$. %in a given coordinate model of
%$\Mod_{0,S}$
% which do not cross normally.
 Since the operation of
blowing up is non-commutative, we have to specify %
that the blow-ups occur along subvarieties in increasing order of
dimension.
%in which order
%to perform them.
% The convention is that we blow up in decreasing order of
%codimension.
A useful intuitive picture of the polytopes
$\overline{X}_{S,\delta}$ can be obtained by blowing up the unit
hypercube $[0,1]^\ell$
 along the  divisors $x_i=\ldots= x_j=1$ for $1\leq i<j\leq \ell$. The set of real points  in the blow-up can be visualised by
truncating the unit hypercube along the hyperboloids $x_i\ldots
x_j=1-\varepsilon$, where $i<j$ for some fixed
$\varepsilon>0$
 which is sufficiently small (see fig. 10).
Alternatively, one could truncate  the simplicial model of
$X_{S,\delta}$  to obtain another explicit construction of
$K_{n-1}$ (see [Dev2]).  This involves a greater number of
truncations, however.
\end{rem}

% (Let us
%denote by $V(\overline{X}_{S,\delta})$ the set of all such
%vertices, which will serve as a base point for the fundamental
%group.)

\subsection{The  compactification $\overline{\Mod}_{0,S}$ and its divisors at infinity}\label{sect26}
%We now consider the full compactification $\overline{\Mod}_{0,S}$
%of $\Mod_{0,S}$ and the structure of its divisors at infinity.
The
set of all cross-ratios defines an embedding:
\begin{equation} \label{fullcomp}
\{[i\, j\, |k\, l] \}: \Mod_{0,S} \To
\Mod_{0,\{i,j,k,l\}}^{\binom{n}{4}} \cong (\Pro^1\backslash\{0,1,\infty\})^{\binom{n}{4}}\
.
\end{equation}
The coordinates $\{[i\, j|k\, l]\}$  satisfy identities $(\ref{S4})$ and $(\ref{mult})$. These
identities define a projective
scheme we denote
$$\overline{\Mod}_{0,S}\subset (\Pro^1)^{\binom{n}{4}}\ ,$$
which is defined over $\Z$. This representation of
$\overline{\Mod}_{0,S}$ is degenerate since some coordinates are
the same, but it is clearly invariant under the action of the
symmetric group. % $ \sigma \big[ij\,|\,kl\big] = \big[\sigma(i)\,
%\sigma(j) \,|\, \sigma(k) \,\sigma(l)\big] $ for all $\sigma \in
%\Sym(S)$.
 Now for every dihedral structure $\delta$ on $S$, there
is an embedding
$$j_{\delta}: \Modf_{0,S}\rightarrow \overline{\Mod}_{0,S}$$
given by lemma $2.1$, which expresses every cross ratio $[ij|kl]$
as a product of dihedral coordinates. Letting $\delta$ vary, we
obtain in this way an affine covering of $\overline{\Mod}_{0,S}$.

\begin{lem}\label{lemfullcompactif} The  compactification $\overline{\Mod}_{0,S}$ is
covered by affine charts $\Modf_{0,S}$, as $\delta$ ranges over the
set of all dihedral structures on $S$:
\begin{equation}\label{Fullcompactcover}
\overline{\Mod}_{0,S} = \bigcup_{\delta \in \Sym_n/D_{2n}}
j_{\delta} \big( \Modf_{0,S} \big) \ .\end{equation}
\end{lem}

\begin{proof}
The right-hand side is clearly contained in the left. But one can
show in a similar manner to the proof of lemma $\ref{lemmadense}$
that $\Mod_{0,S}$ is dense in $\overline{\Mod}_{0,S}$ as defined above. The point is that %, by  lemma $\ref{lemmacrossratio}$,
 $j_{\delta}(\Mod_{0,S}) \subset \{[i\,j|k\,l]\neq 0\}\subset \overline{\Mod}_{0,S}$ for all
 $\delta$. Any cross-ratio $[i\,j|k\,l]$ is a dihedral coordinate $u_{ab}$ for some dihedral structure $\delta_0$.
 Therefore $[i\,j|k\,l]=0$ is in the closure  of $j_{\delta_0}(\Mod_{0,S})$ by
 lemma $\ref{lemmadense}$. The same holds for the divisors $[i\,j|k\,l]=1,\infty$ by $(\ref{S4})$.
 This proves that both sides  are equal.
\end{proof}
Theorem $(\ref{partialissmooth})$ implies the following corollary.
\begin{cor} $\overline{\Mod}_{0,S}$ is smooth and
 $\overline{\Mod}_{0,S}\backslash \Mod_{0,S}$ is a normal
crossing divisor.
\end{cor}
\noindent The irreducible components  at infinity of $\overline{\Mod}_{0,S}\backslash
\Mod_{0,S}$ can be described as follows. %by applying the symmetric group to

\begin{lem} Let $\delta$ denote a dihedral structure on $S$, and let $\{p,q\} \in \chi_{S,\delta}$. The chord $\{p,q\}$
partitions the set $S$, viewed as edges of the $n$-gon $(S,\delta)$,
into two sets $P_1\cup P_2$. Then the divisor
$j_{\delta}(D_{pq}) \subset \overline{\Mod}_{0,S}\backslash \Mod_{0,S}$
is determined by the equations $[i\,j|k\,l]=0$ for all distinct indices $i,j,k,l$ such that
$$\qquad \{i,k\} \subset P_1 \hbox{ and } \{j,l\} \subset P_2\ ,$$
$$\hbox{ or } \quad \{j,l\} \subset P_1 \hbox{ and } \{i,k\} \subset  P_2 \ .$$
\end{lem}
\begin{proof} On the chart  $j_{\delta}(\Modf_{0,S})$, these equations imply in particular that $u_{pq}=0$, and  therefore determine the divisor
$j_{\delta}(D_{pq})$.
That the remaining cross-ratios also vanish on $j_{\delta}(D_{pq})$ follows from lemma $\ref{lemmacrossratio}$.
\end{proof}
It follows that two divisors $j_{\delta_1}(D_1)$ and $j_{\delta_2}(D_2)$ coincide on $\overline{\Mod}_{0,S}\backslash \Mod_{0,S}$
if and only if the corresponding partitions of $S$ agree.
\begin{defn}
A partition $P_1\cup P_2 =S$ is  \emph{stable} if
$|P_1|\geq 2$ and $|P_2|\geq 2$.  A dihedral structure $\delta$ on $S$ is
\emph{compatible} with $P_1\cup P_2$ if the elements of each set $P_1$ and $P_2$ are consecutive with respect to $\delta$.
A divisor $D \subset \overline{\Mod}_{0,S}\backslash \Mod_{0,S}$ is said to be \emph{at finite distance} with respect to a
dihedral structure $\delta $, if $D \subset j_\delta(\Modf_{0,S})$.
\end{defn}

Observe that $D \subset j_\delta(\Modf_{0,S})$ if and only if $D\cap j_\delta(\Modf_{0,S})\neq \emptyset.$

\begin{prop} \label{divisorsatinfinity} There is  a bijection between the irreducible
components of the  divisors at  infinity of
$\overline{\Mod}_{0,S}\backslash \Mod_{0,S}$, and stable partitions $S=P_1\cup P_2.$
   The component $D$ corresponding to this partition  is canonically isomorphic to
$$\overline{\Mod}_{0,P_1\cup\{e\}}\times
\overline{\Mod}_{0,P_2\cup\{e\}}\ ,$$
where $e$ is a symbol.
A divisor is at finite distance with respect to a dihedral structure $\delta$ if and only if $\delta$ is compatible
with the corresponding partition of $S$.
\end{prop}

\begin{proof}
%Let $P_1\cup P_2$ be any such partition of $S$. There is a  permutation  $\sigma \in \Sym_n$ such that
% $\sigma(P_1)= \{k+1,\ldots, n\}$ and $\sigma(P_2)= \{1,\ldots, k\}$.
%Viewing $\sigma(P)$ as the set of edges of an $n$-gon, this  partition is induced by cutting along the chord $\{k,n\}$ and
%corresponds to the divisor $D_{kn}$.
%Thus the partitition $P_1\cup P_2$ corresponds to the divisor $\sigma^{-1} (D_{kn}) \in \sigma^{-1} j_{\delta}( \Modf_{0,S})$, which
%decomposes as a product by lemma $\ref{lemmadecomp}$.
The bijection between stable partitions and divisors follows immediately from the previous remarks and the covering
$(\ref{Fullcompactcover})$. The last statement of the proposition holds by definition. It remains to prove the decomposition.
Suppose that we are given a stable partition $S=P_1\cup P_2$, and let $D$ denote the corresponding
divisor. There is a bijection between the set of  dihedral structures $\delta$ which are
 compatible with $P_1\cup P_2$, and the set of  induced dihedral structures $\delta_1$, $\delta_2$
on the sets $P_1\cup\{e\}$ and $P_2\cup\{e\}$ (compare fig. 3). It follows from lemma  $\ref{lemmadecomp}$ that:
$$D \cap j_\delta(\Modf_{0,S}) \cong
  \begin{cases}
    j_\delta(\Mod^{\delta_1}_{0,P_1\cup\{e\}} \times \Mod^{\delta_2}_{0,P_2\cup\{e\}}) &\text{if } \delta \text{ is compatible with } P_1\cup P_2 \ , \\
    \emptyset & \text{otherwise}\ .
  \end{cases}
$$
If we identify $j_\delta(\Mod^{\delta_1}_{0,P_1\cup\{e\}} \times \Mod^{\delta_2}_{0,P_2\cup\{e\}})$ with  $j_{\delta_1}(\Mod^{\delta_1}_{0,P_1\cup\{e\}})
 \times j_{\delta_2}( \Mod^{\delta_2}_{0,P_2\cup\{e\}})$ in
$\overline{\Mod}_{0,P_1\cup\{e\}}\times \overline{\Mod}_{0,P_1\cup\{e\}}$,
we obtain:
$$D= D\cap \overline{\Mod}_{0,S} = D \cap \bigcup_{\delta} j_\delta(\Modf_{0,S})\cong \bigcup_{\delta_1, \delta_2} j_{\delta}
(\Mod^{\delta_1}_{0,P_1\cup\{e\}} \times \Mod^{\delta_2}_{0,P_2\cup\{e\}})$$
$$ \cong \bigcup_{\delta_1} j_{\delta_1}(\Mod^{\delta_1}_{0,P_1\cup\{e\}}) \times
\bigcup_{\delta_2} j_{\delta_2}(\Mod^{\delta_2}_{0,P_2\cup\{e\}})= \overline{\Mod}_{0,P_1\cup\{e\}} \times
\overline{\Mod}_{0,P_2\cup\{e\}} \ .$$
%This completes the proof, since the last statement  holds by definition.
\end{proof}

 We introduce the following notation.  Let $D$ denote the divisor given by a stable partition $S=P_1\cup P_2$.
Then for any pair of
indices $i,j\in S$, we set
\begin{equation} \label{indicatordefn} \I_D(i,j) = \I(\{i,j\} \subset P^1) + \I(\{i,j\}\subset P^2)\ ,
\end{equation}
where $\I(A\subset B)$ is the indicator function which takes the
value $1$ if a set $A$ is contained in $B$, and $0$ otherwise.

\begin{cor} \label{corbasicordformula} Let $D$ denote the divisor corresponding to the stable partition $S=P_1\cup P_2$. The order of vanishing of any
 cross-ratio along $D$ is given
by:
$$\ord_D \,[ij|kl]
= {1\over 2} \Big[ \I_D(i,k) +\I_D(j,l)
-\I_D(i,l)-\I_D(j,k)\Big]\ .$$
%\begin{eqnarray}
%\ord_D \,[ij|kl]
%& =& \ord_D \, {(z_i-z_k)(z_j-z_l) \over
%(z_i-z_l)(z_j-z_k)} \nonumber \\
%&=& {1\over 2} \Big[ \I_D(i,k) +\I_D(j,l)
%-\I_D(i,l)-\I_D(j,k)\Big]\ . \nonumber
%\end{eqnarray}
\end{cor}
\begin{proof} The formula is invariant under the action of $\Sym(S)$ on divisors and cross-ratios.
 We can therefore fix a dihedral structure $\delta$ on $S$
and assume that  $D=D_{2a}$, where $a\in \{4,\ldots, n\}$. The formula is also compatible with $(\ref{S4})$ and  additive with
respect to $(\ref{mult})$. By lemma $\ref{lemmacrossratio}$, it therefore suffices to verify the formula for  $[p\, p+1| q+1\, q] =u_{pq},$
where $\{p,q\} \in \chi_{S,\delta}$. It follows from  $(\ref{id})$ that $\ord_{D_{2a}} \,u_{pq}$  is $1$ if $\{p,q\}=\{2,a\}$ and is $0$
otherwise. The partition corresponding to  $D_{2a}$ is $\{3,4,\ldots, a\}\cup\{a+1,\ldots,n,1\}$, and it is easy to check that the formula holds in this case.
\end{proof}

  A stable  partition $S=P_1\cup P_2$ is conveniently
represented as the union of two circles, joined at a point $e$,
with  marked points corresponding to $P_1$ on the first circle,
and those corresponding to $P_2$ on the other. Taking iterated
intersections of divisors, one obtains an operad of bubble
diagrams (fig. 11) [Dev2]. Such a diagram defines a tree. If we take the  dual graph,
 we obtain a partial decomposition of a polygon. Note that we
can  find  dihedral structures $\delta$ for which the
labellings of the outer edges are in dihedral order with respect
to $\delta$. In this way, any bubble diagram corresponds to an
intersection of divisors at finite distance on a certain number of  affine  pieces
$\Modf_{0,S}$ in $\overline{\Mod}_{0,S}$.

% One can compute the isotropy group of each stratum $D_{\alpha}$ of $\overline{\Mod}_{0,S}$
%under the action of $\Sym_n$ by similar considerations (see [Dev-R], $\S5$).

\begin{figure}[h!]
  \begin{center}
%    \leavevmode
    \epsfxsize=10.0cm \epsfbox{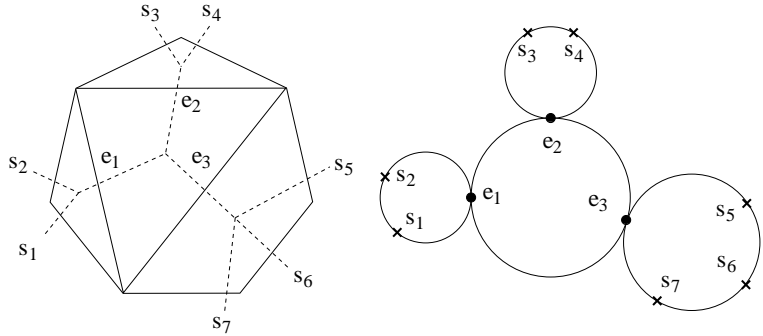}
  \label{Bubbles}
  \caption{A partial decomposition of a heptagon (left), its dual graph (dotted lines), and the corresponding bubble diagram %on
%$\overline{\Mod}_{0,7}$
(right).
 If the tree is rooted at $s_1$, this corresponds to the bracketing $(s_2,((s_3,s_4),(s_5,s_6,s_7)))$.}
  \end{center}
\end{figure}
%\noindent  In the case when there are just two bubbles, then there
%is no ambiguity in specifying the partition $S=P_1 \cup P_2$

\subsection{Product maps} \label{sect27} The projection maps on $\Mod_{0,S}$ defined above  decrease the dimension by one.
 We  will also need to   consider
various  maps between products of moduli spaces $\Mod_{0,T_i}$
which preserve the dimensions. These give rise to special
coordinate systems on $\Mod_{0,S}$ and will be used to  define
products on period integrals.
 %and
%the coordinate systems these give rise to.
 %By proposition  $1.6$, there is a canonical section
%$i_T:\Mod_{0,T} \rightarrow \Mod_{0,S}$ such that $f_T\circ i_T $
%is the identity which identifies $\Mod_{0,T}$ as an intersection
%of divisors .
Given two subsets $T_1,T_2 \subset S$ such that $|T_i|\geq 4$, we consider maps of the form
\begin{equation} \label{productmap}
f= f_{T_1}\times f_{T_2}: \Mod_{0,S} \To
\Mod_{0,T_1}\times \Mod_{0,T_2}\ .\end{equation}
Such a map will be called a \emph{product map} if
\begin{eqnarray}\label{productmapequalities}
|T_1\cap T_2|&=& 3\ , \\
S&=& T_1 \cup T_2  \ . \nonumber
\end{eqnarray}
In this case the dimensions on both sides of $(\ref{productmap})$ are equal, since
the equalities $(\ref{productmapequalities})$  imply that $|S|-3=|T_1|-3+|T_2|-3$. The map $f$ is an embedding, because we
can place the three points in $T_1\cap T_2$ at $0,1$, and $\infty$, and  each remaining marked point $s\in S$ is
then uniquely determined by the map $f_{T_i}$ where $i\in\{1,2\}$ and $s\in T_i$.
 We can  iterate this construction by further decomposing $T_i$ as a  union of sets satisfying
$(\ref{productmapequalities})$. Since the composition of two forgetful maps $f_T$ is itself a forgetful map, we obtain
a family of subsets $T_1,\ldots, T_k \subset
S$ such that $|T_i|\geq 4$, and a  map
\begin{equation} \label{productmapmany}
f=\prod_{i=1}^k f_{T_i}: \Mod_{0,S} \To \prod_{i=1}^k
\Mod_{0,T_i}\ .\end{equation}
This is an embedding by construction. The sets $T_i$ cover $S$, {\it i.e.},
 $S= \cup_{i=1}^k T_i$, and the equality of  dimensions on the left and right hand sides of $(\ref{productmapmany})$
 implies that
\begin{equation} \label{dimcond}
|S|+3\,(k-1) = \sum_{i=1}^k |T_i|\ .\end{equation}
%In this case,
%the dimensions of the left and right hand sides of
%$(\ref{productmap})$ are equal, and $f$ is a birational map.
We
can then regard  $\Mod_{0,S}$  as a %Zariski-
dense open subscheme of
$\prod_{i=1}^k \Mod_{0,T_i}$, and we say that $f$ is a
\emph{non-degenerate coordinate system} on $\Mod_{0,S}$.
%In this case, no
%set $T_i$ can be contained in the union of the other sets
%$\bigcup_{j\neq i} T_j$. For if this were so, we could suppress
%the coordinate $f_{T_i}$, and the resulting map $\prod_{j\neq i}
%f_{T_j}$ would still be an embedding,  contradicting the condition
%on the dimensions $(\ref{dimcond})$.
% The dihedral embedding $(\ref{Embed})$ is an
%example of such a map. To see this, cover $S$ with sets
%$T_{ij}=\{i,i+1,j,j+1\}$ for each chord $\{i,j\}\in
%\chi_{S,\delta}$, and identify each $\Mod_{0,T_{ij}}$ with
%$\Pro^1\backslash\{0,1,\infty\}$ using the dihedral coordinate
%$u_{ij}.$ Similarly,
Any set of vertex coordinates $\{x^\alpha_i\}$ corresponding to a triangulation
$\alpha \in \chi^\ell_{S,\delta}$, when $\alpha$ has no internal triangles, %($\S\ref{sect23}$)
is an example of  a non-degenerate coordinate system (this can be verified by induction).
 More precisely, if $x_i^\alpha= u_{p_i q_i}$, then we can cover $S$ with the sets $T_i=\{p_i, p_i+1, q_i, q_i+1\}$,
 and identify  $\Mod_{0,T_i}$ with $\Pro^1\backslash\{0,1,\infty\}$ using the coordinate $x_i^\alpha$, for $1\leq i\leq \ell$.
 If, however, $\alpha$ has internal triangles, then the map $f$ corresponding to
the set of vertex coordinates $x_i^{\alpha}$ is not an embedding, and therefore cannot be a non-degenerate coordinate system.

Let us fix a dihedral structure $\delta$ on $S$. This is
equivalent, by $\S\ref{sect24}$, to choosing one of the open cells
$X_{S,\delta}$ which cover $\Mod_{0,S}(\R)$. This induces a
dihedral structure  $\delta_i$ on each subset $T_i\subset S$,
which in turn defines a fundamental cell $X_{T_i,\delta_i}$, for
  $1\leq i\leq k$. By construction, $f_{T_i}(X_{S,\delta})
\subset X_{T_i,\delta_i}$, and therefore $f(X_{S,\delta})\subset
\prod_{i=1}^k X_{T_i,\delta_i}$. We define $G_f$ to be the set of
all  dihedral structures on $S$ which are compatible with the
dihedral structures on each $T_i$ induced by $\delta$: {\it i.e.},
\begin{equation} \label{Gdef}
G_f= \{ \gamma \in \Sym(S)/ D_{2n} \hbox{ such that }
\gamma|_{T_i} = \delta_i \}. \end{equation}
 The precise relation between
the domains $\prod_i X_{T_i,\delta_i}$ and $X_{S,\delta}$ is given
by:
\begin{equation} \label{domainrel}
 f^{-1}\Big(\prod_{i=1}^k X_{T_i,\delta_i}\Big) =
\coprod_{\gamma \in G_f} \, X_{S,\gamma}  \ .
\end{equation}
 Any point $x \in X_{S,\gamma}$ maps via $f$ into $\prod_{i}
X_{T_i,\delta_i}$ if and only if $\gamma|_{T_i}=\delta_i$.
Identity $(\ref{domainrel})$ follows because the set of  cells
$X_{S,\gamma}$, for $\gamma\in \Sym(S)/D_{2n}$, cover
$\Mod_{0,S}(\R)$ disjointly, by $\S\ref{sect24}$.

 We consider two examples of such a map $f$, one for which
$G_f$ is trivial, which gives rise to cubical coordinates, and the
other for which $G_f$ is as large as possible, which defines
simplicial coordinates. We will see later in $\S7$ that these
special cases give rise to the stuffle and shuffle relations  for
multiple zeta values, respectively.

We fix a dihedral structure $\delta$ on $S$, and write
$S=\{s_1,\ldots,s_n\}$, as usual.
 Consider first of all the
covering $S= \bigcup_{i=4}^{n} T_{i}$, where
\begin{equation}\label{cubicalcover}
T_i=\{s_2,s_3,s_i,s_{i+1}\} \qquad \hbox{for } \quad 4\leq i\leq
n\ ,
\end{equation}
 and all indices  are modulo $n$, as usual. This defines
a map
 \begin{equation} \label{fcube}  f_\Box: \Mod_{0,S} \To
\prod_{i=4}^n \Mod_{0,T_i}\ ,
 \end{equation}
 which satisfies
condition $(\ref{dimcond})$. One verifies without difficulty that this is a non-degenerate coordinate system as defined above (or
use the fact that $|T_i\cap T_{i+1}|=3$ for  $4\leq i\leq n-1$).
We call $f_\Box$ a system of
\emph{cubical coordinates} on $\Mod_{0,S}$. In this case,
$G_{f_\Box}$ is trivial, since if
 $\gamma $ is a dihedral
structure on $S$ compatible with all the dihedral structures
$s_2<s_3<s_i<s_{i+1}$ on $T_i$ (or $s_2>s_3>s_i>s_{i+1}$), then we must have $s_1<\ldots
<s_n<s_1$ (or $s_1>\ldots>s_n>s_1$). Each moduli space $\Mod_{0,T_i}\cong \Mod_{0,4}$ is
isomorphic to $\Pro^1\backslash\{0,1,\infty\}$ in six natural
ways, corresponding to the six choices of cross-ratio on
$\Mod_{0,4}$. If we  identify $\Mod_{0,T_{i}}$ with
$\Pro^1\backslash\{0,1,\infty\}$ using the coordinate $u_{2\,i}=[2\, 3|i+1\, i]$, for $4\leq i\leq n$,
then, since $x_{i-3}=u_{2\,i}$ by  $(\ref{tcoords})$, we
retrieve the explicit cubical coordinates defined in $\S
\ref{sect21}$. In other words, $(\ref{fcube})$ is just
$$f_{\Box}=(x_1,\ldots, x_{\ell}): \Mod_{0,S} \rightarrow
\big(\A^1\backslash\{0,1\}\big)^{\ell} \ ,$$  and coincides with
$(\ref{ExpCubeCoords})$.
 Each cell $X_{T_i,\delta_i}$ is the unit
interval $(0,1)$ in these coordinates, and $X_{S,\delta}$ maps
under $f_{\Box}$ to $(0,1)^\ell$. In this case, equation
$(\ref{domainrel})$ simply states that a product of $\ell$ unit
intervals is the unit $\ell$-dimensional  hypercube.

  %We can deduce a product
%from the cubical map $f_\Box$ as follows.
Cubical coordinates come from a product map.  If $k\geq 4$, we
set
$$S_1= \bigcup_{i=4}^k T_i = \{s_2,s_3,\ldots, s_{k+1}\} \quad \hbox{ and }
\quad  S_2= \bigcup_{i=k+1}^n T_i = \{s_{k+1},\ldots,
s_n,s_1,s_2,s_3\}\ .$$  Setting $m=k-3$,  we define the
\emph{cubical product map} $m_{\Box}$ to be
\begin{eqnarray}
 \label{cubeprod} m_{\Box} =f_{S_1}\times f_{S_2}: \Mod_{0,S} &\To& \Mod_{0,S_1}\times
\Mod_{0,S_2} \\
(x_1,\ldots, x_\ell )&\mapsto &\big((x_1,\ldots, x_m),
(x_{m+1},\ldots, x_\ell)\big) \nonumber  \ .
\end{eqnarray}

%
%We now assume, by the triple transitivity of the action of
%$\PSL_2(\C)$, that $z_1=1$, $z_2=\infty$, $z_3=0$, and write down
%an explicit system of cubical coordinates.

%and hence by , $X_S\cong \prod_{i=1}^{n-3} X_{T_i}$.
Now, the simplicial case arises by  considering the covering $S=\bigcup_{i=4}^n T_i,$ where
\begin{equation}
T_i=\{s_1,s_2,s_3,s_{i}\} \qquad \hbox{for } \quad 4\leq i\leq n\
.\end{equation}
 This defines a map
\begin{equation}\label{ftriangdef}
f_\triangle: \Mod_{0,S} \To \prod_{i=4}^n \Mod_{0,T_i}\ ,
\end{equation}
 which satisfies condition $(\ref{dimcond})$ and is a non-degenerate coordinate system for the same reasons as above (namely
 $|T_i\cap T_j|=3$ for all $i\neq j$).
We call $f_{\triangle}$ a \emph{system of simplicial coordinates}
on $\Mod_{0,S}$.
 It is easy
to check that   $G_{f_\triangle}$ is  bijective to the symmetric group
on $\ell$ letters
$$G_{f_\triangle} = \Sym(\{s_4,\ldots, s_n\})\ .$$
As above, we obtain an explicit set of simplicial coordinates by
choosing the coordinate $t_i=[i+3\, 1| 3\, 2] : \Mod_{0,T_{i+3}}
\cong \Pro^1\backslash\{0,1,\infty\}$, for $1\leq i\leq \ell$.  Thus $(\ref{ftriangdef})$
 can be written
$$f_\triangle=(t_1,\ldots, t_{\ell}) : \Mod_{0,S}
\rightarrow \big(\A^1\backslash\{0,1\}\big)^{\ell}\ ,$$ and we retrieve the
isomorphism $(\ref{ExpSimpCoords})$.
% and gives an isomorphism
%$$f_\triangle: \Mod_{0,S} \cong \{ t_1,\ldots, t_\ell \in
%(\C\backslash\{0,1\})^{\ell} : t_i\neq t_j \hbox{ for all } i\neq
%j\}\ .$$
 %A simple calculation shows that $u_{13}=1-t_1$, and
%$u_{2n}=t_\ell$, and for $4\leq i+3\leq n-1$,
%\begin{equation} \label{tcoords}
% u_{1 \, i+3} = {1-t_{i+1} \over
%1-t_i},  \quad \hbox{and} \quad t_{2 \,i+3} = {t_i \over t_{i+1}}
%\ .
%\end{equation}
 As before, the domains $X_{T_i,\delta_i}$ map to  unit
intervals $(0,1)$ under $t_{i-3}$, and $X_{S,\delta}$ maps
bijectively under $f_\triangle$  to the unit simplex
$\{0<t_1<\ldots<t_{\ell}<1\}$. In this case, equation
$(\ref{domainrel})$ states that
$$\coprod_{\sigma \in
\Sym(\{\sigma_4,\ldots, \sigma_n\})}  f_{\triangle} (X_{S,\sigma}) = \prod_i X_{T_i}\ ,$$
{\it i.e.},  the unit cube $(0,1)^{\ell}$ is tesselated with
$\ell!$ copies of the unit simplex.

Now let  $k\geq 4$, $m=k-3$,  and set
$$S_1= \bigcup_{i=4}^k T_i = \{s_1,s_2,s_3,\ldots, s_{k-1},s_k\} \quad \hbox{ and }
\quad  S_2= \bigcup_{i=k+1}^n T_i = \{s_1,s_2,s_3, s_{k+1},\ldots,
s_n\}\ .$$ We define the \emph{simplicial product map}
$m_{\triangle}$ to be
 \begin{eqnarray}
\label{simpprodmap} m_{\triangle}=f_{S_1}\times f_{S_2} :
\Mod_{0,S} &\To& \Mod_{0,S_1}\times \Mod_{0,S_2} \\
(t_1,\ldots, t_\ell) & \mapsto & \big((t_1,\ldots, t_m),
(t_{m+1},\ldots, t_\ell)\big) \nonumber \ .
\end{eqnarray}
In this case, the set $G_{m_{\triangle}}$ is exactly the set
$\Sym(m,\ell-m)$ of all possible ways of shuffling together the
sets $\{s_4,\ldots, s_k\}$ and $\{s_{k+1},\ldots, s_n\}$ whilst
preserving the orderings $s_4<\ldots <s_k$ and $s_{k+1}<\ldots
<s_n$. In explicit simplicial coordinates, which involves setting
$s_1=1$, $s_2=\infty$, $s_3=0$, and $s_{i+3}=t_i$ for $1\leq i\leq
\ell$, equation $(\ref{domainrel})$  is the well-known
formula for the decomposition of a product of simplices:
\begin{equation}\label{simplexdecomp}
\{0<t_1<.\,.<t_m<1\} \times \{0<t_{m+1}<.\,.< t_\ell <1\}\cong
\end{equation}$$
\!\!\!\!\!\!\coprod_{\sigma \in \Sym(m,\ell-m)}
\!\!\!\!\!\!\{0<t_{\sigma(1)} <.\,. < t_{\sigma(\ell)}<1\} \ .
$$

\newpage

\section{The reduced bar
construction and Picard-Vessiot theory}

The main tool for computing the periods of moduli spaces is a
triviality theorem for the cohomology of a variant of the bar construction on the de Rham complex of
$\Mod_{0,S}$. Although many of the results below hold in considerably greater generality,
we %restrict to  the case of
consider the  complement of an affine hyperplane arrangement $M$,
which is more than adequate for our purposes. We first show that
the reduced bar construction on $\Omega^\star(M)$ defines a
Picard-Vessiot extension of its ring of regular functions. This is
an abstract algebraic analogue of the ring of iterated integrals
over $M$. Then, by showing that the bar construction decomposes as
a tensor product over a fibration, we prove that the cohomology of
the bar construction is trivial for fiber-type
arrangements. %of hyperplanes.
This result is also proved for quadratic arrangements in the
appendix. Our point of view, using differential Galois theory, is
different from classical approaches to this subject [Ha1-2,Ch1-2,
Ko]. The main technical idea is the notion of unipotent extensions
of differentially simple algebras, which is developed in $\S
\ref{sect33}$. The example $\Mod_{0,5}$ is discussed in
$\S\ref{sect38}$.

\subsection{Shuffle algebras and non-commutative formal power series} \label{sect31}
Let $R$ be a commutative unitary ring. Let $k\geq 1$, let
$A=\{a_1,\ldots,a_k\}$ denote an alphabet with $k$ symbols, and
let $A^*$ denote the free non-commutative monoid generated by $A$,
{\it i.e.}, the set of all words $w$ in the symbols $a_i$, along
with the empty word $1$. Let
 $R\langle A \rangle $ be the
free non-commutative $R$-algebra generated by $A$. %{\it i.e.,}  the
%free $R$-module generated by all words $w$ in the symbols $a_i$
%along with empty word
% $1$.
 If  $V_1$ is the free $R$-module with basis
 $\{a_1,\ldots, a_k\}$,  and if we set $V_m= V_1^{\otimes m}$ and $V_0=R$, then clearly
$$R\langle A\rangle = \bigoplus_{m\geq 0} V_m\ .$$
 It is well-known [Bou] that $R\langle A\rangle $
  can be given the structure of a cocommutative
graded Hopf algebra. The multiplication law on
$R\langle A \rangle$ is given by  concatenation of words,
%(or the
%tensor product)
 and the coproduct $ \Gamma : R\langle A\rangle
\rightarrow R\langle A\rangle\otimes R\langle A\rangle$ is
defined to be the unique coproduct for which the elements of $A$
are all primitive:
\begin{eqnarray}
\Gamma(a_i) = a_i \otimes 1 + 1\otimes a_i \ . \nonumber
\end{eqnarray}
The counit $\varepsilon:R\langle A \rangle \rightarrow R$ is
given by projection onto the unit word $1$. If $|w|$ denotes the
number of symbols occuring in a word $w\in A^*$, then the antipode
map is defined by  $w\mapsto (-1)^{|w|} \widetilde{w}$, where the
mirror map $w\mapsto \widetilde{w}$ reverses the order of the
symbols in each word. One verifies that this defines a graded Hopf
algebra structure.% on $R\langle A\rangle$.

Let $V_1^\vee$ denote the $R$-module dual to $V_1$, and let
$A'=\{a_1',\ldots, a_k'\}$ denote the basis dual to $A$. Then
$R\langle A'\rangle$, the free tensor algebra over $V_1^\vee$, is
the
    graded dual of $R\langle A\rangle$,
 and inherits a commutative Hopf algebra
   structure by duality. The multiplication law is now given by the shuffle product $\sha:
 R\langle A'
\rangle\otimes R\langle A'\rangle \rightarrow  R\langle A' \rangle $ which
 is defined  recursively by the formulae: $w \sha 1 = 1\sha w=  w,$ and
\begin{eqnarray}\label{sha}
a'_i w_1 \sha a'_j w_2& =& a'_i ( w_1\sha a'_j w_2) + a'_j(a'_i w_1 \sha
w_2) \ ,\end{eqnarray} for all words $w_1,w_2\in A'^*$, and all
$a'_i,a'_j\in A'$.
 This is a commutative, associative product with no zero divisors.
 The algebra $(R\langle A'\rangle, \sha)$, will be called the \emph{free shuffle algebra} on the generators
 $a'_1,\ldots, a_k'$.
The coproduct is defined by the map
\begin{eqnarray}
\label{coprod} \Delta: R\langle A' \rangle& \rightarrow &
R\langle A' \rangle \otimes R\langle A' \rangle \\
\Delta(w ) &= &\sum_{uv=w} u\otimes v \nonumber \ ,
\end{eqnarray}
 and the antipode is given by the map $w\mapsto (-1)^{|w|} w$. The
counit $\varepsilon: R\langle A'\rangle\rightarrow R$ is given
by projection onto the graded part of  weight $0$, as previously.

Let $R\langle \langle A \rangle \rangle$ and $R\langle \langle
A' \rangle \rangle$ denote the completions of the graded algebras
defined above with respect to the augmentation ideals $\ker
\varepsilon$. These are just the algebras of
formal power series in $A$, $A'$ respectively. The Hopf algebra
structures $\Delta, \Gamma$, and  $\varepsilon$  extend in the
natural way to the completed algebras, and we shall denote them by
the same symbols.

 In addition, we introduce
$k$ truncation operators $\partial_{a'_i}$ for $1\leq i\leq k$:
\begin{eqnarray}
\partial_{a'_i} : R\langle A' \rangle &\rightarrow& R\langle A'
\rangle \\
\partial_{a'_i} (a'_j w) &=& \delta_{ij} w\ , \nonumber
\end{eqnarray}
for all  $a'_j\in A'$, $w\in A'^*$,  where $\delta_{ij}$ is the
Kronecker delta. It is easy to verify that the $\partial_{a'_i}$
are derivations for the shuffle product, and furthermore, that
this  determines the shuffle product uniquely if we assume that $1$ is the unit. The
operators $\partial_{a'_i}$ %can also be defined by projecting
%graded pieces of  the coproduct, and
 are dual to the  operators $w \mapsto a_i w:R\langle A\rangle \rightarrow R\langle A\rangle $ which affix the letter $a_i$ to the left of words $w\in A^*$.
That $\partial_{a'_i}$ is a derivation is equivalent to the fact that  $a_i$ is primitive for the coproduct $\Gamma$ by duality.

 %.which affix the letter  $a_i$ for the concatenation product.

\subsection{Arrangements of hyperplanes and the bar construction}\label{sect32}
Consider an arrangement of $N$ hyperplanes $H_1,\ldots, H_N $ in
affine space  $\A^\ddim$. Let $k$ denote a field of
characteristic $0$ over which the arrangement is defined. For each $1\leq i\leq N$,
choose a linear form $\alpha_i\in k[x_1,\ldots,x_\ddim]$ such that
$H_i$ is the divisor of zeros of $ \alpha_i$. Let
$$\Or_M= k \big[ x_1,\ldots, x_\ddim,  \{\alpha_i^{-1}\}_{1\leq i\leq N} \big]$$
denote the ring of regular functions on the complement $M=\A^\ddim_k
\backslash \bigcup_i H_i$. We set
$$d = \sum_{i=1}^\ddim {\partial \over \partial x_i} dx_i \ .$$
Consider the de Rham complex of $\Or_M$:
\begin{equation}
0\To \Or_M\overset{d}{\To} \Omega^1 (\Or_M) \overset{d}{\To}  \Omega^2 (\Or_M) \overset{d}{\To}
\ldots \overset{d}{\To} \Omega^\ddim (\Or_M) \To 0\ ,
\end{equation}
 where $\Omega^r (\Or_M)= \bigoplus_{1\leq i_1<\ldots<i_r\leq N} \Or_M
dx_{i_1}\wedge\ldots\wedge dx_{i_r}$ is placed in degree $r$.
 Let $H^i(\Or_M)$, for $0\leq i\leq \ell$, denote the corresponding cohomology groups. These are $k$-vector spaces.
Since $M$ is affine, it follows that $H^i(\Or_M)$ coincides with
the de Rham hypercohomology of $M$ [Gr]. Consider the set of
algebraic $1$-forms:
\begin{equation}
\omega_i = d\log \alpha_i \ ,\quad \hbox{ for } 1\leq i\leq N\
.\end{equation}
The following theorem is due to Arnold and Brieskorn [O-T, \S5.4].

\begin{thm}  \label{thmarnold}
The cohomology ring $H^\star(M)$  is isomorphic to the graded
$k$-algebra $A$ generated by the forms $\omega_i$, for $1\leq i\leq N$.
\end{thm}

\noindent In particular, the cohomology classes of the forms  $\omega_1,\ldots, \omega_N \in
\Omega^1(\Or_M)$ are a  $k$-basis for $H^1(\Or_M)$. In this section, all tensor products will be taken over the field $k$ unless specified
otherwise.  Let $N$ denote the
kernel of the exterior
 product
$$ N= \ker \big(\wedge: H^1(\Or_M) \otimes H^1(\Or_M) \To H^2(\Or_M)\big) \ .$$
We will not require the full strength of theorem $\ref{thmarnold}$, only the following corollary.
\begin{cor}
If a form $\omega \in A$ is a coboundary $d\phi$, then it is zero.
\end{cor}
It follows that $N$
  is also the kernel of the
map
$$\wedge :\bigoplus_{1\leq i,j\leq N} k\, \omega_{i}\otimes \omega_j \To \Omega^2(\Or_M)\ .$$
For each positive integer $m\geq 2$, the vector space $V_m(\Or_M)$ of
integrable words in the forms $\omega_i$ of weight $m$ is defined
to be
\begin{equation} \label{Vmdef}
V_m(\Or_M) = \bigcap_{i+j=m-2} H^1(\Or_M)^{\otimes i} \otimes N \otimes
H^1(\Or_M)^{\otimes j} \ .\end{equation} This is just the intersection of the  kernels of
the maps $\wedge_i$ for $1\leq i\leq m-1$:
\begin{eqnarray}
\wedge_i: H^1(\Or_M)^{\otimes m} & \To & H^1(\Or_M)^{\otimes i-1} \otimes
H^2(\Or_M) \otimes
H^1(\Or_M)^{\otimes m-i-1}\ , \\
\eta_1\otimes\ldots \otimes  \eta_m & \mapsto &
\eta_1\otimes\ldots\otimes (\eta_{i} \wedge \eta_{i+1}) \otimes
\ldots \otimes \eta_m \nonumber \ .
\end{eqnarray}
%By the theorem above, i
Its elements can be written as  linear
combinations of symbols %involving only the forms $\omega_i$:
$$\sum_{I=(i_1,\ldots,i_m)} c_I [\omega_{i_1}|\omega_{i_2}|\ldots
|\omega_{i_m}]\ ,$$ where $1\leq i_j \leq N$, and $c_I \in k$,
which satisfy the integrability condition:
\begin{equation} \label{intcond}
\sum_{I=(i_1,\ldots,i_m)} c_I\, \omega_{i_1}\otimes\ldots\otimes
\omega_{i_{j-1}}\otimes(\omega_{i_{j}}\wedge
\omega_{i_{j+1}})\otimes\omega_{i_{j+2}}\otimes\ldots
\otimes\omega_{i_m}= 0\ , \end{equation} %in   $H^1(\Or_M)^{\otimes k-1}\otimes H^2(\Or_M) \otimes  H^1(\Or_M)^{\otimes m-k-2},$
 for each
$1\leq j\leq m-1.$ We set $V_0(\Or_M) = k,$  and  $V_1(\Or_M)= H^1(\Or_M)=\bigoplus_{i=1}^N k\, \omega_i,$
 and define
\begin{equation}\label{BigVdefn}
V(\Or_m) = \bigoplus_{m\geq 0} V_m(\Or_M)\ .\end{equation}
 The vector space of homotopy-invariant iterated integrals is then
defined to be
\begin{equation} \label{bardef}
 B(\Or_M) =  \,\Or_M\otimes V(\Or_M)\ ,
\end{equation}
with the obvious grading. This is  similar to the zero$^\mathrm{th}$
cohomology group of Chen's reduced bar complex on $\Or_M$, which is
usually written $H^0(B(\Omega^\bullet \Or_M))$, with the difference that it is made
up  of closed $1$-forms only (see [Ch1, Ha1, Ha2, Ko]). %Since we shall only consider integrable
%words, we will omit the $H^0$ to simplify notations.

In order to define a differential on $B(\Or_M)$, we let
$$\Omega^i B(\Or_M)= \Omega^i(\Or_M) \otimes_{\Or_M} B(\Or_M)= \Omega^i(\Or_M)\otimes_k V(\Or_M)\ ,$$
and define $d: \Omega^i B(\Or_M)\rightarrow \Omega^{i+1} B(\Or_M)$ by the formula
\begin{eqnarray} \label{bardifdef}
d\, \sum_{I=(i_1,\ldots,i_m)} \phi_I\otimes
[\omega_{i_1}|\omega_{i_2}|\ldots |\omega_{i_m}]  & =&
\sum_{I=(i_1,\ldots,i_m)} (-1)^{\deg \phi_I} \phi_I \wedge \omega_{i_1}\otimes  [ \omega_{i_2}|\ldots
|\omega_{i_m}] \nonumber  \\
& + & \sum_{I=(i_1,\ldots,i_m)} d\phi_I \otimes [ \omega_{i_1}|\omega_{i_2}|\ldots
|\omega_{i_m}]\ ,
\end{eqnarray}
where $\phi_I \in \Omega^i(\Or_M)$. It follows from
  the integrability condition $(\ref{intcond})$ that
 $d^2=0$.  We can therefore consider the following complex
\begin{equation} \label{BardeRham}
 0\To B(\Or_M)\overset{d}{\To} \Omega^1 B(\Or_M)\overset{d}{\To} \Omega^2
B(\Or_M) \overset{d}{\To} \ldots \overset{d}{\To} \Omega^{\ddim} B(\Or_M) \To 0 \ ,
\end{equation}
 where $\Omega^i B(\Or_M)$ is placed in degree $i$. Its cohomology  will be denoted $H^i_\DR(B(\Or_M))$.
 By definition $(\ref{BigVdefn})$, $ V(\Or_M)$ is contained in the free $\Or_M$ shuffle algebra generated by $V_1(\Or_M)$, which
 is a commutative graded Hopf algebra:
$$V(\Or_M) \subset k \langle \omega_1,\ldots, \omega_N\rangle\ .$$
 The product $\sha$ is the shuffle product defined in $(\ref{sha}),$
and the coproduct $\Delta$ was defined in $(\ref{coprod})$. One can verify that $
V(\Or_M)$ is preserved by $\sha$ and $\Delta$, and is therefore a
graded Hopf subalgebra of $\Or_M\langle \omega_1,\ldots, \omega_N
\rangle$.
\begin{cor}
$B(\Or_M)$  is a commutative graded algebra for the shuffle product $\sha$, and has a natural coproduct
 $\Delta: B(\Or_M)\rightarrow B(\Or_M)\otimes_{\Or_M} B(\Or_M)$.
\end{cor}

\subsection{Unipotent extensions of differentially simple algebras}\label{sectPV}
%In order to state the main properties of the reduced bar construction $B(\Or_M)$,
%       we need some facts about unipotent extensions of differentially simple algebras.

%In this section we develop the theory of unipotent differential algebras over differentially simple rings
%Example to bear in mind is bar construction.

Let $k$ be a field of characteristic 0, and let $R$ denote a commutative, unitary  $k$-algebra with $\mdim$ commuting derivations
$\partial_1,\ldots, \partial_{\mdim}$.
Its  de Rham  complex begins as follows:
$$0 \To R {\To} \bigoplus_{1\leq i\leq \mdim} R \To \bigoplus_{1\leq i<j\leq \mdim} R \To \ldots \ ,$$
where the first map is given by $f\mapsto (\partial_if)_i$, and the second map sends $(f_1,\ldots, f_\mdim)$
to $(\partial_i f_j -\partial_jf_i)_{i<j}$.
The ring of constants of $R$  is the $k$-algebra:
$$ H^0(R) = \bigcap_{i=1}^\mdim \ker \partial_i\ .$$

\begin{defn}
We say that $R$ is  \emph{differentially simple} if  $H^0(R)=k$, and if $R$ is a simple module over its ring of
differential operators
$R[\partial_1,\ldots, \partial_\mdim]$.
\end{defn}

Recall that a differential ideal of  $R$ is
an ideal $I\subset R$ such that $\partial_i I \subset I$ for all $1\leq i\leq \mdim$.
 It is immediate that $R$ is differentially simple if and only if it has
  no  differential ideals apart from $0$ and
$R$. An equivalent condition is that for every non-zero $r\in R$, there exists an operator $D_r\in R[\partial_1,\ldots, \partial_\mdim]$
such that $D_r \, r=1$.
 This is the analogue of the notion of a  field in
differential algebra.

Now let us assume that $R$ is differentially simple.
Let $B$ be a differential $k$-algebra containing  $R$, with differentials we also denote by $\partial_1,\ldots, \partial_\mdim$.
\begin{defn}\label{defnunipotent}
We say that $B$ is
\emph{unipotent} if  $H^0(B)=k$, and if there exists a filtration by $R[\partial_1,\ldots, \partial_\ell]$-subalgebras $W^i B$ of $B$:
$$R=W^0B \subset W^1 B \subset \ldots \subset W^{i+1} B \subset \ldots \subset B\ ,$$
%where $W^iB$ are $R[\partial_1,\ldots, \partial_\ell]$ subalgebras of $B$
 such that $B=\bigcup W^i B$, and $W^{i+1}B $ is generated, as an algebra over $W^i B$, by finitely many elements $y$ such that $\partial_1y,\ldots ,\partial_\mdim y \in W^{i}B$.
\end{defn}
In other words, $B$ is obtained by adding successive primitives to $R$ with respect to the operators $\partial_1,\ldots, \partial_\mdim$.
The following lemma is a variant of a well-known
result concerning extensions of differential fields by adjoining
primitives.

\begin{lem}
Let $R$ be a  differentially simple $k$-algebra, and let $r_1,\ldots, r_\mdim \in R$ such that
$\partial_i r_j = \partial_j r_i$ for all  $1\leq i,j\leq \mdim.$
On the polynomial ring  $R[y]$, we  extend the derivations $\partial_1,\ldots, \partial_\mdim$ by setting
$$\partial_i y = r_i \in R \qquad \hbox{for } 1\leq i\leq \mdim\ . $$
The extended operators $\partial_i$ commute and are unique. Suppose that no
element $u\in R$ satisfies $\partial_i u= r_i$ ({\it i.e.}, the class of $(r_1,\ldots, r_\ell)$ is non-zero in $H^1(R)$). Then
%$y$ is transcendent over $R$, and
$R[y]$ is differentially simple.
\end{lem}
\begin{proof}
Let $I$ be a differential ideal in $R[y]$, and suppose that
$f(y)\in I$ is a  polynomial in $y$  of
minimal degree $n\geq 1$:
$$f(y)=a_ny^n + a_{n-1}y^{n-1} + \ldots + a_0 \in I\ ,$$ where $a_i
\in R$, $a_n\neq 0$. Since $R$ is differentially simple, there
exists an operator $D\in R[\partial_i]$ such that $D\,
a_n=1$. After applying this operator to the equation above, we may
assume that $a_n=1$. On applying $\partial_i$, we obtain
$$(n r_i +\partial_i a_{n-1})y^{n-1} + \ldots + (a_1 r_i+ \partial_i a_0)\in I\ . $$
By the minimality of $f(y)$, this polynomial is identically $0$, so the set of
equations $
\partial_i u = r_i$, for $1\leq i\leq \mdim$, already has a solution $u=-a_{n-1}/n \in
R$. This contradicts the assumption, and proves that
%
%$y$ is
%transcendent. Now suppose that $I \subsetneq R[y]$ is a non-trivial
%differential ideal, and let $f(y) \in I$ of minimal degree.
%Exactly the same computation as above shows that the set of
%equations $\partial_i u=r_i$ must already have had a solution in
%$R$. It follows that
 $R[y]$ has no non-trivial differential
ideals.
\end{proof}

\begin{rem}
Since $R[y]$ is differentially simple, it has no non-trivial
quotients. For  any differential $R$-algebra $R[\eta]$, where
$\partial_i \eta \in R$ for $1\leq i\leq \ell$, and $\eta$
satisfies the conditions of the lemma, the element $\eta$ is
therefore transcendental.
\end{rem}
\begin{cor} \label{corunipispolynomial} Let $B$ denote a unipotent extension of $R$, where $R$ is differentially simple.
 Then $B$ is a polynomial algebra, and
every differential $R$-subalgebra of $B$ is
 differentially
simple.
\end{cor}
 \begin{proof} Let $A$ denote a differential $R$-subalgebra of $B$.  We can formally add primitives $y_1,\ldots,
y_p,\ldots$ to $R$, where $y_p\in A,$ to obtain a sequence of
differential algebras
$$R\subset R[y_1]\subset R[y_1,y_2]\subset \ldots \subset A=R[y_1,\ldots, y_p,\ldots ]\ .$$
We can assume that each inclusion is strict, {\it i.e.}, $y_{p+1}
$ is not in $R[y_1,\ldots,y_p]$ for each $p\geq 0$. Let
$$\partial_i y_{p+1}=r_{p+1,i} \in R[y_1,\ldots,y_p] \ .$$
 Since the ring of
constants of $B$ is $k$, it follows that the
primitive $y_{p+1}$ is the unique solution to the equations
$\partial_iu =r_{p+1,i}$ for $1\leq i\leq \mdim$ in $B$,
up to some constant in $k$. Applying the previous lemma
inductively, we deduce that $R[y_1,\ldots, y_p]$ is differentially
simple and pure transcendent for all $p\geq 1$. It follows that
$A$ is differentially simple, and that $A$ is a polynomial
algebra.
 \end{proof}

\begin{defn}
Let $B$ denote a unipotent extension of a differentially simple $k$-algebra $R$.
We say that $B$ is a \emph{unipotent closure} of $R$ if%, furthermore, it satisfies
$$H^0(B)=k \ , \quad \hbox{ and  } \quad H^1(B)=0\ .$$
\end{defn}
A unipotent closure is closed under the operation of taking 1-primitives: for all $f_1,\ldots,f_\mdim\in B$
such that $\partial_i f_j = \partial_j f_i$ for all $1\leq i,j\leq \ell$, there exists a primitive $F\in B$ such that
$\partial_1 F=f_1$, $\ldots$ , $\partial_\ell F=f_\ell$. % for all $1\leq i\leq d$.
%$$ \forall f_i \in B \hbox{ such that } \partial_i f_j= \partial_j f_i \hbox{ for all } 1\leq i,j\leq d, \hbox{ there exists }$$

\begin{defn}  A \emph{pointed differential $k$-algebra} $(R,\varepsilon)$  is a differential $k$-algebra $R$ and a $k$-linear homomorphism
of algebras $\varepsilon: R\rightarrow k$. Now suppose  that $R$ is differentially simple. We define
$\Up(R,\varepsilon)$ to be  the category of unipotent pointed extensions of $(R,\varepsilon)$.
Its objects are $(B,\varepsilon')$, where $B$ is a unipotent extension of $R$, such that the composition
$R\rightarrow B \overset{\varepsilon'}{\rightarrow }k $ coincides with $\varepsilon:R \rightarrow k$.
%, and $\varepsilon:B\rightarrow k$
%satisfies $\varepsilon=\varepsilon'\circ i: R\rightarrow k$.
 A morphism $\phi$
from $(B_1,\varepsilon_1)$ to $(B_2,\varepsilon_2)$,  is given by a commutative diagram:
$$\begin{array}{ccc}
     & R &   \\
   \swarrow & \,\,| &    \searrow \\
\!\!\!\! \phi: B_1 \quad &\To&  \quad B_2 \\
 \varepsilon_1 \searrow  &\!\!\!\varepsilon \downarrow & \!\!\!\!\varepsilon_2\swarrow   \\
    & k &
\end{array}
$$

\end{defn}

Any object $B\in \Up(R,\varepsilon)$ is differentially simple by the previous corollary. It follows that morphisms in
$\Up(R,\varepsilon)$ are necessarily injective.

\begin{lem}\label{lemmorphismsunique} Morphisms in $\Up(R,\varepsilon)$ are unique.
\end{lem}
\begin{proof}
Consider two morphisms $\phi,\phi': (B_1,\varepsilon_1)\rightarrow (B_2,\varepsilon_2)$ of pointed unipotent algebras over $(R,\varepsilon)$.
If $\partial_i$, for $1\leq i\leq \ell$, are the differentials on $R$, we denote their extensions to $B_1$ and $B_2$ by the same symbols.
Let $W^\bullet \, B_1$ denote a filtration on $B_1$ as in  definition $\ref{defnunipotent}$, and suppose by induction that
$\phi=\phi'$ on $W^p B_1$. Let $y\in W^{p+1}B$ such that $\partial_i y\in W^p B_1$ for all $1\leq i\leq \ell$. Then
$$\partial_i(\phi-\phi')(y)=(\phi-\phi') (\partial_i y) = 0\ ,\qquad \hbox{ for all } 1\leq i\leq \ell\ ,$$
and therefore $\phi(y)-\phi'(y)\in H^0(B_2)=k.$ Since $\varepsilon_2\phi(y)=\varepsilon_2\phi'(y)=\varepsilon_1 y$, it follows that
$\phi(y)=\phi'(y)$. Thus $\phi=\phi'$ on $W^{p+1}B$ and the uniqueness follows by induction.
\end{proof}

\begin{prop} \label{propPVproperties}
Let $(R,\{\partial_i\}_{1\leq i\leq \mdim},\varepsilon)$ and $(R', \{\partial'_i\}_{1\leq i\leq \mdim},\varepsilon')$ denote two differentially simple pointed
$k$-algebras, and let $\phi:(R,\varepsilon) \To (R',\varepsilon')$
be  a non-zero differential homomorphism.  Let $(U, \varepsilon')$
be a unipotent closure of $(R', \varepsilon')$, and let $(B, \varepsilon)$ be any unipotent extension of $(R,\varepsilon)$.
Then there is a unique morphism of differential algebras  $\phil : B \To U$  which extends $\phi$ and which  is necessarily injective,
 such that the following diagram commutes:
$$
\begin{array}{ccc}
  R & \xrightarrow{\,\,\phi\,\,} & R' \\
  \downarrow &  & \downarrow \\
  B & \xrightarrow{\,\,\phi_\star\,\,} & U \\
  \quad \varepsilon\searrow &  & \swarrow\varepsilon' \quad \\
   & k &
\end{array}
$$
 The map $\phil$
preserves any given filtrations on $B$ and $U$, {\it i.e.}, $\phi(W^p B)\subset W^p U$ for all $p\geq 0$.
 If, furthermore, $H^1(B)=0$ and $\phi$ is an isomorphism, then
$\phil$ is also an isomorphism.
\end{prop}
\begin{proof}
Suppose by induction that $\phil$ has been defined on $W^p B$. Since $B$ is unipotent, $W^{p+1} B$ is generated by elements
$y$ such that $\partial_i y \subset W^pB$. For such  $y$,
$$\partial_j' \phil(\partial_i y) = \phil (\partial_j \partial_i y) = \phil (\partial_i \partial_j y)
= \partial_i'\phil (\partial_j y)\ .$$
Since $H^1(U)=0$, there exists $f\in U$ such that $\partial_i'f =\phil(\partial_i y)$ for all $1\leq i\leq \mdim$.
We extend the definition of $\phil$ by setting $\phil(y) = f+ k_y$, where the constant of integration $k_y\in k$ is chosen such that $k_y +\varepsilon'(f)= \varepsilon(y)$.
By the previous lemma, $y$ is transcendent, and therefore $\phil$ is well-defined.
We obtain a   map $\phil$ on the whole of $B$ by induction. The previous corollary implies that $B$ is differentially simple.
It follows that $\phil$ is
injective because its kernel is a differential ideal in $B$ not equal to $B$ itself. The fact that $\phil$ preserves the  filtrations is clear from the construction.

Now suppose that $\phi$ is an isomorphism and that $H^1(B)=0$. Applying the same construction to $\phi^{-1}$, we obtain a map
$(\phi^{-1})_*: U\rightarrow B$. Because of the uniqueness of morphisms in $\Up(R,\varepsilon)$, $ \phil  (\phi^{-1})_*$ is the identity,
and therefore $\phil$ is an isomorphism.
%
%Suppose by induction that $\phil$ is surjective on elements of weight $p$,
%{\it i.e.}, $W^p U \subset  \phil(B)$. For any element $y\in W^{p+1}U$ such that $\partial_1'y,\ldots \partial'_\mdim y\in W^p U
% \subset \phil(B)$,
%there exist elements $x_1,\ldots, x_\mdim\in B$ such that
%$$\phil(x_i) = \partial'_i y \qquad \hbox{ for all } 1\leq i\leq \mdim\ .$$
%But $\phil(\partial_j x_i) = \partial_j' \partial_i' y = \partial_i' \partial_j' y= \phil(\partial_i x_j)$, and therefore
%$\partial_jx_i=\partial_ix_j$ for all $1\leq i,j\leq \mdim$ by the injectivity of $\phil$. Since, by assumption, $B$ contains 1-primitives, there exists
%$f\in B$ such that $\partial_i f =x_i$, and hence
%$$\partial_i' \phil(f) = \partial'_i y \qquad \hbox{for all } 1\leq i\leq \mdim\ .$$
%It follows that $\phil(f)-y\in k$, since primitives are unique in $U$ up to constants, and therefore $y\in  \phil(B)$.
%Such elements $y$ span $W^{p+1} U$ because $U$ is unipotent, and it follows by induction that $\phil:B\rightarrow U$ is
%an isomorphism.
\end{proof}

\begin{cor} \label{corBunipclos} A pointed unipotent closure $(U,\varepsilon)$ over $(R,\varepsilon)$ is a
 final object in $\Up(R,\varepsilon)$, {\it i.e.},
every unipotent extension $(U,\varepsilon)\rightarrow (B,\varepsilon')$ is an isomorphism.
\end{cor}
We can therefore speak of the  unipotent closure $U$ of a pointed differentially simple ring $(R,\varepsilon)$ whenever it exists.
Since $U$ is a union of polynomial algebras (corollary $\ref{corunipispolynomial}$),
 it necessarily has   a $k$-valued point over $\varepsilon$.

%\begin{rem} \label{derivext}
% Let $(R,\varepsilon)$ be a differentially simple $k$-algebra  and  let $(U,\varepsilon)$ denote its unipotent closure.
% The  proof of the proposition also  shows that if $d$ is a $k$-linear derivation on $R$ which commutes with the differential
%operators
%$\partial_1,\ldots, \partial_\mdim$ and satisfies $\varepsilon d=0$, then $d$ has a  unique extension to $U$ such that  $\varepsilon d=0$.
%\end{rem}

%It follows from the proposition that a unipotent closure $U_R$ of $R$ is functorial in $R$, and is unique up to isomorphisms.
%Furthermore, $U_R$ is a direct limit:
 %of all unipotent extensions of $R$:%
%$$U_R = \lim_{\rightarrow} B\ ,$$
%where $B$ ranges over the set of  all unipotent extensions of $R$, ordered by inclusion.
%  It follows that every unipotent extension of $U_R$ is trivial, {\it i.e.}, $U_R$ is the maximal unipotent extension of
%  $R$.

%\begin{cor} Let $B$ be a unipotent extension of $R$. Any
% $R$-differential homomorphism $\phi: U_R \rightarrow B$ is an isomorphism.
%\end{cor}

\begin{defn} If $U$ is the unipotent closure of a differentially simple $k$-algebra $R$, let  $\Gal(U/R)$
 be the group  of  differential automorphisms $\phi:U \rightarrow U$ over $R$.
\end{defn}

It follows from the definitions that $\Gal(U/R)$ is a pro-unipotent group. Now let $\epsilon:R\rightarrow k$ denote a $k$-valued point on
$\Spec R$. The  set of $k$-valued points $ \{ \phi\in \Hom_k(U,k) \,: \, \phi|_R = \varepsilon\}
$
on $\Spec U$ lying above $\varepsilon$, is a principal homogeneous space over $\Gal(U/R)$.
There is thus a complete analogy between the theory of unipotent differentially simple extensions and the theory of covering spaces.
%pi_1^\DR (R,\varepsilon)=

\subsection{Base points at infinity}\label{sect32b}
We need to repeat the theory of unipotent closures in the case where the base points are at infinity.
In order to do this, we need to generalise the notion of a $k-$valued point for certain differential  algebras.
\begin{defn} Let $k$ be a field. We define
\begin{equation} \label{kepsilonsmany}
\kes =k[[\epsilon_1,\ldots, \epsilon_\ell]] \Big[ {1\over
\epsilon_1},\ldots, {1\over \epsilon_\ell} \Big]    \
,\end{equation} to be the differential $k-$algebra of Laurent
series in $\epsilon_i$, equipped with $\ell$ commuting
differentials $\partial_{\epsilon_i}$, for $1\leq i\leq \ell$. Now
define the extension
\begin{equation} \label{Uepsilonsmany}
\Ues = \kes[L_{\epsilon_1},\ldots, L_{\epsilon_\ell}] \
,\end{equation} where $L_{\epsilon_i}$ is the formal logarithm of
$\epsilon_i$, {\it i.e.},  $\partial_{\epsilon_i} L_{\epsilon_i} =
\epsilon_i^{-1}$ for $1\leq i\leq \ell$.
\end{defn}
The ring of constants of $\kes$ is $k$,  and
the extension $\Ues$ is easily verified to be  a unipotent  closure of $\kes$, since
 $H^0(\Ues)=k$, and $H^1(\Ues)=0$.
%since every element in $\Ue$ has a primitive.
% One could
% impose restrictions on the convergence of  series in $\ke$, but this will not be required here.
\begin{defn} \label{defngeneralisedpoint} Let $R$ be a differentially simple $k$-algebra with $\ell$ commuting differentials $\partial_1,\ldots, \partial_\ell$.
%A $k$-valued point is a $k$-linear homomorphism $R\rightarrow k$.
We define a $\kes$-point on $R$ to be a $k-$linear homomorphism
$$p: R \To  \kes\ ,$$
which satisfies
$$ p \, \partial_i = \partial_{\epsilon_i} \, p \ ,\qquad \hbox{ for all }\quad  1\leq i\leq \ell\ .$$
 A $\kes$-point $p: R\rightarrow \kes$
defines an ordinary $k$-valued point if it factorises through $R\rightarrow k[[\epsilon_1,\ldots, \epsilon_\ell]]$:
$$\begin{array}{ccccc}
  R& \To  & k[[\epsilon_1,\ldots, \epsilon_\ell]] & \xrightarrow{\epsilon_1=\ldots=\epsilon_\ell=0} & k \\
   & \searrow  &  \downarrow &  &  \\
   & & \kes & &
\end{array}
$$
At the other extreme, we say that the point $p$ is \emph{at infinity} if $\epsilon_1^{-1},\ldots, \epsilon_\ell^{-1}\in \Image p$.
\end{defn}

\begin{example}
Consider the case where $k=\Q$, and $R= \Q[x,1/x, 1/(1-x)]$ with
differential $\partial/\partial x$. This corresponds to the
projective line minus three points $\Pro^1\backslash
\{0,1,\infty\}=\A^1\backslash \{0,1\}$. The set of $k$-valued
points on $R$ is the set $k\backslash \{0,1\}$. Each $\ke$-point
$$p: \Q\Big[x,{1\over x}, {1\over 1-x} \Big] \To \ke $$
satisfies $\partial_\epsilon p(x)= p(\partial_x x)=1$, and
  takes $x$ to $\epsilon + c$, where $c\in k$. The set of $\ke$-points is therefore the set $k$.
 In this case, there are just two points at infinity, given by the maps
$p_\lambda:R \rightarrow \ke$, where $\lambda =0,1$; defined as follows:
\begin{eqnarray}
x& \mapsto & \epsilon+\lambda \nonumber \\
{1\over x-\lambda} & \mapsto & {1\over \epsilon } \nonumber \ .
\end{eqnarray}
More generally, every corner of the Stasheff polytope
$\overline{X}_{S,\delta}\subset \Modf_{0,S}$ defines a base-point
at infinity on the ring $\Or(\Mod_{0,S})$. Given a triangulation
$\alpha \in \chi^\ell_{S,\delta}$ of the $n$-gon $(S,\delta)$, a
set of vertex coordinates $x^\alpha_1,\ldots, x_\ell^\alpha$
(\S\ref{sect23}) gives rise to  a map
$$\big(\Or(\Mod_{0,S}), \partial/\partial x^{\alpha}_i\big) \To \kes$$
which sends $x^\alpha_i$ to $\epsilon_i$ for $1\leq i\leq \ell$.
% (see, for example,
%equation $(\ref{pointatinfform05})$).
\end{example}
A point at infinity corresponds to a point which is the
intersection of a number of normal crossing divisors, and will
play the role of a tangential base point.

\begin{defn} Let $R$ denote any differentially simple $k$-algebra, with derivations $\partial_1,\ldots, \partial_\ell$.
 We define
 a \emph{logarithmic Laurent expansion} to be a homomorphism
%$  R \To    \Ues  $
 of differential $k$-algebras:
$$
%\begin{array}{ccc}
  \phi:R  \To     \Ues  %\\
%   & \searrow  &   \downarrow   \\
% & &   \Ues
%\end{array}
$$
There is a natural map
$\lambda: \Ues \To k$ which projects on to the constant coefficient in the logarithmic Laurent series. It factorises through
$$\Ues \To \kes \To k\ ,$$
where the first map sends $L_{\epsilon_i}$ to $0$ for $1\leq i\leq \ell$, and the second map picks out the constant term in the
Laurent expansion
$$\sum_{i_1,\ldots, i_\ell\geq -N} a_{i_1,\ldots, i_\ell} \epsilon^{i_1}_1\ldots \epsilon^{i_\ell}_\ell\quad \mapsto \quad a_{0,\ldots,0}\ .$$
The map $\lambda$ has a certain number of formal properties, which
we will not require explicitly.
% except for the following identity:
%\begin{equation}
%\partial_{\epsilon_i}\lambda(ab) = \lambda((\partial_i a) b) + \lambda(a
%(\partial_i b)) \ , \qquad \hbox{for } 1\leq i\leq \ell\ .
%\end{equation}
 Given a logarithmic Laurent
expansion $\phi: R \rightarrow \Ues$, we define \emph{the map of
constants} of $\phi$ to be the  ($k$-linear, additive) map
 $$\lambda\circ\phi: R \rightarrow k\ .$$
\end{defn}
\begin{lem} \label{lemmapofconstants} Let $R$ be a differentially simple
$k$-algebra, and let $p:R \rightarrow \kes$ denote a $\kes$-point.
Let $B$ denote a unipotent extension of $R$. Consider any
logarithmic Laurent expansion $\phi: B \rightarrow \Ues$ over the
point $p$, {\it i.e.}, such that the following diagram commutes:
 $$\begin{array}{ccc}
  B & \overset{\phi}{\To} & \Ues \\
  \uparrow &  & \uparrow \\
   R & \overset{p}{\To} & \kes  \ .
 \end{array}$$
Then $\phi$ is uniquely determined by its map of constants
$\lambda \circ \phi: B \rightarrow k$.
\end{lem}
\begin{proof}This follows immediately using the method of   proof of proposition $\ref{propPVproperties}$.
%$\ref{lemmorphismsunique}$.
%The map of constants determines the constant of integration at each stage.
\end{proof}

%A logarithmic Laurent expansion amounts to choosing branches of functions in $B$, and
A map of constants amounts to choosing a constant of integration for each successive primitive in a unipotent extension $B$ of $R$.
We can now copy the results of the previous sections for base points at infinity.

\begin{defn}
Let $R$ denote a  differentially simple $k$-algebra, and let
$p:R\rightarrow \kes$ be a $\kes$-point. Let $\Ut(R,p)$ denote the
category of pointed unipotent extensions of $(R,p)$, whose objects
are unipotent $R$-algebras $(B, \phi)$,  where
 $\phi:B\rightarrow \Ues$ is
a logarithmic Laurent expansion (or, equivalently, the corresponding map of constants). Morphisms are defined in a similar manner to the category $\Up$.
 \end{defn}

The proof of lemma $\ref{lemmorphismsunique}$  and proposition $\ref{propPVproperties}$ go through without any difficulty.
\begin{prop}
Morphisms are unique in $\Ut(R,p)$, and a unipotent closure $U$ of
$R$ is a final object in the category $\Ut(R,p)$.
\end{prop}
\noindent If $U$ is the unipotent closure of $(R,p)$, where $p$ is a $\kes$-point, then
\begin{equation}\label{p1derhamatinfinity}
%\pi_1^{\DR}(R, p)=
 \{ \phi:U\rightarrow \Ues\ , \phi|_R=p\}\ . \end{equation}
%It is easy to verify that $\pi_1^{\DR}(R, p)$
 is a principal homogeneous space over $\Gal(U/R)$.
 %In the case of base points at infinity, a new feature arises. First observe
% that the group $G_\ell$ of differential automorphisms
%of $\Ues$ over $\kes$, is%, which we denote $G_\ell=\Gal(\Ues/\kes)$, is a torus:
%$$ G_\ell \cong k^+ \times\ldots \times k^+\ .$$
%The group $G_\ell$ acts naturally on $(\ref{p1derhamatinfinity})$ on the left. This action is trivial when the base point is finite.
%This feature will be discussed in detail in the case of moduli spaces  in $\S\ref{sect67}$ and gives  a splitting
%of the set of iterated integrals into convergent and non-convergent parts.

\subsection{One-dimensional fibrations and their relative unipotent closures} \label{sect32a}
Let $R$ denote a  differentially simple $k$-algebra, with
commuting differentials $\partial_1,\ldots, \partial_\ell$.
Suppose that we are given $N$ elements $f_1,\ldots, f_N \in R$
which satisfy the condition
\begin{equation} \label{Rfibhyp}
{1 \over f_i-f_j} \in R \quad \hbox{ for all } 1\leq i<j\leq N\
.\end{equation} Consider the $R$-algebra
\begin{equation}\label{Rhatdefn}
\widehat{R}=R\Big[y, {1 \over y-f_1},\ldots, {1\over y-f_N} \Big]\
,
\end{equation}
equipped with the derivation $\partial_y$ which is the unique
$R$-linear derivation satisfying $\partial_y
 \,y=1$. Clearly $\partial_i, \partial_y$ commute for all $1\leq i\leq \ell$. Consider the free shuffle algebra
 generated by the symbols $\omega_1,\ldots, \omega_N$ over $\widehat{R}$:
\begin{equation}\label{Ufreeshuffle}
\Urel =   \widehat{R}\otimes_k k\langle \omega_1,\ldots,
\omega_N\rangle \ ,
\end{equation}
and let us extend the definition of $\partial_y$ to $\Urel$ by
setting
\begin{equation} \label{GaussManinderiv}
\partial_y = \partial_y \otimes 1 + \sum_{i=1}^N {1 \over y-f_i} \otimes \partial_{\omega_i}\ ,\end{equation}
where the left truncation operators $\partial_{\omega_i}$ were
defined in $\S\ref{sect31}$.
 This makes $\Urel$ into a differential $R$-algebra. A similar algebra was considered in [Ma],
and a specific case is studied in detail in $\S5$.
 The symbol $\omega_i$ represents the formal logarithm $\log (y-f_i)$, for $1\leq i\leq N$. By analogy with the bar
 construction,
 we will write $[\omega_{i_1}|\ldots |\omega_{i_m}]$ for the tensor $\omega_{i_1}\otimes \ldots \otimes\omega_{i_m}$.
The following proposition states that $\Urel$ is a relative
unipotent closure over the base $R$.

\begin{prop} \label{propRelnocohomology} $H^0(\Urel)= R$ and $H^1(\Urel)=0$.
\end{prop}

\begin{proof}
It is a simple exercise to show that the ring of constants
 of $\Urel$ is $R$. The argument is given in the proof of lemma $\ref{lemhoisk}$, and works in complete generality. The fact that $H^1(\Urel)$ vanishes
is equivalent to the existence of primitives with respect to
$\partial_y$ over the base $R$. The key observation is the
following identity, which
 is valid in $\Urel$, by assumption  $(\ref{Rfibhyp})$:
\begin{equation} \label{partialfractionidentity}
{1\over (y-f_i)(y-f_j)} = {1\over f_i-f_j} \Big({1\over y-f_i}
-{1\over y-f_j}\Big)\ ,\end{equation} Using this, we can decompose
elements of $\widehat{R}$ into partial fractions. It suffices,
therefore, to find primitives of expressions of the form
$${1\over (y-f_i)^n} [\omega_{i_1}|\ldots|\omega_{i_m}]\ ,$$
where $n\in \Z$, and $1\leq i_1,\ldots, i_m\leq N$. If $n=-1$, a
primitive is given by
$$[\omega_i|\omega_{i_1}|\ldots |\omega_{i_m}]$$
by definition. For other values of $n$, we can reduce to this case
by integrating by parts and using induction.
%For all negative $n$, we can integrate by parts and use induction on $n$ to reduce to this case.
%When $n$ is positive, we can rewrite any polynomial in $y$ as a polynomial in $y-f_{i_1}$, and  use induction
%on the number of symbols $m$ by integrating by parts.
It follows that every element in $\Urel$ has a primitive with
respect to $\partial_y$.
\end{proof}
\noindent Note that there is no integrability condition to be
verified because the fibres of the map $\Spec
\widehat{R}\rightarrow \Spec R$ are one-dimensional. We now show
how to differentiate the symbols
$[\omega_{i_1}|\ldots|\omega_{i_m}]$ with respect to the operators
$\partial_1,\ldots, \partial_\ell$  of the base ring $R$
(differentiation under an iterated integral).
 To do this, consider an $R$-linear map
 $$p: \widehat{R} \To R\{\epsilon\}=R[[\epsilon]]\big[{1\over \epsilon}\big] $$
 which satisfies $p \,\partial_y = \partial_\epsilon \,p$. Then
 there is a unique logarithmic Laurent expansion:
$$
\begin{array}{ccc}
  \Urel & \overset{\phi}{\To} & R\{\epsilon\}[L_\epsilon]  \\
  \uparrow & &\uparrow  \\
  \widehat{R}& \overset{p}{\To} & R\{\epsilon\}
\end{array}
$$
such that the map of constants is zero on the generators of
$\Urel$, {\it i.e.},
\begin{eqnarray} \label{Trivialmapofconstants}
\lambda \circ \phi :\Urel &\To& R  \\
\big[\omega_{i_1}|\ldots|\omega_{i_m}\big] &\mapsto & 0\ .
\nonumber
\end{eqnarray}
This follows from the inductive method of proof of proposition
$\ref{propPVproperties}$: if
$w=[\omega_{i_1}|\ldots|\omega_{i_m}]$, and $\phi(\partial_y w)=a$
has already been defined, then $\phi(w)$ is defined to be a
primitive of $a$ with respect to $\partial_\epsilon$. The constant
of integration is normalised by the condition
$\lambda(\phi(w))=0$, since the map
$(\ref{Trivialmapofconstants})$ is $R$-linear.

\begin{prop}\label{GaussManinversion}
 The action of
 the differential operators $\partial_i$ on $R$, for $1\leq i\leq \ell$, extend uniquely
to $\Urel$ in such a way that the $\partial_i$ commute with  each
other, and such that:
$$[\partial_i,\phi]=[\partial_i, \partial_y] =0 \quad \hbox{ for all } 1\leq i\leq \ell\ .$$
For each element $w = [\omega_{i_1}|\ldots |\omega_{i_m}]\in W^m
\Urel$, $\partial_i w \in W^{m-1} \Urel$. It follows that $\Urel$
is a unipotent differential algebra with respect to all the
operators $\partial_1,\ldots, \partial_\ell, \partial_y$.
\end{prop}
\begin{proof}
The map $\phi$ is injective, since $\Urel$ is differentially
simple, and therefore the action of the operators $\partial_i$ on
$\Urel$ are induced from $R\{\epsilon\}[L_\epsilon]$ by
restriction. More precisely, suppose by induction that the action
of the operators $\partial_i$ have already been defined on $W^p
\Urel$, for some $p \geq 0$. Let $w\in W^{p+1}\Urel$ such that
$\partial_y w \in W^p\Urel$. If we view $\Urel$ as a subalgebra of
$R\{\epsilon\}[L_\epsilon]$, then we can write
$$\partial_y \partial_i w = \partial_i \partial_y w \in \Urel\ .$$
The element $\partial_i w$, which is \emph{a priori}  in
$R\{\epsilon\}[L_\epsilon]$, in fact lies in $\Urel.$ This is
because it is a primitive of $\partial_i \partial_y w\in \Urel$,
and  we know that $H^1(\Urel)=0$, and
$H^0(R\{\epsilon\}[L_\epsilon])=R$.  More explicitly, if
$w=[\omega_{i_1}|\ldots |\omega_{i_m}]$, then we define
$$\partial_i w = (1- \lambda \circ \phi) A\ ,$$
where $A$ is any solution in $\Urel$ to  $\partial_y A= \partial_i
\partial_y w$. The fact that the operators $\partial_i$ decrease
the weight of each such element $w$ is easily proved by induction
and is left to the reader.
\end{proof}
\noindent For each $1\leq i\leq N$, there is a unique $R$-linear
map
\begin{eqnarray} \label{pmaponfibre}
p:\widehat{R} &\rightarrow  & R\{\epsilon\}  \\
y-f_i & \mapsto & \epsilon \ ,\nonumber
\end{eqnarray}
such that $\partial_\epsilon p= p\, \partial_y$. It satisfies:
$$p\Big( {1\over y-f_j}\Big) ={1\over f_i-f_j}
\sum_{k\geq 0} { \epsilon^k \over (f_j-f_i)^{k}}\quad \hbox{ for
all } j\neq i\ .$$

\begin{cor}\label{corunipclosureisproductforfinitebasepoint}
Suppose that $U_R$ is the unipotent closure of $R$. Then the
algebra $U_R\otimes_R \Urel$ is the unipotent closure of
$\widehat{R}$.
\end{cor}
\begin{proof}
By choosing any map $p$ given by equation $(\ref{pmaponfibre})$
above, we obtain a differential $\widehat{R}[\partial_1,\ldots,
\partial_\ell,
\partial_y]$-structure on $U_R\otimes_R \Urel$.
 It is
clear that a tensor product of unipotent algebras is unipotent,
and that $H^0(U_R\otimes_R \Urel)=k$. Since the operator
$\partial_y$ is zero on $U_R$, it follows that $H^1(U_R\otimes_R
\Urel)=0$. Concretely, in order to find 1-primitives in this
algebra, first take a primitive with respect to $\partial_y$ and
then adjust the constant of integration in $U_R$ using the fact
that $H^1(U_R)=0$.
\end{proof}

By iterating the previous corollary, we deduce that  any
differentially simple algebra $R$ which is of fiber-type ({\it
i.e.}, an iterated sequence of fibrations) has an explicit
unipotent closure which is a tensor product of shuffle algebras.

\begin{thm}\label{thmRBarsplits}
Let $R$ denote a differentially simple $k$-algebra,  which can be
expressed as a finite
 series of extensions of the type $(\ref{Rhatdefn})$ satisfying $(\ref{Rfibhyp})$:
\begin{equation}\label{Rfilt1}
k=R_0 \subset R_1 \subset \ldots \subset R_{n}=R\ ,
\end{equation}
where
\begin{equation} \label{Rfilt2}
R_{t}= R_{t-1} \Big[ y_t, \Big({1\over y_t-f_{t,i}}\Big)_{1\leq
i\leq N_t}\Big]\ ,\end{equation} and $f_{t, i}- f_{t,j}$ is
invertible in $R_{t-1}$ for all $1\leq i<j\leq N_t$, and all
$t=1,\ldots, n$. Then the unipotent closure  $U_R$ of $(R,p)$
exists, and is isomorphic (as an algebra) to the tensor product of
free shuffle algebras on $N_t$ generators, for $1\leq t\leq n$:
\begin{equation}\label{URtensprod}
U_R \cong R\otimes_k \bigotimes_{t=1}^n k\langle
\omega_{t,1},\ldots, \omega_{t,N_t}\rangle\ .\end{equation}
Its
differential structure  is uniquely determined by such a tensor
decomposition. %and
%It has a canonical logarithmic Laurent expansion $\phi: U_R
%\rightarrow \Ues$ over $p$, whose map of constants maps each
%tensor product of symbols $[\omega_{t_1,{i_1}}|\ldots
%|\omega_{t_1,{i_m}} ]\otimes \ldots
%\otimes[\omega_{t_k,{i_1}}|\ldots |\omega_{t_k,{i_m}} ] $ to zero.
\end{thm}

\begin{proof}
This follows immediately from the previous corollary by induction.
The differential structure is determined by the construction in
proposition $\ref{GaussManinversion}$.
\end{proof}

 We therefore have an explicit description of the algebraic
structure of the unipotent closure of $R$ for any $R$ which is of
fiber type. Note that there may be several natural isomorphisms of
the form $(\ref{URtensprod})$, even after fixing base-points.

%{\it i.e.}, which can be expressed as a finite
%sequence of fibrations $(\ref{Rfilt1})$ satisfying
%$(\ref{Rfilt2})$.
 %but
%this depends on the choice of a base point on $\Spec R$.

%\begin{rem} In the case where $p:R\rightarrow \kes$ is a
%base-point at infinity which is compatible with the filtration.
%\end{rem}

\subsection{Iterated integrals}\label{sect33}

Let $\Or_M$ denote the ring of regular functions on an affine hyperplane arrangement as considered in $\S \ref{sect32}$.
 $\Or_M$ is a differential algebra with $\ddim$
commuting differentials $\partial/\partial x_1,\ldots, \partial/\partial x_\ell$. % It is differentially simple.
Suppose that  $I \subset \Or_M$ is any non-zero differential
ideal. It must contain a polynomial $P\in k[x_1,\ldots,x_\ddim]$, since we can
multiply by suitable powers of the
hyperplane equations $\alpha_i$ to clear denominators. % ,  $I$ must
 It is clear that there exists a
polynomial $D_P$ in the $\partial/\partial x_i$ such that
$D_P\, P=1,$ and therefore $I=\Or_M$. It follows that $\Or_M$ is differentially simple. % The following theorem is
\begin{thm} \label{thmbartriv}
The de Rham cohomology of $B(\Or_M)$ satisfies:% is trivial in degrees $0$ and $1$:
 $$H^0_{\DR}(B(\Or_M))=k \ , \quad \hbox{ and } \quad H^1_{\DR}(B(\Or_M))=0\ ,$$
and $B(\Or_M)$ is  the unipotent closure of $\Or_M$. It follows that
 every differential $\Or_M$-subalgebra of $B(\Or_M)$ is differentially
 simple, and $B(\Or_M)$ is a polynomial algebra.
\end{thm}
The proof of this theorem is postponed until $\S \ref{sect36}$.

\begin{rem}  The theorem in fact holds in much greater generality. Let $F$ denote any differential algebra
such that $H^1(F)\cong \bigoplus_{i=1}^N k \,\omega_i$, where $\omega_i \in \Omega^1(F)$ satisfy
$$\big(\bigoplus_{i,j} k\, \omega_i \wedge \omega_j \big)\cap  d\Omega^1(F) =0 \ .$$
 If $k$ is the field of constants of $F$, and if $B(F)$ is defined as in $\S\ref{sect32}$,
  then it is clear from the proof $(\S\ref{sect36})$, that
 $H^0(B(F))=k$ and $H^1(B(F))=0$. Furthermore, when $F$ is differentially simple,  every
  differential $F$-subalgebra of $B(F)$ is differentially
 simple, and $B(F)$ is a polynomial algebra.
\end{rem}

% In the case $d=1$, when $M$ is the projective plane minus a finite
%number of points, $B(F)$ is isomorphic to the free shuffle
%algebra. It was proved by Radford that

%If we view $B(F)$ as a differential extension of the
%differentially simple ring $F$, then the fact that $H^0(B(F))=k$
%states that the ring of constants of $B(F)$ is equal to the ring
%of constants $H^0(F)=k$ of $F$. We can therefore view  $B(F)$ as a
%kind of infinite Picard-Vessiot extension of $F$. It is in fact
%the unique smallest extension of $F$ which is closed under the
%operation of taking primitives of one-forms.
%% {\it i.e.}, given $f_1,\ldots, f_d\in B(F)$
%%such that $\partial_i f_j=\partial_j f_i$ for all $1\leq i,j$
%We shall therefore call $B(F)$  the \emph{unipotent
%closure} of $F$. Such a differential ring has many special
%properties, as we shall see later.

%\begin{cor}
% Every unipotent differential equation on $M$ has a solution in
% $B(F)$. It is the smallest extension of $F$ with this property.
%\end{cor}

We now recall the definition of Chen's iterated integrals, which will give an
isomorphism of the abstract algebra $B(F)$ with an algebra  of
multi-valued functions.  %$M$ a smooth manifold.
 Let $\widehat{M}$ be a universal covering
for $M$, and let $p: \widehat{M}
 \rightarrow M$ denote the projection map. Let $b \in
 M$ denote a  base point for $M$. Given any smooth path $\gamma:[0,1]\rightarrow M$  beginning at
$b$, and holomorphic 1-forms $\eta_1,\ldots, \eta_m \in \Omega^1(M)$,
 the iterated integral of the word $\eta_m\ldots \eta_1$ (note the reversed order of symbols) along $\gamma$ is defined by
$$\int_\gamma \eta_1\ldots \eta_m= \int_{0<t_1<\ldots<t_m<1} \gamma^*\eta_1(t_1) \wedge \ldots \wedge \gamma^*\eta_m(t_m) .$$
One can show using the calculus of variations [Ch1] that the
iterated integral of  a linear combination of forms $f=\sum_I c_I\,
\omega_{i_1}\ldots \omega_{i_m}$ only depends on the homotopy
class of $\gamma$ if and only if the integrability condition
$(\ref{intcond})$ is satisfied. In this case, an iterated integral varies holomorphically
  as a function of the endpoint $z=\gamma(1)$ of $\gamma$, and therefore  defines
a holomorphic
function on the universal covering $\widehat{M}$.
We can realise $\Omega^*(\Or_M)$ as an algebra of differential forms on $\widehat{M}$ by taking the pull-back along
the covering map $p:\widehat{M}\rightarrow M$.
When we refer to a multi-valued function (or form) on $M$ it will be a linear combination of such iterated integrals
with coefficients in $\Or_M$ (resp. $\Omega^*(\Or_M)$) (compare the multi-valued de Rham complex defined in [H-M]).

\begin{lem}([Ch1, Ha1, Ha2])\,\label{lemitintsprop} Let $\eta_1,\ldots, \eta_l \in \Omega^1(M)$.
\begin{enumerate}
  \item Let $1\leq m\leq l$, and let $\Sym(m,l-m)$ denote the set of $(m,l-m)$-shuffles defined in $\S\ref{sect27}$. Then
the shuffle product formula holds:
$$
 \int_\gamma \eta_1\ldots \eta_m \int_\gamma \eta_{m+1}\ldots \eta_{l} = \sum_{\sigma \in \Sym(m,l-m)}
\int_\gamma \eta_{\sigma(1)}\ldots \eta_{\sigma(l)} \ .$$

  \item Let $\gamma_z:[0,1]\rightarrow M$ denote a smooth family of paths such that $\gamma_z(0)=b$, and $\gamma_z(1)=z\in M$.
If $\sum_{I} c_I\, \omega_{i_1}\ldots \omega_{i_m}$ satisfies the integrability condition $(\ref{intcond})$, we have:
$$ {d\over dz}  \int_{\gamma_z} \sum_I  c_I\omega_{i_m}\ldots \omega_{i_1} =\sum_I c_I \omega_{i_1} \int_{\gamma_z}  \omega_{i_m}\ldots
\omega_{i_{2}} \ .$$
\end{enumerate}
\end{lem}

\begin{defn} Let $L_b(M)$ denote the  $\Or_M$-module generated
by all such homotopy-invariant iterated integrals on $\widehat{M}$.
We write $\Omega^i(L_b(M)) = L_b(M) \otimes_{\Or_M} \Omega^i(\Or_M)$.
\end{defn}
%This is an
%algebra by the shuffle product formula for iterated integrals
%(which we shall reprove in $\S7$), and has the same differential
%structure as the algebra $B(\Or_M)$.
By the previous lemma, $\Omega^*(L_b(M))$ is a differential algebra, and there is a map:
\begin{eqnarray}\label{realis}
\rho_b:\Omega^* B(\Or_M) &\overset{\sim}{\To}&\Omega^* L_b(\widehat{M}) \\
\sum_I \phi_I [\omega_{i_1}|\ldots |\omega_{i_m}]&\mapsto& \sum_I \phi_I
\int_\gamma \omega_{i_m}\ldots \omega_{i_1} \nonumber \ ,
\end{eqnarray}
 which is  a surjective map of differential
algebras by $(\ref{bardifdef})$.  As above, $\gamma$ denotes a smooth path beginning at the point $b\in M$. The  kernel of $\rho_b$ is a differential ideal, and therefore must reduce to zero since
$B(\Or_M)$ is differentially simple. Therefore $(\ref{realis})$ is an isomorphism.

\begin{cor} If $\{e_i\}$ is a basis for $B(\Or_M)$ over $\Or_M$, then
the  functions $\rho_b(e_i)$ are linearly independent over
$\Or_M$. All algebraic relations between the functions $\rho_b(e_i)$
are determined by the shuffle product.
%The set of functions given by the images of a basis of $B(\Or_M)$
%under the map $(\ref{realis})$
%%$$\sum_I f_I \int\omega_{i_1}\ldots \omega_{i_n}$$
%are linearly  independent over $k(M)$.
%Every unipotent
%differential equation in $M$ has a solution in $L(\widehat{M})$.
\end{cor}
\noindent One can  determine a basis for $B(\Or_M)$ in the
fiber-type case (see \S\ref{sect62}).

% It is the unique smallest extension of
%$\Q(\Mod_{0,S})$ with this property.
%\end{cor}

%-------------------------SECTION ON CHEN's Pi1-de Rham theorem ------------------------------------
%---------------------------------------------------------------------------------------------------
%---------------------------------------------------------------------------------------------------

%----------------------------------------------------------------------------------------------------
%---------------------------------------------------------------------------------------------------
%---------------------------------------------------------------------------------------------------
%---------------------------------------------------------------------------------------------------

\subsection{Proof of theorem $\ref{thmbartriv}$.}
\label{sect36} %We prove theorem $\ref{thmbartriv}$ in a series of elementary lemmas.
We first show that the ring of constants of $B(\Or_M)$ is $k$.
 % trivializes the
%$\mathrm{zero^{th}}$ cohomology group of $\Or_M$.
 For any  $\psi\in B(\Or_M)$, $n\geq0$,  we write $\psi_n=\gr_n^w \, \psi$ for its graded part of weight $n$.

\begin{lem} \label{lemhoisk} $H^0_{\DR}(B(\Or_M))=k$.
\end{lem}

\begin{proof}
 Let $\psi\in  B(\Or_M)$ of weight $m\geq1$  such that $d\psi=0$.
 We write
 $$\psi_r = \sum_{I=(i_1,\ldots,i_r)} f_I \, [\omega_{i_1}|\ldots |\omega_{i_r}] \quad \hbox{
 for } 0\leq r\leq m\ ,$$
where each $f_I \in \Or_M$.
Then the graded weight $m$ part of $d\psi$ is zero:
$$(d\psi)_m = \sum_{|I|=m} df_I\, [\omega_{i_1}|\ldots |\omega_{i_m}]=0\ .$$
Therefore $df_I=0$ and so $f_I\in H^0_{\DR}(\Or_M)= k$ for all ordered sets $I$ such
that $|I|=m$. The weight $m-1$ part of $d\psi$ is also zero:
$$(d\psi)_{m-1} =\sum_{|I|=m}f_I\,\omega_{i_1} \,[\omega_{i_2}|\ldots |\omega_{i_r}]
+ \sum_{J=(i_2,\ldots,i_m)} df_J [\omega_{i_2}|\ldots
|\omega_{i_m}] =0\ ,$$ which implies that
$$f_{i_1,i_2,\ldots, i_m} \omega_{i_1} +
df_{i_2,\ldots, i_m}=0 \qquad \hbox{for all } i_2,\ldots, i_m\ .$$  Therefore the forms $f_{i_1,\ldots, i_m} \omega_{i_1}$
are exact for all $i_1,\ldots, i_m$.
 But since we have shown that
$f_{i_1,\ldots, i_m}\in k$ is constant,  this  can only occur if $f_{i_1,\ldots, i_m}=0$. This
implies that the weight of $\psi$ is at most $m-1$, which
contradicts the initial assumption. Therefore, any $\psi$ such
that $d\psi=0$ is of weight $0$, and lies in $\Or_M$. Hence $\psi\in
H^0_{\DR}(\Or_M)=k$.
\end{proof}

\noindent The  following lemma  states that we can replace a
closed $1-$form in $B(\Or_M)$ with an element in its cohomology class
of strictly lower weight.
% This gives an inductive algorithmic procedure for
%computing iterated integrals.
\begin{lem}
Let $\psi\in \Omega^1(B(\Or_M))$ be an element of weight $m$ such
that $d\psi=0$. Then there exists $\theta\in B(\Or_M)$
%(of weight $m+1$)
 such that $\kappa = \psi- d\theta$ is of weight at
most $m-1$.
\end{lem}

\begin{proof} Let
$$\psi =  \sum_{r=0}^m \sum_{I=(i_1,\ldots,i_r)} \phi_I [\omega_{i_1}|\ldots |\omega_{i_r}]\ ,$$
where $\phi_I \in \Omega^1 \Or_M$ for all indexing sets $I$. Since  $d\psi=0$,
 we deduce that
$$0 = \sum_{|I|=m} d\phi_I \, [\omega_{i_1}|\ldots |\omega_{i_m}]
-\sum_{|I|=m} \phi_I \wedge \omega_{i_1} [\omega_{i_2}|\ldots
|\omega_{i_m}] +\sum_{r=0}^{m-1}d(\psi_{r})\ .$$ This implies
firstly that $d\phi_I=0$ for all sets $I$ with $|I|=m$, and secondly
that
\begin{equation} \label{pf1}
\sum_{r=0}^{m-1}d(\psi_{r})  -\sum_{I=(i_1,\ldots, i_m)} \phi_I
\wedge \omega_{i_1} [\omega_{i_2}|\ldots |\omega_{i_m}]= 0 \ .
\end{equation}
Taking the graded part of this equation of weight $m-1$, we deduce
that
$$ -\sum_{I=(i_1,\ldots,i_m)} \phi_I \wedge \omega_{i_1}
[\omega_{i_2}|\ldots|\omega_{i_m}] + \sum_{i_2,\ldots,i_m}
d\phi_{i_2, \ldots, i_{m}} [\omega_{i_2}|\ldots |\omega_{i_m}]=0\
,$$ and so
\begin{equation} \label{pf2}
 \sum_{i_1} \phi_{i_1,\ldots, i_m} \wedge \omega_{i_1}=d\phi_{i_2,\ldots,i_m} \ ,
\end{equation}
  for all
$I=(i_1,\ldots, i_m)$. We have shown that $\phi_I$ is closed for $|I|=m$,
so we can write
\begin{equation} \label{pf3}
\phi_I = \sum_{j} \alpha_{I,j}\, \omega_{j} + dg_I \ .
\end{equation} where $\alpha_{I,j}\in k$, and $g_I \in
 \Or_M$.
 Substituting into $(\ref{pf2})$ above, we have
$$\sum_{i_1,j} \alpha_{i_1,\ldots, i_m,j}\,  \omega_{j}  \wedge \omega_{i_1} + \sum_{i_1} dg_{i_1,\ldots, i_m}     \wedge
\omega_{i_1} = d\phi_{i_2,\ldots, i_m}\ ,$$
  for all
$i_2,\ldots, i_m$. The corollary to theorem $\ref{thmarnold}$ implies that any
linear combination of  exterior products of forms $\omega_i$
which is exact, is necessarily zero. Using the fact that $dg_{i_1,\ldots,i_m}
\wedge \omega_{i_1} = d(g_{i_1,\ldots,i_m} \wedge \omega_{i_1})$ is exact, we
have
\begin{equation}\label{pf6}
\sum_{i_1,j} \alpha_{i_1,\ldots,i_m,j}\, \omega_{j} \wedge \omega_{i_1} =0 \ ,
 \qquad \hbox{for all } i_2,\ldots, i_m\ .\end{equation}
 Let
$$\theta_1=
\sum_{I=(i_1,\ldots,i_m)} \sum_{j} \alpha_{I,j} \, [\omega_{j} | \omega_{i_1}|\ldots |\omega_{i_m}]\
.$$ Since the integrability condition $(\ref{intcond})$ is
homogeneous with respect to the weight, the integrability of
$\psi$ implies the integrability of $\psi_m=\sum_{|I|=m} \phi_I
[\omega_{i_1}|\ldots |\omega_{i_m}]$. This is equivalent to a
number of linear equations of the form $\sum_{|I|=m} \lambda_I \phi_I
= 0$, where $\lambda_I \in k$. Using the decomposition
$(\ref{pf3})$, and the fact that   $\Image
\big(\bigoplus_{ i<j} k\,
\omega_{i}\wedge \omega_{j} \rightarrow
\Omega^1(\Or_M) \big) $ and $d\Or_M$  are complementary spaces (this follows from theorem $(\ref{thmarnold})$), we deduce that
%$\sum_{|I|=m} \lambda_I dg_I=0$,
%and thus
\begin{equation}\label{pfint}
\sum_{I=(i_1,\ldots,i_m)} \sum_{j} \alpha_{I,j}\, \omega_{j} [\omega_{i_1}|\ldots |\omega_{i_m}]\ .
\end{equation}
is integrable, as is $\sum_{|I|=m} dg_I [\omega_{i_1}|\ldots|
\omega_{i_m}]$. By adding constants, we can assume that the primitives $g_I$ of $dg_I$ satisfy the same linear equations  $\sum_{|I|=m}
\lambda_I g_I=0$. This ensures that
$$\theta_2=\sum_{|I|=m} g_I [\omega_{i_1}|\ldots |\omega_{i_m}]  $$
satisfies the integrability criterion also. The integrability of $\theta_1$ follows from
$(\ref{pfint})$ and $(\ref{pf6})$.
 We set
$\theta=\theta_1+\theta_2\in B(\Or_M)$.
 By construction, we have
\begin{eqnarray}
 d\theta-\psi
&=& d\,\big( \sum_{I=(i_1,\ldots,i_m)} \sum_{j} \alpha_{I,j} \, [\omega_{j} | \omega_{i_1}|\ldots |\omega_{i_m}] +
 g_I \,[\omega_{i_1}|\ldots|\omega_{i_m}]\big)-\psi
 \nonumber \\
&= &\sum_{|I|=m} \big(\sum_{j} \alpha_{I,j}\, \omega_j +dg_I \big)\, [ \omega_{i_1}|\ldots |\omega_{i_m}]
+  g_I\wedge \omega_{i_1}
[\omega_{i_2}|\ldots |\omega_{i_m}] -\psi\ ,
\nonumber \\
&= &
 \sum_{|I|=m} g_I\wedge \omega_{i_1}
[\omega_{i_2}|\ldots |\omega_{i_m}] -(\psi_0+\ldots +\psi_{m-1}) \nonumber \ ,\end{eqnarray} which is of
weight at most $m-1$, since all terms of weight $m$ cancel by $(\ref{pf3})$.
\end{proof}

%\noindent It follows that any element $\psi\in \Omega^1 B(\Or_M)$  such that $d\psi=0$ can be replaced by a class of
%strictly lower weight without changing its cohomology class.

%The previous lemma, however, gives an  algorithm for finding the primitive of
%any element $\psi \in \Omega^1(\Or_M)$ which satisfies $d\psi=0$.
Given  a closed form $\psi\in \Omega^1(\Or_M)$ of weight $m$, we defined  an explicit
$\theta\in B(\Or_M)$ such that $\psi= d\theta + \psi_1$,
and  $\psi_1$ is of weight $\leq m-1$. In fact,  $\theta$ is of weight at most $m+1$.
Applying the lemma
repeatedly, we obtain a series of forms $\psi_1,\ldots, \psi_m \in
\Omega^1(\Or_M)$ and $\theta_1,\ldots, \theta_m \in B(\Or_M)$,
where $\psi_i$ is of weight at most $m-i$, such that $$\psi_i=
d\theta_i + \psi_{i+1}\ .$$ At the final stage, $\psi_m=
d\theta_m$. Thus $\psi= d(\theta+\theta_1+\ldots +\theta_m)$, and
 $\theta+\theta_1+\ldots +\theta_m$ is a
primitive of $\psi$ of weight at most $m+1$.

As remarked earlier, the argument in the proof of the lemma can be both generalised and simplified using spectral sequence arguments (see the appendix).

% In order to implement this, it suffices to be
%able to compute the decomposition $(\ref{pf3})$, {\it i.e.},
%computations in the cohomology of $B(\Or_M)$ reduce to computations in
%the cohomology of $\Or_M$.

\begin{cor}
$H^1_{\DR}(B(\Or_M))=0$.
\end{cor}

This completes the proof of theorem $\ref{thmbartriv}$. The fact that every $\Or_M$-subalgebra of
 $B(\Or_M)$ is differentially
simple, and the fact that $B(\Or_M)$ is a polynomial algebra, follows from the results of $\S\ref{sectPV}$.

\subsection{Fibrations of hyperplane arrangements} \label{sect38} We  recall necessary and sufficient conditions for an affine hyperplane arrangement to
decompose as a fibration over an arrangement of smaller dimension [O-T]. We deduce
from the results of $\S\ref{sect32a}$ that the reduced bar construction has trivial cohomology
for fiber-type arrangements.

%In this situation the reduced bar construction
%splits as a tensor product.
%This situation will enable us to use  inductive arguments.

 Let $\Hy=\{H_1,\ldots, H_N\}$
denote any affine hyperplane arrangement in $\A^\ddim$.
Choose any affine subspace $W\cong \A^e $ contained in $\A^{\ddim}$
and let $V_0\subset \A^\ddim$ denote a complementary subspace such that
$$\A^\ddim \cong V_0 \oplus W\ .$$
For  each $z\in \A^e\cong W$, let $V_z=V_0+z$ denote  the affine space
parallel to $V$ passing through the point $z\in W$. The spaces $V_z$
define a vertical direction normal to the base $W$. We define
the set of vertical hyperplanes to be
$$\Hy^v = \{ H \in \Hy : H \hbox{ contains } V_z \hbox{ for some }
z\in W\}\ ,$$ and let $\Hy^h$ denote the set of all remaining
hyperplanes. There is a decomposition
$$\Hy=\Hy^v\sqcup \Hy^h\ ,$$
and it is clear that every horizontal hyperplane $H\in \Hy^h$
intersects each $V_z$ properly. Consider the complements
$$M = \A^\ddim\backslash \bigcup_{H\in \Hy} H\ , \quad\hbox{ and  }\quad  M'= W \backslash
\displaystyle{\bigcup_{H\in \Hy^v}} H \cap W \ .$$ The linear
projection  $\A^\ddim \rightarrow W$ with kernel $V_0$ induces a  surjective map $p:M\rightarrow M'$.
%\begin{equation} \label{pfib}  p: M \To M'\ .\end{equation}
\begin{lem} The map $p$ is a fibration if and only if the
following condition holds: for all $H,H'\in \Hy$ such that $H\cap
H'\neq \emptyset$, there exists $H''\in \Hy^v$ such that
$$H'' \supseteq H \cap H'\ .$$
The fibre over $z\in M'$  is  the complement $
V_z\backslash \displaystyle{\cup_{H\in \Hy^h}} (H\cap V_z)$.
\end{lem}
The proof is left as an exercise.

%\begin{defn} We say that the fibration $p:M'\rightarrow M$ is \emph{at infinity} in
% the special case where  $W$ is an intersection of hyperplanes
%$$W= \bigcap_{i=1}^{\ddim-e}H_i\ ,$$
%where $H_1,\ldots, H_{\ddim-e} \in \Hy$  intersect%
%properly  ({\it i.e.}, $W$ is of dimension exactly $e$).
%

%\end{defn}

\begin{figure}[h!]
  \begin{center}
%    \leavevmode
    \epsfxsize=8.0cm \epsfbox{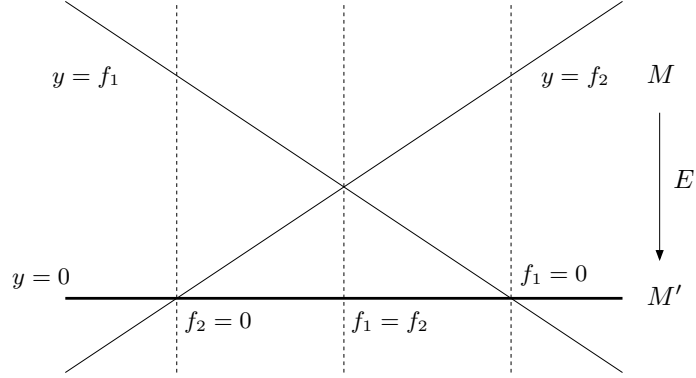}
  \label{Fibration}
\put(-122,17){\small{$f_1=f_2$}}
\put(-185,17){\small{$f_2=0$}}
\put(-58,35 ){\small{$f_1=0$}}
\put(-235,110){\small{$y=f_1$}}
\put(-50,110){\small{$y=f_2$}}
\put(-250,34){\small{$y=0$}}
\put(-10,110){$M$}
\put(0, 70){$E$}
\put(-10,26){$M'$}

  \caption{An arrangement $M$ in $\A^\ddim$ which fibers  over $M'\subset \A^{\ddim-1}$ (the thick  line at the bottom). A vertical
hyperplane (dashed) passes through every point where  horizontal
hyperplanes intersect. }
  \end{center}
\end{figure}

\begin{defn}\label{defnsupersolv} An affine hyperplane arrangement is said to be of \emph{fiber-type} if it can be expressed as an iterated sequence of
linear fibrations whose fibers are of dimension $1$.  Thus there is a sequence of fibrations
\begin{equation} \label{supersolv}
M\overset{E_1}{\To} M_1\ , \quad \ldots \ ,\quad M_{\ddim-2} \overset{E_{\ddim-1}}{\To} M_{\ddim-1}\ ,\end{equation}
where $E_1,\ldots, E_{\ddim-1}$, and $E_\ddim=M_{\ddim-1},$ are of dimension $1$.
%We say that an arrangement is of
%\emph{fiber-type at infinity} if every fibration is a fibration at infinity.
\end{defn}

We consider in greater detail the case where the dimension of the fibres is $1$. Then each fibre is isomorphic to  $\A^1$ minus a finite number of points.
Let us choose coordinates compatible  with the direct sum decomposition $\A^\ddim  = W
\oplus V_0$.  In other words, let $x_1,\ldots, x_{\ddim-1}$ denote coordinates on $\A^{\ddim-1}=W$, and let $y$ denote the vertical coordinate
on $\A^1=V_0$. %We set
%$$d_{M'}= \sum_{i=1}^{\ddim-1} {\partial \over \partial x_i} dx_i\ ,\quad \hbox{
%and } \quad d_E= {\partial \over \partial y} dy\ .
%$$
%Let $d= d_{M'}+d_E$. The operators $d_E$, $d_{M'}$ clearly
%commute.
 Let $\Or_M$ and
$\Or_{M'}$  denote the rings of regular functions on
the affine schemes $M$ and $M'$ respectively.
Let us write the equations of all horizontal hyperplanes in the form
$y-f_i=0$, where $f_i \in \Or_{M'}$, and $1\leq i\leq N_h$, for some integer $N_h$.
Thus $H_i = \ker (y-f_i)$ for $1\leq i\leq N_h$. By the previous lemma, the fact that $M$ is a fibration is equivalent
to the equations (see figure 12):
\begin{equation} \label{finbase}
{1\over f_i -f_j } \in \Or_{M'} \qquad \hbox{for all}\quad 1\leq i,j\leq N_h\ .
\end{equation}
We have
\begin{equation} \label{Mfibrings}
 \Or_M = \Or_{M'} \Big[ y, {1 \over y-f_1},\ldots, {1\over y-f_{N_h}}\Big]\ .\end{equation}
We have already shown that the rings $\Or_M$ and $\Or_{M'}$ are differentially simple over $k$. We are therefore
in the situation considered in $\S\ref{sect32a}$ (compare $(\ref{Rfibhyp})$ and $(\ref{Rhatdefn})$).

% Let $E$ be defined as in the lemma, and
%We will assume
%without loss of generality
% that we have chosen coordinates which
%are compatible with the direct sum decomposition $\A^d  = \A^e
%\oplus V$. In other words, let $\A^e$ have coordinates
%$x_1,\ldots, x_e$ and let $x_{e+1},\ldots, x_d$ be coordinates on
%$V_0$. We set
%$$d_{M'}= \sum_{i=1}^e {\partial \over \partial x_i} dx_i\quad \hbox{
%and } \quad d_E=\sum_{i=e}^d {\partial \over \partial x_i} dx_i\ .
%$$
%Then clearly $d= d_{M'}+d_E$, and the operators $d_E$, $d_{M'}$
%commute.

\begin{defn} Let us write
$\beta_i = dy/(y-f_i)$  for $1\leq i\leq N_h.$
 We define
the \emph{relative bar construction} of  $M$ over the base $M'$ to  be the free $\Or_M$-shuffle algebra:
\begin{equation}
B_{M'}(E) = \Or_{M} \langle \beta_1,\ldots, \beta_{N_h}\rangle\ .
\end{equation}
\end{defn}
%\begin{rem}
%The symbol $[\beta_{i_1}|\ldots |\beta_{i_r}]$ represents the iterated integral (\S5):
%$$\int {dy \over y-f_{i_r}} \ldots {dy \over y-f_{i_2}}   {dy \over y-f_{i_1}}$$
%\end{rem}

The relative bar construction is a differential $\Or_M$-algebra
with respect to the operator $\partial/\partial y$.
  Note that since  $E$ is of dimension $1$, there
is no integrability condition. Proposition
$\ref{propRelnocohomology}$ gives:
\begin{equation} \label{lemfiberprimitives} H^0(B_{M'}(E))=\Or_{M'}\ , \hbox{ and  } H^1(B_{M'}(E))=0\ .
\end{equation}
%Corollary $\ref{corunipclosureisproductforfinitebasepoint}$  then implies the following result.

%For each $1\leq i\leq N_h$, there is a natural map
%$$p : \Or_M
%\rightarrow \Or_{M'}[[\epsilon]][1/\epsilon]$$ which satisfies $p
%\, \partial_y = \partial_\epsilon p$, and which sends $y-f_i$ to
%$\epsilon$.
 Proposition
$\ref{GaussManinversion}$ and its corollary imply the following
result.
\begin{cor} \label{corfib1}  For each  $\Or_{M'}$-linear map  $p : \Or_M
\rightarrow \Or_{M'}[[\epsilon]][1/\epsilon]$
  which satisfies $p \, \partial_y =
\partial_\epsilon p$, there is a natural action of the operators $\partial/\partial x_1,\ldots, \partial/\partial x_{\ell-1}$ on $B_{M'}(E)$, such that
the $\partial/\partial x_i$ commute with $p$. As a result,
 $B(M')\otimes_{\Or_{M'}}
B_{M'}(E)$ is the unipotent closure of $\Or_M$. We deduce that
there is an isomorphism of differential $\Or_M[\partial/\partial
x_1,\ldots, \partial/\partial x_{\ell-1}, \partial/\partial y ]$
algebras:
 $$B(M) \cong B(M') \otimes_{\Or_{M'}} B_{M'}(E)\ .$$
\end{cor}

The following theorem follows by induction.
%The following result follows immediately from theorem
%$\ref{thmRBarsplits}$.
\begin{thm}  \label{thmbarfibration}
 Let $M$ be a fiber-type affine hyperplane arrangement with fibrations $(\ref{supersolv})$.
There is a (non-unique) isomorphism of differential algebras
 $$B(M) \cong B_{{M_1}}(E_1) \otimes_{\Or_{M_1}} \ldots \otimes  B_{M_{\ell-1}}(E_{\ell-1}) \otimes_{\Or_{M_{\ell-1}}} B(E_\ell)\ .$$
\end{thm}

\begin{cor}
The de Rham cohomology of the reduced bar construction on a fiber-type affine hyperplane arrangement defined over a field $k$ is trivial:
$$H^0_{\DR}(B(M))
\cong k\qquad \hbox{ and }\qquad H^i_{\DR}(B(M))=0 \quad \hbox{for all } i\geq 1 \ .$$
\end{cor}
The reason why this result is true is essentially because arrangements of fiber type are rational $K(\pi,1)$ spaces (see [F-R], [H-M]).
By $(\ref{ExpSimpCoords})$,  $\Mod_{0,S}$ is a fiber-type affine hyperplane arrangement over $\Q$.

\begin{cor} \label{corBarmodtriv} In the case of moduli spaces $\Mod_{0,S}$ this gives: % has trivial de Rham cohomology:
$$H^0_{\DR}(B(\Mod_{0,S})) = \Q\ , \quad \hbox{ and } \quad H^i_{\DR}(B(\Mod_{0,S})) =0 \quad \hbox{ for all } \quad i\geq 1\ .$$
The primitive of a closed form $f\in W^b \Omega^i B(\Mod_{0,S})$ is of weight at most $b+1$.
\end{cor}

This result can be proved directly using the fact that the
hyperplane arrangement $\Mod_{0,S}$ is quadratic (see appendix 1).
This is equivalent to corollary 8.7 in [H-M], since $\Mod_{0,p+2}=
Y_1^p$ in the notation of that paper. The fact that primitives
increase the weight by at most one is clear from the definition of
the differential $(\ref{bardifdef})$ on $B(\Or_M)$.

%The choice of the base point $\varepsilon:M\rightarrow k$ is somewhat arbitrary.
%In the case where a hyperplane configuration is of fiber type at infinity, we can circumvent this difficulty by taking
%  base points at infinity. Theorem $\ref{Unipfibrationthmatinfinity}$ then implies the following stronger  version of theorem $\ref{thmbarfibration}$.
%\begin{thm} \label{thmbarfibrationinf}
% Let $M$ be a fiber-type affine hyperplane arrangement at infinity $(\ref{supersolv})$.
%Then every generalised base point at infinity
%defines a canonical isomorphism of  differential algebras
% $$B(M) \cong B_{{M_1}}(E_1) \otimes_{\Or_{M_1}} \ldots \otimes  B_{M_{\ell-1}}(E_{\ell-1}) \otimes_{\Or_{M_{\ell-1}}} B(E_\ell)\ .$$
%\end{thm}

In the case of the moduli spaces $\Mod_{0,S}$, we can make the
decomposition of theorem $\ref{thmbarfibration}$ totally canonical
by working in cubical coordinates $(\ref{ExpCubeCoords})$. The
corresponding fibrations are given by the maps $(x_1,\ldots,
x_\ell) \mapsto (x_1,\ldots, x_{\ell-1})$ (\S\ref{sect25}).
Furthermore, there is a base-point at
%
%there is a base point at infinity at each corner of the Stasheff
%polytope $\overline{X}_{S,\delta}$. It can be written down using
%the corresponding vertex coordinates, and gives rise to a
%canonical decomposition of $B(\Mod_{0,S})$ as a product of shuffle
%algebras. We consider the base point at
infinity corresponding to the origin, which is compatible with
this sequence of fibrations. It is given by the map:
\begin{eqnarray}\label{cubicalbasepointatinfinity}
\Or(\Mod_{0,S}) \cong \Q\Big[(x_i^{\pm 1})_{1\leq i\leq \ell},  \Big({1\over 1-x_{i}\ldots x_j}\Big)_{1\leq i\leq j\leq \ell}\Big]
 &\To& \kes \nonumber \\
x_i & \mapsto & \epsilon_i   \ .  \nonumber
\end{eqnarray}
There is a corresponding logarithmic Laurent expansion over this point, whose map of constants is trivial:
\begin{eqnarray}\label{cubicaltangentialbasepoint}
B(\Mod_{0,S})
 &\To& \Ues  \\
\sum_{I=(i_1,\ldots, i_m)} c_I [\omega_{i_1}|\ldots |\omega_{i_m}] & \mapsto & 0  \ .  \nonumber
\end{eqnarray}

%In cubical coordinates, $\Mod_{0,S}$ is no longer the complement of an affine hyperplane arrangement, but the arguments
%of $\S\ref{sect32b}$ still apply.
Because we have fixed a $\kes$-point, the isomorphism in theorem
$\ref{thmbarfibration}$ is unique.

\begin{cor}\label{corcubicaltensoratorigin} In cubical coordinates, there is a canonical isomorphism
$$B(\Mod_{0,S}) \cong \Or(\Mod_{0,S}) \otimes_\Q \bigotimes_{k=1}^\ell \Q\langle [d\log x_k],
[d\log (1-x_i\ldots x_k)]_{1\leq i\leq k}\rangle \ ,$$ where the
algebras on the right are  free shuffle algebras.
\end{cor}

There is  a similar decomposition for any set of vertex
coordinates $x_1^\alpha,\ldots x_\ell^\alpha$, where $\alpha
\in\chi^\ell_{S,\delta}$ does  not contain an internal triangle.
%This is illustrated in the case of $\Mod_{0,5}$ below.

\begin{rem} \label{remarkonstrongprimitives}
In order to compute the periods of $\Mod_{0,S}$,  we shall only require the fact that $H^\ell(B(\Mod_{0,S}))=0$,
where $\ell=|S|-3$. In cubical coordinates, this is equivalent to  finding a primitive to
$$f \,dx_1\ldots dx_\ell \qquad \hbox{ for all }\quad f \in B(\Mod_{0,S})\ .$$
We have in fact proved a much stronger result. Corollary $\ref{corcubicaltensoratorigin}$  implies that we can find $F\in B(\Mod_{0,S})$ such that
$\partial F/\partial x_\ell= f$. The constant term of $F$ is uniquely determined by the map of constants
$(\ref{cubicaltangentialbasepoint})$. In other words, there is a primitive of the form
$$F\, dx_1\ldots dx_{\ell-1}\ ,$$
where the weight of $F$ is at most one more than the weight of $f$. The primitive $F$ constructed in this way
has the advantage that  is unique.

%Notice  that the decomposition $(\ref{multigrading})$ defines a multi-graduation on the algebra $B(M)$. In the case of moduli spaces $M=\Mod_{0,S}$
%this coincides with the depth filtration on multiple polylogarithms (see $\S6)$.
\end{rem}
%Follows directly using the fact that the algebra is quadratic. Works for more general configurations, for example
%Coxeter groups (refer to future paper).
%
\begin{figure}[h!]
  \begin{center}
%    \leavevmode
    \epsfxsize=5.0cm \epsfbox{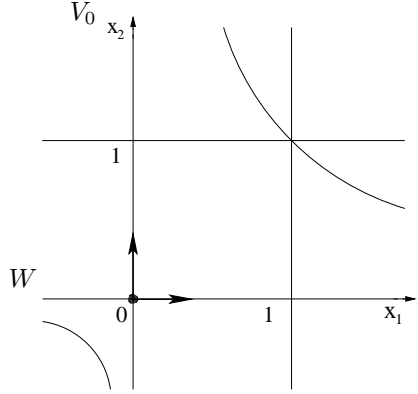}
  \label{basepoints}
%\put(-185,82){$W$}
%\put(-280,140){$V_0$}
\put(-158,39){$W$} \put(-135,140){$V_0$}
 \caption{In cubical coordinates, there is a natural base point at
 infinity on $\Mod_{0,5}$ corresponding to $(0,0)$.
 % It is better to  work in cubical coordinates, shown on the right,
  %which enables us to take the base point at infinity $(0,0)$, since the divisors cross normally there.
  % and project onto $\Pro^1\backslash\{0,1,\infty\}$.
%This is shown on the right in cubical coordinates. %There is a canonical logarithic Laurent  given by the map $B(\Mod_{0,5}) \rightarrow U\{\epsilon_1,\epsilon_2\}$ which
%sends $[d\log x_1]$ to $L_{\epsilon_1}$ and $[d\log x_2]$ to $L_{\epsilon_2}$.
   } \end{center}
\end{figure}
\begin{example}
Consider the fibration $\Mod_{0,5}\rightarrow \Mod_{0,4}$, whose
fibres are isomorphic to $\A^1$ minus 3 points  (fig. 13).
 In
cubical coordinates, we have:
\begin{eqnarray}
\Or_{M'}&\cong &
\Q\big[x,{1\over x}, {1\over 1-x}\big]\ ,\nonumber \\
\Or_{M}& \cong &\Or_{M'}\big[y,{1\over y}, {1\over 1-y}, {x\over
1-xy}\big] \ ,\nonumber
\end{eqnarray}
where  the fibration map is the  projection onto the $x$-axis:
$$(x,y)\mapsto  x : \Mod_{0,5} \rightarrow \Pro^1\backslash \{0,1,\infty\}\ .$$
There is a natural $k\{\epsilon_1,\epsilon_2\}$-point at the origin which sends
\begin{eqnarray} \label{pointatinfform05}
p:\Or_{M}&\To& k\{\epsilon_1,\epsilon_2\} \\
x& \mapsto & \epsilon_1 \nonumber \\
y& \mapsto & \epsilon_2\ , \nonumber
\end{eqnarray}
and which maps,  for example, $x/(1-xy)$ to $  \sum_{i\geq 0} \epsilon^{i+1}_1\epsilon^i_2$.
The differential  algebra $B(M')=B(\Pro^1\backslash\{0,1,\infty\})$ is the universal algebra
of multiple polylogarithms in one variable defined in [Br1], %As an algebra it is
 %the free non-commutative $\Or_{M'}$-algebra on the  symbols ${dx/ x}, {dx/ 1-x}$.
 and $B_{M'}(E)$
is the relative bar construction over $\Or_M'$. As algebras, each
one is the  free non-commutative algebra on two (respectively
three) symbols:
$$B(M') = \Or_{M'} \langle {dx\over x} , {dx \over 1-x} \rangle\ , \quad \hbox{ and }
\quad B_{M'}(E) = \Or_M \langle {dy \over y}, {dy \over 1-y}, {xdy \over 1-xy}\rangle  \ .$$
 Corollary $\ref{corcubicaltensoratorigin}$ gives a canonical isomorphism
$$ B(M')\otimes_{\Or_{M'}} B_{M'}(E) \overset{\sim}{\To} B(M)\ ,$$
and enables us to write down a basis of integrable words in $B(M)$. However, the map $B_{M'}(E) \rightarrow B(M)$ is far from trivial.
 For
example, it gives
$$ 1\otimes \big[ {dy \over 1-y} \big| {x dy \over 1-xy}\big] \mapsto \big[{dy \over 1-y} -{dx \over 1-x} -{dx \over x}\big|
{xdy+ydx \over 1-xy} \big]
+ \big[ {dx \over 1-x} \big|  {dy \over 1-y}\big]\ .$$
A similar formula was given in [Zh] and [Ha3].
It is obvious that the left-hand side is integrable, the right hand side not so. The left-hand coding can be retrieved from the one
on the right by formally setting $dx=0$. The map $B_{M'}(E) \rightarrow B(M)$  is canonically normalised in such a way that,
 apart from all terms of the form $[{dy \over y}]\sha\ldots \sha [{dy\over y}]$, its image vanishes on setting $dy=y=0$.
The logarithmic Laurent expansion of this example %(regularised at the origin)
 is given by the multiple logarithm  (see $\S\ref{sect54}$): %which also vanishes along $y=0$:
$$\Li_{1,1}(x,y)=\sum_{0<k<l} {x^k y^l \over kl}\ ,$$
where we have written $x,y$ instead of $\epsilon_1,\epsilon_2$.
 The coding on the left-hand side of the equation above
only takes into account the differential equations which $\Li_{1,1}(x,y)$ satisfies with respect to the variable $y$ (which are very simple), the right-hand side
encodes the differential equations with respect to both variables $x$ and $y$. The coproduct  of $\Li_{1,1}(x,y)$ can be read off the right-hand coding directly:
$$\Delta \Li_{1,1}(x,y) = \Li_{1,1}(x,y)\otimes 1 + (\log(1-y)-\log(1-x)+\log x) \otimes \log(1-xy) $$
$$+ \log(1-x)\otimes \log(1-y) + 1\otimes \Li_{1,1}(x,y) \ .$$
%Using the results of [H-Z1,2] one can
One can  compare this  with the coproduct  for the motivic multiple polylogarithms defined by Goncharov [Go1].

We therefore have two different points of view on $B(\Mod_{0,S})$.
On the one hand, there is a direct definition in terms of
hyperplane configurations, from which the differential structure
and the action of the symmetric group are evident. The problem is
that the complexity of the set of
 integrable words grows rapidly, and the algebraic structure is obscured.
 On the other hand, using the fibration map above, we have  a description of $B(\Mod_{0,S})$ as
a product of free shuffle algebras, from which its algebraic structure is completely evident. But this point of view breaks the symmetry  and
 only part of the differential structure is visible.
 By exploiting both points of view, one can deduce a lot
 of
 information about the structure of $B(\Mod_{0,S})$. In
 particular, by regarding it as a representation of the symmetric
 group, one obtains many interesting  functional relations between
 multiple polylogarithms.

%), which are certain local systems
%on moduli spaces $\Mod_{0,n}$.
% The two points of view are very close, by the results of [H-Z1,2].% on unipotent variations of mixed Hodge structures.
\end{example}

\newpage

\section{Manifolds with corners and Fuchsian differential
equations}

Let $X$ denote a real analytic manifold with corners. We  consider
functions on $X$ which have logarithmic singularities along the
boundary of $X$, and we  define  the regularised limit of such a
function along components of the boundary $\partial X$. Next, we
state and prove a generalised Fuchs theorem in many variables, and
show that, in the unipotent case, we obtain solutions on $X$ which
are precisely of this type, {\it i.e.}, which have logarithmic
singularities along $\partial X$. Finally,  we state a version of
Stokes' theorem in the case when  $X$ is compact. This requires
some regularity results which allow the integration of functions
with logarithmic divergences along the boundary of $X$. The
example to bear in mind throughout this section is when
$X=\overline{X}_{S,\delta}$ is the closed Stasheff polytope
defined in $\S\ref{sect24}$.
% The analytic part of the proof of the main theorem is
%contained in a Stokes' formula for manifolds with corners.

\subsection{Manifolds with corners} \label{sect41}
A manifold with corners $X$ is a differentiable manifold whose
charts are diffeomorphic to sets of the form
$$U_{p,q} = \R^p \times \R^q_+\ ,$$
where $\R_+= \{x\in \R: x\geq 0\}$, and $p,q \geq 0$ [B-S]. If $q\geq 1$,  the boundary
of $U_{p,q}$ is
\begin{equation} \label{Uboundarystrat}
\partial U_{p,q} = \bigcup_{i+j=q-1} \R^p \times \R_+^{i} \times
\{0\} \times \R_+^j\ ,
\end{equation}
 which is a union of sets diffeomorphic to
$U_{p,q-1}$, and is empty if $q=0$. Let $\partial^i U_{p,q}$
denote the successive submanifolds with corners obtained by
iteration. There is  a stratification
$$U_{p,q} \supseteq
\partial U_{p,q} \supseteq \ldots \supseteq \partial^q U_{p,q}\ ,$$
which has the combinatorial structure of a face  of a hypercube.
There are many different ways to define maps between charts
depending on how rigid we wish to make the manifold $X$. We require  that %directional
derivatives of maps between charts do not vanish along boundary components, and
% The function
%$x\mapsto x^2: \R_+ \rightarrow \R_+$, for example, is ruled out.
 in order for logarithmic regularisation to be well-defined,
we  must  rule out maps  of the form $x\mapsto k  x: \R_+ \rightarrow \R_+$,  where $k\neq1$.

\begin{defn} \label{defnmwc}

Let $n=p+q\geq 1$, and let
$x_1,\ldots, x_n$ be coordinates on $\R^n$ such that
$U_{p,q}=\R^p\times \R^q_+=\{x_1\geq 0,\ldots, x_q\geq 0\}.$ Let $\Sym_q$ denote the symmetric group on $q$
letters which  permutes the indices $1,\ldots, q$.
We define $\Hom_{\an}(U_{p,q}, U_{p,q})$ to be
  the ring of analytic isomorphisms ({\it i.e.}, whose Jacobian does not  vanish anywhere along $U_{p,q}$):
% which extend to give an  analytic  map in an open neighbourhood of $U_{p,q}\subset \R^{p+q}$ which
%is everywhere invertible):
$$\phi=(\phi_1,\ldots, \phi_n): U_{p,q} \To
U_{p,q}$$ which  permute the components  of  the boundary  $\partial U_{p,q}$, {\it i.e.},
$$\phi_i\big|_{x_{\sigma(i)}=0}=0\ , \qquad \hbox{ for } 1\leq i\leq q\ ,$$
 where $\sigma \in \Sym_q$, and which satisfy
$${\partial \phi_i \over \partial x_{\sigma(i)}}\Big|_{x_{\sigma(i)}=0} = 1\ . $$
 In other words, $\phi=(x_{\sigma(1)} f_1,\ldots, x_{\sigma(q)} f_q, \phi_{q+1},\ldots, \phi_n)$,
where $f_i$ are analytic functions such that  $f_i$  is
identically equal to $1$ along the boundary component  $
\{x_{\sigma(i)}=0\}\subset \partial U_{p,q}$, for $1\leq i\leq q$.
For example, if $n=2$ and  $U_{0,2}=\{(x_1,x_2):x_1,x_2\geq 0\}$,
the map $\phi(x_1,x_2) = (x_1+x_1^2x_2, x_2+x_2^2)$ is in
$\Hom_{\an}(U_{0,2},U_{0,2})$.

 We define  an \emph{analytic
manifold with corners} to be  a manifold with corners,  all of whose
transition maps lie in $\Hom_{\an}(U_{p,q},U_{p,q})$.
\end{defn}

% and
%furthermore, that the restriction $\phi|_{\partial U_{p,q}}
%\partial U_{p,q}\rightarrow
%\partial U_{p,q}$ is also in $\Hom_{an}(U_{p,q}, U_{p,q})$ for all $1\leq i\leq q$.

It follows from this definition that  any $\phi \in
\Hom_{\an}(U_{p,q}, U_{p,q})$ preserves the boundary
stratification   of $U_{p,q}$ $(\ref{Uboundarystrat})$. Any analytic manifold with corners
therefore admits a global stratification
$$X=X_0   \supseteq X_1 \supseteq X_2 \supseteq \ldots \supseteq X_n\ ,$$
where each $X_i$ is a manifold with corners, and $X_{i+1}=
\partial X_i$ is the union of the boundary   components of $X_i$.
%We shall write
%$X_i^\circ=X_i \backslash X_{i+1}$ to denote the interior of the
%stratum  $X_i$.

%From here to the end of section - tie in with definitions of first
%section, and state lemma ; stasheff is a mwc. Give charts in nhd
%of each partial triangulation. Use fact that normal crossing.
%Think about explicit corners.
%\\

 Consider  the closed Stasheff polytope $X=\overline{X}_{S,\delta}$ contained in
$\Modf_{0,S}(\R)$. Then $X$ is a manifold with corners whose
stratification is given by $(\ref{Firststrat})$. To see this,
%we
%define an explicit atlas of $\overline{X}_{S,\delta}$ as follows.
let $0<\varepsilon \ll 1$ denote a small constant, and let $e\in
\chi^{k}_{S,\delta}$  denote a $k$-decomposition of the regular
$n$-gon $(S,\delta)$.
%Then $F_e$ is a codimension $k$ face of
%$\overline{X}_{S,\delta}$.
 This can be completed, in a non-unique way, to a full
triangulation $\alpha\in \chi^{\ell}_{S,\delta}$. Then $F_\alpha$
is a corner contained in the face $F_e=\{u_{ij}=0:\,\, \{i,j\}\in
e\}$. Define
\begin{equation}\label{Uchartdefn}
U_e(\varepsilon) = \{0\leq u_{ij}< \varepsilon  \hbox{ for }
\{i,j\} \in e\ , 0<u_{kl}<1  \hbox{ for } \{k,l\} \in \alpha
\backslash e \}\subset \Modf_{0,S}(\R)\ .
\end{equation}
 Since we know that
$\{u_{ij},\, \, \{i,j\} \in \alpha\}$ defines a local coordinate
system  on $\Modf_{0,S}(\R)$ (proposition  $\ref{prop211})$, $U_e(\varepsilon)$ is diffeomorphic to
a chart $U_{\ell-k,k}=\R^{\ell-k}\times \R_+^k$ when $\varepsilon$ is sufficiently small. We have (\emph{c.f.}, $(\ref{wierdstrat})$):
$$\overline{X}_{S,\delta} =  X_{S,\delta}\cup \bigcup_{k\geq 1}\bigcup_{e\in \chi^k_{S,\delta}}
U_e(\varepsilon)\ , \qquad \hbox{  for some } \varepsilon >0 \ .$$
This proves that $\overline{X}_{S,\delta}$ is indeed an analytic manifold with corners, %according to the previous definition,
since all transition maps between  boundary components of charts are given by permutations of coordinates. The action of the
dihedral group of symmetries on $\overline{X}_{S,\delta}$ is a morphism of analytic manifolds with corners.

%In fact, one can find a positive
%$\varepsilon$ for which the sets $U_e(\varepsilon)$ completely
%cover $\overline{X}_{S,\delta}$.

% we have
%$$X^1 = \bigcup_{\{i,j\}\in \chi_n} D_{ij},\,\ldots,\, X^\ell =
%\bigcup_{v\in V^{\delta}} v\ ,$$  and in general, % a component of
%the codimension
%$q$ boundary $X^q$ is a union of spaces
%$$\bigcap_{k=1}^q D_{i_k, j_k} = \{u_{i_1\,j_1}=\ldots=
%u_{i_q\,j_q}=0\}\ ,$$ where $\{i_1,j_1\},\ldots, \{i_q,j_q\}\in
%\chi_n$ are  $q$ non-intersecting chords in an $n$-gon. It will be
%useful to define  explicit %atlas
%charts on $\overline{X}_{S,\delta}$ in the neighbourhood of each
%vertex. Let $v\in V^\delta$, and define for each $0<t\leq 1$,
%$$U_v(t)= \{ 0\leq u_{i\,j} < t : \quad \{i,j\} \in F_v\} \ .$$
%Each set $U_v(t)$ is diffeomorphic to $\R^\ell_+$, and we deduce
%the following algebraic description of $X_{S,\delta}$ as a
%manifold with corners.
%\begin{lem} There exists $t_n<1$ such that for  all $t\geq t_n$, \quad
%$\overline{X}_{S,\delta} = \bigcup_{v\in V^\delta} U_v(t). $
%\end{lem}

\subsection{Logarithmic singularities and regularisation} \label{sect42}
%-------------NB_------------------------------
%-----I have changed complex valued functions to
%------Real valued functions ------------------
%---------------------------------------------------
 We define three sheaves of functions on an analytic
manifold  with corners $X$ which have  singularities along its
boundary $\partial X$. They are:
$$\Fa^{\an} \subset \Fa^{\log} \subset \Fa^{\log}_p\ ,$$
where $\Fa^{\an}$ denotes the sheaf of analytic functions on $X$,
$\Fa^{\log}$ denotes the sheaf of  functions with logarithmic
singularities along $\partial X$, and $\Fa^{\log}_p$ denotes the
sheaf of functions with both logarithmic singularities and
ordinary poles along $\partial X$.

More precisely, let $p,q\geq 0$ where $n=p+q\geq 1$, and let
$x_1,\ldots, x_n$ be coordinates on $\R^n$ such that
$U_{p,q}=\R^p\times \R^q_+=\{x_1\geq 0,\ldots, x_q\geq 0\}. $ Then
we define
$$\Fa^{\an}(U_{p,q}) \subset  \R[[x_1,\ldots, x_n]]\ ,$$
to be the ring of convergent Taylor series in the variables
$x_1,\ldots, x_n$. Next, we define
\begin{eqnarray} \label{Sheafdef}
\Fa^{\log} (U_{p,q})& =&
\Fa^{\an}(U_{p,q})[\log x_1,\ldots, \log x_q ]\ ,  \\
\Fa^{\log}_p (U_{p,q}) &=&  \Fa^{\an}(U_{p,q})[x_1^{-1},\ldots,
x_q^{-1}, \log x_1,\ldots, \log x_q ]\ , \nonumber
\end{eqnarray}
where $\log x_i$  is the principal branch of the logarithm  along
$\R_+$. It follows by a monodromy argument that the functions
$\log x_i$ are linearly independent over the ring
$\Fa^{\an}(U_{p,q})$. Similar rings of functions in one variable
(polynomials in $\log x$ with analytic coefficients) were
considered in [BM].

For each $1\leq i\leq n$, let $v_i$ denote the valuation map on
$\Fa^{\an}(U_{p,q})$ which associates to any function the order of
its vanishing along $x_i=0$. It extends to a valuation
 $$v_i :\Fa_p^{\log}(U_{p,q}) \To \Z \ ,$$
once we have adopted the convention that $v_i(\log x_i)=0$.

\begin{lem} Let $X$ denote an analytic manifold with corners. Then
 $\Fa^{\an}$, $\Fa^{\log}$, and $\Fa^{\log}_p$ define sheaves on $X$,
and for each boundary component $D$ of $\partial X$, the valuation
map $v_D $ on $ \Fa^{\log}_p$ is well-defined.
\end{lem}

\begin{proof} Let $\phi \in \Hom_{\an}(U_{p,q},U_{p,q})$. It
suffices to check that the composition with $\phi$ preserves
$\Fa^{\log} (U_{p,q})$. Let $\phi=(\phi_1,\ldots, \phi_n)$. By definition $\ref{defnmwc}$, and by permuting the coordinates if necessary, we have
\begin{equation}\label{phival}
\phi_i (x_1,\ldots, x_n) = x_i f_i(x_1,\ldots, x_n)\ , \qquad \hbox{ for } 1\leq i\leq q\ ,
\end{equation}
 where $f_i \in \Fa^{\an}(U_{p,q})$. This implies that $f_i(x_1,\ldots, x_n)\geq 0$ for all $(x_1,\ldots, x_n)\in U_{p,q}$, and
 furthermore, $v_i(f_i)=0$. It follows that
$$\log (\phi_i(x_1,\ldots, x_n))= \log (x_i) + \log (f_i(x_1,\ldots, x_n))\ , \hbox{ for } 1\leq i\leq q\ ,$$
where $\log f_i\in \Fa^{\an}(U_{p,q})$ is analytic. It follows
that
 $\phi^*  \Fa^{\log}(U_{p,q}) \subset \Fa^{\log} (U_{p,q})$, and,  similarly,
  $\phi^* \Fa_p^{\log}(U_{p,q}) \subset \Fa_p^{\log}(U_{p,q})$.
The fact that the valuations are well-defined along the components
of $\partial X$ follows immediately from $(\ref{phival})$.
\end{proof}

We can define the regularised value of a function along boundary components of $X$ by formally
 setting the functions $\log x_i$ to $0$, for  $1\leq i\leq q$,  on each chart $U_{p,q}$ of $X$. %and goes back to Hadamard.
\begin{defn} \label{defnlogreg} Let $f\in \Fo^{\log}(U_{p,q})$, and let $1\leq l\leq q$. We can  write
$$f= \sum_{I=(i_1,\ldots, i_l)\in \N^l } f_I  \, \log^{i_1} x_1\ldots \log^{i_l} x_l\ , \quad
\hbox{ where}\quad f_I \in \Fo^{\log}(U_{p+l,q-l}) \ ,$$ and $f_I$
is zero for all but finitely many indices $I$.
%The \emph{regularized coefficient} of $f$ is the function
%$\Reg\, f = f_{(0,\ldots, 0)}(z_1,\ldots, z_n) \in A$.
 The \emph{regularized value of $f$ along }
$D=\{(x_1,\ldots,x_n): x_1=\ldots=x_l=0\}\subset \partial^l
U_{p,q}$  is defined to be:
%
%
%Let $C= \{(x_1,\ldots,x_q)\in U_{p,q}: \, x_1=\ldots =x_q=0\}=\partial^q U_{p,q}$ denote the corner of $U_{p,q}$.
$$\Reg (f, D) = f_{(0,\ldots,0)} (0,\ldots, 0,x_{l+1},\ldots,x_q, x_{q+1},\ldots, x_n)\ ,$$
viewed as a function on $D\cong  U_{p,q-l}$. By construction,
$\Reg(f,D) \in \Fo^{\log}(D)$.
\end{defn}

\begin{propdef}
Let $X$ denote an analytic manifold with corners, and let $D \subset \partial^l X$ denote any boundary component of $X$.
Then there is a well-defined regularisation map along the component $D$:
$$\Reg(\bullet, D): \Fo^{\log}(X) \To \Fo^{\log} (D)\ .$$
\end{propdef}
\begin{proof} The transition maps are compatible with logarithmic regularisation by definition $\ref{defnmwc}$. Let $\phi \in \Hom_{\an}(U_{p,q},U_{p,q})$.
Then, up to permuting coordinates, $\phi$
is of the form $\phi=(x_1(1+x_1 g_1), \ldots , x_q(1+x_q g_q), \phi_{q+1},\ldots, \phi_n)$. It follows that
$$\log \phi_i = \log x_i + \log (1+ x_i g_i)\ , \quad \hbox{for } 1\leq i\leq q\ ,$$
and the term $\log(1+x_ig_i)$ vanishes at $x_i=0$. It follows that  logarithmic regularisation along $x_i=0$ is well-defined for $1\leq i\leq q.$
By regularising with respect to one variable at a time, it follows that
regularisation along an arbitrary boundary component $D\subset \partial^l X$
is well-defined also.
\end{proof}

\begin{rem} \label{remregofpolestoo}
Clearly we can extend the regularisation map  for polar
singularities
$$\Reg(\bullet, D): \Fo_p^{\log}(X) \To \Fo_p^{\log} (D)\ ,$$ by
mapping all negative powers of coordinates $x_1,\ldots, x_l$ to
zero also (this is just a map of constants as defined in
$\S\ref{sect32b}$).
\end{rem}

\subsection{Fuchsian differential equations in several complex variables.} \label{sect43}

%The boundary of a manifold with corners  locally  resembles a connected component in the set of real points of the complement of a
%normal crossing divisor.

Consider the  open complex affine space obtained by complexifying $U_{p,q}$:
$$\C^{p+q}
\backslash \{z_1\ldots z_q=0\}\ ,$$ and let $V_{p,q}$ denote an
open polydisk neighbourhood of the origin contained in
$\C^{p+q}\backslash \{z_1\ldots z_q=0\}$.
%We define a ring
%$$\Fa^{\an}(V_{p,q}) = \Fa^{\an}(\C^n) \subset \C[[z_1,\ldots, z_n]]\ ,$$
% of convergent power series in the complex
%variables $z_1,\ldots, z_n$, and set
%\begin{eqnarray}
%\Fa^{\log} (V_{p,q})& =& \C[\log z_1,\ldots, \log z_q ] \otimes
%%\Fa^{\an}(V_{p,q})\ ,  \\
%\Fa^{\log}_p(V_{p,q}) &=& \C[z_1^{-1},\ldots, z_q^{-1}, \log
%z_1,\ldots, \log z_q ] \otimes \Fa^{\an}(V_{p,q})\ . \nonumber
%\end{eqnarray}
%In this case we have to choose branches $\log z_1,\ldots,
%\log z_q$ of the logarithm. Thus we obtain sheaves $\Fo^{\an} $,
%$\Fo^{\log}$ and $\Fo_p^{\log}$ of multi-valued functions on any manifold which is covered
%by charts of the form $V_{p,q}$. The manifold $\Mod_{0,S}$ is such    an example, by lemmas $\ref{lemUalphascover}$
%and $\ref{lemfullcompactif}$.
%Given any analytic manifold with corners $X$, we typically obtain
%functions in the sheaf $\Fa^{\log}(X)$ as solutions of
%differential equations with singularities in $\partial X$. To this
%end,
We require a
 generalised Fuchs' theorem  which we solve locally on the spaces $V_{p,q}$.
Let $m\geq 1$, and consider the
%following matrix
 differential equation:
\begin{equation} \label{fuchs}
d F = \Omega F\ ,
\end{equation}
 where $F$ takes values in the set of $m\times m$ complex matrices $M_m(\C)$, and
where
\begin{equation}\label{fuchsomega}
 \Omega = \sum_{i=1}^n N_i \, {dz_i \over z_i} + A_i
dz_i  \ .\end{equation} Here,  $N_i\in M_m(\C)$ are constant
matrices, and each $A_i$ is a holomorphic  function on $V_{p,q}$,
which takes values in $M_m(\C)$. Assume that $\Omega$ is
integrable, {\it i.e.},
$$d\Omega= \Omega\wedge \Omega\ .$$
This implies, in particular,
that the matrices $N_i$ commute:
\begin{equation} \label{Ncomm}
[N_i,N_j] = 0 \qquad \hbox{ for all } \quad 1\leq i,j\leq n\ .
\end{equation} Let us write $n=p+q\geq 1$, and suppose that $N_i= 0$ for all $
q+1\leq i\leq n$. The form $\Omega$ is continuous on $V_{p,q}$.
Let us fix branches of the logarithm $\log z_i$ for $i=1,\ldots,
q$ on $V_{p,q}$. In practice, we will choose a real subspace
$\R^p\times \R^q_+\cong U_{p,q} \subset V_{p,q}$ and take the
unique branches of the logarithm which are real-valued on
$U_{p,q}$. The function
$$D= \exp(\sum_{i=1}^q N_i \log z_i)\ ,$$
is a well-defined multi-valued function on $V_{p,q}$ because the matrices $N_i$ commute.
%$U_{p,q}=
%\R^p\times \R^q_+$.
%Let $\log z_i$ denote the branch of the logarithms fixed
%earlier.
The following  result is a  generalized Fuchs' theorem in many
complex variables.
Similar situations have been considered in [De3,Yo]. %but we give a complete proof here.

\begin{thm} \label{thmfuchs} %(Fuchs' theorem in many complex variables.)
Suppose that for each $1\leq i\leq n$, no pair of eigenvalues
of the matrix $N_i$  differ by an non-zero integer. Let $H_0\in M_m(\C)$ be
any constant matrix. Then $(\ref{fuchs})$ has a unique solution
$$F= H \, D\ ,$$
where  $H:V_{p,q}\rightarrow M_m(\C)$ is  holomorphic  and takes
the value $H_0$ at the origin.

\end{thm}

\begin{proof} The matrix $D$ is invertible, and is a solution to
the differential equation
$$d D =\big( \sum_{i=1}^n N_i {dz_i \over z_i}\big)\, D\ .$$
It follows that $F=HD $ is a solution of $(\ref{fuchs})$ if and only
if \begin{equation}\label{Hfuchseq} dH = \sum_{i=1}^n\big[  N_i ,
H\big]{dz_i \over z_i} + A_i H\, dz_i \ .
\end{equation}
If we write $\partial_k = \partial/\partial z_k$, then  this is equivalent to the set of equations
\begin{equation}\label{Hfuchpartials} \partial_k \, H = \big[  N_k ,
H\big]{1\over z_k} + A_k H \ ,\qquad \hbox{ for } 1\leq k\leq n\ .
\end{equation}
 A solution $H$ is
holomorphic on $V_{p,q}$ %in the neighbourhood of $0$
 if and only if it can be
written as a power series
\begin{equation} \label{Hcomdev}
H= \sum_{0\leq i_1,\ldots, 0\leq i_n} H_{(i_1,\ldots, i_n)}
z_1^{i_1}\ldots z_n^{i_n} \ , \end{equation}
 where the coefficients $H_{(i_1,\ldots,
i_n)}\in M_m(\C)$ satisfy a growth condition. By substituting
such a power series expansion into $(\ref{Hfuchpartials})$ and
considering the coefficient of $z_1^{i_1}\ldots z_n^{i_n}$, we obtain
the following recurrence relations:
\begin{equation} \label{Hrecur}
\big(i_{k} + 1 - \ad(N_k)\big) \, H_{(i_1,\ldots,i_k+1,\ldots,
i_n)}= \sum_{0\leq j_1\leq i_1,\ldots ,0\leq j_n\leq i_n} (A_k)_{(j_1,\ldots, j_n)}  \,
H_{(i_1-j_1,\ldots, i_n-j_n)}
\end{equation}
for each $1\leq k\leq n$,  where $(A_k)_{(j_1,\ldots, j_n)}$
are the coefficients in the power series expansion of $A_k$. Now consider a matrix $M\in M_m(\C)$.
If we denote the  eigenvalues of  $M$ by
$\alpha_1,\ldots, \alpha_m$, then the eigenvalues of $\ad \, M$
are $\alpha_{i}-\alpha_{j}$. The assumption on the eigenvalues of
$N_i$ is therefore equivalent to the invertibility of the
operators
$$\big(m - \ad(N_k)\big)\qquad \hbox{ for all } m \in\N \ ,$$
and for each $1\leq k\leq n$. The operator on the left hand side of $(\ref{Hrecur})$ is therefore invertible, so
 we can  solve $(\ref{Hrecur})$
iteratively, provided that these equations  are
compatible. This means that we must show that the two different ways of
obtaining $H_{( i_1,\ldots, i_k+1,\ldots, i_l+1,\ldots, i_n)}$
by applying $(\ref{Hrecur})$   first for $k$ and then for $l$,
or the other way round, both lead to the same answer. This is
equivalent to the integrability of the form $\Omega$. In order to see this, write
$\Omega=\sum_i \Omega_i dz_i$, where $\Omega_i = R_i + A_i$, and $R_i= N_i/z_i$, for $1\leq i\leq n$.
 The integrability of $\Omega$ and the commutativity of the $R_i$ implies the following equations for all $1\leq i,j\leq n$:
 \begin{eqnarray}
\partial_j \Omega_i &=& \partial_i \Omega_j + [\Omega_j,\Omega_i]\ , \nonumber\\
\partial_j R_i &=& \partial_i R_j + [R_j,R_i] \ .\nonumber
 \end{eqnarray}
 It follows that the expression
$$\phi_{ij}(M)=(\partial_j \Omega_i + \Omega_i \Omega_j)M - \Omega_i M R_j - \Omega_j M R_i + M (R_jR_i-\partial_j R_i) $$
is symmetric in $i, j$ for all matrices $M\in M_m(\C)$. But equations $(\ref{Hfuchpartials})$ are precisely the set of equations
$\partial_i H= \Omega_i H - H R_i$ for $1\leq i\leq n$. It follows, on applying $\partial_j$ to each equation,  that
$$\partial_{j} \partial_i \, H = \phi_{ij} H\qquad \hbox{ for all } 1\leq i,j\leq n\ .$$
One can check by differentiating a truncated power series
expansion for $H$  that  the compatibility of the equations
$(\ref{Hrecur})$ up to a given weight is  a consequence of  the
symmetry of the operators $\phi_{kl}$, for all $1\leq k,l\leq n$.
We can therefore solve $(\ref{Hrecur})$ recursively to obtain a
solution of $(\ref{Hfuchseq})$ of the form $(\ref{Hcomdev})$. It
remains to check that the function $H$ defined in this manner is
holomorphic on $V_{p,q}$. Since the series $A_k$ for $1\leq k\leq
n$ are holomorphic on the polydisk $V_{p,q}$, there exist
constants $r_1,\ldots, r_n>0$ and a constant $c>0$  such that
\begin{equation} \label{Acoeffbound}
||(A_k)_{(i_1, \ldots, i_n)}||\leq c\,  r_1^{i_1}\ldots
r_n^{i_n}\qquad \hbox{ for all } 1\leq k\leq n\ .\end{equation}
For $m\geq 1$, let $\varepsilon_m= \sup_{1\leq k\leq n} ||(m- \ad
N_k) ||^{-1}$. By the assumption on the eigenvalues of $N_k$ and
the remarks above, $\varepsilon_m\rightarrow 0$ as $m\rightarrow
\infty$. It follows from $(\ref{Hrecur})$ that
\begin{equation} \label{Hproof1}|| H_{(i_1,\ldots,i_k+1,\ldots,
i_n)}||\leq  \varepsilon_{i_k+1}  \sum_{0\leq j_l\leq i_l}
||(A_k)_{(j_1,\ldots, j_n)}|| \, ||  \, H_{(i_1-j_1,\ldots,
i_n-j_n)}||\ .\end{equation} Now let  $s_k>r_k$, for $1\leq k\leq
n$,  and
 let $m$ be sufficiently large
such that $\varepsilon_{m}\, c  \prod_{i=1}^n ({s_i\over s_i-r_i})
<1$. Set
$$e = \sup_{0\leq i_1,\ldots, i_n \leq m_r} {|| H_{(i_1,\ldots,i_n)}|| \over s_1^{i_1}\ldots s_n^{i_n}}<\infty \ .$$
Let $M\geq m$, and suppose by induction that  $||H_{(i_1,\ldots,
i_n)}||\leq e\, s_1^{i_1}\ldots s_n^{i_n}$ for all $0\leq
i_1,\ldots, i_n\leq M$. This is true when $M=m$ by the definition
of $e$. Then, by  applying $(\ref{Acoeffbound})$, we deduce from
$(\ref{Hproof1})$ that
$$ || H_{(i_1,\ldots,M+1,\ldots,
i_n)}||\leq  \varepsilon_{M+1}   \sum_{0\leq j_1\leq i_1,\ldots,
0\leq j_n\leq i_n} c\,e\, \Big({r_1\over s_1}\Big)^{j_1} \ldots
\Big({r_n\over s_n}\Big)^{j_n} s_1^{i_1}\ldots s_n^{i_n}
$$
$$
\leq \varepsilon_{M+1} c \,e \prod_{i=1}^n \Big({s_i \over
s_i-r_i} \Big) \,s_1^{i_1}\ldots s_n^{i_n} \leq
 e\, s_1^{i_1}\ldots s_n^{i_n}
\ .$$ By induction we deduce that $|| H_{(i_1,\ldots, i_n)}||\leq
e\, s_1^{i_1}\ldots s_n^{i_n}$ for all $(i_1,\ldots, i_n)$. This
holds for any set of constants $s_1,\ldots, s_n$ satisfying
$s_i>r_i$, which proves that $H$ is holomorphic on $V_{p,q}$, as
required.
%
%
% to find a holomorphic function $H$ satisfying $(\ref{Hfuchseq})$.
% Now suppose by induction that we have solved $(\ref{Hrecur})$ for all coefficients
%up to total weight $m$. Set $M= \sum H $( weight $m-1$), and we get::
\end{proof}

We will be interested in the case where the matrices $N_i$ are all
nilpotent. It then follows that the matrix
$$D=\exp \big( \sum_{i=1}^q N_i \log z_i\big)$$
has coefficients which are polynomials in $\log z_i$. Since all
eigenvalues of $N_i$ are $0$, the condition of the previous theorem is
satisfied,  and therefore there exists a matrix solution $F$ to equation
$(\ref{fuchs})$ whose  entries $F_{ab}$ are polynomials in $\log z_1,\ldots, \log z_q$ whose coefficients
are convergent Taylor series in $z_1,\ldots, z_n$:
$$F_{ab} \subset \C[[z_1,\ldots, z_n]][\log z_1,\ldots, \log z_q] \ .$$

\begin{defn} \label{defnFuchsUnipType}
Let $X$ denote an analytic manifold with corners. An integrable
1-form $\Omega$ defined on $X$ is \emph{unipotent of Fuchs' type}
if, locally on each chart of the form $U_{p,q}$, $\Omega$
restricts to a 1-form of type $(\ref{fuchsomega})$, where the
matrices $N_i$ are nilpotent.
\end{defn}
As remarked in $\S\ref{sect42}$, there are canonical branches of the functions $\log z_i$ on local charts of $X$. The
solutions to $(\ref{fuchs})$ will therefore be real-valued on $X$.

\begin{cor} \label{corfuchs}
Let $\Omega$ be a real-valued unipotent integrable 1-form of
Fuchs' type on $X$. Suppose that $X$ is simply connected.
%every boundary component $\partial^i X$ in
%the stratification $X\supseteq \partial X\supseteq \ldots
%\supseteq \partial^i X\supseteq \ldots $ is simply connected.
Then any solution $(F_{ab})$ to $(\ref{fuchs})$ defined in the
neighbourhood of any point $x\in X$ extends over the whole of $X$. This gives a  global  solution %$F_{ij}$
of $(\ref{fuchs})$  whose coefficients satisfy $F_{ab}\in
\Gamma(X,\Fa^{\log}).$
%Then given any point $x \in X$, and a matrix $(F_x)_{ij}$, there
%is a unique solution $(F)_{ij}$ to equation $(\ref{fuchs})$
%defined globally on $X$ such that $F_{ij} \in
%\Gamma(X,\Fa^{\log}),$ and ...
\end{cor}

\subsection{Stokes' theorem with logarithmic singularities}
The key argument in our proof of the main theorem is to apply a version of Stokes' theorem to the manifold with corners
$\overline{X}_{S,\delta}$. This requires integrating functions which have logarithmic singularities along the boundary.

% manifold with corners. The variables $z_i$ will be denoted $x_i$ as above.
% which corresponds locally to
%the set of real points in $V_{p,q}$.
 %In this case, there are canonical choices of  branches for the logarithm functions  $\log x_i$, for $1\leq i\leq q$.
\begin{lem} \label{lemmanopolesint} Let $X$ denote a compact analytic manifold with corners of dimension $n$. Let
$\psi\in \Omega^n(X)$ denote an $n$-form on $X$ whose coefficients lie in $\Fo^{\log}_p$. Then $\psi$ is absolutely integrable on
$X$ if and only if $\psi$ has no poles along $\partial X$.
\end{lem}
\begin{proof} If $\psi$ has a pole of order $k\geq 1$ along some component of $\partial
X$, then there is a chart on $X$ of the form $U_{p,1}$ such that
$\psi=f\, dx_1\ldots dx_n$, where $f$ can be written
$$f(x_1,\ldots, x_n) ={1\over x^k_1}\sum_{i= 0}^N   f_i(x_2,\ldots, x_n)\log^i x_1 +{1\over x_1^{k-1}}
\sum_{i=0}^M g_i(x_1,\ldots, x_n)\log^i x_1\ ,$$ where $f_i, g_i\in \Fo^{\an}(U_{p,q})$  are analytic on $x_1>0,\ldots, x_n>0$, and $f_N$ is not identically zero. Since the
term $(\log x_1)^N$ dominates the other powers of $\log x_1$ near
 $x_1=0$, it follows by continuity that there is a small box
$$B(\varepsilon)=\{(x_1,\ldots, x_n) : x_1\in [0,\varepsilon],\,\,
x_2-\alpha_2,\ldots, x_n-\alpha_n \in [-\varepsilon, \varepsilon]\}\ ,$$ where $\alpha_2,\ldots, \alpha_n>0$, and a
constant $c>0$  such that
$$|f| \geq {c\over x^{k}_1}  |\log x_1|^N\ ,$$
 for
all $ (x_1,x_2,\ldots, x_n) \in B(\varepsilon)$ whenever
$\varepsilon>0$ is sufficiently small. It follows that
$$\int_X |f|\, dx_1\ldots dx_n \geq c \, (2\varepsilon)^{n-1}\int_{0}^{\varepsilon} {1\over x^{k}_1}  |\log
x_1|^Ndx_1 = \infty\ ,$$ and therefore $\psi$ is not absolutely
integrable.

Now suppose that $\psi$ has no poles along $\partial X$. Then in
each small chart of the form $U_{p,q}$, we can write $\psi
=f(x_1,\ldots, x_n)dx_1\ldots dx_n$, where
$$f(x_1,\ldots,x_n) = \sum_{I=(i_1,\ldots, i_q) }  (\log x_1)^{i_1} \ldots (\log
x_q)^{i_q} f_I (x_1,\ldots, x_n)\ ,$$ where $f_I(x_1,\ldots,
x_n)\in \Fo^{\an}(U_{p,q})$, and almost all $f_I$ are identically
zero. But the function $\log x$ is integrable on any interval
$[0,t)$, where $t>0$, and since sums and products of integrable
functions are integrable, it follows that $f$ is integrable
locally. Since $X$ is compact, we can find a finite partition of
unity on $X$, and deduce that $f$ is absolutely  integrable over
the whole of $X$.
\end{proof}

We can therefore integrate functions which have at most
logarithmic singularities. The following lemma implies  that
primitives of functions on $X$ which have at most logarithmic
singularities  extend continuously to $\partial X$. The essential
point is that the 1-form $\log x \,dx$ on $\R_+$ has a logarithmic
singularity at $0$, but its primitive,
 $x\log x -x +c$, is continuous at $x=0$.

\begin{lem}\label{lemnopolesextendscont} Let $X$ be an analytic manifold with corners.
Let $\psi\in \Omega^n(X)$ have at most logarithmic singularities along $\partial X$, and let
 $\Psi\in \Omega^{n-1}(X)$ denote a primitive of $\psi$ which has no poles along $\partial X$. Then $\Psi$
is continuous on the interior of $\partial X$.
\end{lem}

\begin{proof} It suffices to prove the result on each chart of $X$ isomorphic to $U_{p,q}$
with coordinates $x_1,\ldots, x_n$ as above.
 Let $\psi=f dx_1\ldots dx_n$, where $f\in \Fo^{\log}(U_{p,q})$. We  write
 $\Psi= \sum_{i=1}^n (-1)^{i-1} F_i dx_1 \ldots \widehat{dx}_i \ldots dx_n$,
  where $F_i \in\Fo^{\log}(U_{p,q})$ for $1\leq i\leq n$. Let
 $$F_i = \sum_{k \geq 0 } \log^k x_i\, F_{i,k} \ ,$$
 where $F_{i,k}\in \Fo^{\log}(U_{p,q})$ is analytic in the coordinate $x_i$ and is zero for all but finitely many
 indices $k$. Since $\sum_{i=1}^n {\partial F_i /\partial x_i} =f$,
  we have
 $$\sum_{i=1}^n \sum_{k\geq 1} {k\log^{k-1} x_i \over x_i} F_{i,k}(x_1,\ldots,x_{i-1}, 0,x_{i+1}, \ldots, x_n)
  \in \Fo^{\log}(U_{p,q}) \ .$$
This implies that $F_{i,k} (x_1,\ldots, x_{i-1},0, x_{i+1},\ldots ,x_n)=0 $ for all $1\leq i\leq n$, $k\geq 1$, and therefore
$F_i\in \C[\log x_1,\ldots, x_i\log x_i,\ldots, \log x_q][[x_1,\ldots, x_n]]$. It follows that
$$\Psi \Big|_{x_i=0} = (-1)^{i-1} F_i\Big|_{x_i=0} dx_1\ldots \widehat{dx}_i \ldots dx_n$$
is continuous for all $1\leq i\leq n$. Thus $\Psi$ is continuous
along the interior of $\partial X$.
\end{proof}

%\begin{defn} \label{defncompletelyreg} Let $M$ denote an analytic manifold modelled on charts  $V_{p,q}$. We will say that
%a sheaf
%$\Fo$ on $M$ is \emph{completely regularisable} if locally on each chart $V_{p,q}$
%there exists $A\subset \Fo^{\an}(V_{p,q})$ which satisfies the condition $(\ref{condA})$, such that
%$$\Fo(V_{p,q}) \cong A \otimes \C[z_1^{-1},\ldots, z_q^{-1}, \log z_1,\ldots, \log z_q]= A^{\log}_p\ ,$$
%There is an analogous definition  if $X$ is an analytic manifold with corners.
%\end{defn}
%A completely regularisable sheaf allows both types of regularisation.

We can now state the
following version of Stokes' theorem.

\begin{thm}\label{corStokesform} Let $X$ denote a compact analytic manifold with corners of dimension $n$.
Let $\psi\in \Omega^n(X)$ be an $n$-form such that $\psi$ has no
poles along $\partial X$, and let $\Psi\in \Omega^{n-1}(X)$ be a
primitive of $\psi$ such that $\Psi$ has no poles along $\partial
X$ either. Then $\Psi$ extends continuously to $\partial X$, and
$$\int_X \psi = \int_{\partial X} \Psi\ ,$$ where both integrals
are finite.
\end{thm}

\begin{proof} Let $U_{p,q}(\varepsilon) = \R^p\times
\R^q_\varepsilon$ where $\R_\varepsilon= \{x\in \R: x\geq
\varepsilon\}$. By lemma $\ref{lemmanopolesint}$, $\psi$ is integrable on $X$.
%By lemma $\ref{lemremovepoles}$, we can find a primitive $\Psi$ of $\psi$ which has no poles on
%$\partial X$ and which is therefore  continuous  on $\partial X$, and hence integrable.
We know that $\Psi$ extends continuously to $\partial X$ by the previous lemma.
On each small chart of $X$ we can apply Stokes' theorem:
$$\int_{U_{p,q}} \psi = \lim_{\varepsilon\rightarrow 0}
\int_{U_{p,q}(\varepsilon)} \psi= \lim_{\varepsilon\rightarrow 0}
\int_{\partial U_{p,q}(\varepsilon)} \Psi= \int_{\partial U_{p,q}}
 \Psi\ , $$
and all terms are finite. Since $X$ is compact, we can find a
finite partition of unity and apply the above identity locally.
The result then follows in exactly the same way as the usual proof
of Stokes' theorem.
\end{proof}
%This result is more general than the usual version of Stokes'
%theorem because the form $\psi$ may be unbounded on the boundary
%$\partial X$. %, and the manifold $X$ is not necessarily compact.

In the case which interests us, when $X=\overline{X}_{S,\delta}$,
 we can define the following
 exhaustion of the polytopes $\overline{X}_{S,\delta}$.  For all
small $\varepsilon>0$, we set
$$
\overline{X}_{S,\delta}^{\,\varepsilon} = \{ u_{ij}\geq \varepsilon,
\quad \{i,j\} \in \chi_{S,\delta}\}\ .$$
The required version of Stokes'  theorem is then immediate:
%Given forms $\psi$ and $\Psi$ where $dF=\psi$, as in the previous theorem,
$$\int_{X_{S,\delta}} \psi = \lim_{\varepsilon\rightarrow 0}
\int_{\overline{X}_{S,\delta}^{\,\varepsilon}} \psi=
\lim_{\varepsilon\rightarrow 0} \int_{\partial
\overline{X}^\varepsilon_{S,\delta}} \Psi= \int_{\partial
\overline{X}_S}
 \Psi\ . $$

\newpage

\section{Hyperlogarithms}

We give an explicit description of the  constructions  in the
previous two sections  when the dimension is $1$, {\it i.e.}, when
$M$ is  the affine line $\A^1$ minus $N+1$ fixed points $\sigma_0,\ldots, \sigma_N$. However, we need to consider iterated integrals whose path of integration
has endpoints at one of the removed points $\sigma_i$, and so do not necessarily converge. This requires a regularisation procedure  which can be solved
for all iterated integrals simultaneously by considering their generating series.

\subsection{Hyperlogarithms and differential equations}
 Let
$N\geq 1$, and let $A=\{a_0,\ldots,a_N\}$ be an alphabet with
$N+1$ letters.
 We fix any injective map of sets $j: A\hookrightarrow \C$,
and set $\sigma_0=j(a_0)$, \ldots, $\sigma_N=j(a_N)$. Let $\Sigma$
denote the  set $j(A)\cup\{ \infty\}$, and let $D=
\Pro^1(\C)\backslash \Sigma $ denote the complex plane with the points
$\sigma_k$ removed.
Consider the following formal differential equation:
\begin{equation}\label{hstar} {\partial \over \partial z}F(z)=
 \sum_{i=0}^N {a_i \over
z-\sigma_i}  F(z)\ ,\end{equation} which is an equation of Fuchs
type, whose singularities are simple poles in $\Sigma$. Let $F(z)$
be a solution on $D$ taking values in $\C\ld A\rd$. If we write
$$F(z)= \sum_{w\in A^*} F_w(z)\,w\ ,$$
 then $(\ref{hstar})$ is equivalent to the system of equations
\begin{equation}\label{star'}
{\partial \over \partial z} F_{a_kw}(z) = {F_{w}(z) \over
z-\sigma_k}\ , \end{equation}
 for all $0\leq k \leq N$ and all  $w \in A^*$, together with the  initial equation $\partial F_1(z)/\partial  z =0$, where $1$ denotes the empty
 word in $A^*$. The  term $F_1(z)$ is therefore constant.

One can construct explicit holomorphic solutions $L_w(z)$ to
$(\ref{star'})$ on a certain  domain $U$ obtained by cutting $\C$.
These functions  extend by analytic continuation to multi-valued
functions on the punctured plane $D$, and can equivalently be
regarded as holomorphic functions on a universal covering space
$p:\widehat{D}\rightarrow D$. Since no confusion arises, we shall
always denote these functions by the same symbol $L_w(z)$.
For each $0\leq k\leq N$, choose  closed half-lines
$\ell({\sigma_k})\subset \C$ starting at $\sigma_k$, such that no
two intersect. Let $U= \C\backslash \bigcup_{\sigma_k\in\Sigma}
\ell(\sigma_k)$ be  the simply-connected open subset of $\C$
obtained by cutting along these half-lines. Fix  a branch  of % the logarithm
 $\log (z-\sigma_0)$ on $\C\backslash \ell(\sigma_0)$.

\begin{prop} \label{prophstarsoln}
Equation $(\ref{hstar})$ has a  unique solution
$L(z)$ on $U$ such that
$$L(z) =  f_0(z) \exp(a_0 \log (z-\sigma_0))\ ,$$  %\quad \hbox{ as } \quad |z-\sigma_0|\rightarrow 0,
where $f_0(z)$ is a holomorphic function on $ \C\backslash
\bigcup_{k\neq 0} \ell(\sigma_k)$ which satisfies
$f_0(\sigma_0)=1$. We write this $L(z) \sim (z-\sigma_0)^{a_0}$ as
$z\rightarrow \sigma_0$. Furthermore, every solution of
$(\ref{hstar})$ which is holomorphic on $U$ can be written $L(z)\,
C,$ where $C\in \C\ld X\rd$ is a constant series ({\it i.e.},
depending only on $\Sigma$, and not on $z$).
\end{prop}

The proposition can be deduced from theorem $\ref{thmfuchs}$,  and a direct solution is given in Gonzalez-Lorca's
thesis [GL]. We use  another   approach here, since we  require an
explicit formula for the functions $L_w(z)$ which is originally
due to Poincar\'e and Lappo-Danilevsky [P, LD].
First, let $A_c^*$ denote the subset of all words in $A^*$ which
do not end in the letter $a_0$, and let $\C\langle A_c\rangle
\subset \C\langle A\rangle$ denote the sub-vector space they generate. It is easy to verify that $\C\langle A_c\rangle$ is preserved by the shuffle product.  If $w \in A_c^*$ and $w\neq 1$, the limiting condition given in the proposition is just $\lim_{z
\rightarrow \sigma_0} L_w(z)=0$. If we write $w=a_0^{n_r}a_{i_r}
a_0^{n_{r-1}}a_{i_{r-1}}\ldots a_0^{n_1} a_{i_1},$
 where $1\!\leq
\! i_1,\ldots,i_r\!\leq \!N$, then  $L_w(z)$ is defined in  a
neighbourhood of $\sigma_0$ by the formula \begin{equation}
\label{convdef}
 \sum_{1\leq m_1 <\cdots<m_r} {(-1)^r \over
m_1^{n_1+1} \!\!\!\ldots m_r^{n_r+1}} \Big({z\!-\!\sigma_0 \over
\sigma_{i_1}\!-\!\sigma_0}\Big)^{m_1}
 \Big({z\!-\!\sigma_0 \over
\sigma_{i_2}\!-\!\sigma_0} \Big)^{m_2-m_1}
 \!\!\!\! \ldots
  \Big(  {z\!-\!\sigma_0 \over
\sigma_{i_r}\!-\!\sigma_0}\Big)^{m_r-m_{r-1}} \end{equation}
 which converges
absolutely for $|z-\sigma_0| <
\inf\{|\sigma_{i_1}-\sigma_0|,\ldots, |\sigma_{i_r}-\sigma_0|\}$.
One can easily check that this defines a family of  holomorphic
functions satisfying the  equations $(\ref{star'})$ in this open
disk, and  that the limiting condition is trivially satisfied.

The functions $L_w(z)$ extend analytically to the whole of $U$ by
the recursive  integral formula:
\begin{equation} \label{int} L_{a_k w}(z)= \int_{\sigma_0}^z {L_{w}(t) \over
t-\sigma_k }\,dt\  , \end{equation} which is valid for all $0\!\leq
\! k\!\leq \!N$ and all $w\in A_c^*$. Since iterated integrals are homomorphisms for the shuffle product (lemma $\ref{lemitintsprop}$), we also have
\begin{equation}\label{shufflesforhypers}
L_w(z) L_{w'}(z) = L_{w\sha w'}(z) \qquad \hbox{for all }\quad  w,w'\in A_c^*\ ,\end{equation}
where $L$ is extended by linearity to all words $w\in \C\langle A_c\rangle$.
  It follows from the definition of the shuffle product
that any word in $A^*$ can be uniquely written as a linear combination of
shuffles of $a_0^n$ with  words in $A_c^*$:
$$w = \sum_{n\geq 0} a_0^n \sha v_n\ ,\hbox{ where } v_n\in  \C\langle A_c\rangle\ .$$
  We  can therefore set
$$L_{a_0}(z)=\log(z-\sigma_0)\ ,$$
 and extend the  definition of $L_w(z)$ to all words $w\in A^*$ by demanding that $L_w(z)$ satisfy the
shuffle relations $L_w (z) L_w'(z) = L_{w\sha w'}(z)$  for all
$w,w'$ in $ A^*.$  One verifies that the functions $L_w(z)$ can be written in the form $(\ref{int})$ for all words
$w\in A^*$, for $w\neq a_0^n$, and are   solutions to $(\ref{star'})$.
In order to prove that  $f_0(z)= L(z) \exp(-a_0 \log(z-\sigma_0))$ is holomorphic at
$z=\sigma_0$, we  use the following lemma.
\begin{lem}\quad
$\sum_{i=0}^n (-1)^i \,w a_0^{n-i} \sha a_0^i  \equiv 0 \mod \C\langle
A_c^*\rangle \qquad \hbox{for all } w \in A_c^*\ .$
\end{lem}
\begin{proof}
Let $\widetilde{\partial}_{a_0}$ denote the truncation operator
with respect to the letter $a_0$ defined in $\S3.1$, but which
acts by truncation on the right, {\it i.e.}, $\widetilde{\partial}_{a_0} w a_{i}= \delta_{0i} w $, where $\delta_{0i}$ is the Kronecker delta. It is a derivation with
respect to $\sha$. If we apply it to the left-hand side of the equation, we
obtain zero, by the Leibniz formula. This implies that the left hand side is a linear
combination of  words not ending in $a_0$.
\end{proof}
\begin{rem}
The operators $\widetilde{\partial}_{a_i}$ are related to the
`d\'erivations \'etrang\`eres' defined by Ecalle [E].
\end{rem}

Using the fact that  $a_0^{\sha i}= i!\, a_0^i$, we have
$$f_0(z)=L(z) \exp(-a_0 \,\log(z-\sigma_0))=\sum_{w\in A^*} L_w(z)\, w \, \sum_{i\geq 0} (-1)^i a_0^i\, L_{a_0^i}(z)\ .$$
\noindent  It follows from the  previous lemma and the shuffle relations for the functions $L_w(z)$, that the coefficient of
each word $w a_0^n$, where $w\in A^*_c$ and $n\geq 0$,  is a linear combination of  $L_{w'}(z)$, where
$w'\in A^*_c$. These are holomorphic at $z=\sigma_0$ by  construction, and this
 proves the regularity condition for $f_0(z)$.

In order to prove the  uniqueness statement in the proposition, let $K(z)$ be any other
solution of $(\ref{hstar})$ which is holomorphic on $U$. The series
$L(z)$ defined above is invertible, as its leading coefficient is
the constant function $1$.  Let $F(z)= L(z)^{-1}K(z)$. On
differentiating the equation $K(z)=L(z)F(z)$, we obtain
$$\sum_{i=0}^N{a_i \over z-\sigma_i} K(z)= \sum_{i=0}^N{a_i \over z-\sigma_i} L(z) F(z)+ L(z)
F'(z)\ ,$$
 by
$(\ref{hstar})$, and therefore $L(z)F'(z)=0$. Since $L(z)$ is invertible,
$F'(z)=0$, and so $F(z)$ is  constant. This completes the proof of the proposition.

\begin{rem} The functions $L_w(z)$ are known as \emph{hyperlogarithms} and were originally defined by Poincar\'e and Lappo-Danilevsky.
 They were recently resurrected
by Aomoto [Ao1-3], Ecalle [E], and Goncharov [Go1-3].
 It is clear that $L_1(z)=1$,  and %for $n\in \N$,
\begin{eqnarray}
L_{a_i^n}(z)&=&{1\over n!} \log^n \Big({z-\sigma_i \over
\sigma_0-\sigma_i}\Big)\quad \hbox{ if } \quad
i\geq1, \nonumber \\
L_{a_0^n}(z) &= &{1\over n!} \log^n (z-\sigma_0)\ ,  \nonumber
\end{eqnarray}
for all  $n\in \N$.
Note that  $L_{a_0^n}(z)$ depends on the choice of branch of
$\log(z-\sigma_0)$ which was fixed previously, but that the
functions $L_{a_i^n}(z)$ do not. They are the unique branches
which satisfy the limiting condition  $L_{a_i^n}(\sigma_0)=0$.
\end{rem}

% prove the regularity condition at $z=\sigma_0$, let
%$\Mo_{\sigma_0}$ denote the monodromy operator obtained by
%analytic continuation around a small loop encircling the point
%$\sigma_0$. Since $\Mo_{\sigma_0}$ commutes with $d$, we deduce
%that $\Mo_{\sigma_0} L(z)= L(z) M$ for some series $M= \sum_{w\in
%X^*} M(w)\,w \in \C \ld X \rd$. For each $1\neq w\in X^*$ not
%ending in $x_0$, the formula $(\ref{convdef})$ shows that $L_w(z)$
%is holomorphic at $\sigma_0$, and so  $M(w) =0$. Since $M$
%satisfies the shuffle relations, it is completely determined by
%$M(x_0)= 2i\pi$, the monodromy of $\log(z-\sigma_0)$. It follows
%that $\Mo_{\sigma_k} L(z) = L(z) \exp(2 i \pi x_0)$, which proves
%that $f_0(z)= L(z) \exp(-x_0 \log(z-\sigma_0))$ has trivial
%monodromy at $\sigma_0$, and therefore extends to a holomorphic
%function at $\sigma_0$.

Given  a branch  of $\log (z-\sigma_k)$ on $\C\backslash
\ell(\sigma_k)$  for each $1\leq k\leq N$, we obtain by symmetry a
solution to $(\ref{hstar})$ corresponding to each  singularity.
\begin{cor} For every $0\leq k \leq N$, there exists a
unique  solution $L^{\sigma_k}(z)$ of equation $(\ref{star'})$ on
$U$ such that
$$L^{\sigma_k}(z)= f_k(z) \exp(a_k \log (z-\sigma_k))\ ,$$
where $f_k(z)$ is holomorphic on $ \C\backslash \bigcup_{i\neq k}
\ell(\sigma_i)$ and satisfies $f_k(\sigma_k)=1$.
\end{cor}
%Likewise, %by the  change of variables
%$z'-\sigma_0=(z-\sigma_0)^{-1}$,
% one also obtains  a  solution
%$L^{\infty}(z)$ to $(\ref{star})$ on $U$ which corresponds to the
%point at infinity.  If we fix any branch  of $\log z$ on $U$, then
%one can show that $L^\infty(z) \sim \exp(a_{\infty} \log z)$ where
%$a_{\infty}=\sum_{0\leq k\leq N} a_k$.

\noindent
The quotient of any two such solutions is a constant non-commutative series known as a regularised zeta series. Using
these series, one can  determine the monodromy of hyperlogarithms explicitly  ([Br2]).

\subsection{The bar construction on $\Pro^1\backslash \Sigma$} \label{sect52} In this situation, the variant of the bar construction
defined in $\S3$ is very easy to describe.
%and coincides with the
%usual reduced bar construction defined by  Chen.
 Let $k$ denote any subfield of $\C$ which contains $\sigma_0,\ldots, \sigma_N$.
 The ring of
regular functions on $\Pro^1\backslash \Sigma$ is simply
$$\Or_{\Sigma} =k\Big[ z, \Big({1\over z-{\sigma_j}}\Big)_{0\leq j\leq N}\Big]\ .$$ Since $\Pro^1\backslash \Sigma$ is of
dimension one,  the integrability condition  is trivially satisfied.
Let $A^\vee=\{\psi_0,\ldots, \psi_N\}$, where $ \psi_i = d\log(z-\sigma_i)$, for $0\leq i\leq N$.
The cohomology classes  of the forms $\psi_i$ form a $k$-basis for $H^1(\Pro^1\backslash \Sigma)$.
Clearly $\psi_i\wedge \psi_j=0$ for all $0\leq i,j\leq N$, and therefore $B(\Pro^1\backslash \Sigma)$ is a shuffle algebra
%$B(\Pro^1\backslash \Sigma) = \Or_\Sigma \otimes k\langle A^\vee\rangle$.
% In this  case,  there
%is a canonical  identification $\Z\langle A\rangle \cong \Z\langle
%A^\vee \rangle,$ which maps $a_i$ to $d\log(z-\sigma_i)$, so we can
%in fact define $B(\Pro^1\backslash\Sigma)$ to be
\begin{equation}B(\Pro^1\backslash\Sigma) = \Or_\Sigma
\otimes_k k\langle A^\vee\rangle\ ,\end{equation} equipped with the derivation
\begin{equation} \label{derivddef}
d={ d \over dz} \otimes 1 +\sum_{i=0}^N
\Big({1\over z-\sigma_i}\Big) \otimes \partial_{\psi_i}\ ,
\end{equation}
where   the
truncation operators $\partial_{\psi_i}$ were defined in $\S \ref{sect31}$.
Let $L(\Pro^1\backslash \Sigma)$ denote the $\Or_\Sigma$-algebra generated by the coefficients of a solution $L$ to $(\ref{hstar})$. The analogue of
the map $(\ref{realis})$ is the differential homomorphism:
\begin{eqnarray} \label{hyperlogrealis}
\rho: B(\Pro^1\backslash\Sigma) &\To& L(\Pro^1\backslash \Sigma)  \\
w & \mapsto & L_w(z)\ ,\nonumber
\end{eqnarray}
which is the identity on $\Or_\Sigma$.
 Theorem $\ref{thmbartriv}$
implies that this map is an isomorphism.%following statement.

\begin{cor}
The functions $L_w(z)$, for $w\in A^*$, are linearly independent
over $\Or_\Sigma$. Every function  in %$F(z) \in
$L(\Pro^1\backslash
\Sigma)$ has a primitive which is unique up to a constant.
\end{cor}
%We will refer to $B(\Pro^1\backslash \Sigma)$ as the universal%
%algebra of hyperlogarithms.

 The construction of the functions $L_w(z)$ used a decomposition
of $B(\Pro^1\backslash \Sigma)$ into convergent and non-convergent
parts. This used the fact that the map
$$ ( u\otimes v\mapsto u\sha v): \Z\langle A_c\rangle \otimes \Z\langle a_0\rangle\rightarrow \Z\langle A\rangle\ ,$$
is an isomorphism of algebras. We can therefore define
\begin{equation}\label{Bconv}
B_{\sigma_0}(\Pro^1 \backslash \Sigma) = \Or_\Sigma\langle A_c\rangle
\ ,\end{equation}
to be the sub-algebra of convergent iterated integrals (indexed by words not ending in $a_0$). It is a differential algebra for the derivation $d$ defined in $(\ref{derivddef})$. We have
$$B(\Pro^1\backslash \Sigma) \cong B_{\sigma_0}(\Pro^1\backslash \Sigma)\otimes_{k[z,1/(z-\sigma_0)]} B(\Pro^1 \backslash \{\sigma_0,\infty\})\ .$$
There is a corresponding decomposition
 $$L(\Pro^1\backslash\Sigma))=L_{\sigma_0}(\Pro^1\backslash\Sigma) \otimes_{k[z,1/(z-\sigma_0)]} L(\Pro^1 \backslash \{\sigma_0,\infty\})\ ,$$
  where
 $L(\Pro^1\backslash \{\sigma_0,\infty \} )\cong k[z,1/(z-\sigma_0), \log (z-\sigma_0)]$.
 Correspondingly, the hyperlogarithm realisation $(\ref{hyperlogrealis})$ decomposes as a product $\rho=\rho_{\sigma_0} \otimes \rho'$,
 where
 $$\rho_{\sigma_0}(w) = \int_{\sigma_0}^z w\ , \qquad \hbox{ for all } w \in A_c\ .$$
 This is a convergent iterated integral, even though the base point $\sigma_0$ does not lie in the space $\Pro^1\backslash \Sigma$.
 The logarithmic divergences are completely determined by the realisation $\rho':B(\Pro^1\backslash\{\sigma_0,\infty\})\rightarrow
 L(\Pro^1\backslash\{\sigma_0,\infty\})$, where  $\rho'(a_0)= \log (z-\sigma_0)$.
\begin{rem}
In general, the points $\sigma_0,\ldots, \sigma_N$ will not be
arranged symmetrically. In this case, one needs to do a genuine
analytic continuation of the  functions $L_w(z)$, since the formula
$(\ref{convdef})$ is not valid outside its radius of convergence.
Lappo-Danilevsky described a technique for dealing with this
situation, which is described in [Br2]. This extra complication will not arise in the present context.
%The  set of moduli spaces
%$\Mod_{0,S}$ form a  special case of fiber-type hyperplane
%arrangements for which this problem fortunately never occurs.
\end{rem}

\subsection{Quotients of the hyperlogarithm equation}
Now we shall consider the case  where the coefficients $a_i$ in $(\ref{hstar})$ satisfy relations. Therefore, let $A=\{a_0,\ldots, a_N\}$
be an alphabet with $N+1$ letters as before, and consider an ideal
$$I\subset \C\langle a_0,\ldots, a_N\rangle\ .$$
 Typically, $I$ will be generated by commutators  of the form $[a_i,a_j]$ for $i\neq j$.
It defines a closed ideal we also denote by $I$ in the completed algebra $\C\ld A\rd$.
 Let
$$\pi: \C\ld A \rd \To \C\ld A \rd /I$$
denote the  quotient map.  Consider the analogue of equation $(\ref{hstar})$:
\begin{equation} \label{pihstar} {\partial \over \partial z}F(z)=
 \sum_{i=0}^N {\pi(a_i) \over
z-\sigma_i}  F(z)\ , \end{equation}
where,  this time, $F$ takes values in the quotient ring $\C\ld A\rd/I$. An equation of this type will be called
a \emph{hyperlogarithm quotient}  equation.
\begin{cor} There exists a unique solution $F$ to the hyperlogarithm quotient equation  $(\ref{pihstar})$ with solutions in $\C\ld A\rd/I$ such that
$F(z) \sim (z-\sigma_0)^{\pi(a_0)}$ as $z\rightarrow \sigma_0$.

\end{cor}
\begin{proof} The existence follows immediately from proposition $\ref{prophstarsoln}$, on applying $\pi$ to a solution
of $(\ref{hstar}).$ The uniqueness is proved in the same way.
\end{proof}
Now let $L(I)$ denote the $\Or_{\Sigma}$-module of functions generated by the coefficients of a solution $F$ to  $(\ref{pihstar})$.
It is a differential submodule of $L(\Pro^1\backslash \Sigma)$. More precisely,
$$L(I) \cong \Or_{\Sigma} \otimes \big(\C\langle A \rangle /I \big)^\vee \subset
\Or_{\Sigma} \otimes \big(\C\langle A \rangle \big)^\vee \cong L(\Pro^1\backslash \Sigma)\ .$$
It follows that the coefficients of solutions to $(\ref{pihstar})$ are linear combinations of hyperlogarithms.
If $I$ is a Hopf ideal, {\it i.e.}, $\Gamma I \subset 1 \otimes  I + I\otimes 1$,  then  $L(I)$ is an algebra by duality. Theorem
 $\ref{thmbartriv}$ immediately implies the following corollary.
\begin{cor} Suppose that  $I$ is a Hopf ideal. In this case,
$L(I)$ is a unipotent extension of $\Or_\Sigma$. In particular, it is a differentially simple polynomial algebra over $\Or_{\Sigma}$ whose ring
of constants is $k$. %extension of $\Or_\Sigma$.
\end{cor}

 As an example,
consider the equation:
$${d F \over dz}  = \Big( {a_0\over z} + {a_1 \over z-1}\Big) F$$
on $\Pro^1\backslash \{0,1,\infty\}$, and let $I\subset \C\langle a_0,a_1\rangle$ denote the ideal generated by
$[a_0,a_1]$. Then
$F= \exp (a_0\log z +a_1 \log (z-1))$ is the unique solution satisfying $F\sim \exp (a_0 \log z)$ as $z\rightarrow 0$. The differential algebra
$L(I)$ is just  $\C[z,1/z,1/(z-1), \log z, \log (z-1)]$.
% If, instead, we consider the ideal $I$ generated by the relations
%$a_1a_0=a_1^2=0$ (which is not a Hopf ideal) then $L(I)$ is the free differential module over $\C[x,1/x,1/1-x]$ generated by the functions $\Li_n(x)$.

\subsection{Multiple polylogarithms and hyperlogarithms}\label{sect54} We recall the definition of the multiple polylogarithm functions,
which were defined by Goncharov [Go1-4]. Let $n_1,\ldots, n_r\in \N$, and consider the power series
\begin{equation}\label{multpolydefn}
\Li_{n_1,\ldots, n_r} (z_1,\ldots, z_r) = \sum_{0<k_1<\ldots <k_r} { z_1^{k_1} \ldots z_r^{k_r} \over k_1^{n_1}\ldots k_r^{n_r}}\ ,
\end{equation}
which converges absolutely for $|z_i|\leq 1$ if $n_r\geq 2$ and for $|z_i|<1$ in general.

Now let $\ell\geq 2$, $x_1,\ldots, x_{\ell-1} \in \C$, and set $\Sigma=\{\sigma_0,\ldots, \sigma_{\ell}, \infty \}$, where
%$$ \sigma_0=0\ , \qquad \sigma_1=1 \ , \qquad \sigma_2 = x_{\ell-1}^{-1} \ ,\quad  \ldots \ ,\quad  \sigma_{\ell}=(x_1\ldots x_{\ell-1})^{-1} \ .$$
$$ \sigma_0=0\ , \qquad \sigma_1=1 \ , \quad \hbox{ and } \quad  \sigma_{i}=(x_{\ell-i+1}\ldots x_{\ell-1})^{-1}\quad \hbox{ for }\quad  2\leq i\leq \ell \ .$$
Let $A=\{a_0,\ldots, a_{\ell}\}$ as previously, and let
$w= a_0^{n_r-1}a_{i_r}\ldots a_0^{n_1-1}a_{i_1}\in \C \langle A \rangle $, where $1\leq i_1,\ldots, i_r\leq \ell$.
We suppose that the points $\sigma_i$ are distinct and finite (compare $(\ref{ExpCubeCoords})$). %, that is $x_1,\ldots, x_{\ell-1}$ are distinct and not equal to $0$ or $1$.
Let us consider the points $x_1,\ldots, x_{\ell-1}$ as being fixed, and let $x_\ell \in \Pro^1\backslash \Sigma$ denote  a free variable.
By $(\ref{convdef})$, the coefficients of the corresponding hyperlogarithm function  with respect to $x_\ell$, are given near $x_\ell=0$ by the formula
\begin{eqnarray} \label{multpolysandhyp} L_w(x_\ell)&=&\!\!\! \sum_{1\leq m_1<\ldots <m_r}{(-1)^r \over m_1^{n_1} \ldots m_r^{n_r}}( x_{j_1}.\, .\, x_\ell)^{m_1}
(x_{j_2}.\,.\, x_{\ell})^{m_2-m_1}
.\,.\,(x_{j_r}.\,.\, x_{\ell})^{m_r-m_{r-1}} \nonumber\\
& = &
(-1)^r \,\Li_{n_1,\ldots, n_r} \Big( {x_{j_1}\ldots x_{\ell} \over x_{j_2}\ldots x_\ell}, \ldots, {x_{j_{r-1}}\ldots x_{\ell}
 \over x_{j_r}\ldots x_\ell}
,x_{j_{r}}\ldots x_{\ell}\Big)\ , \nonumber
\end{eqnarray}
where we have set $j_k = \ell-i_k+1$ for $1\leq k\leq r$.
It follows that such a  multiple polylogarithm,
% after substituting
%arguments $z_i$ of the form $(x_{s_i}x_{s_i+1}\ldots x_{t_i})^{\pm 1}$ for $1\leq i\leq r-1$, and  $z_r = x_{s_r}\ldots x_\ell$,
 considered as a function of the single variable $x_\ell$, is a hyperlogarithm function on $\Pro^1\backslash \Sigma$.
The relation between the multiple polylogarithm viewed as a hyperlogarithm in $x_\ell$, and the multiple polylogarithm viewed as  a function of all its variables,
 is given by the
fibration sequence between moduli spaces $\Mod_{0,n}$  ($\S\ref{sect65}$).

\subsection{Multiple zeta values and $\Pro^1\backslash\{0,1,\infty\}$} \label{sect55}
 In the case where  $\sigma_0=0$ and  $\sigma_1=1$,
$D$ is the projective line minus three points.
Since $\Pro^1\backslash\{0,1,\infty\}$ also coincides with $\Mod_{0,4}$,  it
is natural to
 make a change of sign and define $X=\{x_0,x_1\} $, where $x_0=a_0$ and $x_1=-a_1$.
  Let $\log z$ denote the  branch of the logarithm which is real for $z\in \R$, $z>0$. By  proposition $(\ref{prophstarsoln})$, the equations
\begin{eqnarray} \label{P1eqn}
{d L(z) \over d z} &= &\Big( {x_0 \over z} + {x_1 \over 1-z} \Big) L(z)  \\
L(z) &\sim& \exp( x_0 \log z) \nonumber
\end{eqnarray}
have a unique solution $L(z)\in \C\ld X\rd$, known as  the
 generating series of multiple polylogarithms in one variable.
Its
  coefficients  are written $\Li_w(z)$, for $w\in X^*$. We have
$\Li_{x_0}(z)= \log z$ and  $\Li_{x_1}(z) = -\log(1-z)$. Now consider a word $w\in x_0X^*x_1$ which
begins in $x_0$ and ends in $x_1$. It can be written
 $$w =x_0^{n_r-1} x_1 x_0^{n_{r-1}-1}x_1\ldots x_0^{n_1-1}x_1\ ,$$
 where $n_r\geq 2$. Equation $(\ref{convdef})$ therefore gives a power series expansion:
 \begin{equation}\label{multpolinonevar}\Li_w(x) = \Li_{n_1,\ldots, n_r}(1,\ldots,1, z)=\sum_{0<k_1<\ldots<k_r} {z^{k_r}
 \over k_1^{n_1}\ldots k_r^{n_r}} \ ,\end{equation}
which is regular at $z=1$.  The numbers $\Li_w(1)$ satisfy the shuffle relations by $(\ref{shufflesforhypers})$.

 \begin{defn}
Let $w \in x_0X^*x_1$ as above. %Then $w=x_0^{n_r-1} x_1 x_0^{n_{r-1}-1}x_1\ldots x_0^{n_1-1}x_1$ where $n_r\geq 2$.
 The
\emph{multiple zeta value}
 of weight $n_1+\ldots+n_r$ and depth $r$ is the real number defined by the convergent sum:
$$\zeta(w)= \zeta(n_1,\ldots, n_r) =  \Li_w(1)= \sum_{0<k_1<\ldots<k_r} {1 \over k_1^{n_1} \ldots k_r^{n_r}}\ ,\qquad n_r\geq 2 \\ .$$
The function $\zeta$ extends by linearity to the $\Q$-vector space spanned by $x_0X^*x_1$.
We define $\MZV$ to be the  $\Q$-module generated by the set of all multiple zeta values:
\begin{equation} \label{MZVdefn}
\MZV= \Q[\zeta(w): \, w\in x_0X^*x_1]\ .\end{equation}
 Because the multiple zeta values satisfy the shuffle relation $\zeta(w\sha w')=\zeta(w)\zeta(w')$, and because $\Q[w: w\in x_0X^*x_1]$
 is stable under the shuffle product, $\MZV\subset \R$ is an algebra. It is naturally filtered by the weight [Wa].
 \end{defn}
It is not difficult to verify that every word $w\in X^*$ is a linear combination of shuffles of $x_0$, $x_1$ and words $\eta \in x_0X^*x_1$.
The map $w\mapsto \zeta(w)$ extends to a unique function on $\Z\langle X\rangle $ which satisfies
$$\zeta_{\sha} (x_0)= 0\ ,  \qquad \zeta_{\sha}(x_1)= 0\ ,$$
$$ \zeta_{\sha}(w \sha w') = \zeta_{\sha}(w) \zeta_{\sha}(w')\ , \qquad \hbox{for all } w,w'\in X^*\ .$$
\begin{defn}
The \emph{Drinfeld associator} [Dr] is the non-commutative  series
$$Z^{0,1} = \sum_{w\in X^*} \zeta_{\sha}(w) w \in \MZV\ld X \rd\ .$$
\end{defn}
\noindent
  It follows that
Drinfeld's associator is precisely the regularised value of $L(z)$ at $1$:
\begin{equation} \Reg(L(z),1) = Z^{0,1}\ .\end{equation}
%Essentially, it is possible to regularise  $L(z)$ at $z=0$ and $z=1$ simultaneously because there is a decomposition
%$\Z\langle X\rangle = \Z\langle x_0 \rangle \sha \Z\langle x_0X^*x_1 \rangle \sha \Z\langle x_1 \rangle$ in this case.

%$$\Li_m(x)=\sum_{n\geq 1}
%{x^n \over n^m}$$ is  the classical $m^{\text{th}}$ polylogarithm
%function.

%______________________SUPERSOLVABLE BLURB_-------------------------
%-------------------------------BLOB blob ------------------------
%--------------TO INTRODUCTION
%\subsection{} Now consider the general case of affine hyperplane
%arrangements in any dimension ($\S3.2$) whose complement we denote
%$M$. In general, it is difficult to construct words in $B(M)$
%because they rapidly become very complicated. However, we say that
%$M$ is \emph{fiber-type} if it can be expressed as an iterated
%sequence of fibrations $F_0,\ldots, F_d$ whose fibres are of
%dimension one. It follows from theorem $3.10$ that
%$$B(M) = \bigotimes_{i=1}^n B(F_i)\ ,$$
%where $F_i$ is the projective line with a certain number of points
%removed. Since $B(F_i)$ is a free tensor algebra on a certain
%alphabet, we deduce the following corollary.%

%\begin{cor}
%If $M$ is the complement of an fiber-type hyperplane
%arrangement, then $B(M)$ is isomorphic to a tensor product of free
%tensor algebras. Futhermore, every iterated integral on $M$ is a
%product of hyperlogarithms.
%\end{cor}
%-------------------------------------------------------------------
%-------------------------------------------------------------------
%-------------------------------------------------------------------

\newpage

\section{The universal algebra of polylogarithms on $\Mod_{0,n}$}
In this section, we give an explicit construction of the algebra of
all homotopy-invariant iterated integrals on $\Mod_{0,S}$ in terms
of multiple polylogarithms.
%We show that this defines a differential
%algebra extension of $\Q(\Mod_{0,n})$ of Picard-Vessiot type which
%is universal for unipotent connections on $\Mod_{0,n}$.
 By
decomposing this algebra as a tensor product of hyperlogarithm
algebras, we compute its monodromy  in terms of  multiple zeta
values.

\subsection{The cohomology ring of $\Mod_{0,S}$.}
%We first compute the cohomology  of  $\Mod_{0,S}$.
Recall that
$\Mod_{0,S}$ was defined as the quotient of the configuration
space of $n=|S|$ distinct points $(z_s)_{s\in S}\in(\Pro^1)^S$, modulo
the action of $\PSL_2$. Let $i,j,k,l\in S$. The cross-ratio $[i\, j|k\, l]:(\Pro^1)_*^{S}\rightarrow \Pro^1$
defines a function on $(\Pro^1)_*^{S}$, and since $[i\, j|k\, l]=1- [i\, k|j\, l]$, we have:
\begin{equation} \label{crosswedge}
d \log [i\, j|k\,l]\wedge d \log [i\,k|j\,l]=0\ .\end{equation} We
introduce the notation
\begin{equation}\label{Deltaij}
\Delta_{ij} = d \log (z_i -z_j)={dz_i -dz_j\over z_i-z_j}\ , \quad
\hbox{ for all } \quad 1\leq i<j\leq n\ .
\end{equation}
 where $\Delta_{ij}=\Delta_{ji}$, and
$\Delta_{ii}=0$, for all $1\leq i,j\leq n$. Equation $(\ref{crosswedge})$ gives a quadratic relation between the $\Delta_{ij}$,
which can be simplified as follows.
 Since $\PSL_2(\C)$ acts transitively
on the projective line $\Pro^1(\C)$, and since the cross-ratio is invariant under its action, we can  place the point $z_1$ at infinity, and
it follows that
\begin{equation} \label{confquot}
\Mod_{0,S}(\C) = \C_*^{n-1}/B\ ,\end{equation} where $\C^{n-1}_*$ denotes
the set of distinct $n-1$-tuples $z_2,\ldots, z_{n}\in \C$, and
$B\cong \C^\times  \ltimes   \C$ is the subgroup of $\PSL_2(\C)$ which
stabilizes  $\infty$. The projection map $\C_*^{n-1}\rightarrow
\Mod_{0,S}(\C)$ is a trivial fibration with fibres isomorphic to $B$,
and it follows that
\begin{equation} \label{cohomprod}
H^\star(\C_*^{n-1}) \cong H^\star (\Mod_{0,S}(\C)) \otimes H^\star (B)\ .
\end{equation}
We can therefore deduce the cohomology of $\Mod_{0,S}(\C)$ from the
structure of $H^\star(\C_*^{n-1})$, which can be described as follows.
We apply $(\ref{crosswedge})$ with $l=1$. Using the fact that $z_1=\infty$, we deduce that $d\log [i\,j|k\,l]=\Delta_{ik}-\Delta_{jk}$,
and $d\log [i\,k|j\,l]= \Delta_{ij}-\Delta_{kj}$,
viewed as  1-forms on $\C_*^{n-1}$. Then $(\ref{crosswedge})$ yields Arnold's relation:
\begin{eqnarray} \label{Arnold}
\Delta_{ij} \wedge \Delta_{jk} + \Delta_{jk}\wedge \Delta_{ki} +
\Delta_{ki} \wedge \Delta_{ij} = 0\ ,
\end{eqnarray}
for any distinct indices $2\leq i,j,k\leq n$. % By symmetry, it holds for all $1\leq i,j,k\leq n$.

\begin{thm}(Arnold [Ar]). \,\,
$H^\star(\C_*^{n-1})$ is the  quotient of the free exterior algebra generated
by $\Delta_{ij}$ for $2\leq i,j\leq n$, by the quadratic relations $(\ref{Arnold})$.
\end{thm}

Now let us fix a dihedral structure $\delta$ on $S$. In $\S2$ we defined  1-forms
$$\omega_{ij} = d\log u_{ij}\ ,   \quad \hbox{ for } \qquad \{i,j\} \in \chi_{S,\delta}\ .$$
Their cohomology classes $[\omega_{ij}]$ form a basis for $H^1(\Mod_{0,S}(\C))$.
%
%Recall from \S2 that% $H^1(\Modf_{0,n})=1$,
%\begin{equation} \label{H1}
%H^1(\Mod_{0,S}) = \bigoplus_{\{ij\}\in S_n} \Z \, [d\log u_{ij}]
%%= \bigoplus_{\{ij\}\in S_n} \Z \,{du_{ij} \over u_{ij}} \
%.\end{equation}
 Recall the definition of $N$  ($\S\ref{sect32}$) as the kernel of the
exterior product:
$$N = \ker \, \big(\wedge: H^1(\Mod_{0,S}(\C)) \otimes H^1(\Mod_{0,S}(\C))
\To H^2(\Mod_{0,S}(\C))\big)\ .$$
% Note that although the forms $\omega_{ij}$ are linearly independent over $\Q$, the same is not true for
%$du_{ij}$.

\begin{prop} \label{Ndescriptiondihed} $N$ is spanned by the following elements:
\begin{equation} \label{Modcohomrel}
\Big( \sum_{\{i,j\}\in A} [\omega_{ij}] \Big) \otimes \Big(
\sum_{\{k,l\}\in B} [\omega_{kl}]\Big) \ ,
%\Big( \sum_{\{i,j\}\in A} [\omega_{ij}] \Big) \wedge \Big( \sum_{\{k,l\}\in B} [\omega_{kl}]\Big) =0\ ,
\end{equation}
 where $A,B\subset \chi_{S,\delta}$ are any two sets of chords which cross completely $(\S\ref{sect22})$. The cohomology ring $H^\star(\Mod_{0,S}(\C))$ is isomorphic
to the free exterior algebra generated by $[\omega_{ij}]$, for
$\{i,j\} \in \chi_{S,\delta}$, modulo the image of elements of the
form $(\ref{Modcohomrel})$.
\end{prop}
\begin{proof}
First we regard each dihedral coordinate $u_{ij}$ as a function on
$(\Pro^1)_*^{n}$. By the defining equations  ($\ref{id})$,  we
have $1-\prod_{a\in A} u_a = \prod_{b\in B}u_b$ for all sets of
chords $A, B \subset \chi_{S,\delta}$ which cross completely. This
implies that $d\log \prod_A u_a\wedge d\log \prod_B u_b=0$, which
is precisely $$\Big( \sum_{\{i,j\}\in A} [\omega_{ij}] \Big)
\wedge \Big( \sum_{\{k,l\}\in B} [\omega_{kl}]\Big) =0\ .$$
 Furthermore, every instance of
$(\ref{crosswedge})$ occurs in this way, since each cross ratio
$[i\,j|k\,l]$ can be written as  a product $\prod_{A} u_a$ or its
inverse,  by lemma $\ref{lemmacrossratio}$. We now place
$z_1=\infty$ as above, and view the corresponding relations on
$\C^{n-1}_*$. By Arnold's theorem, this implies that
$(\ref{Modcohomrel})$ generates the set of all  relations on
$\C^{n-1}_*$. In particular, $(\ref{Modcohomrel})$ generates $N$,
since by $(\ref{cohomprod})$, $H^\star(\Mod_{0,S}(\C))\subset
H^\star(\C^{n-1}_*)$ and so any relation satisfied by the
$\omega_{ij}$ is also satisfied in $H^\star(\C^{n-1}_*)$.

Now, since $B\cong \C^\times \ltimes \C$ is homotopy equivalent to
a circle, $H^*(B)$ is the exterior algebra generated by a single
cohomology class which we denote  $\beta\in H^1(B)$. It follows
from $(\ref{cohomprod})$
 that $H^*(\Mod_{0,S}(\C))$ is the subalgebra of $H^*(\C^{n-1})$ of  degree $0$ in $\beta$.
 We deduce  from Arnold's theorem that $H^*(\Mod_{0,S}(\C))$ is the quotient of the free exterior algebra generated by
 a basis of $H^1(\Mod_{0,S}(\C))$, modulo $N$.
%
%
%
%
% By Arnold's theorem, $H^\star(\C_*^{n-1})$
%is the exterior algebra generated by the forms $\Delta_{ij}$, so
%on taking the quotient modulo $\beta$, it follows from
%$(\ref{cohomprod})$ that $H^\star(\Mod_{0,S})$ is the exterior
%algebra generated by the forms $d\log u_{ij}$, $\{i,j\} \in
%\chi_n$. It suffices to compute the relations between the $d\log
%u_{ij}$ implied by $(\ref{Arnold})$.
% By definition, we have
%\begin{equation} \label{duexp}
%d\log (u_{ij}) = \Delta_{i\, j+1} + \Delta_{i+1\, j} -
%\Delta_{i\,j} -\Delta_{i+1\, j+1}\ , \end{equation} for all
%$\{i,j\} \in \chi_n$, such that $i,j\notin \{n,1\}$. Otherwise, if
%$i \notin \{2,n\}$, and $j\notin \{n-1, 1\}$,
%\begin{eqnarray} \label{duexpsmall}
%d\log (u_{i\,1})& = & \Delta_{2\,i} - \Delta_{2\, i+1}\ ,
%\\
%d\log (u_{n\,j})& = & \Delta_{n\,j+1} - \Delta_{n\, j}\ .
%\nonumber
%\end{eqnarray}
%These equations imply that $H^1(\Mod_{0,S})$ is the augmentation
%ideal in $H^1(\C_*^{n-1})$, and we can take $\beta= \sum_{i,j}
%\Delta_{ij}$.
\end{proof}

%The relations  $(\ref{Modcohomrel})$ are not independent, as a
%simple  counting argument shows.  This relation is more
%conveniently written $(d^\times u_{ik} + d^\times u_{km})\wedge
%d^\times u_{jl}=0$, and is equivalent to the existence of the
%dilogarithm function, as we shall see below.

Similar results have been obtained by Getzler [Ge].

\begin{rem}
The quadratic relations $(\ref{Arnold})$ are  equivalent to the existence of the dilogarithm function, in the following sense.
Let $f=[i\, j\, |k\, l]$. Then
identity $(\ref{Arnold})$ is precisely the integrability of the element
$$[ d\log\,f \big| d\log(1-f)] \in W^2 \, B(\Mod_{0,S})\ .$$
 The iterated integral ($\S\ref{sect33}$) corresponding to this element is the function $\Li_2(f)$.
\end{rem}

\subsection{The universal algebra of polylogarithms on $\Mod_{0,S}$}\label{sect62} Recall that in simplicial coordinates
$(\ref{ExpSimpCoords})$, the space $\Mod_{0,S}$ is the open
complement of an affine hyperplane arrangement. Its ring of
regular functions is
$$\Or(\Mod_{0,S})=\Q[u_{ij}, u_{ij}^{-1}]\cong \Q\Big[\big({t_i}\big)_{1\leq i\leq \ell}, \Big({1\over t_i}\Big)_{1\leq i\leq \ell}, {1\over (1-t_i)}_{1\leq i\leq \ell}, {1\over (t_i-t_j)}_{1\leq i<j\leq \ell}\Big]\ ,$$
which  is a differential algebra with respect to  the partial
differential operators $\partial/\partial t_i$.
 We defined the abstract algebra of homotopy-invariant iterated integrals on $\Mod_{0,S}$ using
the reduced bar construction  in  $\S\ref{sect32}$.
\begin{defn}
 The \emph{universal algebra of polylogarithms} on
$\Mod_{0,S}$ is the differential  graded %Hopf
 algebra
$ B(\Mod_{0,S}) = B\big(\Or(\Mod_{0,S})\big).$
\end{defn}
Recall that $B(\Mod_{0,S})$ is the unipotent closure of
$\Or(\Mod_{0,S})$, and that its de Rham cohomology is trivial,
{\it i.e.}, $H_{\DR}^0(B(\Mod_{0,S}))\cong \Q,$ and
$H_{\DR}^i(B(\Mod_{0,S}))=0$ for all $i\geq 1$. The structure of
$B(\Mod_{0,S})$ is particularly rich: it has  a natural Hopf
algebra structure over $\Or(\Mod_{0,S})$, and also carries an
action of the symmetric group $\Sym(S)$ by functoriality. The
graded pieces of the set of indecomposable elements  in
$B(\Mod_{0,S})$ of fixed weight yield very interesting
finite-dimensional representations of $\Sym(S)$. Correspondingly,
there is an action by  the subgroup of dihedral  symmetries
$D_{2n}$ of the $n$-gon $(S,\delta)$. This action is evident from
the symmetric description of $H^1(\Mod_{0,S})$ and $N$ in terms of
the forms $\omega_{ij}$, for $\{i,j\}\in \chi_{S,\delta}$, given
in proposition $\ref{Ndescriptiondihed}$.

Now if we pass to  cubical coordinates, %(and, ultimately, any set of vertex coordinates),
we can split $B(\Mod_{0,S})$ as a tensor product of shuffle
algebras, and subsequently decompose it into convergent and
non-convergent pieces. First, recall that we defined a base point
at infinity $(\ref{cubicalbasepointatinfinity})$, corresponding to
the origin $x_1=\ldots =x_\ell=0$, which is locally a normal
crossing divisor. The base point is given by the map
$\Or(\Mod_{0,S}) \rightarrow \kes$ which maps $x_i$ to
$\epsilon_i$, for $1\leq i\leq \ell$.  The projection map
$(x_1,\ldots, x_\ell) \mapsto (x_1,\ldots, x_{\ell-1})$ defines a
linear fibration  $\Mod_{0,\{s_1, \ldots, s_n\}}\rightarrow
\Mod_{0,\{s_2, \ldots, s_{n}\}}$, which forgets the  point marked
$s_1$. In the notations of $\S\ref{sect25}$, it corresponds to the
choice of sets $T_1=\{s_n,s_1,s_2,s_3\}$, $T_2=\{s_2,s_3,\ldots,
s_{n-1}\}$, and $T_1\cap T_2=\{s_2,s_3\}$. By iterating in this
manner, we obtain
 a sequence of fibrations:
$\Mod_{0,\{s_i,\ldots, s_n\}} \To \Mod_{0,\{s_{i+1},\ldots,
s_n\}}, $ obtained by forgetting the marked point $s_i$, for
$i=1,n,n-1,\ldots, 4$. By applying theorem $\ref{thmbarfibration}$
to these fibrations, we deduce that there is a canonical
isomorphism:
\begin{equation}\label{Btdecompforcubicals}
 B(\Mod_{0,S})
\cong   \bigotimes_{i=1}^{n-3} B_{\Mod_{0,\Sigma_i}} (\Pro^1\backslash \Sigma_i)\ ,\end{equation}
where  $\Sigma_i=\{s_{2},s_3\ldots, s_{n-i+1}\}.$ Each algebra $B(\Pro^1\backslash \Sigma_i)$ is a universal algebra of hyperlogarithms,
and is a free shuffle algebra on $n-i-1$ generators by  $\S\ref{sect52}$.
\begin{cor}
$B(\Mod_{0,S})$ is isomorphic, as a $\Or(\Mod_{0,S})$-algebra, to the tensor product  of the free shuffle algebras on $2,3,\ldots, n-2$ generators.
\end{cor}
\noindent
Using results of Radford, one can write down a basis of any free shuffle algebra in terms of Lyndon words (see [Rad]).
The corollary implies that a basis of $B(\Mod_{0,S})$ is given by  tensor products of Lyndon words.

We can now decompose each component of
$(\ref{Btdecompforcubicals})$ into convergent and non-convergent
parts. We can define the subalgebra
$B'_{\Mod_{0,\Sigma_i}}(\Pro^1\backslash \Sigma_i)\subset
B_{\Mod_{0,\Sigma_i}}(\Pro^1\backslash \Sigma_i) $ of convergent
words  in a
 similar manner to $(\ref{Bconv})$.
%
%  $\S5$ we defined the set
%of convergent iterated integrals $B_{s_{3}}
%(\Pro^1\backslash\Sigma_i)$, regularised at $s_{3}$, where
%$s_{3}\in \Sigma_i$.
 We can then define
$$B_0(\Mod_{0,S}) = \Or(\Mod_{0,S}) \otimes \bigotimes_{i=1}^{n-3} B'_{\Mod_{0,\Sigma_i}}(\Pro^1\backslash \Sigma_i)\ .$$
Then $B(\Mod_{0,S})$ decomposes as a commutative  tensor product
 $$B(\Mod_{0,S}) \cong B_0(\Mod_{0,S}) \otimes_\Q \Q \big[[d\log x_1],\ldots, [d\log x_\ell]\big] \ .$$
  The algebra on the right
 is the free commutative
(polynomial) algebra on generators $[\omega_{2\,i+3}]=[d\log x_i]$ for $i=1,\ldots,\ell$.
\begin{lem} The subalgebra $B_0(\Mod_{0,S})\subset B(\Mod_{0,S})$ is generated as  a vector space by the set of integrable words,
no element of which ends in a symbol $\omega_{2k}$, for $4\leq k\leq n$.
\end{lem}
\begin{proof}Let $A$ denote the $\Or(\Mod_{0,S})$-subalgebra of $B(\Mod_{0,S})$ generated by  the set of all integrable words
$$\sum_I c_I [\omega_{i_1 j_1}|\ldots | \omega_{i_r\, j_r}]\ ,\qquad c_I \in \Q\ ,$$
where $\{i_r,j_r\}\notin \{\{2,4\}, \ldots, \{2,n\}\}$. It is
clear that $A\subset B_0(\Mod_{0,S})$ is a differential
subalgebra. Furthermore, one easily checks that every element
$a\in \Omega^1(\Mod_{0,S})\otimes_{\Or_{\Mod_{0,S}}} A$ of weight
at least $ 1$ has a primitive in $A$. This follows from the proof
of theorem $\ref{thmbartriv}$ or the argument given in the
appendix, since taking primitives involves adding symbols to the
left of each word. Using the techniques of $\S3$, proposition
$\ref{propPVproperties}$, it follows immediately that the map
$A\rightarrow B_0(\Mod_{0,S})$ is surjective.
\end{proof}

Likewise, for every vertex $v\in V^{\delta}$, the set of vertex
coordinates at $v$ defines a base point at infinity, and (by
considering the action of the differential Galois group of $\Ues$
over $\kes$, for example), one defines
   a subalgebra of convergent words $B_{v,\delta}(\Mod_{0,S})$ such that
$$B(\Mod_{0,S}) \cong B_{v,\delta}(\Mod_{0,S}) \otimes_\Q \Q[\omega_{i_1\, j_1}] \otimes_\Q \ldots \otimes_\Q \Q[\omega_{i_\ell\, j_\ell}]\ ,$$
where $\{i_1,j_1\},\ldots, \{i_\ell, j_\ell\}\in F_v$ are the
chords occurring in the triangulation corresponding to $v$. As
above, $B_{v,\delta}(\Mod_{0,S})$ corresponds to the set of all
integrable words which do not terminate in any  symbol
$\omega_{i_k\, j_k}$. The case $B_0(\Mod_{0,S})$ corresponds to
the vertex whose  triangulation is $\{\{2,4\},\ldots, \{2,n\}\}$.
This is just the point $x_1=\ldots =x_\ell=0$ in cubical
coordinates.

%\begin{rem} \label{remnonprojbarsplit}
%More generally, using the projection maps $(\ref{compactgenfib})$, which are not always fibrations,   one can prove using properties of unipotent
%closures that
%$$B(\Mod_{0,S})\subset B(\Pro^1\backslash \Sigma) \otimes_{\Q(\Mod_{0,S})} B(\Mod_{0,T_1})\otimes_{\Q(\Mod_{0,S})} B(\Mod_{0,T_2})\ ,$$
%where $\Sigma$, $T_1$, and $T_2$ are defined in $\S \ref{sect25}$.
%For example, in the case $\Mod_{0,6}$, this expresses
%$B(\Mod_{0,6})$ as a sub-algebra of the tensor product of free shuffle algebras with $2,2$, and $5$ generators.
%Its image can in fact be characterized explicitly (\S\ref{sect64}).
% This will not be
%required in the sequel, but an equivalent result is proved by a different method in $\S \ref{sect64}$.

%\end{rem}

\subsection{The dihedral connection on $\Mod_{0,S}$} \label{sect63} There is a canonical differential equation on $\Mod_{0,S}$
whose solutions can be expressed in terms of multiple
polylogarithms. Let
  $\Z\langle\delta_{ij}\rangle$ denote  the free
non-commutative Hopf algebra generated by  the symbols
  $\delta_{ij}=\delta_{ji}$, for $\{i,j\} \in \chi_{S,\delta}$, where $\delta_{ij}$ is primitive (see \S$\ref{sect31})$. It is convenient to set
$\delta_{ii}= \delta_{i\, i+1}=0$ for all indices $i\in \Z/n\Z$.
   Consider
  the following formal 1-form on $\Mod_{0,S}$:
\begin{equation} \label{KZ}
\Omega_{S,\delta}= \sum_{\{i,j\}\in \chi_{S,\delta}} { \delta_{ij}}{d u_{ij}
\over u_{ij}} \ .\end{equation}
% Note that although the forms $d\log
%u_{ij}$ are linearly independent, the same is not true for $d
%u_{ij}$.
 The form $\Omega_{S,\delta}$ is integrable if and only if
  $d\,\Omega_{S,\delta}=
\Omega_{S,\delta}\wedge \Omega_{S,\delta}\ .$
Since $\Omega_{S,\delta}$ is closed, this reduces to $\Omega_{S,\delta}\wedge \Omega_{S,\delta}=0$.
We define the
\emph{dihedral infinitesimal braid relations} to be the identities:
\begin{equation} \label{braid}
[\delta_{i-1\, j} + \delta_{i\, j-1} -\delta_{i-1 \, j-1} -\delta_{ij},
\delta_{k-1\, l} + \delta_{k\, l-1} -\delta_{k-1 \, l-1} -\delta_{kl}]=0\ ,
\end{equation}
for all $i,j,k,l\in S$.

 For
each $1\leq i,j\leq n$, consider a set of formal symbols $t_{ij}$,
where $t_{ii}=0$ and $t_{ij}=t_{ji}$. The Knizhnik-Zamolodchikov (KZ)
form on $(\Pro^1_*)^{n}$ is the 1-form:
 \begin{equation}  \Omega_{KZ_n} = \sum_{1\leq i<j\leq n}
t_{ij}\, \Delta_{ij}\ ,
\end{equation}
where $\Delta_{ij}=\Delta_{ji}$ is given by $(\ref{Deltaij})$. Let
us assume that
\begin{equation} \label{thomeq} \sum_{k=1}^n t_{kl } = 0\qquad \hbox{for all }  1\leq l\leq n\ .
\end{equation}
This variant of the KZ-equation has been considered by Ihara, amongst others. It corresponds to the usual KZ-equation
on $\C_*^{n-1}$ , except that
 it has an extra set
of symbols at infinity, and one extra relation which kills the center of the braid algebra. One can prove that
$\Omega_{KZ_n}$ is integrable if and only if the following single relation holds:
\begin{equation} \label{finalbraid}
[t_{ij}, t_{kl}]= 0 \qquad \hbox{ for all } i,j,k,l \hbox{ distinct}\ .\end{equation}
One verifies by computing  $[t_{ij}, \sum_{l=1}^n t_{kl}]=0$, that this is equivalent to the usual
infinitesimal braid relations:
\begin{eqnarray}\label{infbraid}
\mathrm{[}t_{ij}, t_{kl} \mathrm{]} &= & 0 \ , \\
 \mathrm{[}t_{ij}, t_{ik}+t_{jk}\mathrm{]}&=& 0\ , \nonumber
\end{eqnarray}
which hold  for all distinct indices $2\leq i,j,k,l\leq n$.
%for all distinct indices $2\leq i,j,k,l\leq n$.
Now, by $(\ref{udef})$, we have $\omega_{ij} = d\log u_{ij} =
\Delta_{i\, j+1} + \Delta_{i+1\, j} - \Delta_{i+1\, j+1} -\Delta_{ij}$, for $\{i,j\}\in \chi_{S,\delta}$. If we write
$$\Omega_{KZ_n} = \Omega_{S,\delta}\ ,$$
then this is equivalent to the identities
\begin{equation}
\label{braidtexp} t_{ij}  =  \delta_{i\,j-1} + \delta_{i-1\, j} -
\delta_{i-1\, j-1} -\delta_{ij}\ ,
\end{equation}
 for all $1\leq i,j\leq n$, as is easily verified. Since $\delta_{ii}=\delta_{i\, i+1}=0$ for
 $1\leq i\leq n$, then $(\ref{braidtexp})$ implies that
\begin{eqnarray} \label{braidtone}
t_{i j} &=& \delta_{i-1\, j} \qquad \qquad\qquad\qquad\qquad\qquad\hbox{ if }\quad j=i+1\ , \\
\label{bradittwo}
t_{i j} & = & \delta_{i-1\,j} - \delta_{i-1\, j-1} -
\delta_{i  j} \qquad \qquad \quad\hbox{ if }\quad j=i+2\ .
\end{eqnarray}
The following lemma implies that
the set of equations $(\ref{braidtexp})$ are invertible over $\Z$.
\begin{lem}
For all $1\leq i<j \leq n$,
\begin{equation} \label{deltaandt}
\delta_{ij} = \sum_{i<a<b\leq j} t_{ab} \ .
\end{equation}
\end{lem}
\begin{proof}
By equation $(\ref{braidtone})$, $\delta_{i-1\, i+1} = t_{i\,
i+1}$ for $1<i<n$. Substituting into $(\ref{bradittwo})$ gives
$\delta_{i-1\, i+2} = t_{i\, i+2} +t_{i\, i+1}+ t_{i+1\, i+2}$.
Let $m\geq 4$, and suppose by induction that $(\ref{deltaandt})$
holds for all $0< j-i\leq m-1 $. Then for $j-i=m-1$,
$(\ref{braidtexp})$ gives
\begin{eqnarray}\delta_{i-1\, j} &=& t_{ij} + \delta_{ij} + \delta_{i-1\, j-1} - \delta_{i\, j-1} \ ,\nonumber \\
&=& t_{ij} + \sum_{i<a<b\leq j} t_{ab}
+ \sum_{i-1<a<b\leq j-1} t_{ab}- \sum_{i<a<b\leq j-1} t_{ab} \ ,\nonumber \\
&=&t_{ij} + \sum_{i<a<b\leq j} t_{ab} +\sum_{i<b\leq j-1} t_{ib}=  \sum_{i-1<a<b\leq j} t_{ab}\ . \nonumber
\end{eqnarray}
This proves $(\ref{deltaandt})$ when $j-i=m$. The result follows by induction.
\end{proof}

Now if we substitute the expressions $(\ref{braidtexp})$ for
$t_{ij}$ and $t_{kl}$ in terms of $\delta_{ab}$ in equation
$(\ref{finalbraid})$, then we obtain $(\ref{braid})$. This  proves
the following result.
%
%
%
%$$\Omega_{KZ_n} = \Omega_{S,\delta} + \sum_{2\leq i,j\leq n}
%t_{ij} \beta$$ where $\beta=\sum_{ij} \Delta_{ij}$.
%
%
\begin{prop} The form $\Omega_{S,\delta}$ is  integrable if and
only if the dihedral braid relations $(\ref{braid})$ hold.
\end{prop}

\begin{lem} The dihedral braid relations imply that
\begin{equation} \label{braidx}
\mathrm{[} \delta_{ij}, \delta_{kl} \mathrm{]}=0\ ,
\end{equation}
for all chords $\{i,j\}$, $\{k,l\}\in \chi_{S,\delta}$ which do not cross.
\end{lem}

\begin{proof}
 Without loss of generality, we can assume that $1\leq i<j<k<l\leq n$.  Then,
by identity $(\ref{finalbraid})$,
$$[\delta_{ij}, \delta_{kl}] = [ \sum_{i<a<b\leq j} t_{ab}, \sum_{k<c<d\leq l} t_{cd} ] =0\ ,$$
since all sets of four indices $\{a,b,c,d\} $ occurring in the summation are distinct.
\end{proof}

\begin{example}
In the case $S=\{1,2,3,4,5\}$, relation $(\ref{braid})$ with $i=2$, $j=4$, $k=3$, $l=5$ implies that
$[\delta_{14}-\delta_{13}-\delta_{24}, \delta_{25}-\delta_{35}-\delta_{24}, ]=0$. By $(\ref{braidx})$ this gives the following
five-term relation:
\begin{equation}
[\delta_{13}, \delta_{24}]+ [\delta_{24}, \delta_{35}]+ [\delta_{35}, \delta_{41}]+ [\delta_{41}, \delta_{52}] + [\delta_{52},\delta_{13}]=0\ .
\end{equation}
 This is dual to the functional equation of the dilogarithm.
 A similar relation in fact holds for any five consecutive indices on $\Mod_{0,S}$.
\end{example}

% One can verify
%that the group $B$ acts on $\Omega_{KZ_n}$  by scalar
%multiplication, so it follows that
%\\

% First observe that since $\Omega_{S,\delta}$ is a sum of exact forms, $d\Omega_{S,\delta}=0$, and so the integrability condition is
%equivalent to $\Omega_{S,\delta}\wedge \Omega_{S,\delta}=0$. Using
%this, one can prove the relations directly using lemma $4.1$.
%Instead, we prove the relations on $\Mod_{0,5}$ using the known
%braid relations.
\begin{defn}
Let $R$ denote a commutative unitary ring, and let
 $I$ denote the ideal in $R\langle\, \delta_{ij}:\{i,j\}\in
\chi_{S,\delta} \rangle$ generated by the dihedral relations
$(\ref{braid})$  above. The \emph{dihedral braid algebra} over $R$ is  the free
non-commutative $R$-algebra
\begin{equation}
\DB_{S,\delta}(R)= R\langle\, \delta_{ij}:\{i,j\}\in \chi_{S,\delta}
\rangle/I  \ . \end{equation}
\end{defn}
This is a co-commutative graded Hopf algebra over $R$ ($\S
\ref{sect31})$, where $\deg \delta_{ij}=1$. The product is the
concatenation product, and the coproduct $\Gamma$ is the unique
coproduct with respect to which the generators $\delta_{ij}$ are
primitive ($I$ is a Hopf ideal because it is generated by
commutators of primitive elements). It is the universal enveloping
algebra of the free Lie algebra generated by the
symbols $\delta_{ij}$, subject to relation $(\ref{braid})$. As in $\S3.1$, its completion %with respect to the $\ker \varepsilon$
 is
 the $R$-Hopf algebra
\begin{equation} \label{BSdeltahatdefn}
\widehat{\DB}_{S,\delta}(R)= R\ld \, \delta_{ij}:\{i,j\}\in \chi_{S,\delta}\ ,
\rd/\widehat{I}
\end{equation}
where $\widehat{I}$ is the closed ideal generated by $I$. %whose coproduct we also denote by $\Gamma$.
%If $R$ is any commutative ring, we denote $\DB_{S,\delta}(R)= \DB_{S,\delta}\otimes R$
%and $\widehat{\DB}_{S,\delta} (R) =\widehat{\DB}_{S,\delta}\otimes R$ for the extension of scalars to $R$.
It follows from the previous calculations that $\DB_{S,\delta}(R)$
is just the free non-commutative $R$-algebra generated by the
symbols $t_{ij}$, for $1\leq i,j\leq n$, which satisfy
$(\ref{thomeq})$, modulo the relations $(\ref{finalbraid})$.
 This is isomorphic to  the ordinary infinitesimal braid
algebra modulo its center. The difference here is that we have fixed a set of generators for
this algebra which depend on the dihedral structure $\delta$.
 %We write $e=1$ for the empty word
%or unit element in $U_n$, and define the weight of a word $w$ to
%be its length, {\it i.e.} the number of symbols occuring in $w$.

%--------------COMMENTS ON THE WEIGHT. I'M NOT REALLY SURE THAT IT MAKES SENSE,
%--------------SO I HAVE AXED IT --------------------------------------------
%----------------------------------------------------------------------------
%Since the dihedral braid relations are homogeneous of weight two,
%one can verify that any element $u\in \DB_{S,\delta}$ has a
%representation as a linear combination of words of minimal length.
%We define the weight of $u\in \DB_{S,\delta}$ to be the weight of
%the word of maximal length occuring in this representation. This
%defines a grading on
%$\DB_{S,\delta}$.
%--------------------------------------------------------------------------------
%---------------------------------------------------------------------------------

% whose filtration we write $W^\cdot DB_S$.% be denoted $W^\cdot U_n$, with graded pieces

%$\gr_W^\cdot $. %Depth filtration.......
%We write $\DB_{S,\delta}^{>0}$ to denote the set of all words of
%weight $\geq 1$, which is just the kernel of the augmentation map
%$\varepsilon$. We have
%$$W^n \DB_{S,\delta} \cong \DB_{S,\delta} /
%\big(\DB_{S,\delta}^{>0}\big)^{n-1} $$

 Let $\widehat{\Mod}_{0,S}$ be a universal covering space for
$\Mod_{0,S}$, and let $p: \widehat{\Mod}_{0,S}
 \rightarrow \Mod_{0,S}$ denote the projection map.
A multi-valued function on $\Mod_{0,S}$ is defined to be a
holomorphic function on $\widehat{\Mod}_{0,S}$.
 Since the integrability conditions are satisfied in
$\DB_{S,\delta}(\C)$ we can  consider the following formal
differential equation on
$\widehat{\Mod}_{0,S}$:
\begin{equation}\label{star}
d L = \Omega_{S,\delta}\, L\ . \end{equation} A solution  $L$
takes values in $\widehat{\DB}_{S,\delta}(\C)$. Its coefficients are
multi-valued functions on $\Mod_{0,S}$.
%We will write $L(\uij)$
%to indicate the dependence on the dihedral coordinates $\uij=\{u_{ij}\}$,
%bearing in mind that they are not independent.
 We can fix a
solution to $(\ref{star})$ by specifying its value at a point of
$\Mod_{0,S}$, or its limiting value at an intersection of boundary divisors. It suffices to define solutions
at  intersections of divisors of maximal codimension. Therefore, we define
$$V^\delta= \{ \alpha\in \chi^{\ell}_{S,\delta}\}$$
to be the set of all triangulations of the $n$-gon. By  $\S
\ref{sect24}$, each such triangulation determines a unique vertex
of the associahedron $\overline{X}_{S,\delta}$. For each vertex
$v\in V^{\delta}$, let $F_v=\{\{i,j\}\in \chi_{S,\delta}:
u_{ij}(v)=0\}$ denote the set of faces of the associahedron
$\overline{X}_{S,\delta}$ which meet at $v$. Let $\log(u_{ij})$,
for  $\{i,j\}\in \chi_{S,\delta}$, denote the principal branch of
the logarithm on $u_{ij}>0$ (see \S\ref{sect42}).

\begin{thm} \label{thmKZhomog} Let $v\in V^{\delta}$.
 There exists
a unique solution $L_{v,\delta}$ to $(\ref{star})$ such that in a
neighbourhood of $v$,
$$L_{v,\delta} (\uij) = f_{v,\delta}(\uij) \prod_{\{i,j\}\in F_v} \exp( \delta_{ij} \log
(u_{ij}))\ ,$$ where $f_{v,\delta}(\uij)\in \widehat{\DB}_{S,\delta}(\C)
$ extends to a  holomorphic function in the neighbourhood of $v\in
\Modf_{0,S}$,  and takes the value $1$ at  $v$. The function $f_{v,\delta}(\uij)$ extends holomorphically
to an open neighbourhood of the interior of every face $F$ meeting $v$.
%Here,
%$\log(u_{ij})$ is the branch of the logarithm which is real for $u_{ij}>0$.
\end{thm}
\begin{rem} The product  $\prod_{\{i,j\}\in F_v} \exp(\delta_{ij} \log (u_{ij}))$
is well-defined, because by $(\ref{braidx})$, the symbols $\delta_{ij}$ and
$\delta_{kl}$ commute whenever $\{i,j\}$ and $\{k,l\}$ do not
cross, and no two chords $\{i,j\}, \{k,l\} \in F_v$ can cross because
$F_v$ is a triangulation of the $n$-gon $(S,\delta)$.
%the corresponding divisors $D_{ij}$ and $D_{kl}$ have non-empty
%intersection.
\end{rem}

\begin{proof}
Let  $\DB_{S,\delta}^{>0}(\C)\subset \DB_{S,\delta}(\C)$ denote the
kernel of the counit $\varepsilon: \DB_{S,\delta}(\C) \rightarrow \C$. For each integer $N\geq 1$,
define
$$W_N= \DB_{S,\delta}(\C)/ \big(\DB_{S,\delta}^{>0}(\C)\big)^{N+1}\ .$$
%This is isomorphic to $W^N \DB_{S,\delta}$ as a vector space.
If
we write $\delta_{ij}$ for the map which acts by left
multiplication by the symbol $\delta_{ij}$, for each $\{i,j\} \in
\chi_{S,\delta}$, then each $\delta_{ij}$ is a nilpotent operator
on the space $W_N$.

In  $\S\ref{sect41}$ we showed that $\overline{X}_{S,\delta}$ is
a manifold with corners by constructing a specific atlas
$\{U_e(\varepsilon)\}$. We will show that $\Omega_{S,\delta}$ defines a unipotent equation of Fuchs' type on each chart (definition \ref{defnFuchsUnipType}), and apply the
results of $\S\ref{sect43}$.
Therefore, let $\alpha \in \chi^{k}_{S,\delta}$ denote a partial decomposition of the $n$-gon, where $1\leq k\leq \ell$. To $\alpha$
corresponds the face $F_\alpha$ of  $\overline{X}_{S,\delta}$. Choose any complete triangulation  $\alpha'\in \chi^{\ell}_{S,\delta}$
which contains $\alpha$.
By proposition $(\ref{prop211})$, the vertex coordinates
$$\{x_i^{\alpha'}: 1\leq i\leq \ell\}=\{u_{ij}: \,\, \{i,j\} \in \alpha'\}$$
form a system of normal  coordinates in a neighbourhood of $F_\alpha$. We can therefore write
$$\Omega_{S,\delta} = \sum_{\{i,j\} \in \alpha'} \delta_{ij} {d u_{ij} \over u_{ij} } + A_{ij} d u_{ij}\ ,$$
where $A_{ij}$ are holomorphic functions in a neighbourhood of
$F_{\alpha}$. Since the operators $\delta_{ij}$ are nilpotent on
$W_N$, and since the open neighbourhoods  of every face $F_\alpha$
(including $F_{\emptyset}=X_{S,\delta}$) cover
$\overline{X}_{S,\delta}$, it follows that $(\ref{star})$ is
unipotent of Fuchs' type, as required.

By theorem $\ref{thmfuchs}$, we can
 find a local solution $L_{v,\delta}^{(N)}$ to $(\ref{star})$ with values in $W_N(\C)$,  which satisfies the
asymptotic condition stated above. By corollary $\ref{corfuchs}$, this
solution extends globally over the whole Stasheff polytope
$\overline{X}_{S,\delta}$. The theorem follows on taking the limit
as $N$ tends to infinity, since $\widehat{\DB}_{S,\delta}(\C)=\displaystyle{\lim_{\leftarrow }}\, \DB_{S,\delta}(\C)/ (\DB^{>0}_{S,\delta}(\C))^N$.
\end{proof}

\noindent
Any such solution $L_{v,\delta}$ to $(\ref{star})$ extends by analytic continuation to give a multi-valued function on the whole of $\Mod_{0,S}$.
By construction,  the theorem defines  a unique  real-valued branch on the interior
of the associahedron $X_{S,\delta}$. It is convenient to write the asymptotic boundary condition
$$L_{v,\delta} \sim
\prod_{\{i,j\}\in F_v} u_{ij}^{\delta_{ij}}\quad\hbox{ near } v\
.$$

\begin{figure}[h!]
  \begin{center}
%    \leavevmode
    \epsfxsize=5.0cm \epsfbox{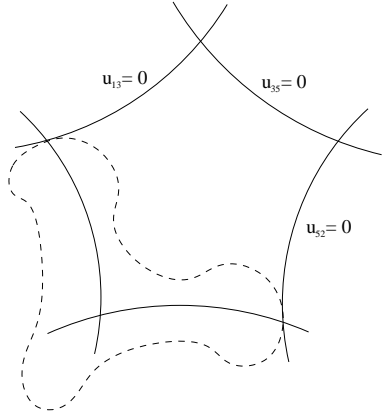}
  \label{Butterfly}
  \caption{A local picture of $\overline{\Mod}_{0,5}(\R)$ and a cell $X_{5,\delta_i}$. The dashed region depicts
  the set of real points of a  domain of holomorphy for the regularised function $f_{v,\delta}$ given in theorem $\ref{thmKZhomog}$.}
  \end{center}
\end{figure}

\begin{rem}
The formal equation $(\ref{star})$ is a homogeneous version of
the Knizhnik-Zamolodchikov equation on $\C_*^{n-1}$. Drinfeld studied solutions to the KZ equation on $\C_*^3, \C_*^4$  with prescribed asymptotics
in certain zones [Dr], which were subsequently generalised by Kapranov [Ka]. Such a zone is determined by a permutation on $n-1$ letters,
plus a bracketing on the set with $n-1$ letters.
Combinatorially, a permutation on $n-1$ letters corresponds to a cyclic structure on $n$ letters, and
 a bracketing corresponds to a triangulation of an $n$-gon, {\it i.e.},
a vertex in the associahedron of dimension $n-3$ (fig. 11).
Kapranov   interpreted   each zone as
the region near a corner at infinity in a certain two-fold
cover $\widetilde{S}$ of the compactified real moduli space:
$$\widetilde{S} \rightarrow \overline{\Mod}_{0,S}(\R)\ .$$
In the previous theorem, we have constructed a canonical solution $L_{v,\delta}$ near each corner on $\overline{\Mod}_{0,S}(\R)$.
%This follows from lemma $\ref{lemmatess}$, which implies that the set of corners are in bijection with the set of pairs $(v,\delta)$.
Therefore, each solution $L_{v,\delta}$ corresponds to exactly two  Drinfeld-Kapranov zones. One way to see this is that
the manifold $\widetilde{S}$ is obtained by gluing together  a set of  associahedra which are parametrized by the set of all cyclic, rather than dihedral
structures. The number of cyclic structures on $S$, where $|S|=n$ is exactly $(n-1)!$. As an example, let us consider the simplest
case $\Mod_{0,4}$. Then equation $(\ref{star})$ reads:
\begin{eqnarray}
{dL \over d x} = \Big( {\delta_{13} \over x} + {\delta_{24} \over x-1} \Big) L\ ,
\end{eqnarray}
Consider the solution $L_{0^+}(x)$ which satisfies $L_{0^+}(x)
\sim \exp (\delta_{13} \log x) $ as $x\rightarrow 0^+$.
Remembering that $x=u_{13}=[12|43]$, and after placing the point
$z_4$ at $\infty$, we obtain $x=(z_2-z_3)/(z_1-z_3)$. If $0<x<1$,
then either $z_1<z_2<z_3$ or $z_3<z_2<z_1$. These are the two
cyclic structures on $\{z_1,z_2,z_3,\infty\}$ which map to the
dihedral structure $\delta$ corresponding to the cell
$[0,1]=\overline{X}_{S,\delta}$. The zone $0<x \ll 1$ therefore
corresponds to the pair of zones $z_1(z_2z_3)$ and $(z_3z_2)z_1$
in Kapranov's notations.

\end{rem}

\subsection{Functoriality with respect to projection maps} \label{sect64}
It is well-known that  solutions to the KZ-equation decompose as
products of hyperlogarithm equations by considering fibration maps between configuration spaces.
% and differentiating with
%respect to a single variable at a time
%([GL]).
 One obtains a more general decomposition result
for the homogeneous equation $(\ref{star})$, by using  the
projection maps
 considered in $\S\ref{sect25}$. First, we fix a dihedral structure $\delta$
on $S$ and choose a chord $\{i,j\}\in \chi_{S,\delta}$. Recall
from \S\ref{sect25}, that if we set $T_1=\{s_{j+1},\ldots,
s_{i},s_{i+1}\}$, and $T_2=\{s_i,s_{i+1},\ldots, s_j\}$, and
denote the  induced dihedral structures by $\delta_1$, $\delta_2$,
then $T_1\cap T_2=\{s_i, s_{i+1}\}$, and  there is a projection
map:
$$
f_{T_1}\times f_{T_2}:  \Modf_{0,S} \To
\Mod^{\delta_1}_{0,T_1}\times \Mod^{\delta_2}_{0,T_2} \ .$$ By
lemma $\ref{lemmadecomp}$, this map has a section whose image is
the divisor $D_{ij}=\{u_{ij}=0\}$:
$$i: \Mod^{\delta_1}_{0,T_1}\times
\Mod^{\delta_2}_{0,T_2} \overset{\sim}{\To} D_{ij} \subset \Modf_{0,S}\ .$$
 In order to  fix solutions of $(\ref{star})$, let $v\in
V^{\delta}$ be a vertex of the polytope $\overline{X}_{S,\delta}$ such that
 $v\in D_{ij}$, {\it i.e.,} $u_{ij}(v)=0$.
 %Since
%$\phi$ extends to the partial compactifications
%$(\ref{compactgenfib})$
By projecting down, we obtain  vertices in
$\overline{X}_{T_k,\delta_k}$:
$$v_k= f_{T_k} (v) \in \overline{X}_{T_k,\delta_k} \qquad \hbox{ for } k=1,2\ .$$
If $v$ is given by a triangulation $\alpha \in
\chi^\ell_{S,\delta}$ of the $n$-gon $S$, then $v_1$, $v_2$ are
given by the restrictions $\alpha_1, \alpha_2$ of this
triangulation to $T_1$ and $T_2$ respectively (compare fig. 4). As
sets of chords, we have $\alpha = \alpha_1 \sqcup \alpha_2 \sqcup
\{i,j\}.$
% Since we are only considering this particular kind of projection,
 We need an extra technical condition that the only
chord in $\alpha$ emanating from the vertex $(j)$ is the chord
$\{i,j\}$ (the dihedral coordinate corresponding to such a chord
would not be preserved since $f_{T_1}$ and $f_{T_2}$ contract one
of the edges $j$ or $j+1$). Let
$$i_{T_1} : \Mod^{\delta_1}_{0,T_1}\To \Mod^{\delta_1}_{0,T_1}\times
\Mod^{\delta_2}_{0,T_2} \To D_{ij}$$
denote the map which sends $x$ to  $i(x,v_2)\in D_{ij}$.  Define a map $i_{T_2}$ similarly, and set
\begin{equation}
\pi_k = i_{T_k}\circ f_{T_k}\qquad \hbox{for } k=1,2\ .
\end{equation}
Then $\pi_k$ is a projection map
$$\pi_k: \Modf_{0,S} \rightarrow D_{\alpha_k\sqcup \{i,j\}}\qquad \hbox{ for } k=1,2\ ,$$ because $f_{T_k} \circ i_{T_k}$ is the identity.
We can write these maps explicitly in vertex coordinates $x_1^{\alpha},\ldots, x_\ell^{\alpha}$, which
form a local system of coordinates on $\Modf_{0,S}(\C)$
(see $\S\ref{sect23}$).  We can choose an  ordering on
$\alpha$ such that $D_{\alpha_2} =\{x^{\alpha}_1=\ldots =x^{\alpha}_{m-1}=0\}$, $u_{ij}= x^{\alpha}_m$, and  $D_{\alpha_1}=\{x^{\alpha}_{m+1}=\ldots=x^{\alpha}_\ell=0\}$ for
some $m$. In that case, we have
$$f_{T_1}\times f_{T_2}: (x^{\alpha}_1,\ldots, x^{\alpha}_\ell) \To \big((x^{\alpha}_1,\ldots, x^{\alpha}_{m-1}),
(x^{\alpha}_{m+1},\ldots, x^{\alpha}_\ell)\big)\ ,$$
and
\begin{eqnarray}
\pi_1:(x^{\alpha}_1,\ldots, x^{\alpha}_\ell)&\mapsto& (x^{\alpha}_1,\ldots, x^{\alpha}_{m-1},0,\ldots, 0) \nonumber \\
\pi_2:(x^{\alpha}_1,\ldots, x^{\alpha}_\ell)&\mapsto& (0,\ldots,0 ,x^{\alpha}_{m+1},\ldots, x^{\alpha}_{\ell}) \nonumber \ .
\end{eqnarray}
%In the example given in $\S\ref{sect25}$, we had $i=2$, $m=j-3$. If we work in cubical coordinates, and suppose that the vertex
%$v$ is given by the point $\{x_1=\ldots=x_\ell=0\}$ corresponding to the triangulation , then  $u_{ij}=x_m$ and
We shall  write $\partial/\partial u_{ij}$ to denote partial
differentiation with respect to the vertex coordinate
$x_m^\alpha=u_{ij}$ %which is well-defined
in a neighbourhood of
$v\in \Modf_{0,S}(\C)$.

Let $L_{v,\delta}$ denote the unique solution to $(\ref{star})$ on $\Mod_{0,S}$ given by theorem $\ref{thmKZhomog}$. We define
$$L_k = L \circ \pi_k\ ,\qquad \hbox{for } k=1,2\ .$$
Then the functions $L_k:\widehat{\Mod}_{0,S} \rightarrow
\widehat{\DB}_{S,\delta}(\C)$ satisfy the differential equation
\begin{equation} \label{Ltildeeqn}
dL_k = \Omega_k L_k\ , \qquad \hbox{for } k=1,2\ ,\end{equation}
where $\Omega_k$ %viewed as a 1-form on $D_{\alpha_k} \subset \Modf_{0,S}$
 is given by:
$$\Omega_k= (\pi_k)_* \Omega_{S,\delta}= \sum_{\{a,b\} \in \chi_{T_k,\delta_k} } {\delta_{ab} } {d u_{ab} \over u_{ab}}\ ,\qquad \hbox{ for } k=1,2\ .$$
By construction, the solutions $L_k$ satisfy the asymptotic condition
\begin{equation} \label{Lkasympt}
 L_k = f_k\big( \sum_{\{a,b\} \in \alpha_k} \delta_{ab} \log u_{ab}\big) \ ,\qquad \hbox{for } k=1,2\ ,\end{equation}
where $f_k$ is holomorphic in a neighbourhood of $v_k$ on $D_{\alpha_k}$, where it takes the value $1$.
It is clear that $\Omega_k$ and $L_k$ only involve the symbols $\delta_{ab}$ where $\{a,b\}\in \chi_{T_k,\delta_k}$.
 Since no chord
in $\chi_{T_1,\delta_1}$ crosses any  chord in
$\chi_{T_2,\delta_2}$, it follows from $(\ref{braidx})$ that
$L_1$ and $L_2$  commute. Likewise, we have
$[L_1, \Omega_2]=[L_2, \Omega_1]=[\Omega_1,\Omega_2]=0$.
 The three
series $L_{v,\delta}, L_1, L_2$ are all   formal
power series in $\widehat{B}_{S,\delta}(\C)$ whose coefficients
are multi-valued functions on  $\Mod_{0,S}$. They are related as
follows.

% Proof hugely simplified since we alre

\begin{prop} \label{propKZsplit} Let  $\{i,j\}\in \chi_{S,\delta}$ be any chord, and let $v\in V^{\delta}$ such that
$u_{ij}(v)=0$. With the notations above, there is
a decomposition
\begin{equation} \label{Lspliteqn}
L_{v,\delta}= h  \,\,L_1 L_2\ ,
\end{equation}
where $h$ is the unique solution in $\widehat{B}_{S,\delta}(\C)$   to the hyperlogarithm  quotient equation
\begin{equation} \label{hsubhyplogeqn}
{\partial h \over \partial u_{ij}} = \Big( %{\delta_{ij} \over u_{ij}} +
 \sum_{\{k,l\} \in \chi_{S,\delta}}   %\x \{i-1 ,j\} }
 { \delta_{kl}}\,{\partial\,  \log \,u_{kl} \over \partial
u_{ij}} \Big)h\ ,
\end{equation}
 which satisfies the boundary condition
\begin{equation} \label{hsubhyplogbound}
h = g \exp( \delta_{ij} \log u_{ij})\ ,\end{equation}
where $g$ is a holomorphic function of $u_{ij}$  in a
neighbourhood of $0$, and  $g|_{u_{ij}=0}$ is the constant
function $1$.
\end{prop}

\begin{proof}
Define a formal power series $h\in \widehat{B}_{S,\delta}(\C)$ by the equation $h = L_{v,\delta} \,
L_2^{-1} \, L_1^{-1} .$  If we
differentiate this equation,  we deduce that
$$\Omega_{S,\delta} L_{v,\delta}= dh \,\, L_1\, L_2 + h
\,\Omega_1\,L_1 \,L_2 +  h
\,\Omega_2\,L_1 L_2 \ , $$
where we have used the fact that $[\Omega_2,
L_1]=0$. It follows that
\begin{equation} \label{heqn}
dh= \Omega_{S,\delta} \,h  - h \, \big(\Omega_1+
\Omega_2\big)\ .
\end{equation}
By definition, the functions $L_k$
do not depend on the variable $u_{ij}$, {\it i.e.}, $\partial L_k
/\partial u_{ij}=0$, for $k=1,2$.
 Therefore $\partial h/ \partial u_{ij} =\partial L_{v,\delta} /\partial u_{ij} \,L_2^{-1} \,L_1^{-1}
  $, and $(\ref{star})$ implies that
$$ {\partial h\over \partial u_{ij}}=
\Big(\sum_{\{k,l\}\in \chi_{S,\delta} }% \x \{i-1,j\} }
 { \delta_{kl} \over u_{kl}} \,{\partial u_{kl} \over \partial
u_{ij}} \Big) \, h\ .$$
%This is just the coefficient of $d u_{ij}$ in
%$\Omega_{S,\delta}$.
%
% In order to prove the asymptotic boundary
%condition, let $F_v,F_1,F_2$ denote the sets of chords   in
%$\chi_{S,\delta}$, $\chi_{T_1,\delta_1}$, $\chi_{T_2,\delta_2}$
%which meet the vertices $v,v_1,v_2$ respectively. Each defines a triangulation, and  it follows from
%the decomposition $(\ref{chidecomp})$ that
%$F_v = F_1 \cup F_2 \cup \{i,j\}$.
%By theorem $\ref{thmKZhomog}$,
%$$ \widetilde{L}_k = f_k \exp \big( \sum_{\{a,b\} \in F_k} \delta_{ab} \log
%u_{ab}\big) \qquad \hbox{ for } k=1,2\ ,$$ where $f_k$ is
%holomorphic in the neighbourhood of $v_k$, for $k=1,2$.
By definition of the solution $L_{v,\delta}$ and equations $(\ref{Lkasympt})$, we have
$$h = f_{v,\delta} \exp \big(\!\!\! \sum_{\{a,b\} \in \alpha} \!\!\!\delta_{ab} \log
u_{ab}\big) \,  \exp \big( -\!\!\!\sum_{\{a,b\}\in \alpha_2}\!\!\!
\delta_{ab} \log u_{ab}\big) \, f_2^{-1} \exp \big(
-\!\!\!\sum_{\{a,b\} \in \alpha_1}\! \!\!\delta_{ab} \log u_{ab}\big)
f^{-1}_1\ ,  $$ where $f_{v,\delta}$ is holomorphic in a neighbourhood of $v$, and $\alpha_1,\alpha_2,\alpha$ are the
triangulations of $T_1,T_2,S$ corresponding to $v_1,v_2,v$ respectively.
Since $L_1$ and $L_2$ commute,
$$ h = f_{v,\delta} \exp(\delta_{ij} \log u_{ij}) \,f_2^{-1}\, f^{-1}_1 = f_{v,\delta}\, f^{-1}_2\,
f^{-1}_1 \exp(\delta_{ij} \log u_{ij})\ .$$ Let $g =  f_{v,\delta} \, f^{-1}_2\,
f^{-1}_1=f_{v,\delta} \, f^{-1}_1\, f^{-1}_2$, which is holomorphic in the neighbourhood
of $v$. In order to complete the proof, it suffices to show that
the function  $g|_{u_{ij}=0}$ is the constant function $1$. This, along with the differential equation for $h$,
will determine $h$ uniquely. Let
$G$ denote the restriction of $g$ to the divisor
$D_{ij}=\{u_{ij}=0\}$. We already know by construction that $G(v)=1$. Since $g$ is holomorphic in the
neighbourhood of $u_{ij}=0$,  $G$ satisfies a differential equation which is obtained by projecting  $(\ref{heqn})$ onto $u_{ij}=0$,
which amounts to pulling back $\Omega_{S,\delta}$ by $(\pi_1\times \pi_2)_*$.
%By $(\ref{id})$, $u_{ij}=0$ implies that $u_{kl}=1$ for all chords $\{k,l\}$ which cross $\{i,j\}$.
% Therefore in the expression
%$$\Omega_{S,\delta} = \delta_{ij} \, d \log u_{ij} + \sum_{\{k,l\}\x
%\{i,j\} }\delta_{kl}\, d \log u_{kl} + \widetilde{\Omega}_1 +
%\widetilde{\Omega}_2\ , $$
%all terms are zero along $u_{ij}=0$ except $\widetilde{\Omega}_1 + \widetilde{\Omega}_2$.
By definition,
$$\Omega_{S,\delta}\Big|_{D_{ij}} = \Omega_1 + \Omega_2\ .$$
Equation $(\ref{heqn})$ therefore restricts to give the following differential equation for  $G$:
\begin{equation} \label{Gequation}
d G = [\Omega_1 + \Omega_2
,G]\ ,\end{equation} where $G(v)=1$.  This equation only has constant solutions. To see this,  consider the conjugate
$H=(L_1 L_2)^{-1} G\, L_1 L_2$. Substituting into $(\ref{Gequation})$ gives  $dH=0$.
Therefore $H$ is the constant function $1$, and so the same is true of  $G$, which
completes the proof.
\end{proof}

%The differential equation for $h$ can be written very simply in
%cubical coordinates. We can assume without loss of generality that
%$\{i,j\}=\{2,j\}$ and that the fibration map is given by the
%projection $x_m =0$, where $m=j+3$, as in the proof of lemma
%$2.16$.

%\subsection{}
%In order to compute the holonomy of the moduli spaces
%$\Mod_{0,S}$, we need to
% break the dihedral symmetry and
%express the functions on $\Mod_{0,S}$ as products of
%hyperlogarithms.
% This amounts to regarding the iterated integrals on $\Mod_{0,S}$
% as functions of one variable at a time.
%
%\begin{prop} There is an isomorphism of graded  algebras
%$$ h_z:  B(\Mod_{0,S\backslash\{z\}})\otimes B(\Pro^1\backslash \{s_1,\ldots, s_{n-1}\})\overset{\sim}{\rightarrow} B(\Mod_{0,S})\ .$$
%\end{prop}
One can  verify from the definitions that the map $f_{T_k}$ induces a map
$$(f_{T_k})_*: \, \DB_{S,\delta}(\C) \To \DB_{T_k,\delta_k}(\C)\ , \qquad \hbox{for } k=1,2\ ,$$
which sends $\delta_{ab}$ to zero for all chords $\{a,b\}$ which are not in $\chi_{T_k,\delta_k},$ and is the identity on $\delta_{ab}$ for
all chords $\{a,b\} \in \chi_{T_k,\delta_k}$.
This  also follows immediately from the fact that $\Omega_k= (\pi_k)_* \Omega_{S,\delta}$ is
integrable. We can then consider
$$(f_{T_k})_* L_k : \Mod_{0,T_k} \To \widehat{\DB}_{T_k,\delta_k}(\C) \, \qquad \hbox{for } k=1,2\ .$$
By $(\ref{Lkasympt})$ and the uniqueness part of theorem $\ref{thmKZhomog}$, we conclude that
$$L_{v_k,\delta_k} =(f_{T_k})_* L_k  \qquad \hbox{ for } k=1,2\ . $$
In conclusion, a solution $L_{v,\delta}$ to $(\ref{star})$ on
$\Mod_{0,S}$ is equivalent to a pair of solutions
$L_{v_k,\delta_k}$ on $\Mod_{0,T_k}$ for $k=1,2$ plus a solution
to the hyperlogarithm quotient equation $(\ref{hsubhyplogeqn})$.
Note that    many of  the terms in the differential equation
$(\ref{hsubhyplogeqn})$ vanish.
\begin{example} Consider the case $\Mod_{0,6}$, and let $\{i,j\} = \{2,5\}$ (see fig. 3).
Then $T_1=\{2,3,4,5\}$ and $T_2=\{6,1,2,3\}$. We shall work in cubical coordinates, and write
$(x_1,x_2,x_3)=(x,y,z)$. Then $u_{24}=x$, $u_{25}=y$, and $u_{26}=z$. The map
$\Mod_{0,6}\rightarrow D_{25}\cong \Mod_{0,T_1}\times \Mod_{0,T_2}$ is given by projecting onto the divisor $y=0$.
Therefore $L_1(x,y,z)= L(x,0,0)$, and $L_2(x,y,z)=L(0,0,z)$ and we have:
\begin{eqnarray}
d L_1 &= &\Big(  \delta_{24}{dx \over x} + \delta_{35} {dx\over x-1}\Big) L_1\ , \nonumber\\
d L_2 &= &\Big( \delta_{26} {dz\over z} + \delta_{15} {dz\over z-1}\Big) L_2 \ .\nonumber
%\widetilde{L}_1 & \sim & \exp (\delta_{26} \log x)
\end{eqnarray}
Thus $L_1, L_2$ are generating series of multiple polylogarithms in one variable (but note that there is a difference
in sign in [Br1]).
On the other hand, $h$ is the unique solution to
\begin{eqnarray} \label{hinsixcase}
{\partial h\over \partial y}&=&  \Big({\delta_{25}\over
y}+{t_{56}\over y-1} + {t_{46} \over y-x^{-1}} + { t_{15} \over
y-z^{-1}} +
{ t_{14} \over y-(xz)^{-1}}\Big)\,h\ , \\
h& \sim& \exp(\delta_{25} \log y) \qquad \hbox{ as}\quad  y \rightarrow 0\ ,\nonumber
\end{eqnarray}
where, according to ($\ref{braidtexp}$),  $t_{56}=\delta_{46}$,
$t_{15}=\delta_{14}-\delta_{15}-\delta_{46}$,
$t_{46}=\delta_{36}-\delta_{46}-\delta_{35}$, and $t_{14}=
\delta_{13}+\delta_{46}-\delta_{14}-\delta_{36}$. By
$(\ref{infbraid})$, there are commutation relations
$$[t_{14}, t_{56}]=0, \quad \hbox{ and }\quad [t_{15}, t_{46}]=0\ .$$
Therefore $(\ref{hinsixcase})$ is  a hyperlogarithm  quotient
equation on the punctured affine line $\Pro^1\backslash
\{0,1,\infty, x^{-1}, z^{-1}, (xz)^{-1}\}$. Compare remark
$\ref{remmakefib}$. % and  $\ref{remnonprojbarsplit}$.
\end{example}

By applying the proposition repeatedly, we obtain an
 explicit decomposition of the generating series $L_{v,\delta}$ as products of
hyperlogarithms. Let us apply the proposition in the case where
$\{i,j\}=\{2,n\}$. In cubical coordinates, $u_{2n}=x_\ell$, and
one can check that $h$ is the unique function satisfying
\begin{eqnarray} \label{hhypcubeeqn} {\partial h \over \partial x_\ell }& =& \Big(
{\delta_{2n} \over
x_\ell} + {\delta_{1\, n-1} \over x_\ell-1} + %{ \delta_{n-2\, 1} -
%\delta_{n-2\, n} -\delta_{n-1\, 1} \over x_\ell-x_{\ell-1}^{-1}}+
\sum_{i=1}^{\ell-1} { \delta_{i+3 \, n} + \delta_{i+2\, 1} -
\delta_{i+2\, n} -\delta_{i+3\, 1} \over x_\ell- (x_i\ldots
x_{\ell-1})^{-1}} \Big)\, h\ , \\
h & \sim & \exp(\delta_{2n} \log x_\ell) \qquad \hbox{ as } \,\,x_\ell
\rightarrow 0\ . \nonumber
\end{eqnarray} where, as usual, $\delta_{ii}= \delta_{i-1 i}=0$ by
convention. The function $\log x$ is the unique branch satisfying $\log 1=0$. %(the computation is carried out in lemma $7.8$).
 This
defines a multi-valued function on $\Pro^1\backslash \Sigma$ where
$\Sigma=\{\sigma_0,\ldots, \sigma_\ell\}$, with
$$\sigma_0=0\ , \quad \sigma_1=1\ , \quad \sigma_2=x_{\ell-1}^{-1}\ ,\quad \ldots \ ,\quad  \sigma_{\ell-1} = (x_1\ldots
x_{\ell-1})^{-1}\ .$$
 The
series  $h \exp(- \delta_{2n} \log x_\ell) $ is holomorphic in the
neighbourhood of $x_\ell=0$.  In this case, the projection map $f_{T_1}\times f_{T_2}$ is a fibration. It follows  that
$h$ is a hyperlogarithm equation ({\it i.e.}, there are no relations between the coefficients in
  $(\ref{hhypcubeeqn})$). Notice also that the coefficient $\delta_{i+3 \, n} + \delta_{i+2\, 1} -
\delta_{i+2\, n} -\delta_{i+3\, 1}$ is just $t_{1\, i+3}$.

 By substituting the above values of $\sigma_i$ into the formula
$(\ref{convdef})$, we deduce that  the coefficients of the formal power series $h$ are the
multiple polylogarithms (\S\ref{sect54}):
\begin{equation}\label{multiplepoly}
\Li_{n_1,\ldots, n_r} \Big( {x_{j_1}\ldots x_{\ell} \over x_{j_2}\ldots x_\ell}, \ldots,
{x_{j_{r-1}}\ldots x_{\ell} \over x_{j_r}\ldots x_\ell}
,x_{j_{r}}\ldots x_{\ell}\Big)\ ,
\end{equation}
 where $1\leq j_1,\ldots, j_r\leq \ell$ are any
indices. By applying the proposition inductively, we obtain an
explicit decomposition of $L_{v,\delta}$ in terms of
hyperlogarithm generating series $h$.

\begin{cor} \label{corproductsofmultiplepolylogs}
$L_{v,\delta}$ is a product of hyperlogarithm generating series. Its coefficients are sums of products of multiple polylogarithms
 of the form $(\ref{multiplepoly})$ with the functions $\log x_1,\ldots, \log x_\ell$.
\end{cor}
% A similar result holds more generally for any set of vertex
%coordinates $x_1^{\alpha}$,\ldots, $x_\ell^\alpha$.
\subsection{Regularised zeta series and monodromy } \label{sect65}%Regularization, Zeta series, monodromy, holonomy
The monodromy of the KZ equation was first computed by Drinfeld
[Dr].  We shall follow the argument given in [H-P-V] for
$\Mod_{0,4}\cong \Pro^1\backslash\{0,1,\infty\}$ (see also [GL]).
Consider two vertices $u,v\in V^{\delta}$. By theorem
$\ref{thmKZhomog}$, each vertex defines a generating series of
multi-valued functions $L_{u,\delta}$, $L_{v,\delta}$ on
$\Mod_{0,S}$. The ratio of any two solutions to $(\ref{star})$ is
a constant series.

\begin{defn}
The \emph{regularised zeta series} corresponding to  $u,v\in
V^{\delta}$ is
$$Z^{u,v}=   (L_{u,\delta}(x))^{-1} \,L_{v,\delta}(x) \in
\widehat{\DB}_{S,\delta}(\C)\ ,$$ for any $x\in X_{S,\delta}$,
{\it i.e.}, $x=(u_{ij})$, where $0<u_{ij}<1$.
\end{defn}
\noindent
Since  $L_{v,\delta}$ is real-valued on $X_{S,\delta}$,  $Z^{u,v}\in
\widehat{\DB}_{S,\delta}(\R)$  has real coefficients.  Clearly
\begin{equation} \label{productforzetaseries}
Z^{u,v}Z^{v,w}=Z^{u,w}
\end{equation}
 for all $u,v,w\in V^\delta$, and in
particular, $Z^{u,v}Z^{v,u}=1$. The zeta series describe the
limiting behaviour of a solution of $(\ref{star})$ near the
boundary of $\overline{X}_{S,\delta}$.

\begin{lem}\label{Zseriesaslimit}  For all $u,v\in V^{\delta}$,
$$Z^{u,v} =\displaystyle{ \lim_{x\rightarrow u}}
\prod_{\{i,j\}\in F_{u}} \exp(-\delta_{ij} \log u_{ij})  \,  L_{v,\delta}(x)
 \ ,
$$
where $x=(u_{ij}) \in X_{S,\delta}.$
\end{lem}

\begin{proof} Let $x=(u_{ij})$ and let $x\rightarrow u$
along a path in $\overline{X}_{S,\delta}$.
   By theorem $\ref{thmKZhomog}$,
\begin{equation} \label{asympexp}
L_{v,\delta}(x) =  L_{u,\delta}(x)  Z^{u,v} =  f_{u,\delta} \Big(\prod_{\{i,j\}\in F_{u}}
\exp(\delta_{ij} \log u_{ij} )\Big)\,  Z^{u,v} \ ,
\end{equation}
which implies that
$$ Z^{u,v} = \lim_{x\rightarrow u}\Big(\prod_{\{i,j\}\in F_{u}}
\exp(-\delta_{ij} \log u_{ij} )\Big)\, f_{u,\delta}^{-1} \,L_{v,\delta}(x)\ .$$
But $f_{u,\delta}^{-1}$ is a non-commutative series which is holomorphic in a neighbourhood of $u$ where it takes the value $1$.
We can write $f_{u,\delta}^{-1}= 1 + g(x)$, where $g$ is holomorphic and vanishes at $x=u$.
%coefficient of any non-empty word therefore vanishes at $x=u$.
Since $z\log^n z\rightarrow 0$ as $z\rightarrow 0$ for all $n\in \N$,
 and since $L_{v,\delta}(x)$ has at most logarithmic singularities at $x=u$, we deduce that only the constant term $1$ in $f_{u,\delta}^{-1}$
  gives  a non-zero contribution in the limit, which proves the result.
  \end{proof}

\begin{thm} \label{thmcoeffsofzaremzvs}
 The coefficients of the series $Z^{u,v}$ are multiple zeta values:
 $$Z^{u,v} \in \MZV \ld\, \delta_{ij}: \, \{i,j\}\in \chi_{S,\delta}\rd  /I
 \quad \hbox{ for all }u,v\in V^\delta \ ,$$
 where $I$ denotes the closed ideal generated by the dihedral braid relations $(\ref{braid})$.
\end{thm}
\begin{proof}
By  the relations $(\ref{productforzetaseries})$, it suffices to compute the coefficients of $Z^{u,v}$, where $u,v$ are adjacent
corners of $\overline{X}_{S,\delta}$.
In other words, $u,v$ are given by triangulations $\alpha,\beta\in \chi^\ell_{S,\delta}$ which differ by one chord only.
Let us write $\{a,a'\} = \alpha\backslash \alpha\cap \beta$ and $\{b,b'\}= \beta\backslash \alpha\cap \beta$.
Since $D_{\alpha\cap \beta}$ is of dimension 1, there is an isomorphism
%There is a subset of four elements $T\subset S$ with the induced  dihedral structure $\delta'$ such that the natural map
$$i_{uv}: \Mod^{\delta'}_{0,T} \overset{\sim}{\To} D_{\alpha\cap \beta} \subset \Modf_{0,S}\ ,$$
where $|T|=4$, which
maps the cell $\overline{X}_{T,\delta'}$ onto the 1-dimensional face of $\overline{X}_{S,\delta}$ which connects $u$ and $v$.
%More precisely, if  $\alpha\cap\beta=\{\{i_1,j_1\},\ldots, \{i_{\ell-1},j_{\ell-1}\}\}$, then the image of $i_{uv}$ is identified with the  1-dimensional subspace $D_{\alpha\cap\beta}$ by
% setting $u_{i_1\,j_1}=\ldots = u_{i_{\ell-1}\, j_{\ell-1}}=0$.
%$T$ is the  set of vertices of the quadrilateral (inscribed in the regular $n$-gon whose vertices are labelled by $S$) which contains the two remaining chords  in $\alpha\cup \beta$.
%
Consider the solution $L_{u,\delta}$ to $(\ref{star})$ given by theorem $\ref{thmKZhomog}$.  We can identify
$\Mod_{0,T}$ with $\Pro^1\backslash \{0,1,\infty\}$ and $X_{T,\delta'}$ with the interval $(0,1)$ in such a way that $i_{uv}(0)=u$ and $i_{uv}(1)=v$.
 The pull-back  $F=(i_{uv})_* L_{u,\delta}$ then satisfies
the differential equation:
\begin{eqnarray}
{dF\over dx} &= &\Big( {\delta_u \over x} + {\delta_v \over x-1}\Big) F\ ,\\
F &\sim &\exp(\delta_u \log x) \qquad \hbox{ as } x\rightarrow 0\ , \nonumber
\end{eqnarray}
on $\Pro^1\backslash\{0,1,\infty\}$.  This follows from proposition $\ref{propKZsplit}$. Here, $\delta_u, \delta_v$ are the dihedral symbols corresponding to the chords
$\{a,a'\}$,  and $\{b,b'\}$. We have
\begin{equation}\label{Associatorformula}
Z^{u,v} = \Reg(L_{u,\delta}, v) = \Reg(F, 1) =Z^{0,1}(\delta_u,\delta_v)\ .
\end{equation}
It follows from the calculations in  \S\ref{sect55} that the coefficients of $Z^{u,v}$ lie in $\MZV$.
\end{proof}

As an example, consider the case $\Mod_{0,5}$. Then $X_{5,\delta}$ is a pentagon with vertices $v_1,v_2,..,v_5$ in order.  It follows from $(\ref{productforzetaseries})$ that
 $Z^{v_1v_2}Z^{v_2v_3}Z^{v_3v_4}Z^{v_4v_5} Z^{v_5v_1}=1$. Applying $(\ref{Associatorformula})$, we deduce the pentagonal relation
  due to Drinfeld [Dr]:
 $$Z^{0,1}(\delta_{25},\delta_{14}) Z^{0,1}(\delta_{24},\delta_{13})
Z^{0,1}(\delta_{14},\delta_{35})
Z^{0,1}(\delta_{13},\delta_{25})
Z^{0,1}(\delta_{35},\delta_{24})
=1 \in \widehat{\DB}_{5,\delta}(\MZV)\ .$$

\begin{rem}
One can  prove the previous theorem directly using corollary  $\ref{corproductsofmultiplepolylogs}$.
In cubical coordinates, the coefficients of $L_{v,\delta}$ are sums of products of logarithms with the multiple polylogarithms
%$(\ref{multiplepoly})$:
$\Li_{n_1,\ldots, n_r} \Big( {x_{j_1}\ldots x_{\ell} \over x_{j_2}\ldots x_\ell}, \ldots,
{x_{j_{r-1}}\ldots x_{\ell} \over x_{j_r}\ldots x_\ell}
,x_{j_{r}}\ldots x_{\ell}\Big). $
%where $v_0$ is the vertex corresponding to the triangulation
%$\{\{2,4\},\{2,5\},\ldots, \{2,n\}\}$ which has cubical coordinates
 %$x_1=\ldots= x_\ell=0$.
  By taking suitable limits in such coordinate systems, one can see directly that the coefficients of each
  regularised zeta series are multiple zeta values.
\end{rem}

We can compute the monodromy of a solution $L_{v,\delta}(x)$ to
$(\ref{star})$ explicitly in terms of the zeta series defined
above. First, let us define
\begin{equation}
\pi_1^\delta (\Mod_{0,S}) = \pi_1(\Mod_{0,S}, X_{S,\delta})\end{equation}
to be the fundamental group of $\Mod_{0,S}$ relative  to the  set  $X_{S,\delta}$, which can be taken
as a base point because it is contractible.
%just like a base point.
%It can also  be useful to take the set of all vertices $V^{\delta}$ to be a base point,
% but this will not be required here.
For each $\{i,j\} \in \chi_{S,\delta}$, let
$$\gamma_{ij} \in \pi_1^\delta(\Mod_{0,S}) $$
denote a small path which winds once around the face
$D_{ij}=\{u_{ij}=0\}$ in the positive direction, {\it i.e.}, such
that
$$\int_{\gamma_{ij}} {du_{ij} \over u_{ij}} = +2\pi i\ .$$
For each $\{i,j\}\in \chi_{S,\delta}$, let $\Mo_{ij}$ denote the monodromy
operator given by analytic continuation of functions along a loop
which is homotopy equivalent to  $\gamma_{ij}$. The operators
$\Mo_{ij}$ commute with multiplication and differentiation. It
follows that the $\Mo_{ij}$, for  $ \{i,j\}\in \chi_{S,\delta}$, %necessarily
 act on  $L_{v,\delta}(x)$ by
right multiplication by   constant series.
\begin{prop} Let $\{i,j\} \in \chi_{S,\delta}$, and let $v\in V^{\delta}$ denote any  vertex of $\overline{X}_{S,\delta}$.
Choose any vertex $w\in V^{\delta}$ which lies on the face
$D_{ij}$, {\it i.e.}, $u_{ij}(w)=0$. Then for all $x\in
X_{S,\delta}$,
$$\Mo_{ij} \, L_{v,\delta}(x) = L_{v,\delta}(x)\, Z^{v,w} \,  e^{2i\pi
\delta_{ij}} \, Z^{w,v}\ .$$
\end{prop}
\begin{proof}
By theorem $\ref{thmKZhomog}$,
$$L_{w,\delta}(x) = f_w(x) \prod_{\{k,l\} \in F_w}
\exp(\delta_{kl} \log u_{kl})\ ,$$ where $f_w(x)$ is holomorphic
in a neighbourhood of $w\in \Modf_{0,S}(\C)$ which contains the interior of the face $D_{ij}$. By analytic
continuation along a small loop $\gamma_{ij}$ which is contained in this neighbourhood
and winds once around  $D_{ij}$,
 we deduce that
$$\Mo_{ij} L_{w,\delta}(x) = L_{w,\delta}(x) \,\exp(2i \pi \, \delta_{ij} )\ .$$
It follows from the  definition of the zeta series that  $L_{v,\delta}(x) =
L_{w,\delta}(x) \, Z^{w, v}$ for all $x\in X_{S,\delta}$. Since the
quotient $(L_{v,\delta}(x))^{-1} L_{w,\delta}(x) \, Z^{w, v}$ is the constant
function $1$, which is single-valued, the same equation must also
hold for all $x$ in the universal covering space of $\Mod_{0,S}(\C)$.
Therefore,
$$\Mo_{ij}\, L_{v,\delta}(x) = \Mo_{ij} \,L_{w,\delta}(x) \, Z^{w,v} = L_{w,\delta}(x) \, e^{2i\pi  \delta_{ij}}\, Z^{w,v}=
 L_{v,\delta}(x) \,Z^{v,w}\, e^{2i\pi \delta_{ij}} \, Z^{w,v}\ .  $$
%
%The quotient is the constant function $1$ on $X_{S,\delta}$, which
%extends to $\Modf_{S}$. It follows that the same equation holds
%for all $x\in \Modf_{S}$.
\end{proof}
The previous lemma holds for any pair of vertices $w,w'\in
V^{\delta}$ which meet $D_{ij}$. We immediately deduce that the
following identity holds  in $\widehat{\DB}_{S,\delta}(\C)$:
\begin{equation} \label{Zrel}
 Z^{v,w}\,
e^{2i\pi \delta_{ij}} \, Z^{w,v}=Z^{v,w'}\, e^{2i\pi \delta_{ij}}
\, Z^{w',v}    \ .
\end{equation}
This identity in fact follows from the commutation relation $(\ref{braidx})$.
It follows from the previous theorem that the monodromy of
$\Mod_{0,S}$ can be completely expressed in terms of multiple zeta
values, and the constant $2\pi i$.
\begin{cor}
 The monodromy ring of
$\Mod_{0,S}$ is  $\MZV[2\pi i]$.
\end{cor}

%\begin{cor} The  ring  of
%$\Mod_{0,S}$ is  $\MZV$.
%\end{cor}

\subsection{Regularisation of polylogarithms on $\Mod_{0,S}$} \label{sect66}
\begin{defn}For any $v\in V^{\delta}$, let
$L^{v,\delta}(\Mod_{0,S})$
denote the  $\Or(\Mod_{0,S})$-module generated by the coefficients of the solution $L_{v,\delta}$ to $(\ref{star})$ given
by theorem $\ref{thmKZhomog}$.
 Now let us define
\begin{equation}\label{LZdefn}
L^{\delta}_\MZV(\Mod_{0,S}) = L^{v,\delta}(\Mod_{0,S}) \otimes_\Q
\MZV\ .\end{equation} It does not depend, up to isomorphism, on
the choice of the vertex $v$  by lemma $\ref{Zseriesaslimit}$  and
theorem $\ref{thmcoeffsofzaremzvs}$. We write $\Omega^k
L^{v,\delta}(\Mod_{0,S})=L^{v,\delta}(\Mod_{0,S})\otimes_{\Or(\Mod_{0,S})}
\Omega^k(\Mod_{0,S})$ and $\Omega^k
L^{\delta}_\MZV(\Mod_{0,S})=\Omega^k
L^{v,\delta}(\Mod_{0,S})\otimes_\Q \MZV$, for $k\geq 0$.
%which are differential algebras in the obvious way.
\end{defn}

Since $\MZV$ is filtered by the weight, we deduce a natural weight filtration on $L^\delta_{\MZV}(\Mod_{0,S})$ which we denote by $W^\centerdot$.
It follows immediately that any dihedral
symmetry $\sigma\in D_{2n}$ of the $n$-gon $(S,\delta)$ induces an isomorphism of filtered algebras
\begin{equation} \label{dihedralactiononL}
\sigma_*: L^\delta_\MZV(\Mod_{0,S})\overset{\sim}{\rightarrow} L^\delta_\MZV(\Mod_{0,S})\ .
\end{equation}
%Each choice of vertex $v\in V^{\delta}$  defines a different $\Z$-structure on $L^\delta_\MZV(\Mod_{0,S})$ (generated by the coefficients
%of $L_{v,\delta}$).
 Each algebra $L^{v,\delta}(\Mod_{0,S})$ is in fact a graded Hopf algebra (\S\ref{sect67}), although we lose
the grading when we pass to $L^\delta_\MZV(\Mod_{0,S})$, because it is not yet known whether $\MZV$ is graded by the weight.
 %Indeed, this algebra is completely functorial with
%respect to the various maps we have defined between moduli spaces.
The following theorem shows that the function theory of multiple
polylogarithms is dictated by the geometry of the Stasheff
polytopes $\overline{X}_{S,\delta}$.
%We write $\Reg$ to mean the regularisation map which sends not only logarithmic, but also
%polar terms to zero.

\begin{thm} \label{proprestrictionoffunctions} Let $\{i,j\} \in \chi_{S,\delta}$. For any function $f\in L^\delta_\MZV(\Mod_{0,S})$, let
$\mathrm{Reg}(f,D_{ij})$ denote the regularised restriction of $f$
to the divisor $D_{ij}$, which maps not only logarithmic, but also polar singularities to zero. %(definition $\ref{defnlogreg}$)
 As in lemma $\ref{lemmadecomp}$, let $T_1 \cup T_2 = S$ denote the
partition corresponding to the chord $e=\{i,j\}$, such that
 $$D_{ij}\cong \overline{\Mod}^{\delta_1}_{0,T_1\cup e}\times
\overline{\Mod}^{\delta_2}_{0,T_2\cup e}\ .$$  Then there is an
isomorphism of filtered algebras:
$$\mathrm{Reg} \, \big(L^{\delta}_{\MZV} (\Mod_{0,S}), D_{ij}\big) \cong L_{\MZV}^{\delta_1}(\Mod_{0,T_1}) \otimes_{\MZV}
L^{\delta_2}_{\MZV}(\Mod_{0,T_2}) \ .$$
\end{thm}

\begin{proof}
Let us choose any vertex $v\in V^\delta $ such that $u_{ij}(v)=0$.
The algebra $L^{v,\delta}(\Mod_{0,S})\otimes_\Q \MZV$ is generated
by the coefficients of the generating series $L_{v,\delta}(x)$
over $\MZV$. By proposition $\ref{propKZsplit}$, there is a
decomposition
$$L_{v,\delta} = h \,L_1 L_2\ ,$$
where $L_1$, $L_2$ can be viewed as solutions
of $(\ref{star})$ on $\Mod_{0,T_1}$ and $\Mod_{0,T_2}$
respectively, and do not depend on $u_{ij}$. Since the series $h
\, \exp(-\delta_{ij} \log u_{ij})$ is holomorphic in $u_{ij}$ and
is the constant function $1$ along $D_{ij}=\{u_{ij}=0\}$, we have $\mathrm{Reg} (h,D_{ij}) =1$.
Therefore
$$\mathrm{Reg}(L_{v,\delta}, D_{ij}) = L_1 L_2\ .$$
Likewise, for any coefficient $f$ of $L_{v,\delta}$, and any $k\in \Z$,
$$\mathrm{Reg}\Big({f \over u_{ij}^k}, D_{ij}\Big) \in
L^{v_1,\delta_1}(\Mod_{0,T_1})\otimes_\Q
L^{v_2,\delta_2}(\Mod_{0,T_2})\ ,$$ where $v_1,v_2$ are the images
of $v$ defined in $\S6.5$, and $\delta_1,\delta_2$ are the induced
dihedral structures on $T_1,T_2$. This proves that  there is an
isomorphism of filtered algebras $\mathrm{Reg}
(L^{v,\delta}(\Mod_{0,S}), D_{ij}) \cong
L^{v_1,\delta_1}(\Mod_{0,T_1})\otimes_\Q
L^{v_2,\delta_2}(\Mod_{0,T_2})$.  On taking the tensor product
with $\MZV$, we obtain the statement of the theorem.
\end{proof}
\noindent
The theorem states that if we restrict a multiple polylogarithm of weight $m$ to the divisor $u_{ij}=0$, then
we obtain a linear combination of products of multiple zeta values and multiple polylogarithms such that the total weight is at most $m$.

\subsection{The regularised realisation of polylogarithms}  \label{sect67}
Let $S=\{s_1,\ldots, s_n\}$ with the obvious dihedral structure
$\delta$, and let $v\in V^\delta$. We can now define a realisation
of $B(\Mod_{0,S})$ which is regularised at $v$. Let us first
suppose that $v$ corresponds to the vertex $x_1=x_2=\ldots
=x_\ell=0$ in cubical coordinates, in order to exploit the
decomposition of $B(\Mod_{0,S})$ as a product of shuffle algebras.
The corresponding triangulation of the $n$-gon $(S,\delta)$
consists of all chords $\{\{2,4\}, \ldots, \{2,n\}\}$. The
projection map onto $x_\ell=0$ gives a fibration
$$\Mod_{0,\{s_1,\ldots, s_n\}} \To \Mod_{0,\{s_2,\ldots, s_n\}}\ .$$
Correspondingly, we proved that there is a decomposition
$$B(\Mod_{0,S}) = B(\Mod_{0,S'}) \otimes_{\Or(\Mod_{0,S'})} B_{\Mod_{0,S'}}(\Pro^1\backslash \Sigma)\ ,$$
where $S'=\{s_2,\ldots, s_n\}$, and $\Sigma$ is given in $\S\ref{sect54}$.
%The vertex $v$ corresponds to  a triangulation $\alpha \in \chi^{\ell}_{S,\delta}$ of the regular $n$-gon. It necessarily
%has at least two free vertices (\S\ref{sect23}). To each free vertex corresponds a short chord $\{i,i+2\} \in \alpha$. Recall
%that the projection map onto $D_{i\, i+2}$ obtained by setting $u_{i\, i+2}=0$ is a fibration (\S\ref{sect25}).
%  Without loss
%of generality, we can assume $s_n$ is a free vertex, and set $i=n-1$. We proved that
%there is a decomposition
%$$B(\Mod_{0,S}) = B(\Mod_{0,S'}) \otimes  B_{\Mod_{0,S'}}(\Pro^1\backslash \Sigma)\ ,$$
%where $S'=\{s_1,\ldots, s_{n-1}\}$ and $\Sigma=S'$ corresponds to the set of chords emanating from the vertex $n$.
Now let $v'$ denote the vertex corresponding to the restricted
triangulation $\alpha'$ of $S'$ with induced dihedral structure
$\delta'$. By  proposition $\ref{propKZsplit}$, there is a
decomposition $L_{v,\delta}=h\, L_{v',\delta'}$, where $h$ is a
hyperlogarithm equation on $\Pro^1\backslash \Sigma$ in the
variable $x_\ell$. We deduce that
$$L^{v,\delta}(\Mod_{0,S}) = L^{v',\delta'}(\Mod_{0,S'}) \otimes_{\Or(\Mod_{0,S'})}  L_{\Mod_{0,S'}}(\Pro^1\backslash \Sigma)\ ,$$
where  $L_{\Mod_{0,S'}}(\Pro^1\backslash \Sigma)$  denotes the
$\Or(\Mod_{0,S'})$-algebra generated by the coefficients of $h$. From the   realisation
$(\ref{hyperlogrealis})$ we obtain a realisation:
\begin{equation} \label{fibreisom}
\rho_{S'} : B_{\Mod_{0,S'}}(\Pro^1\backslash\Sigma) \overset{\sim}{\To} L_{\Mod_{0,S'}} (\Pro^1\backslash \Sigma)\ ,\end{equation}
which is regularised at  $x_\ell=0$. It is an isomorphism of graded $\Or(\Mod_{0,S'})[\partial/\partial x_\ell]$-algebras.
If we iterate this argument, we obtain two analogous decompositions
\begin{eqnarray} \label{tensordecompforregreal}
B(\Mod_{0,S}) & =& \bigotimes_{1\leq i\leq \ell} B_{\Mod_{0,S_i}} (\Pro^1\backslash \Sigma_i) \ , \\
L^{v,\delta}(\Mod_{0,S}) & = & \bigotimes_{1\leq i\leq \ell} L_{\Mod_{0,S'}} (\Pro^1 \backslash \Sigma_i) \ . \nonumber
\end{eqnarray}
for some subsets $S_1 \subsetneq S_2 \subsetneq \ldots \subsetneq
S_{\ell} \subsetneq  S$ where $|S_1|=3$, and $\Sigma_i\cong S_i$.
Taking the tensor product of the fibre-wise isomorphisms
$(\ref{fibreisom})$, we obtain a map
$$\rho_{v,\delta}: B(\Mod_{0,S}) {\To} L^{v,\delta}(\Mod_{0,S})\ .$$

\begin{thm} \label{thmstructureofL} The map $\rho_{v,\delta}$ is an isomorphism of differential graded algebras. It follows that
every $\Or(\Mod_{0,S})$-differential subalgebra of $L^{v,\delta}(\Mod_{0,S})$ is differentially simple,
and that $L^{v,\delta}(\Mod_{0,S})$ is a polynomial algebra. Furthermore,
$$H^0(L^{v,\delta}(\Mod_{0,S}))=\Q\quad \hbox{ and }\quad H^i(L^{v,\delta}(\Mod_{0,S}))=0 \quad \hbox{for all } i\geq 1\ .$$
The primitive of a closed form $f\in W^b\,\Omega^k L^{v,\delta}(\Mod_{0,S})$ is of weight  at most $b+1$.
\end{thm}

\begin{proof}
The proof of theorem $\ref{thmbarfibration}$ implies that
 the differential structure of the algebras $B(\Mod_{0,S})$ and $L^{v,\delta}(\Mod_{0,S})$ are uniquely
 determined from the tensor decompositions $(\ref{tensordecompforregreal})$, since we have a fixed base point at infinity corresponding
 to $v$.
 It follows that $\rho_{v,\delta}$ is a map of differential graded algebras. The fact that it is an isomorphism then follows immediately from
 corollary $\ref{corBunipclos}$. The rest of the theorem is a consequence of theorem $\ref{thmbartriv}$ and corollary $\ref{corBarmodtriv}$.
\end{proof}
%\begin{cor} $L^{v,\delta}(\Mod_{0,S})$ is isomorphic to the free tensor algebras
%\end{cor}

We obtain a similar decomposition for every vertex $v\in V^\delta$.
\begin{cor}
For each $v\in V^\delta$, there is a canonical realisation
$$\rho_{v,\delta}: B(\Mod_{0,S}) {\To} L^{v,\delta}(\Mod_{0,S})\ ,$$
which is regularised at the vertex $v$.
\end{cor}

%\begin{rem} To any triangulation of an $n$-gon we have associated a filtration of $S$, which in turn gives rise to
%a decomposition of the form $(\ref{tensordecompforregreal})$. There are many ways to do this.
%  For example, if we consider the triangulation
%$\{2,4\}, \{4,6\}, \{6,2\}$ on the hexagon with vertices labelled 1 to 6,  there is the filtration:$$\{2,4,6\} \subset \{1,2,4,6\} \subset \{1,2,3,4,6\} \subset \{1,2,3,4,5,6\}\ .$$
%along with five others. Conversely, any such filtration determines a triangulation in the obvious way.
%\end{rem}

 The map $\rho_{v,\delta}$  has to be defined directly as follows. Recall from $\S\ref{sect62}$ that there is
 a decomposition $B(\Mod_{0,S}) \cong B_{v,\delta}(\Mod_{0,S}) \otimes \Q[[\omega_{i_1\,j_1}],\ldots, [\omega_{i_\ell\, j_\ell}]]$ into convergent and non-convergent
 words, where $\{i_1,j_1\},\ldots, \{i_\ell, j_\ell\}$ are the set of chords in the  triangulation of the $n$-gon
 corresponding to $v$. Then $\rho_{v,\delta}$
  is the unique homomorphism such that
\begin{eqnarray}
\rho_{v,\delta} ([\omega_{i_k\, j_k}] )& =& \log u_{i_k\, j_k} \qquad \hbox{ for all } 1\leq k\leq \ell\ ,\nonumber \\
\rho_{v,\delta} \big(\sum_I c_I [\omega_{i_1}|\ldots |\omega_{i_n}]\big)&= & \sum_I
c_I \int_\gamma \omega_{i_n}\ldots \omega_{i_1} \nonumber\ ,
\end{eqnarray}
for all $\sum_I c_I [\omega_{i_1}|\ldots |\omega_{i_n}] \in
B_{v,\delta}(\Mod_{0,S})$, where $\gamma $ is a smooth path such
that $\gamma(0)=v$ and $\gamma(1)=z\in \Mod_{0,S}(\C)$. Such an
iterated integral converges, since  $B_{v,\delta}(\Mod_{0,S})$ is
spanned by the set of integrable words no element of which ever
ends in a symbol $\omega_{i_k\,j_k}$ for $1\leq k\leq \ell$. The
integrability condition $(\ref{intcond})$ ensures that it only
depends on the homotopy class of $\gamma$ and therefore defines a
multi-valued function on $\Mod_{0,S}(\C)$.  Correspondingly, there
is a decomposition
$$L^{v,\delta}(\Mod_{0,S}) = L^{v,\delta}_c(\Mod_{0,S}) \otimes_\Q \Q[\log u_{i_k\, j_k}: \,\, 1\leq k\leq \ell]\ ,$$
where
$L^{v,\delta}_c(\Mod_{0,S})=\rho_{v,\delta}(B_{v,\delta}(\Mod_{0,S}))$
is the algebra generated by the coefficients of $f_{v,\delta}$
(defined in theorem $\ref{thmKZhomog}$). They   are holomorphic in
a neighbourhood
 of $v$.

\newpage
\section{Period integrals on $\overline{\Mod}_{0,n}(\R)$ and generalised shuffle products.}

 Given a regular algebraic $n-3$-form on $\Mod_{0,S}$, we give
necessary and sufficient conditions for its integral over a
fundamental cell $X_{S,\delta}$ to converge. We obtain a formula
for the order of vanishing of any such form along any given
divisor on $\overline{\Mod}_{0,S}\backslash \Mod_{0,S}$. Finally,
we show how the double shuffle relations for multiple zeta values
are a special case of generalised
 multiplicative structures on the set of all period
integrals.

\subsection{}\label{sect71} The set of all  regular algebraic
$\ell$-forms on $\Mod_{0,S}$  can be written %using dihedral coordinates
in terms of a canonical dihedrally-invariant form which we
construct as follows. Let $\delta $ be a fixed dihedral structure
on $S$, and correspondingly, write
$S=\{s_1,\ldots,s_n\}$. %This induces a dihedral structure on the%
%set of marked points $z_s\in \Pro^1$.
First we define  the following form on $(\Pro^1)_*^n$, where the
indices are taken modulo $n$:
\begin{equation} \label{volform}
 \widetilde{\omega}_{S,\delta} = \bigwedge_{j=1}^n {dz_j \over z_j-z_{j+2}} \ .
\end{equation}
The forms $\widetilde{\omega}_{S,\delta}$ are
$\PSL_2(\C)$-invariant, since if we set
$$z_i' = {\alpha z_i + \beta \over \gamma z_i + \delta}\ , \quad \hbox{for } 1\leq i\leq n\ , \quad
\hbox{ where} \quad \begin{pmatrix}
  \alpha & \beta \\
  \gamma & \delta
\end{pmatrix} \in \PSL_2(\C)\ ,$$
then $dz_i'=(\gamma z_i+\delta)^{-2} dz_i$ and $z'_i-z_j'= (\gamma
z_i+\delta)^{-1} (\gamma z_j +\delta)^{-1} (z_i-z_j),$ and
therefore $\widetilde{\omega}_{S,\delta}$ is unchanged on
replacing each $z_i$ by $z_i'$. In order to define a form on
$\Mod_{0,S}$, consider the quotient map $p:
(\Pro_*^1)^n \rightarrow \Mod_{0,S}$, which has %is a fibration with
fibres $\PSL_2$. Let $v$ denote a fixed non-zero algebraic
invariant $3$-form on $\PSL_2(\C)$ which is  defined over $\Q$.
This is uniquely determined up to a non-zero rational multiple.
% so
%we take
%$$v= \big( x_{11} dx_{21} - x_{21}dx_{11}\big) \wedge dx_{12} \wedge
%dx_{22}\ ,  \hbox{ where} \quad \begin{pmatrix}
%  x_{11} & x_{12} \\
%  x_{21} & x_{22}
%\end{pmatrix} \in \PSL_2(\C)\ .$$
Then there exists a unique form $\omega_{S,\delta}$ on $\Mod_{0,S}$
such that
$$ p^*(\omega_{S,\delta}) \wedge v =  \widetilde{\omega}_{S,\delta} \ .$$
The form  $\omega_{S,\delta}$ is  defined over $\Q$, and is
 $D_{2n}$-invariant by construction.
%It is well-defined up to
%multiplication by $\Q^\times$.
 In simplicial
coordinates $(\ref{ExpSimpCoords})$, and using  the
$\PSL_2(\C)$-invariance of $(\ref{volform})$, we can normalise the rational coefficient of $\omega_{S,\delta}$ such that:
\begin{equation}\label{omegasimp}
\omega_{S,\delta} = {dt_1\wedge\ldots \wedge dt_\ell  \over t_2
(t_3-t_1)(t_4-t_2)\ldots(t_{\ell}-t_{\ell-2}) (1-t_{\ell-1}) } \ ,
 \end{equation}
if $\ell\geq 2$, and $\omega_{S,\delta}=dt_1$ if $\ell=1$.
In dihedral coordinates, one has $\omega_{S,\delta}=du_{24}$ if $\ell=1$, and if $\ell\geq 2$, one can write $(\ref{omegasimp})$ using
$(\ref{tcoords})$ as follows:\begin{equation}\label{omeganexpl}
\omega_{S,\delta} = {du_{24}\wedge du_{25}\wedge \ldots \wedge
du_{2 \, n-1} \wedge du_{2n} \over (1-u_{24}u_{25})
(1-u_{25}u_{26}) \ldots (1-u_{2\, n-1}u_{2n}) }\ .
\end{equation}
The latter representation is not unique because of the various
relations between the functions $u_{ij}$ and their differentials.
 The form $\omega_{S,\delta}$ clearly defines a meromorphic form
 on the compactification $\overline{\Mod}_{0,S}$. For any boundary divisor
 $D\subset \overline{\Mod}_{0,S} \backslash \Mod_{0,S}$ we denote
 by
 $\ord_D\, \omega_{S,\delta} $ the order of vanishing of
 $\omega_{S,\delta}$ along $D$.

\begin{lem}\label{omegahasnopoles}
The form $\omega_{S,\delta}$ has neither zeros nor poles on
$\Modf_{0,S}\backslash \Mod_{0,S}$.
\end{lem}

\begin{proof}
In cubical  coordinates, $\omega_{S,\delta}$ has the
representation:
\begin{equation}\label{omegansimpcube}
\omega_{S,\delta} = {dx_1\wedge\ldots \wedge dx_\ell \over
(1-x_1x_2)\ldots (1-x_{\ell-1}x_\ell)}\ .
\end{equation}
It is clear that  $\omega_{S,\delta}$ is not identically zero nor
infinite along the divisors $x_i=0$, for $1\leq i\leq \ell$. In
other words,  the order of vanishing of $\omega_{S,\delta}$ is
zero along the divisor $u_{2\, i}=0$ for each $4\leq i\leq n$. But
since $\omega_{S,\delta}$ is $D_{2n}$-invariant, it follows that
the order of vanishing of $\omega_{S,\delta}$ is zero along all
divisors at finite distance $u_{i j}=0$, where $\{i,j\}\in
\chi_{S,\delta}$.
\end{proof}
In other words, given any fixed dihedral structure $\delta$ on
$S$, we can define $\omega_{S,\delta}$ to be  the unique (up to
multiplication by $\Q^\times$) non-zero volume form on
$\Mod_{0,S}(\R)$ which has no zeros or poles at finite distance. Equivalently, it has no zeros or poles on the boundary
of the closed Stasheff polytope $\overline{X}_{S,\delta}$.

It follows from $(\ref{regring})$ that every algebraic volume form
on $\Modf_{0,S}(\R)$ can be written as a linear combination of
forms
\begin{equation} \label{genform} \prod_{\{i,j\}\in \chi_{S,\delta}} u_{ij}^{\alpha_{ij}}\,
\omega_{S,\delta} \ , \quad \hbox{ where } \alpha_{ij}\in \Z \quad
\hbox{for each } \{i,j\} \in \chi_{S,\delta}\ .
\end{equation} Now suppose that we are
given a collection of coefficients
$\alpha=(\alpha_{ij})_{\{i,j\}\in \chi_{S,\delta}}$ which are all
non-negative. We define the following family of period integrals:
\begin{equation} \label{Idef}
I_{S,\delta}(\alpha_{ij}) = \int_{\overline{X}_{S,\delta}}
\prod_{\{i,j\}\in \chi_{S,\delta}} u_{ij}^{\alpha_{ij}}
\,\,\omega_{S,\delta} \ .\end{equation} The integral is finite
because  each function $u_{ij}$ is continuous and bounded on the
compact set $\overline{X}_{S,\delta}$.
 Since $\omega_{S,\delta}$ is positive on
$\overline{X}_{S,\delta}$ and invariant under the action of
$D_{2n}$, it follows that $I_{S,\delta}(\alpha_{ij})$ is also
positive, and we have a dihedral transformation formula:
\begin{equation} \label{d2ntrans}
I_{S,\delta}(\alpha_{ij}) = I_{S,\delta}(\alpha_{\sigma(i)\,
\sigma(j)})\quad \hbox{for all } \sigma \in D_{2n}\ .
\end{equation}

 These  integrals %$I_{S,\delta}(\alpha_{ij})$
can be written explicitly in simplicial and cubical coordinates.

\begin{lem} In cubical coordinates, we have the following formula:
\begin{equation} \label{Icube}
I_{S,\delta}(\alpha_{ij}) =\int_{ [0,1]^\ell} \prod_{i=1}^\ell
x_i^{a_i} (1-x_i)^{b_i} \prod_{1\leq i<j\leq \ell} (1-x_ix_{i+1}\ldots x_j
)^{c_{ij}}\, dx_1\ldots dx_\ell \ ,
\end{equation}
 where the indices $a_i$, $b_i$, $c_{ij}\in \Z$ are given by:
\begin{eqnarray}\label{cubeindices}
a_i & = & \alpha_{2\,i+3}\ ,  \\
b_i & = &\alpha_{i+2\, i+4}\ , \nonumber \\
c_{i\, i+1} & = &   \alpha_{i+2\, i+5} - \alpha_{i+2\, i+4} -
\alpha_{i+3\, i+5} -1  \ ,\nonumber \\
c_{i j} & =  &     \alpha_{i+3\, j+3} +\alpha_{i+2\, j+4} -
\alpha_{i+2\, j+3} - \alpha_{i+3\, j+4}\ , \qquad \hbox{ if }
j\geq i+2 \ .\nonumber
\end{eqnarray}
\end{lem}
\begin{proof}
In cubical coordinates, the domain of integration is
$\overline{X}_{S,\delta} \cong [0,1]^{\ell}$, and the only factors
 that occur in the denominator of
$\omega_{S,\delta}$ are $(1-x_ix_{i+1})$  by
$(\ref{omegansimpcube})$. Using the definition of the cross-ratios
$u_{ij}$, we can rewrite the function
$$f=\prod_{\{i,j\}\in \chi_{S,\delta}} u_{ij}^{\alpha_{ij}} =
\pm \prod_{1\leq p< q \leq n} (z_p-z_q)^{s_{pq}}\ ,$$ where the
indices $s_{pq}$ are given by  $s_{pq} = \alpha_{p-1\, q}+
\alpha_{p\, q-1} - \alpha_{p-1\, q-1} -\alpha_{p\, q}$, and where
we set $\alpha_{i\, i+1}=\alpha_{i\,i}=0.$ In cubical coordinates,
we have $z_1=1$, $z_2=\infty$,  $z_3=0$, and $z_{i+3} = x_i\ldots
x_\ell$, for $1\leq i\leq \ell$.  If we put the various elements
together, we obtain the formulae for $b_i$ and $c_{ij}$ given above. The formulae for $a_i$ are easily deduced
using the fact that $x_1=u_{24},\ldots, x_\ell = u_{2n}$.
\end{proof}

 A special sub-family of these integrals were
considered in [Zu], [Fi2], [Zl], where it was also conjectured that
they are expressible in terms of multiple zeta values. It is easy
to verify that the change of variables matrix given by
$(\ref{cubeindices})$ is invertible over $\Z$.

%See also [Fi] for necessary and sufficient conditions on the
%$a_i,b_i,c_{ij}$ for the integral $(\ref{Isim})$ to converge. In
%this context, Main theorem, is conjectured by .... [refs].

Similarly, in simplicial coordinates one can verify that
$$
I_{S,\delta}(\alpha_{ij}) =\int_{\Delta} \prod_{i=1}^\ell
t_i^{a'_i} (1-t_i)^{b'_i} \prod_{1\leq i<j\leq\ell}
(t_j-t_i)^{c'_{ij}} {dt_1\ldots dt_\ell \over t_2 (t_3-t_1) \ldots
(t_\ell-t_{\ell-2}) (1-t_{\ell-1}) } \ ,
$$
where  $\Delta=\{ 0<t_1<\ldots< t_\ell<1\}$ denotes the unit
simplex, and where
\begin{eqnarray}
a_i'& =& %s_{3\, i+3} \quad =
\quad \alpha_{3\,i+2} +\alpha_{2\, i+3} - \alpha_{3\,i+3} -\alpha_{2\,i+2}\ , \qquad 1\leq i\leq \ell\nonumber\\
b_i'& = & %s_{1\, i+3}\quad  =
 \quad \alpha_{n\,i+3} + \alpha_{1\, i+2} - \alpha_{n\,i+2} -\alpha_{1\,i+3}\ ,\qquad 1\leq i\leq \ell \nonumber\\
c_{ij}' &=& %s_{i+3\, j+3}\quad  =
 \quad \alpha_{i+2\, j+3} +\alpha_{i+3\, j+2} -\alpha_{i+3\, j+3} -\alpha_{i+2\,j+2}\ , \quad 1\leq i<j\leq \ell\nonumber
\end{eqnarray}
where we set $\alpha_{i\, i+1}=\alpha_{ii}=0$ as above.
%
%
%$c'_{i\,j} = c_{i+1\, j}$ when $j\geq i+2$,
%$c'_{i\, i+1} = \alpha_{i+2\, i+4}$ for all $1\leq i\leq \ell$,
%$b'_i = c_{\ell \, i}$ if $1\leq i\leq \ell-1$, $a'_i=c_{i\, n-1}$
%for all $1\leq i\leq\ell$, and $b'_\ell =\alpha_{1\, n-1} $.
Once
again, it is not difficult to verify that the corresponding change
of variables matrix is invertible over $\Z$ (this is implied  by equation $(\ref{deltaandt})$).

\subsection{Relative periods and mixed Hodge structures} \label{sectMHS}
% We recall the definition of  the  relative periods of  $\Mod_{0,S}$ and the corresponding mixed
%Hodge structures, as defined in [G-M].
 Let $n=|S|=\ell+3$, and let $A$,
$B$ denote two sets of divisors at infinity on
$\overline{\Mod}_{0,S} \backslash \Mod_{0,S}$, where we assume
that $A\cap B$ is of codimension at least $2$, {\it i.e.}, $A$ and
$B$ have no shared irreducible components. Consider the relative
cohomology group \begin{equation} \label{Hodgestr}
H^\ell(\overline{\Mod}_{0,S} \backslash A, B \backslash B\cap A)\
, \end{equation} which has a canonical mixed Hodge structure [De].
Since the divisor $A\cup B$ is globally normal crossing, this can
be computed using the techniques of [G-S], [V], and it is easily
verified that it is of Tate type. %({\it i.e.}, all Hodge numbers
%are equal).
 Goncharov and Manin  construct an object in the
abelian
 category of mixed Tate motives  $\MT(\Q)$ over $\Q$
%, as defined
%by Goncharov and Deligne
[Go, D-G], whose Hodge realisation is the mixed Hodge structure
$(\ref{Hodgestr})$. They then show that this motive is unramified
over $\Z$.
 We shall write the corresponding motive and
mixed Hodge structure with the same symbol. No confusion arises
because the Hodge realisation functor is fully faithful over $\Q$
([D-G], proposition 2.9).
 Suppose that we are given a  relative
homology cycle
$$[\Delta_B] \in H_\ell(\overline{\Mod}_{0,S} ,B)\ .$$
We can assume that this class is represented by  a smooth compact
real submanifold with corners  $\Delta_B$ whose codimension-$k$
boundary is contained in the $k$-stratum of $B$. More precisely,
if $B$ consists of irreducible components $B_i$, for $1\leq i\leq
N$, then
\begin{equation} \label{Bstrat}
\partial^k \Delta_B = \Delta_B \cap \bigcap_{i_1,\ldots,
i_k}  B_{i_1}\cap \ldots \cap B_{i_k}\ .\end{equation}
% must be contained in the union of
%all $k$-fold intersections $B_{i_1}\cap \ldots \cap B_{i_k}$.
 Suppose  that we are
given an algebraic $\ell$-form $\Omega_A$ on $\Mod_{0,S}$ which is
defined over $\Q$ and whose singularities are contained in $A$.
Then the \emph{relative period integral} of $\Omega_A$ along
$\Delta_B$ is defined to be
\begin{equation} \label{relperiod}
\int_{\Delta_B} \Omega_A \in \C \ .\end{equation} By a
higher-dimensional version of Cauchy's theorem, this integral  is
invariant under continuous deformations of $\Delta_B$ relative to
$B$. We can thus assume that $\Delta_B$ is disjoint from $A$, and
therefore the integral is bounded, since $\Omega_A$ is continuous
on $\Delta_B$, which is  compact. Note that the integral  depends
on the  \emph{relative} cohomology classes of $\Omega_A$ in
$H^\ell(\overline{\Mod}_{0,S} \backslash A, B \backslash B\cap
A)$, and $\Delta_B$ in $H_\ell(\overline{\Mod}_{0,S} \backslash A,
B \backslash B\cap A)$.

 %and not  on its cohomology class $[\Omega_A]\in
%H^\ell(\overline{\Mod}_{0,S} \backslash A)$.
%We say that the integral $(\ref{relperiod})$ is a \emph{real
%relative period} if $\Delta_B \subset \overline{\Mod}_{0,S}(\R)$.
%and $\Omega_A$ is defined over $\Q$.%\footnote{This is not to be confused with the
%notion of the real period of a mixed Hodge structure [BD]}

\begin{lem} \label{Cellsspangradedhomology} $\gr_0^W H_{\ell}(\overline{\Mod}_{0,S}, B)$ is spanned by the homology classes of a number  of cells
$X_{S,\delta}$, where  $\delta$ is in a certain set of dihedral
structures which depends upon $B$.
\end{lem}
\begin{proof}
 The relative cohomology group $H_\ell(\overline{\Mod}_{0,S}, B)$ can be computed using the spectral sequence of the complex
%$$\overline{\Mod}_{0,S} \leftarrow \bigcup_i B_i \leftarrow \bigcup_{i,j} B_{ij} $$
$$\overline{\Mod}_{0,S} \leftarrow \bigsqcup_{i}
B_i \leftleftarrows \bigsqcup_{i,j} B_{i,j} \leftleftarrows \ldots
\leftleftarrows \bigsqcup_{|I|=\ell} B_I \ ,
$$
where $B$ is the union of a set of divisors $B_i$, and
$B_I=\bigcap_{i\in I} B_i.$ The spectral sequence degenerates on
the $E^2$ level and  it follows that
$$\gr^W_0 H_\ell(\overline{\Mod}_{0,S},B) \cong \ker \big( \bigoplus_{|I|=\ell} H_0(B_I) \To \bigoplus_{|J|=\ell-1} H_0(B_J)       \big)\ . $$
In the simplicial complex defined by $B$, this is just $\ker (\C^p
\rightarrow \C^e)$ where $p$ is the number of points, {\it i.e.},
$\ell$-fold intersections of divisors, and $e$ is the number of
edges, {\it i.e.}, $\ell-1$-fold intersections. This also computes
the
number of  cells in $\overline{\Mod}_{0,S}(\R)$ bounded by  $B\cap \overline{\Mod}_{0,S}(\R)$. % - {\it i.e.} it coincides with the
%simplicial cohomology $H^\ell(B\cap \R)$ on
%taking the dual. %definition of simplicial cohomology.
Since $\Mod_{0,S}(\R)$ is tesselated by the cells $X_{S,\delta}$
(lemma $\ref{lemmatess}$),
 they must generate $\gr^W_0
H_\ell(\overline{\Mod}_{0,S},B)$.
%  The result follows easily, since
%. An other way to see this is to observe that $H_0(B_I)= H_0(B_I \cap \R)$.
\end{proof}

\begin{lem}\label{lemrealperiodsareIs}
 Every  relative period integral over a union of cells $\overline{X}_{S,\delta_i}$   is a  $\Q$-linear
combination of $I_{S,\delta}(\alpha_{ij})$, where the
$\alpha_{ij}$ are all non-negative.
\end{lem}
\begin{proof} %By lemma , there is a tesselation $\Mod_{0,S}(\R) =
%\coprod_{\delta} X_{S,\delta}$ where $\delta$ ranges over the set
%of all dihedral structures on $S$. A real relative  cycle
%$\Delta_B$ is therefore homology equivalent to  a union of cells
%$\overline{X}_{S,\delta_i}$ for $1\leq i \leq N$.
Fix a dihedral structure $\delta$ on $S$.  Any such  integral $I$
can be written:
$$I=\sum_{i=1}^N \int_{\overline{X}_{S,\delta_i}} \omega_i =
\int_{\overline{X}_{S,\delta}} \sum_{i=1}^N \sigma_i^* (\omega_i)
\ ,
$$
where   $\omega_i\in \Omega^{\ell}(\Mod_{0,S}(\R))$, and  where
$\sigma_i$ is an element of $\Sym(S)$ which maps the dihedral
structure $\delta_i$ onto $\delta$.
%Since the ring of
%regular functions on $\Mod_{0,S}$ is $\Q(u_{ij}, u_{ij}^{-1})$,
The right-hand side can be written
$$I=\int_{\overline{X}_{S,\delta}} f \,\omega_{S,\delta}\ ,$$
where $f\in \Q[u_{ij}, u_{ij}^{-1}]$ is a regular function on
$\Mod_{0,S}$. Note that by lemma $\ref{lemmanopolesint}$ this
integral converges absolutely if and only if $f\,
\omega_{S,\delta}$ has no poles along $\partial
\overline{X}_{S,\delta}$. Since $\partial \overline{X}_{S,\delta}
$ is the union of divisors $D_{ij}=\{u_{ij}=0\}$, and since
$\ord_{D_{ij}} f\, \omega_{S,\delta} = \ord_{D_{ij}} f$ (by lemma
$\ref{omegahasnopoles}$), this implies that $f\in \Q[u_{ij}]$.
Since $f$ is a polynomial in the $u_{ij}$, it can be written as a
linear combination of monomials with positive exponents, or in
other words, $I$ is a finite $\Q$-linear combination  of integrals
$I_{S,\delta}(\alpha_{ij})$, with $\alpha_{ij}$ all non-negative.
\end{proof}

In order to rephrase the above in motivic terms, we need to recall
the notion of framings from [Go1-2], [BGSV]. Let $m\geq 0$ denote
an integer. An $m$-framing on a mixed  Tate motive (or its Hodge
realisation) is given by two morphisms
$$v: \Q(-m) \rightarrow \gr_{2m}^W M \quad \hbox{ and } \quad
f: \Q(0) \rightarrow (\gr_{0}^W M)^\vee \ .$$ Two framed mixed
Tate motives $(M,v,f)$ and $(M',v',f')$ are said to be equivalent
if there is a morphism $M\rightarrow M'$ which respects the
framings.
% Let $[M,v,f]$ denote the equivalence class for this
%relation.
 The framings on the motive $(\ref{Hodgestr})$ were
defined in [G-M] as follows. There are isomorphisms
\begin{eqnarray}
\gr_0^W H_\ell(\overline{\Mod}_{0,S} \backslash A, B \backslash
B\cap A) &\cong& \gr_0^W H_\ell(\overline{\Mod}_{0,S}, B)\ ,
\nonumber \\
 \gr_{2\ell}^W H^\ell (\overline{\Mod}_{0,S} \backslash A, B \backslash
B\cap A) & \cong & \gr^W_{2\ell}
H^\ell(\overline{\Mod}_{0,S}\backslash A) \nonumber \ .
\end{eqnarray}
Therefore, the classes $[\Delta_B] \in \gr_0^W
H_\ell(\overline{\Mod}_{0,S}, B)$ and $[\Omega_A]\in \gr^W_{2\ell}
H^\ell(\overline{\Mod}_{0,S}\backslash A)$ define an
$\ell$-framing on $(\ref{Hodgestr})$. Note that these framings
could be zero.

We introduce a simplified variant of the above motives. Let
$\delta$ denote a fixed dihedral structure on $S$ and let
$D_\delta$ denote the set of divisors at finite distance in
$\Modf_{0,S}$. These are the  affine varieties which bound the
fundamental cell $\overline{X}_{S,\delta}$. Let $\omega\in
\Omega^\ell(\Mod_{0,S})$ denote an algebraic $\ell$-form with no
singularities along $D_{\delta}$, which is defined over $\Q$. Let
us define the $\ell$-framed mixed Tate motive:
\begin{eqnarray} \label{mymotivedef}
m_{S,\delta}(\omega) =  \big(H^\ell( \Modf_{0,S}, D_{\delta} ) ,
[\overline{X}_{S,\delta}], [\omega] \big)\ ,
\end{eqnarray}
 equipped with the framings given
by the  class of the  fundamental cell
$\Delta_B=\overline{X}_{S,\delta}$, and the  class of $\omega$.
%$$\prod_{\{i,j\}\in \chi_{S,\delta}}
%u_{ij}^{\alpha_{ij}} \, \omega_{S,\delta} \in H^\ell
%(\Modf_{0,S})\ .$$
The framed motives $m_{S,\delta}(\omega)$ are more convenient to
work with  because the varieties $\Modf_{0,S}$ are affine, and we
do not need to keep track of the divisor data at infinity.  Lemmas
$\ref{Cellsspangradedhomology}$ and $\ref{lemrealperiodsareIs}$
imply that the framed mixed Tate motive
$(H^\ell(\overline{\Mod}_{0,S}\backslash A, B\backslash B\cap A),
[\Delta_B], [\Omega_A]) $,
%in the case when $\Delta_B\subset
%\overline{\Mod}_{0,S}(\R)$ is real,
 is equivalent to a linear
combination of motives:
$$m_{S, \delta}\Big(  f \, \omega_{S,\delta}           \Big)\ ,\quad \hbox{ where } f =
\prod_{\{i,j\}\in \chi_{S,\delta}} u_{ij}^{\alpha_{ij}}\ .$$ The
equivalence is  given by natural inclusion maps between moduli
spaces, the action of the symmetric group,
%$\Modf_{0,S} \rightarrow \overline{\Mod}_{0,S}\backslash A$
and
the additivity of framed objects with respect to their framings.
%Note, however, that the lemma above is a stronger result, since th
%The
%motives $\Modf_{0,S}$ are more convenient to work with, because we
%c%an ignore all the divisor data at infinity, all the information
%being contained in the framing.

\subsection{Formulae for the divisor of singularities}\label{sect73} In order to compute the divisor of singularities of
an arbitrary form $(\ref{genform})$, it suffices to compute the
order of the canonical form $\omega_{S,\delta}$ along each divisor
at infinity. This is easily done by exploiting the action of the
symmetric group.

\begin{prop} \label{orderofomega} Let $|S|=n=\ell+3$, and let $D$ denote the divisor given by the stable
partition $S^1\cup S^2 =S$ (proposition
$\ref{divisorsatinfinity}$). Then
%$$2 \,\ord_D \,\omega_{S,\delta}=\ell-\#\{1\leq i \leq n : \hbox{such that } \{i,i+2\} \subset S^1 \hbox{ or } \{i,i+2\} \subset S^2 \}\ .$$
$$\ord\, \omega_{S,\delta} = {\ell-1\over 2} - {1\over 2} \sum_{i\in
\Z/n\Z} \I_D(i,i+2)\ ,$$
where the notation $\I_D$ is defined by equation $(\ref{indicatordefn})$ in $\S\ref{sect26}$.
\end{prop}

\begin{proof}
Let $k\geq 2$ denote the number of elements in $S^1$. Let $\sigma$
denote a permutation $\sigma \in \Sym(n)$ such that $\sigma^{-1}
(S^1)= \{1,2,\ldots, k\}$. By $(\ref{volform})$, we have
$$ \sigma^*(\widetilde{\omega}_{S,\delta}) =  \pm \prod_{i\in \Z/n\Z} \Big({z_i -z_{i+2} \over
z_{\sigma(i)} - z_{\sigma(i+2)} }\Big) \,
\widetilde{\omega}_{S,\delta} \ .$$ By passing to the quotient $p:(\Pro^1)_*^S\rightarrow \Mod_{0,S}$, we
have
$$\ord_D\, \omega_{S,\delta} = \ord_{D}\, \sigma^*(\omega_{S,\delta})
- \ord_D f\ ,$$ where the function
$$f=\prod_{i\in \Z/n\Z} \Big({z_i -z_{i+2} \over
z_{\sigma(i)} - z_{\sigma(i+2)} }\Big) $$ is homogeneous and
 $\PSL_2(\C)$-invariant by the remarks in $\S7.1$. It can
therefore be written as a product of cross-ratios and is  a
well-defined function on $\Mod_{0,S}$.
% \,\prod_{i\in \Z/n\Z} \Big({z_i -z_{i+2} \over
%z_{\sigma(i)} - z_{\sigma(i+2)} }\Big)\ .$$
Now $\ord_{D}\, \sigma^*(\omega_{S,\delta}) = \ord_{\{1,\ldots,
k\}}\, \omega_{S,\delta}=0$,  since the divisor given by the
stable partition $\{1,\ldots,k\} \cup \{k+1,\ldots,n\}$ is
$D_{kn}=\{u_{kn}=0\},$ and we know by  lemma $\ref{omegahasnopoles}$ that
$\omega_{S,\delta}$ has no zeros or poles at finite distance.
Therefore
\begin{equation} \label{pf10}
\ord_D\, \omega_{S,\delta} = - \ord_D\, f\ ,\end{equation} and it
suffices to compute the zeros and poles of $f$.
 Recall from corollary  $\ref{corbasicordformula}$
that
$$\ord_D \, {(z_i-z_k)(z_j-z_l) \over
(z_i-z_l)(z_j-z_k)} \nonumber \\
= {1\over 2} \Big[ \I_D(i,k) +\I_D(j,l)
-\I_D(i,l)-\I_D(j,k)\Big]\ .
$$
%\begin{eqnarray}
%\ord_D \,[ij|kl]& =& \ord_D \, {(z_i-z_k)(z_j-z_l) \over
%(z_i-z_l)(z_j-z_k)} \nonumber \\
%&=& {1\over 2} \Big[ \I_D(i,k) +\I_D(j,l)
%-\I_D(i,l)-\I_D(j,k)\Big]\ . \nonumber
%\end{eqnarray}
 We deduce
that
$$\ord_D\, f ={1\over 2} \sum_{i\in \Z/n\Z} \Big(\I_D(i,i+2) -
\I_D(\sigma(i),\sigma(i+2)) \Big)\ .$$
%\begin{eqnarray}
%\ord_D\, f &=&{1\over 2} \#\{1\leq i \leq n : \hbox{such that }
%\{i,i+2\} \subset S^1 \hbox{ or } \{i,i+2\} \subset S^2 \}
%\nonumber \\
%&-&{1\over 2} \#\{1\leq i \leq n : \hbox{such that } \sigma
%\{i,i+2\} \subset S^1 \hbox{ or } \sigma \{i,i+2\} \subset S^2 \}
%\nonumber
%\end{eqnarray}
But $\{\sigma(i),\sigma(i+2)\}\subset S^1$ if and only if
$\{i,i+2\}\subset \{1,\ldots,k\}$. The number of such pairs is
exactly $k-2$. Likewise, the number of $i$ such that
$\{\sigma(i),\sigma(i+2)\}\subset S^2$ is $n-k-2$. It follows that
the second quantity  in the sum directly above is $n-4=\ell-1$.
This completes the proof on substituting into $(\ref{pf10})$.
\end{proof}

\noindent We immediately deduce the following  formula for the
order of vanishing of an arbitrary form along any divisor
$D\subset \overline{\Mod}_{0,S}\backslash \Mod_{0,S}$. Let
 $$f=\prod_{\{i, j\}\in \chi_{S,\delta}} u_{ij}^{\alpha_{ij}}
\ , \quad  \alpha_{ij} \in \Z \ .$$

\begin{cor}
Let $D$ and $f$ be as above. Then %, and let $f= \prod_{\{i, j\}\in
%\chi_n} u_{ij}^{\alpha_{ij}} \omega_{S,\delta}\ .$ Then
\begin{eqnarray}
2\,\ord_D\, f\omega_{S,\delta}&=& \!\!\!\! \sum_{\{i,j\}\in
\chi_{S,\delta}} \!\!\!\alpha_{ij} \big[
\I_D(i,j+1)+\I_D(i+1,j) -\I_D(i,j)- \I_D(i+1,j+1)\big]\nonumber \\
&& +  \,\,(\ell-1)  - \sum_{i\in \Z/n\Z} \I_D(i,i+2) \ . \nonumber
\end{eqnarray}
\end{cor}
\begin{proof}
This follows immediately from the additivity of $\ord_D$ and the
fact that
$$2\,\ord_{D} u_{ij} = 2\,\ord_D [i \,i+1| j+1\, j]=\I_D(i,j+1)+\I_D(i+1,j) -\I_D(i,j)-
\I_D(i+1,j+1)\ .$$
\end{proof}
% An arbitrary regular function on $\Mod_{0,n}$ can be
%written in the form
%$$f=\prod_{\{i,j\} \in S_n} u_{ij}^{\alpha_{ij}}\ , $$
%where $\alpha_{ij}$ are $n(n-3)/2$ non-negative integers.
\noindent Note that along  each divisor at finite distance
$D_{ij}=\{u_{ij}=0\}$, where $\{i,j\}\in \chi_{S,\delta}$, we
clearly have $\ord_D (f) = \alpha_{ij}$.  In total, there are as
many boundary divisors $D\in \overline{\Mod}_{0,S}$ as there are
partitions of $S$ into two sets, each containing at least two
elements. These number $2^{n-1}-n-1$, but there are only $n(n-3)/2$
parameters $\alpha_{ij}$, which implies that there are many
relations between the quantities $\ord_{D} f$, for varying $D$. The following lemma
gives an alternative approach for  computing  the orders of
functions along divisors.

\begin{lem}
Let $D$ be the divisor of $\overline{\Mod}_{0,S}\backslash
\Mod_{0,S}$  corresponding to a stable partition $S_1\cup S_2 =S$.
For each two-element subset $T=\{s_i,s_j\} \subset S_1$, let $D_T$
denote the divisor given by the partition $T$ and its complement
$S\backslash T$. Then for any function $f\in \Q(\Mod_{0,S})$,
\begin{equation}\ord_{D} f = \sum_{T\subset S_1, |T|=2} \ord_{D_T} f\ .\end{equation}
\end{lem}
\begin{proof} It suffices to verify  the formula for the function
$f=u_{ij}$, where $\{i,j\} \in \chi_{S,\delta}$, since it is
compatible with products. By corollary
$\ref{corbasicordformula}$,
$$\ord_D\, u_{ij} =
\begin{cases}
    -1 & \text{if } \quad  i,j \in S_1,  \qquad \hbox{ and  } \quad i+1,j+1 \in S_2\ ,  \\
 1 & \text{if  } \quad i,j+1\in S_1\, \quad \hbox{ and  } \quad
i+1,j \in S_2\ ,
  \end{cases}
$$
and is $0$ otherwise. One can check that the identity holds for
two-element subsets of $S_1\cap \{i,i+1,j,j+1\}$, from which
it follows in general. For example, if $S_1\cap \{i,i+1,j,j+1\} =
\{i,i+1, j+1\}$, then $\ord_{D} u_{ij}= 0$ and this is equal to
$$\ord_{D_{\{i,i+1\}}} u_{ij} + \ord_{D_{\{i,j+1\}}}
u_{{ij}}+\ord_{D_{\{i+1,j+1\}}} u_{ij} = 0 +1-1\ .$$
\end{proof}
%\noindent Applying the lemma to the set $S=\{1,2,3\}$ and its
%two-element subsets $\{1,3\}$, $\{1,2\}$, and $\{2,3\}$, we
%deduce, for example, that $\ord_{D_{\{1,3\}}} f =
%\alpha_{n3}-\alpha_{13}-\alpha_{n2},$ when $n\geq 5$. Using lemma
%$7.3$, one  easily  shows that $ \ord_{D_{\{1,3\}}}
%\omega_{S,\delta}=-1$, and hence
%\begin{equation} \label{auxiltriple}
%\ord_{D_{\{1,3\}}}\Big(\!\! \prod_{\{i,j\}\in
%\chi_{S,\delta}}u_{ij}^{\alpha_{ij}} \,\omega_{S,\delta}\Big) =
%\alpha_{n3}-\alpha_{13}-\alpha_{n2}-1\ .
%\end{equation}
% Likewise, whenever $n\geq
%6$, one can check that
%\begin{equation} \label{auxilquad}
%\ord_{D_{\{1,4\}}}  \Big(\!\!\prod_{\{i,j\}\in
%\chi_{S,\delta}}u_{ij}^{\alpha_{ij}} \,\omega_{S,\delta}\Big) =
%\alpha_{n3}+ \alpha_{n4} -\alpha_{n3}-\alpha_{14}\ .
%\end{equation}
% In general,  restricting the locus of
%singularities of a set of $\ell$-forms on
%$\overline{\Mod}_{0,S}(\R)$ is equivalent to imposing linear
%inequalities on the indices $\alpha_{ij}$. Linear forms such as
%$(\ref{auxiltriple})$ and $(\ref{auxilquad})$ generalise the
%auxiliary parameters defined in [RV1,2].

\subsection{Singularities of the Kontsevich multiple zeta value forms} \label{sect74}
%Following Goncharov and Manin, we can
There is  a special set of $\ell$-forms on $\Mod_{0,S}$
corresponding to the iterated integral representations  of
multiple zeta values due to Kontsevich.  Let $n\geq 5$. We  apply
the previous proposition to compute the divisor of singularities
of each such form and retrieve one of the results of [G-M]. Let
$\underline{\epsilon}= (\epsilon_1,\ldots, \epsilon_\ell)$ where
$\epsilon_1,\ldots, \epsilon_\ell \in \{0,1\}$. We define
\begin{equation} \label{gammaindexdefn}
\gamma_i = 3-2\epsilon_i\in \{1,3\}\ , \qquad \hbox{for } 1\leq i\leq \ell\ ,
\end{equation}
and set
\begin{equation} \label{omegaeps}
\Omega(\underline{\epsilon}) =\big[5 \,\,n \,|\, 3\,\, 2\big] \,\, \big[2\,\, n|1\,\,3\big]^{\varepsilon_\ell}\,\, \prod_{i=1}^{\ell-1}
 \big[ i+5\,\,\,
\gamma_i\, |\,i+3\, \, \, 2 \big]
 \,\,\,
 \omega_{S,\delta}\ .
\end{equation}
The term in the product corresponding to $i=\ell-1$ requires explanation. We define
$[n+1\, \gamma_{\ell-1}|n\!-\!1\,\, 2]= [1\,3|n\!-\!1\,\, 2]$ if $\gamma_{\ell-1}=3$, and define it to be $1$ if $\gamma_{\ell-1}=1$.
%\begin{equation} \label{omegaeps}
%\Omega(\underline{\epsilon}) =\prod_{i=0}^{\ell-1} \big[ i+5\,\,\,
%\epsilon'_i\, |\,i+3\, \, \, 2 \big]\,\, % \,[5 \,\,n \,|\, 3\,\, 2]\,\,
% \omega_{S,\delta}\ ,
%\end{equation}
%where we set $\epsilon'_i =3- 2\epsilon_i \in \{1,3\}$ for $1\leq
%i\leq \ell-1$ and $\epsilon_0'=n$.
We can write this  expression in explicit  simplicial coordinates $(\ref{ExpSimpCoords})$ by setting $z_1=1$,
$z_2=\infty$ and  $z_3=0$. If we define $t_{\ell+1}=1$,  one can verify
using $(\ref{omegasimp})$ that
\begin{equation}\label{omegasimplicKontdefn}
\Omega(\underline{\epsilon}) = {t_2\over t_\ell} \, \Big({t_\ell \over t_\ell-1}\Big)^{\varepsilon_\ell}
\,
\prod_{i=1}^{\ell-1} \Big( { t_{i+2}-t_i \over
\epsilon_i-t_i} \Big)\,\omega_{S,\delta} = \bigwedge_{i=1}^\ell {dt_i
\over \epsilon_i-t_i}\ .
\end{equation}
 Let $X=\{x_0,x_1\}$ be an alphabet with two letters
 as considered in $\S\ref{sect55}$. Assume that $\epsilon_1=1$ and $\epsilon_\ell=0$,
and define a word $w=x_{\epsilon_\ell}\ldots x_{\epsilon_1}\in x_0X^*x_1$. Let $r_\epsilon = \sum_{i=1}^\ell \epsilon_i$.
 It
follows from $(\ref{Xsimplex})$ and a well-known formula
for $\zeta(w)$  that
\begin{equation}\label{Leibniz-Kontsevich}
\int_{X_{S,\delta}} \Omega(\underline{\epsilon}) =\int_{0<t_1<\ldots <t_\ell<1} \bigwedge_{i=1}^\ell
{dt_i \over \epsilon_i-t_i} =(-1)^{\ell-r}  \zeta(w)\ .
\end{equation} The integral
converges if and only if  $\epsilon_1=1$ and $\epsilon_\ell=0$. It
follows that every multiple zeta value of weight $\ell$ occurs as
a relative period of $\Mod_{0,\ell+3}(\R)$ [G-M].
\begin{lem}
Let $\underline{\epsilon}=(\epsilon_1,\ldots, \epsilon_\ell)$ with
$\epsilon_i \in \{0,1\}$ for all $1\leq i\leq \ell$. Then
$$
2\,\ord_{D}\, \Omega(\underline{\epsilon})=\ell-1 +
\sum_{i=1}^{\ell} \big[\I_D(2,\gamma_i) - \I_D(i+3,\gamma_i)\big]
-\sum_{k \neq 2} \I_D(2, k) - \I_D(1,3)\ ,
$$
%where $E_D= -\I_D(1,3)$ and
where  $\gamma_i \in \{1,3\}$ is defined in $(\ref{gammaindexdefn})$.
\end{lem}

\begin{proof} First we assume that $\epsilon_\ell=0$, and therefore $\gamma_\ell=3$. It follows from
 $(\ref{omegaeps})$  that
 \begin{eqnarray} 2\, \ord_D \,\Omega(\underline{\epsilon}) &=& \sum_{i=1}^{\ell-1}
\I_D(i+3,i+5) + \I_D(2,\gamma_i) - \I_D(2, i+5)- \I_D(
i+3,\gamma_i) \nonumber \\
& + &  \I_D(3,5) + \I_D(2,n) - \I_D(2, 5)- \I_D(
3,n) + \ord_D\,\omega_{S,\delta}\ ,
\end{eqnarray} and the formula
stated above follows on substituting the expression for $\ord_D
\,\omega_{S,\delta}$ given in  proposition $\ref{orderofomega}$.
In the case where $\epsilon_\ell=1$, $\gamma_\ell=1$, a similar formula for
$\Omega(\underline{\epsilon})$ holds except that one must multiply
by an extra  cross-ratio $[2 \, n | 1\, 3]$. This contributes
$$\I_D(2,1) - \I_D(n,1)     - (\I_D(2,3)- \I_D(n,3))$$ in the expression
above, and this is precisely what is required for  the formula to
hold in this case also.
\end{proof}

Let $D\subset \overline{\Mod}_{0,S}\backslash \Mod_{0,S}$ be the
divisor corresponding to a stable partition $S_1\cup S_2$ of $S$.
Then, up to permuting the sets $S_1$ and $S_2$,  $D$ is one of the following four types, where $A\cup B$ is a
partition of  $\{s_4,\ldots, s_n\}$:
\begin{enumerate}
  \item $S_1=\{s_1,s_2,s_3\}\cup A$,\, \,\,\,$S_2=B$.
 \item $S_1=\{s_1,s_3\}\cup A$,\quad\quad $S_2= \{s_2\} \cup B$.
 \item $S_1= \{s_1,s_2\} \cup A$,\quad\quad $S_2= \{s_3\} \cup B$.
\item $S_1= \{s_2,s_3\} \cup A$,\quad\quad $S_2= \{s_1\} \cup
B$.
\end{enumerate}
 \begin{cor} Let $A_0= A\cap \{
s_{i+3}\,\,\hbox{for } 1\leq i\leq \ell \hbox{ such that } \epsilon_i=0\}$, $A_1=A\backslash A_0$, and define $B_0$,
$B_1$, similarly. Then, according to each of the cases above,
$$ \ord_D\, \Omega(\underline{\epsilon}) = \begin{cases}
    |B|-2  & \text{if } D \text{ is as in case (1)}\ , \\
    -1  & \text{if } D \text{ is as in case (2)} \ , \\
     |B_1|-1 & \text{if } D \text{ is as in case (3)} \ , \\
|B_0|-1& \text{if } D \text{ is as in case (4)} \ .
  \end{cases}
$$
\end{cor}

\begin{proof}
In case $(1)$,  the formula stated in the previous  lemma gives,
term by term,
$$2\, \ord_{D}(\Omega(\underline{\epsilon})) = \ell-1 + \ell -
|A|-(|A|+2) -1 \ ,$$ and the formula follows, since
$\ell=|A|+|B|$. In case $(2)$, it gives
$$2\, \ord_{D}(\Omega(\underline{\epsilon})) = \ell-1 + 0 -
|A|-|B| -1 = -2\ .$$ In case $(3)$, we have $S_1= \{s_1,s_2\} \cup
A_0 \cup A_1 $ and $S_2=\{s_3\} \cup B_0 \cup B_1$. The formula in
the previous lemma gives, term by term:
$$2\, \ord_{D}(\Omega(\underline{\epsilon})) = \ell-1 + (|A_1|+|B_1|)
-(|B_0|+|A_1|) -(1+ |A_0|+|A_1|) -0\ ,$$ but since
$\ell=|A_0|+|A_1|+|B_0|+|B_1|$, this is just $2|B_1|-2$, as
required. The formula for case $(4)$ follows by symmetry.
\end{proof}
In case $(1)$, we must have $|S_2|=|B|\geq 2$, otherwise the
partition $S_1\cup S_2$ is not stable, so no singularity ever
occurs along such a divisor. It follows that the divisor of
singularities of $\Omega(\underline{\epsilon})$ are precisely
those divisors of type $(2)$, and those of type $(3)$ (resp.
$(4)$) for which $B_1$ (resp. $B_0$) is empty. Let us set $s_1=1$,
$s_2=\infty$, and $s_3=0$, as usual. Then the divisors of type
$(2)$ correspond to the divisors which are called `type $\infty$'
in [G-M]. The divisors of type $(3)$ for which there is a pole are
partitions of the form $\{1,\infty\}\cup A$ and $\{0\}\cup B$ ,
where $B=B_0$ and hence $B\subset \{s_{i+3}: \, \epsilon_i=0\}$.
These are exactly the divisors of `type $0$' according to [G-M].
Similarly, our type $(4)$ above corresponds to `type $1$' and the
previous result implies proposition $3.1$ of $[$Go-Ma$]$. Note
that the above proof only uses the action of the symmetric group
and does not use any blow-ups.

\subsection{Generalised products and the double shuffle relations}\label{sect75} In $\S\ref{sect27}$  we considered
non-degenerate coordinate systems
$$f= \prod_{i=1}^k f_{T_i} : \Mod_{0,S} \To \prod_{i=1}^k
\Mod_{0,T_i}\ ,$$ where  the
sets $T_i$ cover $S$ and the dimensions satisfy $(\ref{dimcond})$.
Since the dimension of $\Mod_{0,T_i}(\R)$ is $|T_i|-3$, the
K\"{u}nneth formula gives an isomorphism
$$\bigotimes_{i=1}^k
H^{|T_i|-3} (\Mod_{0,T_i}(\R)) \cong H^{|S|-3} (\prod_{i=1}^k
\Mod_{0,T_i}(\R))\ .$$ We deduce %from $(\ref{dimcond})$
  the
existence of a multiplication map for volume forms:
\begin{equation} \label{multforms}
f^*: \bigotimes_{i=1}^k H^{|T_i|-3} \big(\Mod_{0,T_i}(\R)\big) \To
H^{|S|-3} \big(\Mod_{0,S}(\R)\big)\ .
\end{equation}
This in turn gives a product formula for period integrals on the
spaces $\Mod_{0,S}(\R)$. If $S$ has dihedral structure $\delta$,
then it induces dihedral structures $\delta_i$ on $T_i$.
 Recall that the fundamental domains $\prod_{i=1}^k
X_{T_i,\delta_i}$ and $X_{S,\delta}$ are related by the set $G_f$
 defined in $(\ref{Gdef})$ via the formula $(\ref{domainrel})$.

\begin{cor} Let $\omega_i \in H^{|T_i|-3}(\Mod_{0,T_i}(\R))$, for $1\leq i\leq k$.
Then
$$\prod_{i=1}^k \int_{X_{T_i,\delta_i}} \omega_i = \sum_{\gamma \in G_f}
\int_{X_{S,\gamma}} f^*(\omega_1\otimes\ldots \otimes \omega_k) \
.$$
\end{cor}

\noindent It follows that a product of  period integrals on
  real moduli spaces is itself a period of real moduli spaces.

\begin{rem} \label{motiverem}
We know  that  $f$  extends to give a map $\Modf_{0,S} \To
\prod_{i=1}^k \Mod^{\delta_i}_{0,T_i}$.
 The previous corollary therefore implies the following multiplication
formula for the framed mixed Tate motives  defined in
$\S\ref{sectMHS}$:
\begin{equation}
 \bigotimes_{i=1}^k m_{T_i,\delta_i}(\omega_i) = \bigoplus_{\gamma \in G_f}
m_{S,\gamma} (f^*(\omega_1\otimes \ldots \otimes \omega_k))\
.\end{equation}
\end{rem}

 We can apply the product formula above to the set of multiple zeta
forms $\Omega(\underline{\epsilon})$ defined in $\S\ref{sect74}$.
 As in $\S\ref{sect55}$, let $X$ denote an
alphabet with two letters $\{x_0,x_1\}$. Let $w=x_{\epsilon_m}
\ldots x_{\epsilon_1}$ and $w'= x_{\epsilon_{\ell}}\ldots
x_{\epsilon_{m+1}}$ denote two words in $x_0X^*x_1$, where $\epsilon_i\in
\{0,1\}$,  such that $\epsilon_1=\epsilon_{m+1}=1$ and
$\epsilon_m=\epsilon_\ell=0$. Let us write
$\underline{\epsilon}=(\epsilon_1,\ldots, \epsilon_m)$ and
$\underline{\epsilon}'=(\epsilon_{m+1},\ldots, \epsilon_\ell)$.
Recall  the simplicial product map defined in
 $\S\ref{sect27}$:
$$m_{\triangle} : \Mod_{0,S} \To \Mod_{0,S_1} \times
\Mod_{0,S_2}\ ,$$  where $S_1= \{s_1,s_2,\ldots, s_{m+3}\}$ and
$S_2=\{s_1,s_2,s_3, s_{m+4},\ldots, s_n\}$. We deduce that
$$\int_{X_{S_1,\delta_1}} \Omega(\underline{\epsilon}) \int_{X_{S_2,\delta_2}}
\Omega(\underline{\epsilon}')= \sum_{\gamma\in
G_{m_\triangle}}\int_{ X_{S,\gamma}} m_{\triangle}^*
\big(\Omega(\underline{\epsilon}) \otimes
\Omega(\underline{\epsilon}')\big)\ .$$ Recall from $\S2.7$ that
$G_{m_\triangle}$ is the set  of
$(m,\ell-m)$-shuffles, and that,  in simplicial coordinates,
$X_{S,\delta}$ is the unit simplex. We therefore deduce that
$$\int_{0<t_1<\ldots<t_m<1} \bigwedge_{i=1}^m {dt_i\over \epsilon_i-t_i}  \times
\int_{0<t_{m+1}<\ldots<t_\ell<1} \bigwedge_{i=m+1}^\ell {dt_i\over
\epsilon_i-t_i}$$
$$ =
 \sum_{\sigma \in \Sym(m,\ell-m)}
\int_{0<t_{\sigma(1)}<\ldots <t_{\sigma(\ell)}<1}
\bigwedge_{i=1}^\ell {dt_i\over \epsilon_i-t_i}\ ,$$ which, by
$(\ref{Leibniz-Kontsevich})$,  gives the shuffle product formula:
\begin{equation}
\zeta(w) \, \zeta(w') = \sum_{\sigma \in \Sym(m,\ell-m)}
\zeta(x_{\sigma(\epsilon_\ell)}
x_{\sigma(\epsilon_{\ell-1})}\ldots x_{\sigma(\epsilon_1)})
=\zeta(w\sha w') \ .
\end{equation}
 Now let us see what happens in the case of
the cubical product map $(\ref{cubeprod})$:
$$m_{\Box} : \Mod_{0,S} \To \Mod_{0,S_1} \times
\Mod_{0,S_2}\ ,$$  where $S_1= \{s_2,s_3,\ldots, s_{m+4}\}$ and
$S_2=\{ s_{m+4},\ldots, s_n,s_1,s_2,s_3\}$.
 We deduce that
$$\int_{X_{S_1,\delta_1}} \Omega(\underline{\epsilon}) \int_{X_{S_2,\delta_2}}
\Omega(\underline{\epsilon}')= \int_{X_{S,\delta}} m_{\Box}^*
\big(\Omega(\underline{\epsilon}) \otimes
\Omega(\underline{\epsilon}')\big)\ ,$$ since in this case
$G_{m_\Box}$ is the single element $\{\delta\}$. In cubical coordinates, each fundamental
cell is a hypercube, and thus we obtain the formula:
\begin{equation} \label{provestuffle}
\int_{[0,1]^m} \Omega_c(\underline{\epsilon})
\int_{[0,1]^{\ell-m}} \Omega_c(\underline{\epsilon}')=
\int_{[0,1]^{\ell}} \Omega_c(\underline{\epsilon}) \,
\Omega_c(\underline{\epsilon}'\big)\ ,\end{equation} where
\begin{equation}
\Omega_c(\epsilon_1,\ldots, \epsilon_m)= \bigwedge_{i=1}^m
{d(x_i\ldots x_\ell) \over \epsilon_i - x_i\dots x_\ell}\ .
\end{equation} One can write the product $\Omega_c(\underline{\epsilon})
\, \Omega_c(\underline{\epsilon}'\big)$ as a sum of terms
$\Omega_c(\underline{\epsilon}'')$  either using an identity due
to Cartier (see [Zu]) or using a power series expansion due to
Goncharov ([Go3], lemma 9.6). We use the latter approach. Let
$\eta_i = (1,0,\ldots, 0)$ denote a $1$ followed by a  sequence
$i-1$ zeros. Then
$$\Omega_c(\epsilon_1,\ldots, \epsilon_m) =
\Omega_c(\eta_{n_1},\ldots, \eta_{n_r}) = \sum_{0\leq k_1<\ldots <k_r}
y_1^{k_1} \ldots y_r^{k_r} dx_1\ldots dx_m\ ,$$ where
$y_1=x_1\ldots x_{n_1}$, $y_2=x_{n_1+1}\ldots x_{n_2}$,\ldots ,
$y_r=x_{n_r+1}\ldots x_m$. Expanding
$\Omega_c(\underline{\epsilon}')$ in a similar way, we obtain
\begin{eqnarray}
\Omega_c(\underline{\epsilon})\Omega_c(\underline{\epsilon}')& =&
\sum_{0\leq k_1<\ldots <k_r} y_1^{k_1} \ldots y_r^{k_r}
\sum_{0\leq k_{r+1}<\ldots <k_t} y_{r+1}^{k_{r+1}} \ldots
y_{t}^{k_{t}}
 dx_1\ldots
dx_\ell \nonumber \\
& = & \sum_{\sigma \in \overline{\Sigma}_{m,\ell-m}}
\sum_{\sigma_*(k_1,\ldots, k_t)} y_1^{k_1} \ldots y_{t}^{k_{t}}
 dx_1\ldots
dx_\ell\nonumber = \sum_{\sigma \in \overline{\Sigma}_{m,\ell-m}}
\Omega_c(\sigma(\underline{\epsilon} , \underline{\epsilon'})) \
,\end{eqnarray} where $\overline{\Sigma}_{m, \ell-m}$ is the set
of stuffles in the stuffle product for quasi-symmetric power
series [Ho], and $\sigma_*(k_1,\ldots, k_t)$ is the corresponding
domain of summation. Substituting this into
$(\ref{provestuffle})$, we deduce the stuffle product formula (see
[W]):
\begin{equation}
\zeta(w) \zeta(w') = \zeta(w \star w')\ .
\end{equation}
 %Writing out each integral
%in explicit coordinates  we deduce the shuffle product formula
%$\zeta(w)\zeta(w')= \zeta(w\sha w')$, where
%$w=a_{\epsilon_\ell}\ldots a_{\epsilon_0}$ and
%$w=a_{\epsilon'_\ell}\ldots a_{\epsilon'_0}$. If one does the same
%calculation this time using the cubical product map $m_{\square}$,
%then one can deduce the stuffle product formula for multiple zeta
%values (Cartier - see [Zu2]).

%Now explain what happens on framed versions, and refer to
%Goncharov - say that the double shuffles are motivic. (NB not
%regularised versions yet!!).

\begin{rem}
This approach can be used to derive any number of elementary
products between mutiple zeta values. To make such a product explicit, one needs to
fix a rule for decomposing a product of $\ell$-forms into a sum of $\ell-$forms of a preferred type (for example,
$\Omega(\underline{\varepsilon})$ or $\Omega_c(\underline{\varepsilon})$). %$(\ref{omegasimplicKontdefn})$.
 The motivic origin of such a
product formula follows immediately from  remark $\ref{motiverem}$
above (in the case of the shuffle and stuffle formulae, this is
equivalent to an argument due to  Goncharov [Go3]). We see by
looking at the sets $G_{m_\Box}$ and $G_{m_{\triangle}}$, that the
shuffle and stuffle product formulae are extreme cases of a range
of intermediary product formulae, obtained by shuffling together
two subsets of $\{s_3,\ldots, s_\ell\}$ relative to $s_1=1$,
$s_2=\infty$, and $s_3=0$. Such modular products will be studied
elsewhere.
\end{rem}

\subsection{Product formulae for integrals of  generalised polylogarithms} More generally, we can apply the product formulae
to convergent iterated integrals of functions of arbitrary weight, rather than
just regular algebraic forms. If $f$ is a non-degenerate
coordinate system (\S\ref{sect27}),  there is a commutative diagram
$$\begin{array}{ccc}
 f^*:\, \bigotimes_{i=1}^k B(\Mod_{0,T_i}) & \To  & B(\Mod_{0,S}) \\
     {\downarrow}{\wr} &  & {\downarrow}{\wr} \\
 \qquad \bigotimes_{i=1}^k L^{v_i,\delta_i}(\Mod_{0,T_i}) & \To & L^{v,\delta}(\Mod_{0,S})\ , \\
\end{array}
$$
where $f^*$ is  a map of  differential graded algebras. The
vertices $v,v_1,\ldots, v_k$ are chosen such that
$f(v)=(v_1,\ldots, v_k) \in \prod_{i=1}^k
 \Mod^{\delta_i}_{0,T_i}$.
The vertical maps are given by canonical regularisation maps $\rho_{v,\delta}, \rho_{v_i,\delta_i}$ defined in $\S\ref{sect67}$.
Since each $B(\Mod_{0,T_i})$ is differentially simple,
and  since the map $f_{T_i}^*$ is non-zero for each $1\leq i\leq k$, it
follows that $f^*$ is injective. The horizontal map along the
bottom is given by composition and multiplication of functions:
$$(p_1,\ldots,p_k) \mapsto  p_1\circ f_{T_1} \times \ldots \times
p_k \circ f_{T_k}\ .$$ In the same way as $\S\ref{sect75}$, we deduce the
product formula:
\begin{equation} \label{productforpolylogs}
\prod_{i=1}^k \int_{\overline{X}_{T_i,\delta_i}} p_i =
\sum_{\gamma\in G_f} \int_{ X_{S,\gamma}} \prod_{i=1}^k p_i\circ
f_{T_i} \ , \end{equation} where we suppose that all integrals are
convergent. In this way, we can obtain  product formulae for
generalised period integrals (integrals of polylogarithms).

In $\S\ref{sect66}$ we defined an  action of the dihedral group of
symmetries on the space of functions $L^\delta_{\MZV} (\Mod_{0,S})$.
We therefore have the following formula for any function  $f\in L^{v,\delta}(\Mod_{0,S})$ such that
the integral converges:
\begin{equation}\label{dihedralonintegrals}
\int_{\overline{X}_{S,\delta}} f = \int_{\overline{X}_{S,\delta}}
\sigma^* f \qquad \hbox{ for all } \sigma \in D_{2n}\ .
\end{equation}
If we combine  dihedral symmetries with such product formulae, we
have  much freedom for manipulating  integrals of generalised
polylogarithms over $\overline{X}_{S,\delta}$. In particular, we
can replace a period integral over any  given face of the Stasheff
polytope $\overline{X}_{S,\delta}$ with one over a face of fixed
combinatorial type.

\begin{lem} \label{lemremoveproducts}
Let $F_0 \subset \partial \overline{X}_{S,\delta}$ denote a fixed
face of the Stasheff polytope $\overline{X}_{S,\delta}$
% of
%combinatorial type $F_0 \cong \overline{X}_{n-1,\delta'}$.
which corresponds to a short chord $\{i,i+1\}$ in the $n$-gon
$(S,\delta)$. Given any other face $F\subset
\partial \overline{X}_{S,\delta}$, and any  form $\omega\in  \Omega^{\ell-1}  L^{\delta|_F}_\MZV(F)$,
 there exists another form $\omega'\in
\Omega^{\ell-1} L^{\delta|_{F_0}}_\MZV(F_0)$ such that
$$\int_F  \omega = \int_{F_0} \omega'\ ,$$
where the weight of $\omega$ is less than or equal to the weight of $\omega'$.
\end{lem}

\begin{proof}
By using a  product  $(\ref{productforpolylogs})$
%in the case of a
%cubical product $f=m_{\square}$ ,
 we can replace the integral of
$\omega$ over $F\cong \overline{X}_{k,\delta_1}\times
\overline{X}_{n-1-k,\delta_2}$ with an integral over a face of
$\partial \overline{X}_{S,\delta}$ of combinatorial type
$\overline{X}_{n-1,\delta'}$. Since the group of dihedral
symmetries $D_{2n}$ acts transitively on the set of all such
faces, we can replace this with an integral over the face $F_0$ by
applying $(\ref{dihedralonintegrals})$.
\end{proof}
%\begin{lem}
%Any convergent period integral over a boundary component  of the
%Stasheff polytope $\overline{X}_{S,\delta}$ can be written as an
%integral over a face of combinatorial type $\Mod_{0,n-1}$. More
%precisely, let $S=\{s_1,\ldots,s_n\}$ with dihedral structure
%$\delta$, and consider the divisor at finite distance obtained by
%cutting along a chord $e$. This gives rise to a stable partition
%$T_1\cup T_2$  with the induced dihedral structures, where
%$|T_1|=k$ and $|T_2|=n-k$. Set $S'=\{s_1,\ldots, s_{n-1}\}$. Then
%$\overline{X}_{T_1\cup\{e\},\delta}\times\overline{X}_{T_2\cup\{e\},\delta}$
%and $\overline{X}_{S',\delta}$ are both faces of
%$\overline{X}_{S,\delta}$. For any $f_1 \in \Omega^{k-2}
%L(\widehat{\Mod}_{0,T_1\cup\{e\}})$, $f_2 \in \Omega^{n-k-2}
%L(\widehat{\Mod}_{0,T_2\cup\{e\}})$, there exists $f\in
%\Omega^{n-4} L(\widehat{\Mod}_{0,S'})$ such that
%$$\int_{\overline{X}_{T_1\cup\{e\},\delta}}f_1  \times \int_{\overline{X}_{T_2\cup\{e\},\delta}}
%f_2 = \int_{\overline{X}_{S,\delta}} f\ .$$ The weight of the form
%$f$ is the sum of the weights of $f_1$ and $f_2$.
%\end{lem}
%\begin{proof} The lemma follows on applying the cubical product
%map and using dihedral symmetry. Essentially, the left hand side
%can be written in cubical coordinates as a product of integrals
%over the unit cubes $[0,1]^{k-2}$ and $[0,1]^{n-k-2}$
%respectively. Multiplying the integrands gives an integral over
%$[0,1]^{n-4}$, which can be brought back to an integral over
%$\overline{X}_{S,\delta}$ using a suitable set of cubical
%coordinates.%
%\end{proof}

\subsection{Examples of period integrals in small dimensions}
%We give some examples of the  constructions defined above in low
%dimensions.
 First of all, consider the case
$\Modf_{0,S}(\R)$, where $S=\{s_1,\ldots, s_5\}$, with
the obvious dihedral structure which we denote $\delta$. We shall
work in  cubical coordinates $(x_1,x_2)$, which we write $(x,y)$.
 The set of chords $\chi_{S,\delta}$ is $\{13, 24, 35, 41, 52\}$, and  the dihedral coordinates
are
\begin{equation}\label{dihedralcoordsin5case}
u_{13}= 1-xy,\quad u_{24}=x,\quad u_{35}={1-x\over 1-xy}, \quad
u_{41}={1-y\over 1-xy}, \quad u_{52}=y\ .\end{equation}
%The map
%$(\ref{Embed})$ is the anticanonical embedding of a del Pezzo
%surface in $(\Pro^1)^5$.
 The domain $\overline{X}_{S,\delta}$ is bounded by
the five sets of equations $u_{ij}=0$, $u_{i+1\, j+1}=u_{i-1\,
j-1}=1$ for each pair $\{i,j\}\in \chi_{S,\delta}$, and these form
the sides of a pentagon whose interior is $X_{S,\delta}$. In all,
there are ten stable partitions of the set $\{s_1,\ldots, s_5\}$,
which means that there are another five divisors at infinity given by
the five equations $u_{ij}=\infty, u_{i-2\, j-2}=u_{i+2\, j+2}=0
$, where $\{i,j\}\in S$.  These too form a
pentagon. % which is disjoint from the first.

The volume form $\omega_{S,\delta}=d\log u_{13}\wedge d\log
u_{24}=(1-xy)^{-1} dxdy$, and  every real period integral on
$\overline{X}_{S,\delta}$ is a sum of integrals
$$
I_{\chi_{S,\delta}} (\alpha_{ij})  =  \int_{X_{S,\delta}}
u_{13}^{\alpha_{13}}\, u_{24}^{\alpha_{24}}\,
u_{35}^{\alpha_{35}}\, u_{41}^{\alpha_{41}}
\,u_{52}^{\alpha_{52}}\ \omega_{S,\delta}\ , $$ which in cubical
coordinates is just
$$I_5(h,i,j,k,l)   =  \int_{0}^1\int_{0}^1
{x^h (1-x)^i y^k (1-y)^j \over (1-xy)^{i+j-l}} {dxdy \over 1-xy} \
, $$
 where we have set $\alpha_{24}=h$, $\alpha_{35}=i$,
$\alpha_{41}=j$, $\alpha_{52}=k$, $\alpha_{13}=l$.
 This exactly coincides with
the family of integrals first defined  by Dixon, and studied by
Rhin and Viola  [Di, RV1]. The dihedral group $D_{10}$ preserves
the integral and permutes the indices $\{h,i,j,k,l\}$. It is
generated by a cyclic rotation of order five $\tau_5$, and a
reflection of order two $\sigma_5$, where
$$\tau_5(x,y) = \big( {1-xy}, {1-y\over 1-xy}\big)\ , \qquad
\sigma_5(x,y) = (y,x)\ .$$
%The integrand above has poles along
%each divisor at infinity of order $h+i-k+1 ,i+j-l+1, j+k-h+1,
%k+l-i+1, l+h-j+1$. These are the auxiliary variables considered in
%[RV2], and are also permuted by the group $D_{10}$ in the natural
%way.

\begin{rem}
By  combining the action of the dihedral symmetry group $D_{12}$
on $\Mod_{0,6}$ and the product formula for integrals on
$\Mod_{0,4}\times \Mod_{0,5}$ one can  deduce the `hypergeometric
transformation formula' for the integrals above. This remarkable
identity  was discovered by Dixon in 1905, and was exploited by
Rhin and Viola to obtain the best  irrationality measures for
$\zeta(2)$ known to date. It is:
 \begin{equation} \label{hyper} { 1 \over j!\, k!} I(h,i,j,k,l)= { 1 \over
(k+l-i)!(i+j-l)!}
 I(h,i,k+l-i, i+j-l,l)\ .
 \end{equation}

Before proving this identity, first observe that the  real period
integral $I_4$ on $\Mod_{0,4}$ is the following beta integral:
\begin{equation}
\label{Weight1} I_{4} (\alpha_{13},\alpha_{24} ) = \int_0^1
x^{\alpha_{24}} (1-x)^{\alpha_{13}} {dx } = {\alpha_{13}!\,
\alpha_{24}! \over (\alpha_{13}+ \alpha_{24}+1)!} \in \mathbb{Q}\
.
\end{equation}
%\\
Now consider the case $\Mod_{0,6}$. In cubical coordinates $(x,y,z)=(x_1,x_2,x_3)$ we have:
$$ u_{13} =1-xyz\ , \   u_{24}= x \ ,   \ u_{35} = {1\!-\!x\over 1\!-\!xy} \ , \  u_{46} = {(1\!-\!xyz)(1\!-\!y) \over (1\!-\!xy)(1\!-\!yz) }
 \ , \  u_{51}={1\!-\!z\over 1\!-\!yz} \ , \ u_{62} = z\ ,$$
 \begin{equation} \label{dihedralsformo6} u_{14} = {1-yz \over 1-xyz} \ , \quad u_{25} = y \ , \quad u_{36} = {1-xy \over 1-xyz} \ .
 \end{equation}
Let $h,i,j,k,l,m,r,s,t$ be any  nine non-negative integers.
 Let  $I_6(h,i,j,k,l,m\, ; \, r,s,t)$ denote the  period
integral of weight three on $\mathfrak{M}_{0,6}$ which is therefore given by
\begin{equation} \label{Weight3}
\int_0^1\! \!\int_0^1\! \!\int_0^1 {x^h (1-x)^i y^t (1-y)^j z^l
(1-z)^k \over (1\!-\!xy)^{i+j-r}
(1\!-\!yz)^{j+k-s}(1\!-\!xyz)^{r+s-j-m} }
{dxdydz\over(1\!-\!xy)(1\!-\!yz) } \ .
\end{equation}
 The dihedral group  of symmetries  $D_{12}$ for $\Mod_{0,6}$ is generated by a
 cyclic permutation we denote $\tau_6=(h\, i \,j \,k \,l \, m) (r\, s\, t)$, and the reflection $\sigma_6=(h\,
l)(i\,k)(r\,s)$. In the degenerate case where $r+s=j+m$ and
 $r=i+j+1$,  the terms $(1-xy)$ and $(1-xyz)$ vanish in the integrand and
$(\ref{Weight3})$ splits as a product $I_4(h,i)\, I_5(l,k,j,t,s)$.
This is precisely a cubical product. Similarly, if  the terms
$(1-yz)$ and $(1-xyz)$  vanish, then we obtain a different
splitting of the integral. \\

\emph{ Proof of} $(\ref{hyper})$. Let $h,i,j,k,l$ be non-negative
integers. We can assume without loss of generality that $i<l$. Set
$\alpha=k+l-i$, and $\beta=i-l-1$. Consider
$$
I_4(\alpha,\beta)\, I_5(h,i,j,k,l) = I_6(h,i,j,\beta,\alpha,
l+\beta+1\,; \,l, j+\beta+1,k)\ .$$ We apply the cyclic
permutation $\tau^4_6$ and use the fact that $\alpha+\beta+1=k$,
$\alpha+i=k+l$ to obtain a different splitting. %(the terms $(1-yz)$ and
%$(1-xyz)$ disappear).
 This gives
$$I_6(j,\beta,\alpha,
l+\beta+1,h,i\,; \,k,l, j+\beta+1)=
I_4(j,\beta)\,
I_5(h,l+\beta+1,\alpha,j+\beta+1,l)
\ .
$$
Replacing the $I_4$ terms with factorials using $(\ref{Weight1})$,
we obtain the identity $(\ref{hyper})$.
\\

These examples illustrate how many important
identities between multiple zeta values and Euler integrals can
be proved by simple geometric considerations on moduli spaces
$\Mod_{0,n}$. Kontsevich and Zagier have made the very general and
ambitious  conjecture that every identity between periods can be
proved using three elementary operations on integrals: changes of
variables, the linearity of integration, and Stokes' theorem.
 In
our situation, we have %what we would like to call a finite
%Kontsevich-Zagier formalism. This is
 an infinite family of period
integrals, but we have a fixed  set of algebraic  operations which we
can perform on these integrals ({\it e.g.}, the action
of dihedral symmetries and the multiplication rules we defined above).
It would be interesting to see which of the many known identities between multiple zeta values  can be proved  %(for example, the rationality of $\zeta(1-2n)$)
 by using
just these operations.

%One could formalise this notion in terms of
%the action of a certain pro-group on an algebra of symbols. Note
%that one can also
% bring to bear changes of variables coming from
%systems of vertex coordinates ({\it e.g.}, coming from the
%butterfly lemma). This yields  many non-trivial birational
%transformations which preserve the  domain of integration
%$X_{S,\delta}$, and give rise to non-dihedral changes of variables
%such as those considered by Rhin and Viola [RV2, Br3].

%
%Example shows that we have a Kontsevich-Zagier formalism. Group of
%changes of variables. Products, Stokes. Deduce very non-trivial
%identities.  Interesting to study systematically. All motivic -
%rephrase using my motives. Hypergeometric is of modular/motivic
%origin!!

\end{rem}

\newpage

\section{Calculation of the periods of $\Mod_{0,n}$.}

We prove that the integral of a convergent algebraic $\ell$-form
%on $\Mod_{0,S}$
over an associahedron $X_{S,\delta}$ can be written as a linear combination of multiple
zeta values of weight at most $\ell$.
The key to the argument is
the interplay  between logarithmic singularities (which are
permitted), and polar singularities (which are forbidden), along the
boundaries of the associahedron $\overline{X}_{S,\delta}$.

\subsection{Pole-free primitives}
%The key argument in our proof of the main theorem is to apply a
%version  of Stokes' theorem to the manifold with corners
%$\overline{X}_{S,\delta}$. By repeatedly taking primitives, we
%obtain a cascade of integrals of smaller and smaller dimension.

In order to apply Stokes' theorem to the manifold with corners $\overline{X}_{S,\delta}$, we need to
verify that the algebra of generalised polylogarithms on $\Mod_{0,S}^\delta$ satisfies the required properties.

First of all,  it follows from the regularization results of $\S\ref{sect67}$ %immediately
that
 the
coefficients of the generating series of generalized
polylogarithms $L_{v,\delta}(z)$ have at most logarithmic
 singularities along the boundary of the Stasheff polytope
$\overline{X}_{S,\delta}$.  This implies that
\begin{equation}
L^\delta_{\MZV}(\Mod_{0,S}) \subset \Gamma(\overline{X}_{S,\delta},
\Fa_p^{\log})\ .
\end{equation}
Theorem \ref{thmstructureofL} tells us that primitives exist in $L^\delta_{\MZV}(\Mod_{0,S})$.
One difficulty, however, is that primitives of $n$-forms on a
manifold of dimension $n$ are not unique, and we may inadvertently
introduce extra poles, which would give rise to divergent integrals.  We show how to remove these extra poles
below. In order to do this, we define
\begin{equation} \label{lplusdef}
L^{\delta,+}_{\MZV}(\Mod_{0,S})= L^\delta_{\MZV}(\Mod_{0,S}) \cap \Gamma
(\overline{X}_{S,\delta}, \Fa^{\log}) \end{equation} to be the
sub-algebra
 of polylogarithms on $\Mod_{0,S}$ which have at
most logarithmic singularities on the boundary faces of the
Stasheff polytope $\overline{X}_{S,\delta}$. %This is not a differential algebra.
Observe that every generalised polylogarithm
 has a canonical branch on $X_{S,\delta}$, and is therefore  a well-defined, real-valued function.
% Since
%$\overline{X}_{S,\delta}$ is compact,  we know by lemma $\ref{lemmanopolesint}$
% that any function $f \in \Gamma(\overline{X}_{S,\delta},
%\Fa^{\log})$ is integrable on $\overline{X}_{S,\delta}$.

\begin{prop} \label{prop81}
Let $f\in W^k \Omega^\ell  L^{\delta,+}_{\MZV}(\Mod_{0,S})$. There exists a pole-free primitive
$$F
\in W^{k+1}\, \big( \Omega^{\ell-1} L_{\MZV}^{\delta,+}(\Mod_{0,S})\big)
$$ such that
 $dF =f$, which implies that the restriction of $F$ to $\partial \overline{X}_{S,\delta}$ is
continuous. In other words, the conditions of theorem $\ref{corStokesform}$  hold.
\end{prop}

\begin{proof}
We know from theorem \ref{thmstructureofL} that
$L^\delta_\MZV(\Mod_{0,S})$ has trivial de Rham cohomology. It follows
that we can find a primitive  $G\in  \Omega^{\ell-1}
L^\delta_\MZV(\Mod_{0,S})$ for $f$, of weight at most $k+1$, which may have polar
singularities along $\partial \overline{X}_{S,\delta}$.

%We know from corollary \ref{corLiscompletelyreg} that $L^{\delta}_\MZV(\Mod_{0,S})$ is completely regularisable
In order to remove spurious poles in $G$, we work on a single
chart  of $\overline{X}_{S,\delta}$ at a time. Therefore, let
$e\in \chi^q_{S,\delta}$ denote a partial decomposition of the
$n$-gon $(S,\delta)$, and let $\alpha\in \chi^\ell_{S,\delta}$ be
a full triangulation which contains $e$. For every small
$\varepsilon>0$, recall that there is a chart $U_e(\varepsilon)$
(see  $(\ref{Uchartdefn})$), which has  local (vertex) coordinates
$x_1^{\alpha},\ldots, x_\ell^{\alpha}$, which are canonical up to
permutations. Recall that $U_e(\varepsilon)\cong U_{p,q}$ where
$p+q=\ell$ (\S4). It contains the face $F_e=\{u_{ij}=0: \{i,j\}\in
e\}$, and we can assume, by reordering the coordinates if
necessary, that $F_e=\{x_1^{\alpha}=\ldots = x_q^{\alpha}=0\}.$ We
remove polar singularities with respect to each coordinate
$x_1^\alpha$, \ldots, $x_q^\alpha$ in turn.
 First, there is a decomposition $G=G_p+ G'$, where $G'$ has at
 most logarithmic singularities in $x_1^\alpha$,  and $G_p$ is the divergent part of $G$ along $\{x_1^\alpha=0\}$:
$$G_p =
\sum_{a\geq 0, b\geq 1}
 {\log^{a} x_1^{\alpha}  \over (x_1^{\alpha})^{b} } \,
 g_{a,b}(x_2^{\alpha},\ldots, x_\ell^{\alpha})\, \omega_{a,b} + \sum_{c\geq 1} \log^c x^\alpha_1\, h_c(x^\alpha_2,\ldots, x^\alpha_\ell) dx^\alpha_2\ldots dx^\alpha_\ell\ ,$$
where $g_{a,b}(x_2^{\alpha},\ldots, x_\ell^{\alpha}) ,
h_c(x^\alpha_2,\ldots, x^\alpha_\ell) \in
\,\Fo^{\log}_p(U_{p,q-1})$,  and where $\omega_{a,b} $ are any
$\ell-1$ forms  $ \sum_i a_i\, dx^\alpha_1\ldots
\widehat{dx^\alpha_i} \ldots dx^\alpha_\ell$, where $a_i\in \R$.
By differentiating this expression, and using the fact that $dG=f$
has no poles, it is easy to verify that $dG_p=0$ (in other words,
poles can only get worse on differentiating). Therefore $dG'= f$,
and so $G'$ is a primitive of $f$ which has no poles along
$x_1^\alpha=0$. Using the fact that $L^\delta_\MZV(\Mod_{0,S})$ is
closed under differentiation with respect to $x_i^{\alpha}$, and
closed under taking regularised limits at $x_i^{\alpha}=0$, for
$1\leq i\leq \ell$ (theorem $\ref{proprestrictionoffunctions}$),
one can easily check that $G'$ lies in
$L^\delta_\MZV(\Mod_{0,S})$, {\it i.e.}, $G'$ is still a
generalised polylogarithm. Repeating this argument for
$x_1^\alpha$,\ldots, $x_q^\alpha$, in turn, we obtain  a primitive
of $f$ with no poles on the local chart $U_2(\varepsilon)$. The
whole argument can then be repeated on each local chart of
$X_{S,\delta}$, and we end up with a primitive $F$ of $f$ which
has no poles anywhere along $\partial \overline{X}_{S,\delta}$.
This is because, whenever we remove a polar singularity along the
divisor $D_{ij}$, $\{i,j\}\in \chi_{S,\delta}$, no new poles are
created along any other boundary component $D_{kl}$, where
$\{k,l\}\in \chi_{S,\delta}$, and the total weight is  not
increased. This proves the proposition. The fact that the
restriction of $F$ to each component of the boundary is continuous
follows from lemma $\ref{lemnopolesextendscont}$.
\end{proof}

In $\S\ref{sect83}$ we show how to construct canonical primitives which are automatically free of poles along $\partial \overline{X}_{S,\delta}.$

\subsection{Proof of the main theorem}% We can now  prove the following
%theorem.
\label{sect82}
Let $S$ denote a set of order $n=\ell+3$ with a fixed dihedral structure $\delta$.

%\begin{defn} Let $f\in W^{b} \Omega^\ell L(\Mod_{0,S})$.
%Then $f$ is holomorphic on the interior of
%$\overline{X}_{S,\delta}$. We say that $f$ converges if it has no
%poles along the boundary of $\overline{X}_{S,\delta}$. In that
%case the period integral is defined to be
%$$\int_{\overline{X}_{S,\delta}} f \ .$$
%\end{defn}

\begin{thm}\label{thmmain}
For all sets of indices $\alpha_{ij}\geq 0$,
$$I_{S,\delta}(\alpha_{ij})=\int_{\overline{X}_{S,\delta}} \prod_{\{i,j\}\in \chi_{S,\delta}} u_{ij}^{\alpha_{ij}}
\omega_{S,\delta} \in W^{\ell} \MZV\ .$$
\end{thm}

\begin{proof}
The proof is by induction and by repeated application of Stokes'
theorem (theorem \ref{corStokesform}).
% Since
%$\overline{X}_{S,\delta}$ is compact,  we know by lemma $\ref{lemmanopolesint}$
% that any function $f \in \Gamma(\overline{X}_{S,\delta},
%\Fa^{\log})$ is integrable on $\overline{X}_{S,\delta}$.
 We write $S=S_n=\{s_1,\ldots, s_n\}$. First observe that
the regular $\ell$-form
$$f_0 =\prod_{\{i,j\}\in \chi_{S,\delta}} u_{ij}^{\alpha_{ij}}
\omega_{S,\delta} \in W^0( \Omega^{\ell} L^\delta_{\MZV}(\Mod_{0,S}))
$$ has no poles %, and is strictly positive
on the compact set $\overline{X}_{S,\delta}$. Let us write $$S_i =
\{s_1,\ldots, s_{i}\}\quad \hbox{for all }\quad 3\leq i\leq n\ ,$$
and let $\delta_b$ denote the dihedral structure  on $S_b\subset
S$  induced by $\delta$. This is equivalent to choosing a nested
sequence of sub-faces of $\overline{X}_{S,\delta}$ in its
stratification.
 Let $0\leq b\leq \ell$, and suppose by induction
that %we have proved that
there exists an $\ell-b$ form
$$f_b \in W^b( \Omega^{\ell-b}
L^{\delta,+}_{\MZV}(\Mod_{0,S_{n-b}}))\ ,$$
  which
 has no poles %and is positive
 on $\overline{X}_{S_{n-b},
\delta_{n-b}}$, such that
$$I_{S,\delta}(\alpha_{ij}) = \int_{\overline{X}_{S_{n-b}, \delta_{n-b}}} f_b\ .$$
By  proposition \ref{prop81}, there exists  a primitive $P
 \in \Omega^{\ell-b-1} L_{\MZV}^{\delta,+}(\Mod_{0,S_{n-b}})$  of weight at most $b+1$, which
has no poles in $\overline{X}_{S_{n-b},\delta_{n-b}}$,
 and is continuous on the
 interior of
 $\partial \overline{X}_{S_{n-b},\delta_{n-b}}.$
  By the version of Stokes' formula
 stated in  theorem $\ref{corStokesform}$,
$$ I_{S,\delta}(\alpha_{ij}) = \int_{\overline{X}_{S_{n-b},
\delta_{n-b}}} f_b =  \int_{\partial \overline{X}_{S_{n-b}, \delta_{n-b}}} P
\ .$$ By the geometry of the Stasheff polytopes $(\S\ref{sect22})$, we know that
$$\partial \overline{X}_{S_{n-b}, \delta_{n-b}} = \bigcup_{\{i,j\}\in \chi_{S_{n-b},\delta_{n-b}}}
F_{ij}\ ,$$ where $F_{ij}=F_{ij}(\overline{X}_{S_{n-b},
\delta_{n-b}})$ is the face corresponding to the chord $\{i,j\}$,
 and therefore
\begin{equation}\label{bpf1} I_{S,\delta}(\alpha_{ij}) =
\sum_{\{i,j\}\in \chi_{S_{n-b},\delta_{n-b}}}  \int_{F_{ij}}
P\big|_{F_{ij}} \ .\end{equation} Given a chord $\{i,j\} \in
\chi_{S_{n-b},\delta_{n-b}}$, there exists a partition of
$S_{n-b}=T_1\cup  T_2$ such that $F_{ij}\cong
\overline{X}_{T_1\cup\{e\},\delta_1}\times \overline{X}_{T_2\cup\{e\},\delta_2}$, where $e$ corresponds to the chord $\{i,j\}$
 (equation $(\ref{Facestruc})$) and $\delta_1,\delta_2$ are the induced dihedral structures.
 By theorem
$\ref{proprestrictionoffunctions}$, we have
$$P\big|_{F_{ij}} \in  W^{b+1} \Omega^{\ell-b-1}(
L^{\delta_1}_{\MZV}(\Mod_{0,T_1\cup\{e\}} )\otimes
L^{\delta_2}_{\MZV}(\Mod_{0,T_2\cup\{e\}} ))\ . $$ By lemma
$\ref{lemremoveproducts}$, there exists $g_{ij} \in W^{b+1}(
\Omega^{\ell-b-1} L^\delta_{\MZV}(\Mod_{0,S_{n-b-1}} ))$ such that
$$
 \int_{\overline{X}_{T_1\cup\{e\},\delta_1}\times
\overline{X}_{T_2\cup\{e\},\delta_2}} P\big|_{F_{ij}} =
\int_{\overline{X}_{S_{n-b-1},\delta_{n-b-1}}} g_{ij}\ .$$ Thus
each integral in the sum $(\ref{bpf1})$ can be written as an
integral over the fixed face $\overline{X}_{S_{n-b-1},
\delta_{n-b-1}}$ by applying product formulae and using dihedral
symmetries. Since $P|_{F_{ij}}$ is continuous with at most
logarithmic singularities along $\partial F_{ij}$, it follows that
$g_{ij}\in \Omega^{\ell-b-1}(
L^{\delta,+}_{\MZV}(\Mod_{0,S_{n-b-1}}))$. Taking the sum over all
$\{i,j\}\in \chi_{S_{n-b-1},\delta_{n-b-1}}$ in $(\ref{bpf1})$, we
obtain a form $f_{b+1}\in W^{b+1}( \Omega^{\ell-b-1}
L_{\MZV}^{\delta,+}(\Mod_{0,S_{n-b-1}}))$ such that
$$I_{S,\delta}(\alpha_{ij})= \int_{\overline{X}_{S_{n-b-1}, \delta_{n-b-1}}} f_{b+1} \ .$$
This completes the induction step. At the final stage of the
induction, we deduce that $I_{S,\delta}$ is given by evaluating a
multiple polylogarithm in one variable in $W^\ell
L_\MZV^\delta(\Mod_{0,4})$ at a  single point. We conclude (see
$\S\ref{sect55}$) that $I_{S,\delta}(\alpha_{ij})\in W^{\ell}
\MZV\ .$
\end{proof}

\noindent Note that it is not strictly necessary in the course of
the above proof to use the product formula (lemma  $\ref{lemremoveproducts}$).
This  replaces the sum of a product of integrals with a single
integral at each stage, and only  serves to simplify notations.
Lemma $\ref{lemrealperiodsareIs}$ implies the following result.

\begin{cor} Every relative period integral over a union of cells $X_{S,\delta_i}$  can be written as a linear
 combination of multiple zeta values of weight at most $\dim  \Mod_{0,S}(\R)$.
\end{cor}

\subsection{Canonical primitives - an algorithmic approach} \label{sect83}

%The inductive nature of the proof can be simplified to a certain
%extent. At each stage of the integration process we obtain a sum
%of integrals over each face of the Stasheff polytope. Using the
%existence of product structures and the dihedral symmetry, we can
%replace these integrals with a single integral:

%As a result of the lemma, we can always

The existence of primitives uses the fact that $\Mod_{0,S}$ is a
fiber-type hyperplane arrangement. By exploiting the
hyperlogarithm fibration, we can find canonical primitives, as in
remark $\ref{remarkonstrongprimitives}$, which have no spurious
poles.
% In order to do this we
%need to speak of the positivity of an arbitrary real $\ell$-form
%on $\Mod_{0,S}$, which makes sense because
% we have defined a canonical volume form on $\Mod_{0,S}(\R)$.
This gives rise to a simplified series of integrals occurring in
the proof of  theorem $\ref{thmmain}$ above, and yields an effective algorithm
for computing period integrals on $\Mod_{0,S}(\R)$ algebraically.

Let $f\in W^{b} \Omega^{\ell} L_{\MZV}^{\delta,+}(\Mod_{0,S}) $.
Working in cubical coordinates, we can write
$f=g(x_1,\ldots,x_\ell)\, dx_1\ldots dx_\ell$, where $g\in
L^{\delta,+}_\MZV(\Mod_{0,S})$ is of weight at most $b$. Recall
that by remark $\ref{remarkonstrongprimitives}$, there exists a
primitive $F\in W^{b+1} \Omega^{\ell-1} L_\MZV^\delta(\Mod_{0,S})$
such that
$$F= G(x_1,\ldots, x_\ell) \,dx_1\ldots dx_{\ell-1}\ ,$$
where $\partial G/\partial x_\ell =g$. More concretely, let
 $S=\{s_1,\ldots, s_n\}$ and let
$S'=\{s_2,\ldots, s_{n}\}$. Recall from $\S\ref{sect67}$ that the
hyperlogarithm fibration given by projection onto $x_\ell=0$:
$$ \Mod_{0,S} \To \Mod_{0,S'}\ ,$$
gives rise to a
 decomposition of filtered algebras
 $$L^\delta_{\MZV}(\Mod_{0,S})\cong L_{\Mod_{0,S'}}(\Pro^1\backslash\Sigma)\otimes_{\Or(\Mod_{0,S'})}
 L^{\delta'}_{\MZV}(\Mod_{0,S'})\ ,$$
 where  $\Sigma=\{0,1, x_{\ell-1}^{-1}, \ldots, (x_1\ldots x_{\ell-1})^{-1}\}$ as in $\S\ref{sect54}$, and $\delta'$ is
 the induced dihedral structure on $S'$.
 We can therefore write the function $g$ as a  finite sum of products
 $$g(x_1,\ldots, x_\ell) =\sum_i a_i(x_\ell) \, b_i(x_1,\ldots, x_{\ell-1})\ ,$$
 where each $b_i\in L^{\delta'}_\MZV(\Mod_{0,S'})$ is a function of $\ell-1$ variables $x_1,\ldots, x_{\ell-1}$ only,
and each $a_i$, considered as a function of the single variable
 $x_\ell$, is a hyperlogarithm with singularities in $\Sigma\cup
 \infty$. We can
 assume that the $a_i$ are linearly independent.
The weight of each product of $a_i(x_\ell)$ and $b_i(x_1,\ldots,x_{\ell-1})$ is at most
 $b$.
  Now by proposition $\ref{propRelnocohomology}$, each function $a_i(x_\ell)$
 has a primitive (with respect to the variable $x_\ell$) which we denote
 $$A_i (x_\ell) \in   L_{\Mod_{0,S'}}(\Pro^1\backslash\Sigma)\ ,$$
 which is of weight at most one more than the weight of $a_i(x_\ell)$. We can choose the constant of
 integration in such a way that $A_i(x_\ell)$ either vanishes at $0$,
 or is $\log^k(x_\ell)$ for some $k\geq 1$ (see $\S\ref{sect52}$). In the latter case, $a_i(x_\ell)$, and
 hence $g(x_1,\ldots, x_\ell)$ would have a
 pole at the origin, so this cannot occur (this is precisely the argument in the proof of lemma $\ref{lemnopolesextendscont}$). It follows that the
 function
 $$F=\sum_i A_i(x_\ell)\, b_i(x_1,\ldots, x_{\ell-1})\,dx_1\ldots dx_{\ell-1} $$
%
%There is a decomposition of $G$ as a product
%$G(x_1,\ldots, x_\ell)=\sum_i a_i(x_1)\, b_i(x_2,\ldots, x_\ell)$, where $a_i$ is a hyperlogarithm in the variable $x_1$.
%We can choose .
%
is a primitive of $f$, and is identically zero on all faces of
$\overline{X}_{S,\delta}$ except the single face given by
$x_\ell=1$. The primitive $F$ has no poles since this would
contradict the convergence of the integral by lemma
$\ref{lemmanopolesint}$. It is therefore continuous on the
interior of this face, and we can apply Stokes' theorem directly.
This approach to the induction step in the proof of the main
theorem has the advantage that it does not involve any
regularisation, or having to apply a product formula (lemma
$\ref{lemremoveproducts}$).

% The previous proposition actually proves that there
%exist primitives which are non-zero on one face of the Stasheff
%polytope only, so there is no need to apply any product formulae
%at any stage of the integration process.

%-----REMOVE TO SECTION ON BAR CONSTRUCTION FOR MODULI SPACES AND PRODUCTS----------
%--------------------------------------------------------------------------------
\subsection{Taylor expansions of Selberg integrals and multi-beta functions}
The method of  proof of theorem $\ref{thmmain}$ works much more
generally, and enables us to compute  integrals of arbitrary
generalised polylogarithms on $\Mod_{0,S}$, which are allowed
logarithmic singularities along the boundary of the domain of
integration.

\begin{thm}
Let $f\in W^k  L_\MZV^{\delta,+}(\Mod_{0,S})$ denote a generalised polylogarithm on $\Mod_{0,S}$ of weight at most $k$, which has
no poles along $\partial \overline{X}_{S,\delta}$. Then
$$I(f)=\int_{\overline{X}_{S,\delta}} f \, \omega_{S,\delta} \in W^{\ell+k} \MZV\ .$$
\end{thm}
The proof is identical to the proof of theorem $\ref{thmmain}$.
Note that the integrand is always well-defined on the real domain
$X_{S,\delta}$ at each stage of the induction. If we apply this
theorem in the case where $f$ is of weight at most one, {\it
i.e.,} a $\Or(\Mod_{0,S})$-linear combination of logarithms, then
we deduce the following corollary.
\begin{cor}
Let $\{s_{ij}\}$ denote a set of complex parameters.  It follows from the calculations in $\S4$ that
the following integral, viewed as a  function of the variables $s_{ij}$, is holomorphic in the region $\Real\, s_{ij}>-1$:
\begin{equation}\label{multibeta}
\beta_{S,\delta}(\{s_{ij}\})=\int_{\overline{X}_{S,\delta}}
\prod_{\{i,j\}\in \chi_{S,\delta}} u_{ij}^{s_{ij}}
\omega_{S,\delta} \ .\end{equation} The coefficients of its Taylor
expansion (with respect to the variables $s_{ij}$) at any integral
point  $s_{ij} \in \Z$, where $s_{ij}\geq 0$ for all $\{i,j\} \in
\chi_{S,\delta}$,
 are multiple zeta values.
\end{cor}
\noindent Similar kinds of results have been obtained by Terasoma
 [Ter] in certain  cases. The integral $(\ref{multibeta})$
defines a multi-beta function, since in the case $\Mod_{0,4}$ it
reduces to the ordinary beta function. It satisfies many
functional identities coming from the dihedral relations
$(\ref{id})$, the product maps $(\ref{productforpolylogs})$,
 and also the action of the dihedral symmetry group, and would merit further study.
%In the simplest case $\Mod_{0,4}$, this reproves the well-known fact that the coefficients in the Taylor expansion of the beta
%integral
%$$I_{S,\delta}(s_{13}, s_{24}) = \int_{0}^1 (1-x)^{s_{13}} x^{s_{24}} dx $$

%----------------------------------------SECTION ON RESIDUE FORMULAE---------------------------------
%----------------------------sketched calculation of complex periods --------------------------------
%----------------------------------------------------------------------------------------------------
%----------------------------------------------------------------------------------------------------

\subsection{Computation of all relative periods of the moduli spaces $\Mod_{0,S}$.}

Let $A, B$ denote two sets of divisors  in
$\overline{\Mod}_{0,S}\backslash \Mod_{0,S}$ which do not share
any irreducible components ($\S\ref{sectMHS}$). We sketch a proof
of the following result.

\begin{thm}
The periods of $H^\ell(\overline{\Mod}_{0,S}\backslash A, B
\backslash B\cap A)$ are $\Q$-linear combinations of multiple zeta
values and the constant $2i\pi$, of total weight at most $\ell$.
\end{thm}

 Let $\Delta_B\subset
\Mod_{0,S}(\C)$ denote any real smooth compact submanifold with
corners of dimension $\ell$, whose boundary is contained in the
set of complex points of $B$. Let $\omega\in
\Omega^\ell(\overline{\Mod}_{0,S}\backslash A)$. We can assume
that $\Delta_B$ is disjoint from $A$, and that $\Delta_B$ is
stratified by $B$ according to $(\ref{Bstrat})$. Let $A$ be a
union of distinct divisors $A_i$, for $1\leq i\leq N$.

By decomposing the relative homology class $[\Delta_B] \in
H_\ell(\overline{\Mod}_{0,S}\backslash A, B \backslash B\cap A)
 $ into different pieces, we can first consider the case when $\Delta_B$ does not wind
non-trivially around any component $A_i$ of $A$. In this case,
write $X=\Delta_B$, and observe that the argument given in
$\S\ref{sect82}$ goes through as before. In other words, we can
take primitives in the algebra of polylogarithms
$L^\delta_\MZV(\Mod_{0,S})$, and repeatedly apply Stokes' formula
(theorem $\ref{corStokesform}$) to the manifold with corners $X$,
and proceed by induction. Note that, although $X$ is not
necessarily  simply-connected, the argument goes through as long
as the functions we integrate remain single-valued along $X$.
Since $X$ does not wind around $A$, we can always ensure that this
is the case, by taking primitives whose singularities are
contained in $A$. This proves that
$$\int_X \omega \in W^\ell \MZV\ .$$

In the case when $\Delta_B$ winds around some component of $A$, we
apply a residue formula and induction. To make this precise, let
$A_i^c = \bigcup_{j\neq i} A_j$ for all $1\leq i\leq \ell$, and
consider the residue map
$$ H^\ell(\overline{\Mod}_{0,S}\backslash A, B\backslash B\cap A) \To \bigoplus_{i=1}^N H^{\ell-1} (A_i \backslash (A_i\cap A_i^c), (B\cap A_i) \backslash (B\cap A_i\cap A_i^c))\ ,$$
and its dual map
$$H_\ell(\overline{\Mod}_{0,S}\backslash A, B\backslash B\cap A) \longleftarrow \bigoplus_{i=1}^N H_{\ell-1} (A_i \backslash (A_i\cap A_i^c), (B\cap A_i) \backslash (B\cap A_i\cap A_i^c))\ .$$
Suppose that $[\Delta_B]\in H_\ell(\overline{\Mod}_{0,S}\backslash
A, B\backslash B\cap A)$ is  the image of a class $[Y]\in
H_{\ell-1} (A_i \backslash (A_i\cap A_i^c), (B\cap A_i) \backslash
(B\cap A_i\cap A_i^c))$, for some $1\leq i\leq N$, where  $Y
\subset A_i \backslash (A_i\cap A_i^c)$ is a smooth compact
submanifold  with corners of dimension $\ell-1$. Therefore
 $\Delta_B$ resembles a
narrow tube around $A_i$. By  taking the residue along $A_i$, we
get:
$$\int_{\Delta_B} \omega = 2i\pi \int_{Y} \mathrm{Res}\,\, \omega\Big|_{A_i}\ .$$
The corresponding period is therefore $2\pi i$ times a period of
$H_{\ell-1} (A_i \backslash (A_i\cap A_i^c), (B\cap A_i)
\backslash (B\cap A_i\cap A_i^c)).$ Since $A_i$ is itself
isomorphic to a product of moduli spaces,
 we can repeat the argument inductively. We conclude
that, in all cases,
$$\int_{\Delta_B} \omega \in W^\ell \MZV[2i\pi]\ ,$$
where $\MZV[2i\pi]$ has the natural filtration  which gives
$2i\pi$
 weight $1$.

%----------------------------------------------------------------------------------------------------
%----------------------------------------------------------------------------------------------------
%----------------------------------------------------------------------------------------------------

\subsection{Some simple examples}
%We illustrate with some simple examples.
In the following examples, it is  convenient to work in cubical coordinates $x_1,\ldots, x_\ell$. At each stage
we take canonical primitives (as described in $\S\ref{sect83}$) with respect to $x_1 $ or $x_\ell$. This is because
the projection maps onto $x_1=0$ or $x_\ell=0$ are fibrations ($\S\ref{sect25}$), and so we can use the method of partial fractions  to find
primitives. % (by %($\ref{finbase}$) and ($\ref{partialfractionidentity})$).
 At each stage, one can re-confirm (using dihedral coordinates) that the primitives have no poles along
the boundary of the domain of integration
$\overline{X}_{S,\delta}$. First,  assume $|S|=5$, that is
$S=\{s_1,\ldots, s_5\}$. We compute
$$I_1=\int_{\overline{X}_{S,\delta}} \omega_{S,\delta} =\int_{0}^1\int_{0}^1 {dx dy\over 1-xy}\ .$$
Following  $\S \ref{sect83}$, we take the primitive of
$\omega_{S,\delta}$ with respect to the variable $y$. This is $$F=
-\log(1-xy) \,{dx \over x} \ , $$ which vanishes at $y=0$ as
required. In dihedral coordinates $(\ref{dihedralcoordsin5case})$,
this is  $F=-\log u_{13} \,d\log u_{24} $, which has no poles at
finite distance. Then
$$ \int_{\overline{X}_{S,\delta}} \omega_{S,\delta} = \sum_{\{i,j\}\in \chi_{S,\delta}}  \int_{D_{ij}} -\log u_{13} \,d\log u_{24}\ . $$
The only face on which the form does not vanish is the face
$D_{14}$ which is defined by $u_{14}=0$, $u_{25}=u_{35}=1$, which implies that $u_{13}=1-u_{24}$ by $(\ref{id})$. We obtain
$$ \int_{\overline{X}_{S,\delta}} \omega_{S,\delta}= \int_0^1 -\log(1-x)
{dx\over x}\ .$$ Notice that the form  $  \log(1-x) \, dx/x$ is
continuous on the  interval $[0,1)$ but  has a logarithmic
singularity at $x=1$. It has a unique primitive which vanishes at
$0$, namely $\Li_2(x)$, which is now bounded at $x=1$ by lemma
$\ref{lemnopolesextendscont}$. We conclude that
\begin{equation} \label{firstexample} I_1=\int_{\overline{X}_{S,\delta}} \omega_{S,\delta}= \Big[\Li_2(x)\Big]_0^1=\zeta(2) \ .
\end{equation}

Now let $|S|=6$. Consider the following integral on $\Mod_{0,6}$:
$$I_2=\int_{0\leq t_1\leq t_2\leq t_3\leq 1} {dt_1 \over 1-t_1 } {dt_2 \over t_2} {dt_3 \over t_3-t_1} =
\int_{[0,1]^3} {dx \,dy\, dz \over (1-xyz)(1-xy)} =
\int_{\overline{X}_{S,\delta}} u_{14} \,\omega_{S,\delta}\ .$$ The
last formula shows that $I$  converges. Working in cubical
coordinates, we have
$$I_2=\int_{[0,1]^2} \Big[ {-\log(1-xyz) \over xy} \Big]_0^1 \, {dx \, dy \over 1-xy} =
\int_{[0,1]^2} {-\log (1-xy) \over xy} {dx \,dy \over 1-xy}\ .$$
Using partial fractions with respect to the variable $y$,
$$I_2= \int_{[0,1]^2} {-\log(1-xy) \over xy}  - {\log (1-xy) \over 1-xy} dx \,dy =
\int_0^1 {1\over z} \Big[ \Li_2(xy) + {1\over 2}
\log^2(1-xy)\Big]_0^1 \, dx\ .$$ We conclude that
\begin{equation}\label{secondexample} I_2= \int_0^1 {\Li_2(x) \over x} + {\log^2 (1-x) \over 2x} \, dx =
 \Big[ \Li_3(x) +  \Li_{1,2}(x)\Big]_0^1 = \zeta(3)+ {\zeta(1,2)}\ .
\end{equation}
At each stage one can verify that the canonical primitives we have used above do not introduce any new poles along the boundary of the associahedron $\partial \overline{X}_{S,\delta}$. For example, the first primitive can
be written using $(\ref{dihedralsformo6})$:
\begin{equation} \label{examplefirstprim} {-\log(1-xyz) \over xy} {dx \, dy\over 1-xy} = -{\log ( u_{13}) \over u_{24} u_{25} }\, \omega_{5,\delta}\ ,
\end{equation}
where $\omega_{5,\delta}=(1-xy)^{-1}\, dx\, dy$ is the pull-back
of the canonical $2$-form on $ D_{26} \subset  \Modf_{0,S}$ along
the map $(x,y,z)\mapsto (x,y):\Mod_{0,6}\rightarrow \Mod_{0,5}$.
It has no poles at finite distance by definition (lemma
$\ref{omegahasnopoles}$). Furthermore, we  know by $(\ref{id})$
that $1-u_{13} = u_{24}u_{25} u_{26}$, so
$(\ref{examplefirstprim})$ has no poles along
$\overline{X}_{S,\delta}$ as required.

To show how lower weight multiple zeta values can occur, let us also compute
$$I_3=\int_{0\leq t_1\leq t_2\leq t_3\leq 1} {dt_1 \over 1-t_1 } {dt_2 \over t_3-t_1} {dt_3 \over t_3} =
\int_{[0,1]^3} {y\,dx\, dy\, dz \over (1-xyz)(1-xy)} =
\int_{\overline{X}_{S,\delta}} u_{25}\, u_{14}
\,\omega_{S,\delta}\ .$$ By applying a dihedral rotation, and
referring to $(\ref{dihedralsformo6})$, this is just:
$$I_3= \int_{X_{S,\delta}} u_{14} u_{36} \, \omega_{S,\delta}= \int_{[0,1]^3} {dx\, dy\, dz \over (1-xyz)^2}\ .$$
Integrating with respect to the variable $x$, we obtain by
$(\ref{firstexample})$:
\begin{equation} I_3= \int_{[0,1]^2} \Big[ {1\over 1-xyz} \Big]_0^1 {dy\, dz \over yz} = \int_{[0,1]^2} {dy \,dz \over 1-yz} = \zeta(2)\ .
\end{equation}
%The reason why the weight of $I_3$ is two is therefore that the
%cohomology class of the integrand is $0$ in $H^3(\Mod_{0,6})$.

\section{Appendix}
Let $M$ be the complement of an affine hyperplane configuration, as defined in $\S3$, and let $F$ denote its ring of regular functions.
 In the case where the de Rham cohomology ring  $H^\star(F)$ only has quadratic relations, we can prove
directly that all higher cohomology groups of $B(F)$ vanish.
\begin{thm} Let $F$ be the ring of regular functions on an affine  hyperplane arrangement $M$ such that $H^\star(F)$ is a quadratic algebra.
Then
$$H^0_{\DR}(B(F))=k \ ,\quad \hbox{and} \quad H^i_\DR(B(F))=0 \qquad \hbox{ for all }\quad  i\geq 1 \ .$$
\end{thm}

\begin{proof}
Let $A\subset \Omega^*(F)$ denote the  algebra generated by the
$1$-forms $\omega_1,\ldots, \omega_N\in H^1(F)$, where
$\omega_i=d\log \alpha_i$ as defined in $\S\ref{sect32}$.
 Let $V_1(F) =
\bigoplus_{i=1}^N k\, \omega_i$. Let $T$ denote the free tensor
algebra over $V_1$. We can view $A$ as a quotient algebra of $T$
of the form $ A = T / T Q T\ ,$ where $Q$ consists of finitely
many quadratic relations $q_1,\ldots, q_t$ of the form
\begin{equation} \label{aquadrel}
q_l=\sum_{i,j} \lambda^l_{ij} \omega_i\otimes \omega_j  \qquad \hbox{ for  } 1\leq l \leq t\ .
\end{equation}
Now
recall that $\Omega^\star (B(F)) = \Omega^\star(F) \otimes \bigoplus_{m\geq 0} V_m(F)$. The weight filtration on $\Omega^\star (B(F))$,
defined by $W_m \Omega^i B(F) = \Omega^i(F) \otimes \bigoplus_{j=0}^m V_j(F)$, gives rise
 to a spectral sequence with $E_0$ terms
$$E_0^{p,q} (B(F))= \Omega^{p+q}(F)\otimes V_p(F)\ , $$
 which
is bounded below and exhaustive, and therefore converges to the
cohomology of $B(F)$. Consider the differential subalgebra
$A(F)=A\otimes \bigoplus_{m\geq 0} V_m(F)$ of
$\Omega^\star(B(F))$. It also defines a spectral sequence with
$E_1$ terms
$$E_1^{p,q} (A(F))= H^{p+q}(A) \otimes V_p(F) \cong H^{p+q}(F) \otimes V_p(F)=E_1^{p,q}(B(F))\ .$$
It follows that $H^i(B(F)) \cong H^i(A(F))$ for all $i\geq 0$, and it suffices to show that $H^\star(A(F))$ is trivial.
Therefore consider an element
$$f=\sum_{I=(i_1,\ldots, i_m)} \sum_{J=(j_1,\ldots, j_n)}
 \alpha_{I,J}\, \omega_{j_1} \wedge \ldots \wedge \omega_{j_n} [\omega_{i_1}|\ldots |\omega_{i_m}]\in A\otimes V_m(F)\ , $$
such that $df=0$. This implies that
$$ \sum_{i_1,J}  \alpha_{I,J}\, \omega_{j_1} \wedge \ldots \wedge \omega_{j_n} \wedge \omega_{i_1}=0 \quad \hbox{ for all } \,\,i_2,\ldots, i_m\ .  $$
Because $A$ is quadratic, this expression (viewed in the tensor algebra $T$) decomposes as a sum of relations of the form
$w_1 q^l w_2$, where $w_1, w_2\in T$. If $w_2$ is of degree $\neq 0$ in $T$,  the corresponding relation is already zero in $f$. We deduce that $f$
is a sum of terms
$$\sum_{I=(i_2,\ldots, i_m)}\sum_{l=1}^t \sum_{i,j} \lambda^l_{ij}\, \omega_{I,l}\wedge \omega_i\, [ \omega_j | \omega_{i_1}|\ldots |\omega_{i_m}]\ ,$$
where $\omega_{I,l} \in A$. Each such expression has a primitive
$$\sum_{I=(i_2,\ldots, i_m)}\sum_{l=1}^t \sum_{i,j}  \omega_{I,l}\, \lambda^l_{ij}\, [\omega_i| \omega_j | \omega_{i_1}|\ldots |\omega_{i_m}]\ .$$
This is integrable: it maps to 0 under $\wedge_k$ for $k\geq 2$,
because $f$ is integrable, and it maps to $0$ under $\wedge_1$
because of the quadratic relations $(\ref{aquadrel})$.
 It follows that every element $f\in A(F)$ of weight $\geq 1$ such that $df=0$ has a primitive.
 It is easy to write down a primitive of any  element of weight $0$. Thus
$H^i(A(F))=0 $ for all $i\geq 1$, which completes the proof of the theorem.
\end{proof}

\section{Bibliography}
%
%Br1 - CRAS note

%Br2 - preprint on hyperlogarithms%

%Br3 - \emph{Mixed Tate motives and the irrationality of
%$\zeta(3)$}, in preparation.

\begin{itemize}

\item[[An]] Y.\, Andr\'e, {\it Une introduction aux motifs}, Panoramas et Synth\`eses {\bf 17},
S.M.F., Paris (2004)

\item[[Ao1]]  K.\ Aomoto, {\it Fonctions hyperlogarithmes et
groupes de monodromie unipotents}, Journal Fac. Sci. Tokyo, {\bf 25}
(1978), 149-156.

\item[[Ao2]]  K.\ Aomoto, {\it Addition theorem of Abel type for
hyperlogarithms}, Nagoya Math. Journal, {\bf 88} (1982), 55-71.

\item[[Ao3]] K.\ Aomoto, {\it Special values of hyperlogarithms and linear difference schemes}, Illinois J. Math. {\bf 34},
no. 2 (1990), 191-216.

\item[[Ar]] V. I. Arnold, {\em The cohomology ring of the coloured braid group}, Mat. Zametki {\bf 5} (1969), 227-231;
Math Notes {\bf 5} (1969), 138-140.

 \item[[B]]  A.\ Beauville, {\it Monodromie des syst\`emes
diff\'erentiels lin\'eaires \`a p\^oles simples sur la sph\`ere
de Riemann}, S\'em. Bourbaki, (1992-93), no. 765.

\item[[Be]] F. Beukers, {\it A note on the irrationality of
$\zeta(2)$ and $\zeta(3)$}, Bull. London Math. Soc {\bf 11}, (1979),
268-272.

\item[[BM]] L.\ Boutet de Monvel, {\it Polylogarithmes},
\verb"www.math.jussieu.fr/~boutet".

\item[[B-S]] A.\ Borel, J.\ P.\ Serre, {\it Corners and arithmetic groups}, and appendix: A.\ Douady, L.\ H\'erault, {\it
Arrondissement des vari\'et\'es \`a coins}, Comment Math. Helv. {\bf 48} (1973), 436-491.

\item[[Bou]]  N.\ Bourbaki, {\it  Alg\`ebre,} Hermann, Paris.

\item[[Bri]] E. Brieskorn, {\em Sur les groupes des tresses}, in S\'em. Bourbaki 1971/2, Lecture Notes in Math. {\bf 317}, Springer-Verlag
(1973), 21-44.

\item[[Br1]]   F.\ C.\ S.\  Brown, {\it Polylogarithmes
multiples uniformes en une variable}, C. R. Acad. Sci.  Paris,
Ser. I 338 (2004), 527-532.

\item[[Br2]] F. C. S. Brown, {\it Single-valued hyperlogarithms and
unipotent differential equations}, preprint.

\item[[Br3]] F. C. S. Brown, {\it P\'eriodes des espaces des modules $\overline{\Mod}_{0,n}$ et valeurs z\^eta multiples}, to appear
in C. R. Acad. Sci.  Paris, (2006).

  \item [[Ca1]] P. Cartier, {\em D\'eveloppements r\'ecents sur les groupes des tresses}, S\'em. Bourbaki, 1989-90, No. 716, Ast\'erisque (1989),
  17-67.

\item[[Ca2]]  P.\ Cartier, {\it Fonctions polylogarithmes,
nombres polyz\^etas et groupes pro-unipotents}, S\'eminaire
Bourbaki, (2000-01), no. 885.

\item[[Ch1]]  K.\ T.\ Chen, {\it Iterated path integrals}, Bull.
Amer. Math. Soc. {\bf 83}, (1977), 831-879.

\item[[Ch2]]  K.\ T.\ Chen, {\it Extension of $C^{\infty}$
Function Algebra by Integrals and Malcev completion of $\pi_1^*$},
Adv. in Math. {\bf 23}, no. 2 (1977), 181-210.

\item[[dC-P]] C.\ De Concini, C.\ Procesi, {\it Wonderful models
of subspace arrangements.}, Sel. Math., New Ser. 1, {\bf 3}
(1995), 459-494.

\item[[De1]] P. Deligne, {\em Le groupe fondamental de la droite projective moins trois points}, in `Galois groups over $\overline{\Q}$',
Proc. Workshop, Berkeley/CA (1987), Publ. Math. Sci. Res. Inst. {\bf 16} (1989), 79-297.

  \item[[De2]] P. Deligne, {\em Les immeubles des groupes de tresses g\'en\'eralis\'es}, Invent. Math. {\bf 17} (1972), 273-302.

\item[[De3]] P. Deligne, {\em Equations diff\'erentielles \`a points singuliers r\'eguliers},  Lecture Notes in Math. {\bf 163},
Springer-Verlag (1970).

\item[[De4]] P. Deligne, {\em Th\'eorie de Hodge  III}, Publications  Math\'ematiques de l'IH\'ES, {\bf 44} (1974), 5-77  

\item[[D-M]] P.\ Deligne, D.\ Mumford, {\it The irreducibility of
the space of curves of a given genus}, Publ. Math. IHES {\bf 36},
(1969), 75-109.

\item[[D-G]] Deligne, P., Goncharov, A.B.: {\em Groupes
fondamentaux motiviques de Tate mixte}, Ann. Sci. Ecole Norm. Sup., S\'er. IV {\bf 38}, No. 1, (2005), 1-56.

\item[[Dev1]] S. Devadoss, {\it Tesselations of moduli spaces and the mosaic operad}, in Contemp. Math. {\bf 239} (1999), 91-114.

\item[[Dev2]] S. Devadoss, {\it Real compactification of Moduli spaces},
Notices of the A.M.S {\bf 51}, No.6 (2004), 620-628.

%\item[[Dev-R]] S. Devadoss, R. Read, {\it Cellular structures determined by polygons and trees}, ?
\item[[Di]] A. C. Dixon, On a certain double integral, Proc. London Math. Soc. {\bf 2}, no.2 (1905), 8-15.

\item[[Dr]] V. G. Drinfeld, {\it On quasitriangular  quasi-Hopf algebras and on a group that is closely connected with
$\mathrm{Gal}(\overline{\Q}/\Q)$}, Leningrad Math. Journal  {\bf 2} (1991), No. 4, 829-860.

\item[[E]] J.\ Ecalle, {\it Singularit\'es non abordables par la
g\'eom\'etrie}, Ann. Inst. Fourier, Grenoble, 42, 1-2 (1992),
73-164.

%\item[[ENS]]  S\'eminaire de l'ENS (1979-1980), Partie IV, {\it
%Math\'ematique et Physique},  Progress in Math. {\bf 37},
%Birkh\"{a}user, Boston.

  \item[[F-R]] M. Falk, R. Randell, {\em The lower central series of a fiber-type arrangement}, Invent. Math. {\bf 82} (1985), 77-88.

\item[[Fi1]] S. Fischler, {\it Irrationalit\'e de valeurs de z\^eta [d'apr\`es Ap\'ery, Rivoal,..]},
S\'eminaire Bourbaki 2002-2003, Expos\'e 910, Ast\'erisque (2002).

\item[[Fi2]] S. Fischler, {\it Groupes de Rhin-Viola et
int\'egrales multiples}, J. Th\'eorie des Nombres de Bordeaux {\bf 15},
(2003), 479-534.

\item[[G-G-L]] H.\ Gangl, A.\ B.\ Goncharov, A.\ Levin, {\it Multiple polylogarithms, polygons, trees and algebraic cycles}, preprint
(2005) \verb"arXiv:math.NT/0508066".

\item[[Ge]] E. Getzler, {\it Operads and moduli spaces of genus $0$ Riemann surfaces}, Prog. Math {\bf 129} (1995), 199-230.

\item[[Go1]]  A.\ B.\ Goncharov, {\it Multiple polylogarithms and
mixed Tate motives}, preprint (2001),
\verb"arXiv:math.AG/0103059v4".

\item[[Go2]]  A.\ B.\ Goncharov, {\it Periods and
mixed motives},
\verb"arXiv:math.AG/0202154"  (2001).

\item[[Go3]]  A.\ B.\  Goncharov, {\it Multiple $\zeta$-values, Galois
groups and the geometry of modular varieties}, Proceedings of the
third European Congress of Mathematics, (2000).

\item[[Go4]] A.\ B.\ Goncharov, {\it The dihedral Lie algebras and Galois symmetries of $\pi_1^{(l)}(\Pro^1-(\{0,1\}\cup \mu_N))$},
Duke Math. J. {\bf 110}, No. 3, (2001), 397-487.

\item [[G-M]] A. B. Goncharov, Y. I. Manin, {\em Multiple $\zeta$-motives and moduli spaces $\overline{\Mo}_{0,n}$}, Compositio Math.
{\bf 140} (2004), 1-14.

\item[[GL]]   J.\ Gonzalez-Lorca, {\it S\'erie de Drinfel'd,
monodromie et alg\`ebres de Hecke}, Thesis, Ecole Normale Sup\'erieure, (1998).

\item[[G-S]] P.\ Griffiths, W.\  Schmid, {\it Recent developments in Hodge theory}, in `Discrete subgroups of Lie groups and applications
to moduli',  Internat. Colloq. Bombay 1973, Oxford Univ. Press (1975), 31-127.

\item[[Gr]] A.\ Grothendieck, {\it On the de Rham cohomology of
algebraic varieties}, Publ. Math. IHES {\bf 29} (1966), 95-103.

  \item [[Ha1]] R.\ M.\ Hain, {\em The geometry of the Mixed Hodge
Structure on the Fundamental Group}, Proc. Sympos. in Pure Math.
{\bf 46}, (1987), 247-282.
\item[[Ha2]] R.\ M.\ Hain, {\em Classical Polylogarithms}, Proc. Sympos. Pure Math. {\bf 55}, (1994), Part 2.
\item[[H-M]] R.\ M.\ Hain, R.\ MacPherson, {\em Higher logarithms}, Illinois J. Math. {\bf 34}, No. 2 (1990), 392-475.

\item[[H-Z1]] R.\ M.\ Hain, S.\ Zucker, {\em A guide to unipotent variations of Mixed Hodge structures}, Lecture Notes in Math.  {\bf 1246}, Springer (1987), 92-106.
\item[[H-Z2]] R.\ M.\ Hain, S.\ Zucker, {\em  Unipotent variations of Mixed Hodge structure},
Invent. Math {\bf 88}, no. 1 (1987), 83-124.

 \item[[H-P-V]]  M.\ Hoang Ngoc, M.\ Petitot, J.\ Van Den
Hoeven, {\it Polylogarithms and shuffle algebra},  FPSAC '98,
Toronto, Canada, (June 1998).

\item[[Ho]] M.\ E.\ Hoffman, {\it Quasi-shuffle products}, J.
Algebraic Combin. {\bf 11}, 1 (2000), 49-68.

\item[[Ih]] Y.\ Ihara, {\it Profinite braid groups, Galois
representations, and complex multiplication}, Ann. of Math. {\bf
123} (1986), 3-106.

 \item [[Ka]] M. M. Kapranov, {\em The permutoasssociahedron, MacLane's coherence theorem and asymptotic zones for the KZ equation},
 Journal of Pure and Applied Algebra {\bf 85}, (1993), 119-142.

 \item[[Ke]] S.\ Keel, {\em Intersection theory of moduli space of stable $N$-point curves of genus zero}, Trans. AMS {\bf 330}, no.2 (1992), 545-574.
\item[[Kh]] T.\ Kohno, {\em Bar complex of the Orlik-Solomon algebra}, Topology and its applications {\bf 118}, (2002) 147-157.
\item[[Kh2]] T.\ Kohno, {\em S\'erie de Poincar\'e-Koszul associ\'ee aux groupes de tresses pures},
 Invent. Math. {\bf 82}  (1985), 57-75.

\item[[Ko]] E. R. Kolchin, {\em Differential algebra and Algebraic groups}, Academic Press, New York (1976).

%\item[[Ko-Ma]] M. Kontsevich, Y. Manin, {\em Gromov-Witten classes, quantum cohomology, and enumerative geometry}.

\item[[Ku]]  E.\ E.\ Kummer, {\it \"{U}ber die Transcendenten,
welche aus wiederholten Integrationen rationaler Formeln
entstehen}, J. Reine Angew. Mathematik, {\bf 21}, (1840), 74-90,
193-225, 328-371.

\item[[K-Z]]  V.\ G.\  Knizhnik, A.\ B.\ Zamolodchikov, {\it
Current algebras and Wess-Zumino models in two dimensions}, Nucl.
Phys. B {\bf 247} (1984), 83-103.

\item[[Kn]] F.\ F.\ Knudsen, {\it The projectivity of the moduli space of stable curves II. The stacks $\overline{M}_{0,n}$}, Math.
Scand. {\bf 52} (1983), 163-199.

\item[[Ko-Za]] M. Kontsevich, D. Zagier, {\it Periods}, dans
Mathematics unlimited - 2001 and beyond, Ed. Engquist and Schmidt,
pp. 771-808, Springer, 2001.

\item[[LD]]  I.\ Lappo-Danilevskii, {\it M\'emoires sur la
th\'eorie des syst\`emes des \'equations diff\'erentielles
lin\'eaires}, Chelsea, New York (1953).

\item[[Lee]] C. Lee, {\it The associahedron and triangulations of the $n$-gon}, European Journal of Combinatorics {\bf 10} (1989), 551-560.

\item[[Le]]  {\it Structural Properties of Polylogarithms}, ed.
L.\ Lewin, Math. surveys and monographs {\bf 37}, AMS (1991).

\item[[Le-M]]  T.\ Q.\ T.\ Le, J.\ Murakami,
 {\it Kontsevich's integral for the Kauffman polynomial},
Nagoya Math. Journal {\bf 142}, (1996) 39-65.

\item[[Ma]] A. Magid, {\em Lectures on Differential Galois Theory}, Univ. Lecture Series {\bf 7} (1994).

\item[[Man]] Y. Manin, {\em Frobenius manifolds, quantum cohomology, and moduli spaces}, AMS Colloquium Publications {\bf 47}, Providence,
RI (1999).

\item[[Oe]]  J.\ Oesterl\'e, {\it Polylogarithmes}, S\'em.
Bourbaki 1992-93, exp. no. 762, Ast\'erisque  {\em 216} (1993), 49-67.

\item[[O-S]] P. Orlick, L. Solomon, {\em Combinatorics and topology of complements of hyperplanes}, Invent. Math. {\bf 56}, (1980), 167-189.
\item[[O-T]] P. Orlick, H. Terao, {\em Arrangements of hyperplanes}, Grundlehren der Math. Wiss. {\bf 300}, Springer-Verlag (1992).

\item[[Po]]  H.\ Poincar\'e, {\it Sur les groupes d'\'equations
lin\'eaires}, Acta Mathematica {\bf 4}, (1884).

\item[[Ra]] G. Racinet, {\em Double shuffles of multiple polylogarithms at roots of unity}, Publ. Math.  Inst. Hautes Etudes. Sci.
{\bf 95}, 185-231 (2002).

\item[[Rad]]  D.\ E.\ Radford, {\it A natural ring basis for
shuffle algebra and an application to group schemes,} Journal of
Algebra {\bf 58}, (1979), 432-454.

\item[[Re]]  C.\ Reutenauer, {\it Free Lie Algebras}, London
Math. Soc. Mono. {\bf 7}, Clarendon Press, Ox. Sci. Publ., (1993).

\item[[RV1]] G. Rhin, C. Viola, {\it On a permutation group
related to $\zeta(2)$}, Acta. Arith. {\bf 77}, no. 1 (1996),
23-56.
\item[[RV2]] G. Rhin, C. Viola, {\it The group structure for
$\zeta(3)$}, Acta. Arith. {\bf 97}, no. 3 (2001), 269-293.

\item[[S-V]] M.\ F.\ Singer, M.\ Van der Put, {\it Galois theory
of linear differential equations}, Grundlehren der Math. Wiss.
{\bf 328},  Springer (2003).

\item[[St]] J. D. Stasheff, {\it Homotopy associativity of $H$-spaces I}, Trans. Amer. Math. Soc. {\bf 108} (1963), 275-292.

\item[[Ter]] T. Terasoma, {\it Selberg integrals and multiple zeta values},
Compos. Math {\bf 133}, (2002), no. 1, 1-24.

\item[[V]] C.\ Voisin, {\it Hodge theory and complex algebraic geometry II}, Cambridge Studies in Advanced Math. {\bf 77}, (2003).

\item[[W]] M.\ Waldschmidt, {\it Valeurs z\^etas multiples : une introduction}, J. Th\'eorie des Nombres de Bordeaux {\bf 12} (2002), 581-595.

\item[[Yo]] M.\ Yoshida, {\it Fuchsian differential equations}, Aspects of Mathematics E11,  Braunschweig 1987.

\item[[Zh]] J.\ Zhao, {\it Multiple polylogarithms: analytic continuation, monodromy,
and variations of mixed Hodge structures}, in `Contemporary trends in algebraic geometry and algebraic
topology', Tianjin, China (2000),  Nankai Tracts in Math. {\bf 5}, (2002) 167-193.

\item[[Zl]] S.\ A.\ Zlobin, {\it Integrals expressible as linear
forms in generalized polylogarithms}, Math. Zametki [Math. Notes]
{\bf 71} (2002), 782-787 [711-716].

\item[[Zl2]] S.\ A.\ Zlobin, {\it Properties of coefficients of
certain linear forms in generalized polylogarithms}, preprint
(2005), \verb"arXiv:math.NT/0511245v1".

\item[[Zu]] W. Zudilin, {\it Well-poised hypergeometric
transformations of Euler-type multiple integrals}, J. London Math.
Soc. {\bf 70}, no. 2 (2004), 215-230.

%\bibitem{D-G} Deligne, P., Goncharov, A.B.: {\em Groupes
%fondamentaux motiviques de Tate mixte}, Ann. Sci. Ecole Norm.
%Sup., S\'er. IV {\bf 38}, No. 1, (2005), 1-56.

%\bibitem{Vi} C. Viola, {\it Birational transformations
%and values of the Riemann zeta function}, J. Th\'eor. Nombres
%Bordeaux, No. 15 (2003), 561-592 (Actes des 12\`emes Rencontres
%Arithm\'etiques de Caen (29-30 June 2001)).

%\bibitem{Wa} M. Waldschmidt, {\it Multiple Polylogarithms: an Introduction},

%$[$Ho$]$  M.\ E.\ Hoffman, {\it The alebra of multiple
%harmonic series}, Journal of Algebra 194 (1997), 477-495.
%\\

%$[$Rac$]$  G.\ Racinet, {\it Doubles m\'elanges des
%polylogarithmes multiples aux racines de l'unit\'e}, Publ. Math.
%Inst. Hautes \'Etudes Sci. no. 95, (2002), 185-231.
%\\

% $[$W1$]$  Z.\ Wojtkowiak, {\it Mixed Hodge Structures and
%Iterated Integrals I},  Motives, Polylogarithms and Hodge Theory
%(I), F. Bogomolov et L. Katzarkov, International Press, (2002).
%\\

\end{itemize}

\newpage

\section{Index of notations}
We list the most frequently used notations, along with the section in which they are first defined.
The integers $n=\ell+3$, where $\ell\geq 0$, are fixed.

Section 2:
\begin{description}
\item[$\delta$] \qquad A dihedral structure on a set $S$ with $n$ elements.
\item[$(\Pro^1)^n_*$] \qquad The set of $n$ distinct points in the projective line.
  \item[$\Mod_{0,S}$] \qquad The moduli space of curves of genus $0$ with points marked by $S$.
  \item[$\Modf_{0,S}$] \qquad The partial blow-up of $\Mod_{0,S}$ with respect to $\delta$.
   \item[$\overline{\Mod}_{0,S}$] \qquad The full blow-up of $\Mod_{0,S}$.
   \item[$\Or(\Mod_{0,S})$] \qquad The ring of regular functions on $\Mod_{0,S}$.
  \item[$\chi_{S,\delta}$] \qquad The set of all chords in the $n$-gon $(S,\delta)$.
  \item[$\chi^k_{S,\delta}$] \qquad The set of all partial $k$-triangulations of $(S,\delta)$.
  \item[$\{i,j\}\x\{k,l\}$] \qquad  The chords $\{i,j\},\{k,l\} \in \chi_{S,\delta}$ cross.
\item[$\{u_{ij} : \{i,j\}\in \chi_{S,\delta}\}$] \qquad The set of dihedral coordinates on $\Mod_{0,S}$.
\item[$(x_1,\ldots,x_\ell)$] \qquad The set of cubical coordinates on $\Mod_{0,S}$.
\item[$(t_1,\ldots, t_\ell)$] \qquad The set of simplicial coordinates on $\Mod_{0,S}$.
\item[$(x_1^{\alpha},\ldots, x_\ell^{\alpha})$] \qquad The set of vertex coordinates corresponding to $\alpha \in\chi^\ell_{S,\delta}$.
\item[$X_{S,\delta}$]\qquad The open associahedron $X_{S,\delta}\subset \Mod_{0,S}(\R)$.
\item[$\overline{X}_{S,\delta}$] \qquad The closed associahedron $\overline{X}_{S,\delta} \subset \Modf_{0,S}(\R)$.
\item[$F_{ij}$] \qquad The face of $\overline{X}_{S,\delta}$ corresponding to the chord $\{i,j\}\in \chi_{S,\delta}$.
\item[$F_\alpha$] \qquad The intersection of faces $\bigcap_{\{i,j\}\in\alpha} F_{ij}$ corresponding to  $\alpha \in \chi^k_{S,\delta}$.
\item[$f_T$] \qquad The forgetful map $f_T :\Mod_{0,S}\rightarrow \Mod_{0,T}$.
\item[$m_{\Box}$] \qquad The cubical multiplication map.
\item[$m_{\triangle}$] \qquad The simplicial multiplication map.
\end{description}
Section 3:
\begin{description}
  \item[$A$] \qquad An alphabet.
  \item[$\Z\langle A \rangle$] \qquad The free tensor algebra on $A$.
  \item[$\sha$] \qquad The shuffle product.
  \item[$\Delta$] \qquad The coproduct on the shuffle algebra.
  \item[$\varepsilon$] \qquad The counit on the shuffle algebra.
  \item[$\partial_a$]\qquad The  operator acting by  truncation on the left.
  \item[$M=\A^\ell\backslash \cup_{i=1}^N H_i$] \qquad The complement of an affine hyperplane arrangement.
  \item[$\Or_M$] \qquad The ring of regular functions on $M$.
  \item[$B(\Or_M)=B(M)$] \qquad The reduced bar construction on $M$.
  \item[$B_{\Or_{M'}}(E)$] \qquad The relative bar construction on $E$ with coefficients in $\Or_{M'}$.
\item[$\kes$] \qquad The differential $k-$algebra of Laurent series in
$\epsilon_1,\ldots, \epsilon_\ell$.
\item[$\Ues$] \qquad The differential $k$-algebra of logarithmic
Laurent series.
\item[$\Up(R,\varepsilon)$] \qquad The category of unipotent pointed extensions
of $(R,\varepsilon)$.
\item[$\Ut(R,p)$] \qquad The category of unipotent pointed extensions
of $(R,p)$, where $p$ is a base point at infinity.

\item[$\Urel$] \qquad The relative unipotent closure with
respect to a one-dimensional  fibration $R\rightarrow
\widehat{R}$.
% \item[$\Z\widehat{\pi}_1$] \qquad The pro-unipotent completion of the fundamental group algebra of $M$. %$\Z[\pi_1]$.
\end{description}
Section 4:
\begin{description}
  \item[$U_{p,q}$] \qquad The open real complement of  the coordinate hyperplanes in $\R^{p} \times \R^q_+$.
  \item[$V_{p,q}$] \qquad The open complex complement of $q$ coordinate hyperplanes in $\C^{p+q}$.
  \item[$\Fo^{\an}$] \qquad The sheaf of analytic functions.
  \item[$\Fo^{\log}$] \qquad The sheaf of analytic functions with logarithmic singularities.
  \item[$\Fo^{\log}_p$] \qquad  The sheaf of analytic functions with ordinary and logarithmic poles.
\end{description}
Section 5:
\begin{description}
  \item[$\Sigma$] \qquad A set of points $\{\sigma_0,\ldots,\sigma_{N}\} \subset \Pro^1$.
  \item[$\Or_{\Sigma}$] \qquad The ring of regular functions on $\Pro^1\backslash \Sigma$.
\end{description}
Section 6:
\begin{description}
\item[$\C^{n-1}_*$] \qquad The configuration space of $n-1$ distinct points in $\C$.
  \item[$\Delta_{ij}$] \qquad Logarithmic $1$-forms $\Delta_{ij}=d \log (z_i-z_j)$ on  $(\Pro^1)_*^{n}$.
\item[$t_{ij}$]\qquad Generators of the infinitesimal braid algebra, where $1\leq i\leq j\leq n$.
\item[$\Omega_{KZ}$] \qquad The Knizhnik-Zamolodchikov $1$-form on $(\Pro^1)_*^n$.
\item[$\delta_{ij}$] \qquad Infinitesimal dihedral braid elements, indexed by $\{i,j\}\in\chi_{S,\delta}$.
 \item[$\omega_{ij}$] \qquad Logarithmic  1-forms $\omega_{ij}=d\log u_{ij}$.
\item[$\Omega_{S,\delta}$] \qquad The canonical dihedral 1-form on $\Mod_{0,S}$.
\item[$\DB_{S,\delta}$] \qquad The dihedral braid algebra.
\item[$\widehat{\DB}_{S,\delta}$] \qquad The completion of the dihedral braid algebra.
\item[$V^{\delta}$]\qquad The set of vertices of the associahedron $\overline{X}_{S,\delta}$.
\item[$L_{v,\delta}$]\qquad The generating series of polylogarithms on $\Mod_{0,S}$.
\item[$L^{v,\delta}(\Mod_{0,s})$] \qquad The $\Q(\Mod_{0,S})$-algebra of polylogarithms on $\Mod_{0,S}$ whose regularised  value at the vertex $v\in V^\delta$ is zero.
\item[$\rho_{v,\delta}$]\qquad The realisation isomorphism $\rho_{v,\delta}:B(\Mod_{0,S}) \rightarrow L^{v,\delta}(\Mod_{0,S})$.
\item[$\MZV$] \qquad The ring $\Q[\zeta(2),\zeta(3),\ldots]$ generated by all multiple zeta values.
\end{description}
Section 7:

\begin{description}
  \item[$\omega_{S,\delta}$] \qquad The canonical volume form on $\Mod_{0,S}(\R)$.
  \item[$I_{S,\delta}(\alpha_{ij})$]\qquad The period integral over $X_{S,\delta}$.
\item[$m_{S,\delta}(\omega)$] \qquad The framed mixed Tate motive defined by $\omega$.
\item[$\Omega(\underline{\epsilon})$] \qquad The multiple zeta volume form.
\item[$L^\delta_\MZV(\Mod_{0,S})$] \qquad The algebra of polylogarithms on $\Mod_{0,S}$ with $\MZV$ coefficients.
\item[$L^{\delta,+}_\MZV(\Mod_{0,S})$] \qquad The subalgebra of polylogarithms  with  no poles along $\partial X_{S,\delta}$.
\end{description}

%\newpage
%\section{A further example}
%This is not part of the paper. I have included this extra example for the benefit of the informal reader, and for myself when giving
%talks. We work in cubical coordinates on $\Mod_{0,7}$, where  $(x,y,z,w)=(x_1,x_2,x_3,x_4)$. Consider
%$$I_4= \int_{[0,1]^4} {dx\, dy \, dz\, dw \over (1-xyzw)(1-zw)} =\int_{\triangle} {dt_1 \, dt_2 \, dt_3 \, dt_4 \over (1-t_1)t_2 (t_1-t_3)  t_4}
%=\int_{X_{7,\delta}}  u_{51} u_{14} u_{47} u_{73} u_{31} u_{14} \,\omega_{7,\delta}  \ .$$
%Then by taking partial fractions with respect to the fibration variable $w$ this is
%$$ \int_{[0,1]^4} {dx\,dy\,dz\,dw \over 1-xy} \Big({ 1\over 1-zw}- {xy \over 1-xyzw}   \Big)  =
%\int_{[0,1]^3} { dx\,dy\,dz\over (1-xy)} \Big[ -{\log (1-zw) \over z}+{\log( 1-xyzw)\over z}  \Big]_0^1  \ .$$
%Therefore
%$$I_4= \int_{[0,1]^3} {dx\, dy \over 1-xy} \Big( {\log (1-xyz) \over z} - {\log (1-z) \over z} \Big) dz=
%\int_{[0,1]^2} {dx\, dy \over 1-xy} \Big[  \Li_2(z)-\Li_2(xyz) \Big]_0^1\ . $$
%At this stage, we have
%$$I_4= \int_{[0,1]^2} {dx\, dy \over 1-xy} \Big(\zeta(2)-\Li_2(xy)\Big) =
%\zeta(2)I_1 -\int_{[0,1]} {dx \over x} \Big[ \Li_{1,2}(xy)  \Big]_0^1 \ .$$
%But we know that $I_1=\zeta(2)$ by the first example in $\S8$, hence
%$$I_4=\zeta(2)^2- \int_0^1 {\Li_{1,2} (x) \over x} dx =  \zeta(2)^2 - \Big[ \Li_{2,2}(x)\Big]_0^1  = \zeta(2)^2- \zeta(2,2)=\zeta(2,2)+\zeta(4)\ .$$
%

%By applying dihedral rotations, this same integral could also be written
%$$I_4= \int { xy\, dx\, dy\, dz\, dw \over (1-xy)(1-xyz)(1-wxyz)} = \int {yz^2 w \, dx \, dy \, dz\, dw \over (1-wxyz)(1-zw)(1-yz)}$$
\end{document}